\documentclass[a4paper,10pt]{book}

\usepackage[utf8x]{inputenc}
\usepackage{tracefnt,amsmath,tabu,array}
\usepackage{amssymb,graphicx,setspace,amsfonts,amsbsy}
\usepackage{pifont,latexsym,ifthen,amsthm,rotating,calc,textcase,booktabs,cancel,slashed}

\newtheorem{theorem}{Theorem}[chapter]
\newtheorem{lemma}[theorem]{Lemma}
\newtheorem{corollary}[theorem]{Corollary}
\newtheorem{proposition}[theorem]{Proposition}

\newtheorem{remark}[theorem]{Remark}
\newcommand{\filledbox}{\leavevmode
  \hbox to.77778em{%
  \hfil\vbox to.675em{\hrule width.6em height.6em}\hfil}}

\newcommand{\Rm}{{\mathbb R}}

\newcommand{\eps}{\varepsilon}

\begin{document}
%\doublespacing
\tabulinesep=1.0mm
% Enter full title and short title for running headers
\title{On the wave equation with Coulomb potential}

\author{Liang Li, Shenghao Luo and Ruipeng Shen\\
Centre for Applied Mathematics\\
Tianjin University\\
Tianjin, China
}

\maketitle

\frontmatter
 Wave/Schr\"{o}dinger equations with potentials naturally originates from both the quantum physics and the study of nonlinear equations. The distractive Coulomb potential is a quantum mechanical description of distractive Coulomb force between two particles with the same charge. The spectrum of the operator $-\Delta +1/|x|$ is well known and there are also a few results on the Strichartz estimates, local and global well-posedness and scattering result about the nonlinear Schr\"{o}dinger equation with a distractive Coulomb potential. In the contrast, much less is known for the global and asymptotic behaviour of solutions to the corresponding wave equations with a Coulomb potential. 

 In this work we consider the wave equation with a distractive Coulomb potential in dimensions $d\geq 3$. We first describe the asymptotic behaviour of the solutions to the linear homogeneous Coulomb wave equation, especially their energy distribution property and scattering profiles, then show that the radial finite-energy solutions to suitable defocusing Coulomb wave equation are defined for all time and scatter in both two time directions, by establishing a family of radial Strichartz estimates and combining them with the decay of the potential energy. 

\tableofcontents

\mainmatter

\chapter{Introduction} 

In this work we let the space dimension $d\geq 3$ and consider the linear wave equation with distractive Coulomb potential
\begin{equation} \label{linear wave equation with potential}
 \partial_t^2 u -\Delta u + \frac{u}{|x|} = 0,
\end{equation}
as well as the corresponding defocusing non-linear equation 
\begin{equation} \label{defocusing wave equation with potential}
 \partial_t^2 u -\Delta u + \frac{u}{|x|} + |u|^{p-1} u = 0.
\end{equation}
For convenience we may define 
\[
 \mathbf{H} = - \Delta + \frac{1}{|x|},
\]
write them in a unified form 
\begin{equation} \label{unified wave equation} 
 \partial_t^2 u + \mathbf{H} u + \zeta |u|^{p-1} u = 0.
\end{equation}
We will call them Coulomb wave equations for convenience in this work. We start by discussing the background about this the equation and the associated operator. 

\section{Background}

\paragraph{The Coulomb operator} The operator $\mathbf{H}$ originates from a quantum mechanical description of the distractive Coulomb force between two particles with the same charge. Since the potential $V(x) = 1/|x| \in L^2(\Rm^d) + L^\infty(\Rm^d)$, the operator $\mathbf{H}$ is essentially self-adjoint on $C_0^\infty(\Rm^d)$ and self-adjoint on $D(-\Delta)$, as long as $d\geq 3$. In addition, it can be proved by the H-smooth theory that $\mathbf{H}$ has purely absolutely continuous spectrum $\sigma(\mathbf{H})\subseteq [0,+\infty)$ and no eigenvalues. Please refer to Reed-Simon \cite{Simon}, Mizutani \cite{StrichartzSchrCoul}, Mizutani-Yao \cite{Hsmooth1}, for instance. The author would like to mention that if the potential is attractive, then the operator
\[
 -\Delta - \frac{K}{|x|}, \qquad K>0
\]
provides a decent approximation of the hydrogen atom. Its spectrum consists of a continuous part $[0,+\infty)$ and a family of negative, bounded, discrete, finite-multiplicity eigenvalues in $(-\infty,0)$. For more details one may refer to Taylor's book \cite{Taylorbook}. 
Although none of the operators mentioned above come with a non-negative eigenvalue, the radial solutions to the elliptic equation 
\[
 -\Delta u + \frac{u}{|x|} = \lambda u, \qquad x\in \Rm^n \setminus \{0\}
\]
with $\lambda\geq 0$ are still useful. In fact they correspond to the Coulomb wave functions, which are important in some aspect of modern physics, such as nuclear physics. More details about Coulomb wave functions can be found in \cite{handbook, zeros}, for example. 
 
 \paragraph{Schr\"{o}dinger equation} The Schr\"{o}dinger equation with a Coulomb potential has been discussed in a few works. For example, the corresponding Strichartz estimates of the linear equation has been established by Mizutani \cite{StrichartzSchrCoul}. The global existence, blow-up and scattering of solutions to the nonlinear Coulomb Schr\"{o}dinger equation in 3-dimensional space
 \[
  {\mathrm i} u_t + \Delta u + \frac{K}{|x|} = \pm |u|^{p-1} u
 \]
 was discussed in Miao-Zhang-Zheng \cite{SchrCoul} under suitable assumptions on the constants $K$, $p$ and initial data. A similar result for more general repulsive inverse-power potential $V(x) = c|x|^{-\sigma}$ with $\sigma\in (0,2)$ was later given by Dinh \cite{SchrInversePower} in all dimensions $d\geq 3$.  The authors would also mention that there are many works regarding the dispersive estimates, Strichartz estimates or scattering theory of solutions to the Schr\"{o}dinger equation with a potential $V$ which is either inverse-square, i.e. $V(x) = a|x|^{-2}$ or satisfies a similar decay estimates, $V(x) \in L^{n/2}(\Rm^n)$, for example. Please refer to Rodnianski-Schlag \cite{TimedecaySchlag}, Goldberg \cite{Goldberg1, Goldberg2, Goldberg3}, Zhang-Zheng \cite{SchrInverseSquare}, Mizutani \cite{SchrCritical}, for instance. Finally please note that there are also dispersive estimate results for a potential satisfying abstract assumptions. Here we give Barcel\'{o}-Ruiz-Vega \cite{abstractSchr} as an example.  

\paragraph{Wave equation} Much less is known for the wave equation with a Coulomb potential. We start by recalling a few results about the dispersive and Strichartz estimates of linear wave equations with potentials. First of all, dispersive/Strichartz estimates can be established for some potentials with a strong decay at the infinity. Please see Petkov \cite{WPcompact} for compact-supported potentials, Beals-Strauss \cite{WPfastsmall} for small and fast-decaying potentials, Vodev \cite{WPfasthigher} for potentials $|V(x)|\leq C(1+|x|)^{-(n+1)/2}$ in high dimensions $d\geq 4$, Cuccagna \cite{WPfast3d} for potentials with $|x|^{-3}$ decay in dimension 3, Costin-Huang \cite{WPfast1d} and Green \cite{WPfast2d} for 1 and 2 dimensional cases, respectively. There are also decay/Strichartz estimates for potentials $V(x)$ in the Kato class, i.e. 
 \[
  \|V\|_K = \sup_{x\in \Rm^n} \int_{\Rm^n} \frac{|V(y)|}{|x-y|^{n-2}} {\mathrm d} y < +\infty.  
 \]
Please refer to D'ancona-Pierfelice \cite{WPKato1}, Pierfelice \cite{WPKato2} and Bui-Duong-Hong \cite{WPKato3}, for instance. The inverse-square potentials $V(x)=a|x|^{-2}$ and other potentials with a similar decay at the infinity have also been extensively studied. The decay/dispersive estimates were discussed in Donninger-Schlag \cite{WIS1d}, Donninger-Krieger \cite{WIS1d2}, Georgiev-Visciglia \cite{WISGV} and Planchon-Stalker-Tahvildar-Zadeh \cite{WISPST1, WISPST2}. The corresponding Stichartz estimates were given by Burq, Planchon, Stalker and Tahvildar-Zadeh in \cite{WISStrichartz1, WISStrichartz2}. Miao-Zhang-Zheng \cite{WISStrichartzRadial} showed that the range of admissible pairs for Strichartz estimates can be extended if the initial data possess additional angular regularity. The authors would like to mention that there are also dispersive/Strichartz estimates for solutions to the wave equation with magnetic potentials
\[
 u_{tt} + ({\mathrm i}\nabla + A(x))^2 u + V(x) u = 0,
\]
where $V(x)\simeq |x|^{-2}$ is a roughly inverse-square potential. Please see, D'Ancona \cite{WMP1}, D'Ancona-Fanelli \cite{WMP2} and Fanelli-Zhang-Zheng \cite{WMP3}, for example. Finally reversed Strichartz estimate has also been considered in Chen \cite{StrichartzReversed} for a class of $L^2$ potential in $\Rm^3$. 

As for the global behaviour of solutions to nonlinear wave equation with potentials, most previous results are concerning potentials with roughly the same or higher decay than the inverse square potential $|x|^{-2}$. For example, Jia-Liu-Xu \cite{WPcm1} and Jia-Liu-Schlag-Xu \cite{WPcm2, WPcm3} discuss the long time dynamics and center stable manifold for the defocusing wave equation
\[
 u_{tt} - \Delta u + V(x) u + |u|^4 u = 0
\]
with a potential $V\in \mathcal{C}(\Rm^3)$ satisfying $|V(x)|\leq C(1+|x|)^{-\beta}$ for some $\beta>2$. Chen \cite{WPms} constructed multi-soliton to a defocusing energy-critical wave equation with a rapidly decaying potential and discussed their stability. Next we consider the wave equation with the inverse square potential
\[
 u_{tt} - \Delta u + \frac{a}{|x|^2} u = f(u). 
\]
Dai-Fang-Wang \cite{WISsmalldata} gives a long-time existence result for small initial data for $f(u)=|u|^p$. Miao-Murphy-Zheng \cite{WIScritical} considers the energy critical case $f(u) = \pm |u|^{4/(d-2)} u$ in 3 and 4 dimensional spaces and proves the global existence and scattering result for all initial data in the energy space in the defocusing case, and for initial data below the ground state threshold in the focusing case. Miao-Shen-Zhao \cite{WISsubRadial} then gives a scattering result for radial initial data in weighted energy space in the defocusing and energy sub-critical case. Both these two scattering results also make an assumption on the lower bound of the constant $a$. Mizutani \cite{WISnonradial} proves that the following wave operator is well-defined in the energy $\dot{H}^1 \times L^2$ for all $a>-(n-2)^2/4$:
\[
 {\mathrm s} - \lim_{t\rightarrow \pm \infty} \mathbf{S}_0 (-t) \mathbf{S}_a (t).
\]
Here $\mathbf{S}_a$ is the wave propagation operator associated the linear wave equation $u_{tt} - \Delta u + a|x|^{-2} u = 0$. In other words, linear free waves with an inverse square potential can be approximated by a linear free wave of the classic wave equation in the energy space. 

It seems that the inverse square decay is a boundary or critical case, especially for the wave equations. This is because the inverse-square potential has roughly the same strength as the Laplacian $-\Delta$. This can be understood in a few different but related ways:
\begin{itemize}
 \item If we consider the energy or pairing in the $L^2$ space
 \begin{align*}
  &\langle -\Delta u, u \rangle = \int_{\Rm^d} |\nabla u|^2 {\mathrm d} x;& &\langle |x|^{-2} u, u\rangle = \int_{\Rm^d} \frac{|u|^2}{|x|^2} {\mathrm d} x. 
 \end{align*}
 In dimension $d\geq 3$, the Hardy inequality guarantees that the contribution of the potential term can be dominated by the contribution of the Laplacian. This does not hold for any potential with a lower decay rate, in particular, the Coulomb potential. This also explains why the best constant $a = (n-2)^2/4$ in the Hardy inequality plays an important role in the theory with an inverse square potential, as we mentioned above. 
 \item For a typical function $f(x)\simeq |x|^{-\beta}$ with polynomial decay at the infinity, to apply the Laplacian is roughly to multiply by a constant time of $|x|^{-2}$. 
 \item The Laplacian operator $-\Delta$ and the multiplication operator by $|x|^{-2}$ share the same rescaling invariance property 
 \begin{align*}
  &\mathbf{H}_a \left(u(\lambda x)\right) = \lambda^2 (\mathbf{H}_a u)(\lambda x),& & \mathbf{H}_a = -\Delta + a|x|^{-2}. 
 \end{align*}
 \end{itemize}
As a result, a potential with a decay rate higher than the inverse square potential is called a short range potential; a potential with a lower decay rate than the inverse square potential is called a long range potential. The Coulomb potential is a typical long range potential, which imposes a stronger influence on the long-time property of solutions.  We shall give more details about this assertion later in this article. 

None of the results mentioned above cover the case of wave equation with a Coulomb potential. In fact, the authors find almost no literature on the details of asymptotic behaviours of solutions to the wave equation with a Coulomb potential. This work is an attempt to describe the global and asymptotic behaviour of solutions to both the linear and nonlinear wave equation with a distractive Coulomb potential. 

The authors would like to point out that the list of literature about the Schr\"{o}dinger/wave equation with a potential is far from complete. There are many more related works on this topic. 

\begin{remark}
  Also by a rescaling transformation, we may also show the the qualitative property of solutions to the linear homogeneous Coulomb wave equation 
 \[
  \partial_t^2 u - \Delta u + \frac{K}{|x|} u = 0
 \]
 depends on the sign of the constant $K$, but not on the size of $K$. 
\end{remark}

\begin{remark}
 There are also results about the global and asymptotic behaviour of solutions to Klein-Gordon equations with a suitable potential. Please see Bambusi-Cuccagna \cite{KGp1} and Soffer-Weinstein \cite{KGp2}, for instance. 
\end{remark}

\section{Goal of this paper} 

In this work we study the wave equation with a distractive Coulomb potential from the perspective of partial differential equations. More precisely we investigate the following properties of Coulomb wave equation 
\begin{itemize}
 \item As a preparation work, we consider the local and global well-posedness of Coulomb wave equations. Since the well-posedness of linear free Coulomb wave equation is guaranteed by the selfjointness of the operator $\sqrt{-\Delta + 1/|x|}$ and Stone's theorem, we only need to consider the well-posedness of non-linear Coulomb wave equation. For simplicity we mainly focus on the defocusing Coulomb wave equation 
\begin{align*}
 \left\{\begin{array}{l}\displaystyle \partial_t^2 u - \Delta u + \frac{u}{|x|} + |u|^{p-1} u=0, \qquad (x,t) \in \Rm^d \times \Rm;\\
 u|_{t=0} = u_0 \in \mathcal{H}^1(\Rm^d)\cap L^{p+1}(\Rm^d); \\
 u_t|_{t=0} = u_1 \in L^2(\Rm^d);\end{array}\right.
\end{align*}
with $3\leq d\leq 5$ and $p$ between conformal and energy critical exponent $1+\frac{4}{d-1} \leq p \leq 1+\frac{4}{d-2}$. Here the space $\mathcal{H}^s$ is defined in Section \ref{sec: notation}. In fact, this follows a combination of energy conservation law and a fixed-point argument in the energy sub-critical case. While in the energy critical case, this is a direct consequence of the global well-posedness and scattering of solutions to the energy critical wave equation $\partial_t^2 u - \Delta u + |u|^{\frac{4}{d-2}} u = 0$ and a suitable perturbation theory.
\item We are interested in the global and asymptotic behaviour of free Coulomb waves, i.e. how the solutions to 
\[ 
 \partial_t^2 u - \Delta u + \frac{u}{|x|} = 0. 
\]
looks like as the time tends to infinity. For example, we consider the law of energy distribution and conversion in the space-time. This also helps us establish a family of Strichartz estimates in the radial case, which plays an important role in the discussion of nonlinear equations. For another example, we shall show that the scattering profile of a free Coulomb wave can be approximated by that of a solution to the free Klein-Gordon equation via a geometric transformation. This helps us understand the energy dispersion rate and exclude the energy concentration phenomenon in the radial direction. 
\item Finally we consider the scattering of radial solutions to the defocusing Coulomb wave equations. Unlike the free wave equations, whose Strichartz estimates have to satisfy the scaling identity, the radial solutions to the free Coulomb wave equation satisfies a wide family of Strichartz estimates with flexible coefficients. This enable us to prove the scattering of energy solutions to defocusing Coulomb wave equation in both the energy critical and sub-critical case, as long as the potential energy $\int_{\Rm^d} |u(x,t)|{\rm d} x$ vanishes as the time tends to infinity. Please note that a similar argument does not work in the case of classic wave equation with an energy sub-critical exponent. In fact the scattering of solutions to the energy sub-critical wave equation $\partial_t^2 u - \Delta u + |u|^{p-1} u = 0$ is still an open problem, if only the finiteness of energy is assumed, even in the radial case. The authors would like to mention that the Strichartz estimates also give the scattering of small solutions to the intermediate focusing Coulomb wave equation. Please see 
\end{itemize}

\begin{remark} As we mentioned above, it has been proved in \cite{WISnonradial} that the asymptotic behaviour of solutions to homogeneous linear wave equation with inverse square potential 
\begin{equation} \label{inverse square}
 \partial_t^2 u - \Delta u + \frac{\alpha}{|x|^2} u = 0, \qquad \alpha > -\frac{(n-2)^2}{4}
\end{equation} 
is similar to that of solutions to the linear free wave equation. More precisely, if $u$ is a finite-energy solution to \eqref{inverse square}, then there exists a solution $v$ to the free wave equation $\square v = 0$, such that 
\[
 \lim_{t\rightarrow +\infty} \|(v,v_t)-(u,u_t)\|_{\dot{H}^1 \times L^2} = 0. 
\]
On the other hand, the asymptotic behaviour of solutions to the Klein-Gordon equation 
\[
 \partial_t^2 u - \Delta u + u = 0
\]
is much different from that of solutions to the free wave equation. It is reasonable to consider the equation between the wave equation with inverse square potential and the Klein-Gordon equation, namely
\[
 \partial_t^2 u - \Delta u + \frac{u}{|x|^s} =0, \qquad 0<s<2;
\]
and investigate the asymptotic behaviour of solutions to these equations. The wave equation with Coulomb potential is exactly at the middle point, i.e. $s=1$. In this work we will show that the asymptotic behaviour of solutions to this equation exhibits some aspects from both the wave and Klein-Gordon equations, as one may expect. 
\end{remark}

\section{Notations} \label{sec: notation}

We first introduce a few necessary notations before we give the details of our main results. Let us recall
\[
 \mathbf{H} = -\Delta + \frac{1}{|x|}.
\]
The homogeneous Sobolev space corresponding to the operator $-\Delta + 1/|x|$ can be defined by 
\[
 \mathcal{H}^s (\Rm^d) = \left\{u: \left\|\mathbf{H}^{s/2} u\right\|_{L^2} < +\infty \right\} = \left\{u: \left\|\left(-\Delta + \frac{1}{|x|}\right)^{s/2} u\right\|_{L^2} < +\infty \right\}, 
\]
with the norm $\|\mathbf{H}^{s/2} u\|_{L^2}$. In particular we have 
\[
 \mathcal{H}^1 = \left\{u\in H_{\rm loc}^1: \int_{\Rm^d} \left(|\nabla u|^2 + \frac{|u|^2}{|x|}\right){\rm d} x < +\infty \right\},
\]
with 
\[
 \|u\|_{\mathcal{H}^1(\Rm^d)} = \left(\int_{\Rm^d} \left(|\nabla u|^2 + \frac{|u|^2}{|x|}\right){\rm d} x\right)^{1/2};
\]
and 
\begin{align*}
 \mathcal{H}^2 = \left\{u\in H_{\rm loc}^2: \int_{\Rm^d}\left(|\Delta u|^2 + \frac{2|\nabla u|^2}{|x|} + \frac{|u|^2}{|x|^2} + (d-3)\frac{|u|^2}{|x|^3}\right){\rm d} x < +\infty\right\},
\end{align*}
with 
\begin{align*}
 \|u\|_{\mathcal{H}^2(\Rm^d)} & = \left(\int_{\Rm^d}\left(|\Delta u|^2 + \frac{2|\nabla u|^2}{|x|} + \frac{|u|^2}{|x|^2} + (d-3)\frac{|u|^2}{|x|^3}\right){\rm d} x\right)^{1/2}, & & d\geq 4; \\
 \|u\|_{\mathcal{H}^2(\Rm^3)} & = \left(\int_{\Rm^3}\left(|\Delta u|^2 + \frac{2|\nabla u|^2}{|x|} + \frac{|u|^2}{|x|^2} \right){\rm d} x + 4\pi |u(0)|^2 \right)^{1/2}, & & d= 3.
\end{align*}
If $(u_0,u_1)$ are initial data, then the corresponding free Coulomb wave, i.e. the solution to the linear homogeneous wave equation with the Coulomb potential 
\[
 \left\{\begin{array}{l}\displaystyle \partial_t^2 u - \Delta u + \frac{u}{|x|} = 0; \\ (u(0), u_t(0)) = (u_0,u_1);\end{array}\right.
\]
can be given in term of the operator $\mathbf{H}$ by the formula
\begin{align*}
 u (t) & = (\cos t \sqrt{\mathbf{H}}) u_0 + \frac{\sin t \sqrt{\mathbf{H}}}{\sqrt{\mathbf{H}}} u_1; \\
 u_t(t) & = (-\sqrt{\mathbf{H}} \sin t \sqrt{\mathbf{H}}) u_0 + (\cos t\sqrt{\mathbf{H}}) u_1. 
\end{align*}
This immediately gives us the general conservation law 
\begin{equation} \label{general conservation law}
 \|u(\cdot, t)\|_{\mathcal{H}^{s}}^2 + \|u_t(\cdot, t)\|_{\mathcal{H}^{s-1}}^2 = \|u_0\|_{\mathcal{H}^{s}}^2 + \|u_1\|_{\mathcal{H}^{s-1}}^2, \qquad \forall t\in \Rm. 
\end{equation}
Here $s\in \Rm$ is an arbitrary constant. In particular, the case $s=1$ is exactly the energy conservation law. 
\[
 E = \int_{\Rm^d} \left(\frac{1}{2}|\nabla u|^2 + \frac{1}{2}|u_t|^2 + \frac{|u|^2}{2|x|}\right) {\rm d} x = {\rm Const}. 
\]
For convenience we define 
\begin{align*}
 &\vec{\mathcal{H}}^s (\Rm^d)= \mathcal{H}^s (\Rm^d)\times \mathcal{H}^{s-1}(\Rm^d); & & \vec{v} = (v,v_t).
\end{align*}
Thus the conservation law of free Coulomb wave $u$ can be given by
\[
 \|\vec{u}\|_{\vec{\mathcal{H}}^s} = \|(u_0,u_1)\|_{\vec{\mathcal{H}}^s}, \qquad s\in \Rm. 
\]

\paragraph{Linear propagation operators} For convenience we will use the notations $\mathbf{S}_{\mathcal{C}} (t)$ and $\vec{\mathbf{S}}_{\mathcal{C}} (t)$ for the linear propagation operators of the Coulomb wave equation. Namely, if $(u_0,u_1)$ are initial data and $u$ is the corresponding free Coulomb wave, then we define
\begin{align*}
 &\mathbf{S}_{\mathcal{C}} (t) (u_0,u_1) = u(t);& &\vec{\mathbf{S}}_{\mathcal{C}} (t) (u_0,u_1) = (u(t), u_t(t)).
\end{align*}
Similarly we may define the linear propagation operators $\mathbf{S}_{\mathcal{K}}$ for the linear homogeneous Klein-Gordon equation
\[
 \partial_t^2 u - \Delta u + m^2 u = 0;
\]
and $\mathbf{S}_{\mathcal{W}}$ for the classic wave equation. 

\paragraph{Inequality with implicit constant} In this work we the notation $A\lesssim B$ means that there exists a constant $c$ so that $A \leq c B$. We may also add subscript(s) to the notation $\lesssim$ to emphasize that the constant $c$ depends on the subscript(s) but nothing else. In particular, the notation $\lesssim_1$ means that the implicit constant $c$ is an absolute constant. The meaning of notations $\simeq$ and $\gtrsim$ is similar.

\section{Main results}

Now we give three main results of this work. The first two results are concerning the asymptotic behaviour of free Coulomb waves in the energy space; while the last result is about the scattering of radial solutions to the defocusing Coulomb wave equation. 

\begin{theorem}[Asymptotic energy distribution and dispersion] \label{thm general energy distribution}
 Let $u$ be a finite-energy free Coulomb wave, i.e. a solution to the linear homogeneous Coulomb wave equation. Given any $\varepsilon > 0$, there exists two constants $0<c_1<c_2<+\infty$, so that 
 \begin{align*}
  \limsup_{t\rightarrow +\infty} \int_{|x|<t-c_2\ln t} e(x,t) {\rm d} x < \varepsilon; \\
  \limsup_{t\rightarrow +\infty} \int_{|x|>t-c_1\ln t} e(x,t) {\rm d} x < \varepsilon. 
 \end{align*}
 Here $e(x,t)$ is the energy density function defined by
 \[
  e(x,t) = \frac{1}{2}|\nabla u(x,t)|^2 + \frac{1}{2}|u_t(x,t)|^2 + \frac{|u(x,t)|^2}{2|x|}. 
 \]
 Furthermore, given a positive function $\ell(t)$ satisfying the growth condition
 \[
  \lim_{t\rightarrow +\infty} \frac{\ell(t)}{\ln t} = 0, 
 \]
 then we have 
 \[
  \lim_{t\rightarrow +\infty} \left(\sup_{r\geq 0} \int_{r<|x|<r+\ell(t)} e(x,t) {\rm d} x\right) = 0. 
 \]
\end{theorem}

\begin{remark}
 Clearly the energy distribution property of Coulomb free waves is much different from that of the classic wave equation. The majority of energy concentrates in a sphere shell of a constant thickness for the classic wave equation. In other words, the scattering does not happen in the radial direction at all. This also explains the decay rate of $\|u(t)\|_{L^\infty(\Rm^d)} \lesssim t^{-(d-1)/2}$ for sufficiently good initial data, since there are $d-1$ dimensions left if we exclude the radial direction. Let us consider another classic example of Klein-Gordon equation. In this case the scattering also fully happens in the radial direction, thus the decay rate is higher $\|u(t)\|_{L^\infty(\Rm^d)} \lesssim t^{-(d-1)/2}$. The theorem \ref{thm general energy distribution} implies that the Coulomb free waves do scatter in the radial direction, but the dispersion rate is at a very low rate $\ln t$, in comparison with the rate $t$ for the Klein-Gordon equation. 
\end{remark}

Before we introduce the next result, We introduce the following notations. Let $\mathcal{K}$ be the Hilbert space of finite-energy solutions to the one-dimensional homogeneous linear Klein-Gordon equation $v_{tt} - v_{xx} + 2v = 0$, equipped with the norm 
\[
 \|v\|_{\mathcal{K}}^2 = \int_{\Rm} \left(|v_x(x,t)|^2 + |v_t(x,t)|^2 + 2|v(x,t)|^2\right) {\rm d} x. 
\]
We also let $\mathcal{C}$ be the Hilbert space of radial free Coulomb waves with data in $(\mathcal{H}^1 \cap L^2) \times (L^2 \cap \mathcal{H}^{-1})$, equipped with the norm 
\[
 \|u\|_{\mathcal{C}}^2 = \|(u(\cdot,t),u_t(\cdot,t))\|_{\mathcal{H}^1\times L^2}^2 + \frac{1}{2}\|(u(\cdot,t),u_t(\cdot,t))\|_{L^2\times \mathcal{H}^{-1}}^2.
\]
Please note that the right hand sides in the definition of the norms given above are independent of time. 

\begin{theorem}[Scattering profile] \label{main transformation}
 Assume that $d\geq 3$. There exists a linear homeomorphism $\mathbf{T}: \mathcal{K}\rightarrow \mathcal{C}$ preserving the norm (up to a constant)
 \[
  \|\mathbf{T} v\|_{\mathcal{C}}^2 = \sigma_d \|v\|_{\mathcal{K}}^2. 
 \]
 In particular, if the data of $v$ are smooth and compactly supported, then the asymptotic behaviour of the image $\mathbf{T} v$ is similar to that of 
 \[
 u(x,t) = \rho\left(\frac{t-|x|}{t^{1/2}}\right) |x|^{-\frac{d-1}{2}} v\left(\frac{|x|-t+\ln (t+|x|)}{2}, \frac{t-|x|+\ln(t+|x|)}{2}\right). 
\]
Here $\rho: \Rm \rightarrow [0,1]$ is a smooth cut-off function satisfying
 \[
 \rho(s) = \left\{\begin{array}{ll} 0, & s\in [2,+\infty); \\ 1, & s\in (-\infty,1/2].\end{array}\right.
\]
More precisely, we have
 \[
  \lim_{t\rightarrow +\infty} \left\|(\mathbf{T} v - u, \partial_t(\mathbf{T} v - u))\right\|_{(\mathcal{H}^1 \cap L^2) \times (L^2 \cap \mathcal{H}^{-1})} = 0. 
 \]
 \end{theorem}

For the defocusing Coulomb wave equation, we have 

\begin{theorem} \label{main thm scattering}
 Assume $3\leq d\leq 5$ and $1+\frac{4}{d-1} < p \leq 1+\frac{4}{d-2}$. Let $(u_0,u_1)\in \mathcal{H}^1 \times L^2$ be radial initial data. Then the solution to the defocusing Coulomb wave equation 
 \[
  \left\{\begin{array}{l} \displaystyle \partial_t^2 u - \Delta u + \frac{u}{|x|} + |u|^{p-1} u = 0; \\ (u(\cdot,0), u_t(\cdot, 0)) = (u_0,u_1). \end{array}\right.
 \]
 is globally defined for all $t$ and scatters in both two time directions. More precisely, scattering implies that there exist two free Coulomb waves $u^\pm$ with finite energy, such that 
 \[
  \lim_{t\rightarrow \pm \infty} \left\|(u(\cdot,t),u_t(\cdot,t))- (u^\pm(\cdot,t), u_t^\pm(\cdot,t))\right\|_{\mathcal{H}^1\times L^2} = 0. 
 \]
\end{theorem}

\begin{remark} 
 Radial $\mathcal{H}^1$ functions satisfy the Sobolev embedding (see Corollary \ref{radial q0 infinity estimate})
 \[
  \|u\|_{L^{2+\frac{4}{2d-3}}(\Rm^d)} \lesssim_d \|u\|_{\mathcal{H}^1(\Rm^d)}.
 \]
 Combining this with the Sobolev embedding $\mathcal{H}^1 \hookrightarrow \dot{H}^1 \hookrightarrow L^\frac{2d}{d-2}$, we obtain the following Sobolev embedding of radial $\mathcal{H}^1$ functions
 \[
  \|u\|_{L^q(\Rm^d)} \lesssim \|u\|_{\mathcal{H}^1(\Rm^d)}, \qquad 2+\frac{4}{2d-3} \leq q \leq 2 + \frac{4}{d-2} 
 \]
 It follows that any radial initial data $(u_0,u_1) \in \mathcal{H}^1 \times L^2$ come with a finite energy
 \[
  E = \int_{\Rm^d} \left(\frac{1}{2}|\nabla u_0(x)|^2 + \frac{1}{2}|u_1(x)|^2 + \frac{1}{2} \frac{|u_0(x)|^2}{|x|} + \frac{1}{p+1} |u_0(x)|^{p+1}\right) {\rm d} x < +\infty,
 \]
 as long as $1+\frac{4}{d-1}\leq p \leq 1+\frac{4}{d-2}$. This is the reason why we do not need to assume $u_0 \in L^{p+1}$ explicitly in Theorem \ref{main thm scattering}.
\end{remark}

\section{The Structure of this work}

For reader's convenience we give a brief description of each chapter in this work before the end of this chapter. 
\begin{itemize}
\item In Chapter 2 we introduce the well-posedness theory of defocusing Coulomb wave equations and give a regularity result about the radial free Coulomb waves. 
\item In Chapter 3 we discuss the inward/outward energy theory, which gives energy distribution information of solutions to both the free and defocusing Coulomb wave equations. Although the information given by inward/outward energy theory is not detailed enough to prove the main theorems, this information is necessary for us to carry out further argument. 
\item In Chapter 4 we prove the energy retraction property of radial free Coulomb waves. Namely, given any finite-energy radial free Coulomb waves, any forward light cone will eventually contain almost all energy as time tends to $+\infty$. This is much different from a free wave without the Coulomb potential. This energy retraction property is a preparation work to prove the main theorems. 
\item In Chapter 5 we give the homomorphism from the scattering profiles of the Klein-Gordon equations to those of the Coulomb wave equation, thus proves the second main theorem. In addition, we show that the first main theorem is a direct consequence of the second one, if the free Coulomb wave is radial. 
\item In Chapter 6 we utilize the spherical harmonic function decomposition to deduce the energy distribution property of non-radial free Coulomb waves from the corresponding result in the radial case. 
\item In Chapter 7 we combine the results of inward/outward energy theory with the decay estimate of radial $\mathcal{H}^1$ functions to deduce a family of Strichartz estimates for radial free Coulomb waves. 
\item In Chapter 8 we combine the decay estimate given by inward/outward energy theory with the radial Strichartz estimates to prove the scattering of radial solutions to the defocusing intermediate Coulomb wave equation. 
\item Chapter 9 is an appendix, which contains some technical results to be used in the argument of this work. 
\end{itemize}

\begin{remark}
 Throughout this work we consider real-valued solutions to the Coulomb wave equations. Complex-valued solutions can be dealt with in a similar way. 
\end{remark}

\chapter{Well-posedness and regularity} 
In this chapter we consider the well-posedness theory of solutions to the defocusing Coulomb wave equation and the regularity of radial linear Coulomb waves with good initial data. 

\section{Well-posedness of defocusing equation} \label{sec: LWP}

Now we introduce our main result about the local and global well-posedness of solutions to the nonlinear wave equation 
\begin{equation} \label{nonlinear equation pre} 
 \left\{\begin{array}{l} \partial_t^2 u - \Delta u + \frac{u}{|x|} + |u|^{p-1} u = 0, \quad (x,t)\in \Rm^d \times \Rm; \\
  u(0) = u_0,\\ u_t(0) = u_1. \end{array} \right. 
\end{equation}
We focus on the dimension $3\leq d \leq 5$ and the nonlinearity $1+\frac{4}{d-1} \leq p \leq 1+\frac{4}{d-2}$ for simplicity of the argument. The author would like to mention that higher dimensions can be dealt with in a similar way. Our main result of this section is 

\begin{proposition} \label{prop global well}
 Let $(u_0,u_1)\in (\mathcal{H}^1\cap L^{p+1})(\Rm^d) \times L^2(\Rm^d)$ be initial data. Then there exists a unique solution $u$ defined for all time $t\in \Rm$ to the Coulomb wave equation \eqref{nonlinear equation pre} satisfying 
  \begin{align*}
  &(u,u_t) \in \mathcal{C}(\Rm; \mathcal{H}^1\times L^2);& &u \in L_{\rm loc}^\frac{2p}{(d-2)p-d} L^{2p} (\Rm \times \Rm^d);&
 \end{align*}
 and the energy conservation law 
 \[
  E = \int_{\Rm^d} \left(\frac{1}{2} |\nabla u(x,t)|^2 + \frac{1}{2} |u_t(x,t)|^2 + \frac{1}{2} \frac{|u(x,t)|^2}{|x|} + \frac{1}{p+1}|u(x,t)|^{p+1}\right) {\rm d} x = {\rm Const}. 
 \]
\end{proposition}

 We shall prove this in a few steps. We start by considering the local well-posedness theory. The local theory follows a standard fixed-point argument, combined with the Strichartz estimates of the wave equation. The observation that the potential term $u/|x|$ is at an energy sub-critical level plays an important role in the argument.  We first consider the energy sub-critical case $p<1+\frac{4}{d-2}$.
 
\begin{lemma} \label{local well-posedness sub-critical} 
 Assume that $3\leq d\leq 5$ and $1+\frac{4}{d-1} \leq p < 1+\frac{4}{d-2}$. Given any constant $M>0$, there exists a time $T=T(d,p,M)>0$, so that for any initial data $(u_0,u_1)$ with $\|(u_0,u_1)\|_{\dot{H}^1\times L^2} < M$, there is a unique solution $u$ to \eqref{nonlinear equation pre} in the time interval $[0,T]$ satisfying 
  \begin{align*}
  &(u,u_t) \in \mathcal{C}([0,T]; \dot{H}^1\times L^2);& &u \in L^\frac{2p}{(d-2)p-d} L^{2p} ([0,T]\times \Rm^d).
 \end{align*}
\end{lemma}

\begin{proof}
 This follows a classic fixed-point argument. We first recall the solution to the wave equation 
 \[
  \partial_t^2 u - \Delta u = f
 \]
 with initial data $(u_0,u_1)\in \dot{H}^1\times L^2$ satisfies the Strichartz estimates (see Ginibre-Velo\cite{strichartz}, for example)
 \[
  \|u\|_{L^{\frac{2p}{(d-2)p-d}}L^{2p}([0,T]\times \Rm^d)} \leq C_1 \left(\|(u_0,u_1)\|_{\dot{H}^1\times L^2} + \|f\|_{L^1 L^2([0,T]\times \Rm^d)}\right)
 \] 
 and 
 \[
  \|(u,u_t)\|_{\dot{H}^1\times L^2(\Rm^d)} \leq \|(u_0,u_1)\|_{\dot{H}^1\times L^2} +  \|f\|_{L^1 L^2([0,T]\times \Rm^d)}. 
 \]
 We consider the Banach space 
 \[ 
   X(T) = \mathcal{C}([0,T]; \dot{H}^1(\Rm^d)) \cap L^{\frac{2p}{(d-2)p-d}}L^{2p}([0,T]\times \Rm^d)
  \] 
  equipped with the norm 
  \[
   \|u\|_{X(T)} = \max_{t\in [0,T]} \|u\|_{\dot{H}^1(\Rm^d)} + \|u\|_{L^{\frac{2p}{(d-2)p-d}}L^{2p}([0,T]\times \Rm^d)}. 
  \]
  We then introduce a map $\mathbf{T}: X(T) \rightarrow X(T)$ by defining 
 \[ 
  \mathbf{T} u = \mathbf{S}_{\mathcal{W}} (t) (u_0,u_1) + \int_0^t \mathbf{S}_{\mathcal{W}}(t-\tau) (0, -|u|^{p-1} u(\tau) - u(\tau)/|x|) {\rm d} \tau.  
 \]
 It is equivalent to saying that $\mathbf{T} u$ is the solution to the linear wave equation 
 \[
  \partial_t^2 v - \Delta v = - \frac{v}{|x|} -|v|^{p-1} v
 \]
 with initial data $(u_0,u_1)$. It follows the Strichartz estimate above and Hardy inequality that 
 \begin{align*}
  \|\mathbf{T} u\|_{X(T)} & \leq C_1 \left(\|(u_0,u_1)\|_{\mathcal{H}^1\times L^2} + \left\|-|u|^{p-1} u-|x|^{-1} u\right\|_{L^1 L^2 ([0,T]\times \Rm^d)} \right)\\
  & \leq C_1 M + C_1 T^{\frac{d+2-(d-2)p}{2}} \left\|-|u|^{p-1}u\right\|_{L^\frac{2}{(d-2)p-d} L^2([0,T]\times \Rm^d)} +C_2 T \|u\|_{L^\infty ([0,T]; \dot{H}^1)}\\
  & \leq C_1 M + C_1 T^{\frac{d+2-(d-2)p}{2}} \|u\|_{X(T)}^p + C_2 T \|u\|_{X(T)}
 \end{align*}
 and that 
 \begin{align*}
  \|\mathbf{T} u - \mathbf{T} \tilde{u}\|_{X(T)} & \leq C_1 \left\|-|u|^{p-1}u + |\tilde{u}|^{p-1} \tilde{u}\right\|_{L^1 L^2([0,T]\times \Rm^d)} +C_1 \left\||x|^{-1} (\tilde{u}-u)\right\|_{L^1 L^2([0,T]\times \Rm^d)} \\
  & \leq C_1 T^{\frac{d+2-(d-2)p}{2}}  \left\|-|u|^{p-1}u + |\tilde{u}|^{p-1} \tilde{u}\right\|_{L^\frac{2}{(d-2)p-d} L^2([0,T]\times \Rm^d)} \\
  & \qquad + C_2 T \|\tilde{u}-u\|_{L^\infty([0,T]; \dot{H}^1)}\\
  & \leq C_1 p T^{\frac{d+2-(d-2)p}{2}} \|u-\tilde{u}\|_{X(T)} \left(\|u\|_{X(T)}^{p-1} + \|\tilde{u}\|_{X(T)}^{p-1}\right) + C_2 T \|\tilde{u}-u\|_{X(T)}. 
 \end{align*}
 Next we choose a small time $T = T(d,p,M)$ satisfying
 \begin{align*}
  C_1 T^{\frac{d+2-(d-2)p}{2}} (2C_1 M)^p + C_2 T (2C_1 M)& < C_1M;\\
   2 C_1 p T^{\frac{d+2-(d-2)p}{2}} (2C_1 M)^{p-1} + C_2 T & \leq \frac{1}{2};
 \end{align*}
 and let $\mathcal{X} = \{u: \|u\|_{X(T)} \leq 2C_1 M\}$. It is clear that $\mathcal{X}$ is a complete metric space with the distance
 \[
  d(u,\tilde{u}) = \|u-\tilde{u}\|_{X(T)}. 
 \]
 It immediately follows the inequalities above that $\mathbf{T}$ becomes a contraction map form $\mathcal{X}$ to itself. The unique fixed-point $u$ is exactly the solution to \eqref{nonlinear equation pre} with initial data $(u_0,u_1)$. 
\end{proof}

A similar result holds in the energy-critical case as well. Although the time period $T$ depends on not only the norm of initial data $(u_0,u_1)$ but also the profile of them. 

\begin{lemma}
 Assume $3\leq d\leq 5$ and $p=1+\frac{4}{d-2}$. Given any initial data $(u_0,u_1)\in \dot{H}^1 \times L^2$, there exists a time $T>0$, so that there is a unique solution to \eqref{nonlinear equation pre} with initial data $(u_0,u_1)$ satisfying 
 \begin{align*}
  &(u,u_t) \in \mathcal{C}([0,T]; \dot{H}^1\times L^2);& &u \in L^{\frac{d+2}{d-2}}L^{\frac{2(d+2)}{d-2}} ([0,T]\times \Rm^d).
 \end{align*}
\end{lemma}
\begin{proof}
 The argument is similar to the energy sub-critical case. We choose $X_1(T)$ and $X_2(T)$ to be the spaces 
 \begin{align*}
  &X_1(T) = L^{\frac{d+2}{d-2}}L^{\frac{2(d+2)}{d-2}}([0,T]\times \Rm^d);& &X_2(T) = \mathcal{C}([0,T]; \dot{H}^1(\Rm^d)).
 \end{align*}
 We first recall the solution to the wave equation 
 \[
  \partial_t^2 u - \Delta u = f
 \]
 with initial data $(u_0,u_1)\in \dot{H}^1\times L^2$ satisfies the Strichartz estimates
 \begin{align*}
  \|u\|_{X_1 (T)} + \|u\|_{X_2(T)} &\leq C \left(\|(u_0,u_1)\|_{\dot{H}^1\times L^2} + \|f\|_{L^1 L^2([0,T]\times \Rm^d)}\right);\\
  \|u\|_{X_2(T)} & \leq \|(u_0,u_1)\|_{\dot{H}^1\times L^2} + \|f\|_{L^1 L^2([0,T]\times \Rm^d)}. 
 \end{align*}
 Thus the map $\mathbf{T}$ defined by 
 \[
 \mathbf{T} u = \mathbf{S}_{\mathcal{W}} (t) (u_0,u_1) + \int_0^t \mathbf{S}_{\mathcal{W}}(t-\tau) \left(0, -|u|^{\frac{4}{d-2}} u(\tau) - u(\tau)/|x|\right) {\rm d} \tau  
 \]
 satisfies 
 \begin{align*} 
  \|\mathbf{T} u\|_{X_1 (T)} & \leq  \|\mathbf{S}_{\mathcal{W}}(u_0,u_1)\|_{X_1 (T)} + C_1 \|u\|_{X_1 (T)}^\frac{d+2}{d-2} + C_2 T \|u\|_{X_2 (T)}; \\
  \|\mathbf{T} u\|_{X_2 (T)} & \leq  \|(u_0,u_1)\|_{\dot{H}^1\times L^2} + \|u\|_{X_1(T)}^\frac{d+2}{d-2} + T \|u\|_{X_2(T)}; \\
  \|\mathbf{T} u - \mathbf{T} \tilde{u}\|_{X_1(T)\cap X_2(T)}  & \leq C_1 p \|u-\tilde{u}\|_{X_1(T)} \left(\|u\|_{X_1(T)}^{\frac{4}{d-2}} + \|\tilde{u}\|_{X_1(T)}^{\frac{4}{d-2}}\right) + C_2 T \|\tilde{u}-u\|_{X_2 (T)}.
 \end{align*}
 Here we may assume the constants $C_1, C_2$ satisfy $C_1, C_2 > 1$, without loss of generality. We first choose a small constant $\delta > 0$ so that 
 \begin{align*}
  &C_1 (3\delta)^\frac{d+2}{d-2} < \delta; & &2 C_1 p (3\delta)^\frac{4}{d-2} < \frac{1}{2};&
 \end{align*}
 then choose a small time $T>0$ so that 
 \[
  \|\mathbf{S}_{\mathcal{W}}(u_0,u_1)\|_{X_1 (T)} < \delta;
 \]
 and
 \begin{equation*}
  C_2 T (\|(u_0,u_1)\|_{\dot{H}^1\times L^2}+2 \delta) < \delta; \qquad \Rightarrow \qquad C_2 T < \frac{1}{2}.
 \end{equation*}
 Combining all the inequalities above we conclude that $\mathbf{T}$ is a contraction map from the metric space 
 \[
  \mathcal{X} = \{u\in X_1(T) \cap X_2(T): \|u\|_{X_1(T)}\leq 3\delta, \|u\|_{X_2(T)}\leq \|(u_0,u_1)\|_{\dot{H}^1\times L^2} + 2\delta\}
 \]
 with the distance 
 \[
  d(u,\tilde{u}) = \|u-\tilde{u}\|_{X_1(T)\cap X_2(T)} =  \|u-\tilde{u}\|_{X_1(T)} + \|u-\tilde{u}\|_{X_2(T)}. 
 \]
 Again the unique fixed-point $u$ is exactly the solution to \eqref{nonlinear equation pre} with the initial data $(u_0,u_1)$. 
\end{proof}

\begin{remark}
 This local theory works for the following defocusing/focusing equation as well:
 \[
  \partial_t^2 u - \Delta u \pm \frac{u}{|x|} \pm |u|^{p-1} u = 0.
 \]
 Here the signs of the Coulomb potential and the nonlinear term can be chosen independently. A review reveals that the argument above does not depends on the signs of the potential or the nonlinear term. 
\end{remark}

\paragraph{Maximal lifespan} Given a pair of initial data $(u_0,u_1) \in \dot{H}^1 \times L^2$, local theory discussed above guarantees that there is a unique solution in a small time interval. We may further extend its time interval of existence to the maximum and obtain a solution $u$ with a maximal lifespan $(-T_-,T_+)$. In particular, if $T_+ < +\infty$, i.e. the solution blows up in finite time, then we must have 
\[
 \|u\|_{L^{\frac{2p}{(d-2)p-d}}L^{2p}([0,T]\times \Rm^d)} = +\infty. 
\]
This is usually called the finite blow-up criterion. The argument is a standard procedure thus we omit the details here. 

\paragraph{Energy conservation law} Next we show that if the initial data $(u_0,u_1)$ satisfy $(u_0,u_1)\in (\mathcal{H}^1 \cap L^{p+1})\times L^2$, then local solutions given above satisfy the energy conservation law
\[
 E = \int_{\Rm^d} \left(\frac{1}{2} |\nabla u(x,t)|^2 + \frac{1}{2} |u_t(x,t)|^2 + \frac{1}{2}\frac{|u(x,t)|^2}{|x|} + \frac{1}{p+1}|u(x,t)|^{p+1}\right) {\rm d} x = {\rm Const}. 
\]
The proof depends on a smooth approximation technique. Let $\rho_1: \Rm \mapsto [0,+\infty)$ and $\rho_2: \Rm^d \mapsto [0,+\infty)$ be smooth cut-off functions satisfying 
\begin{align*}
 &\rho_1 \in \mathcal{C}_0^\infty(\Rm); & &\int_\Rm \rho_1(t) {\rm d} t = 1;& \\
 &\rho_2 \in \mathcal{C}_0^\infty(\Rm^d); & &\int_{\Rm^d} \rho_2 (x) {\rm d} x = 1. &
\end{align*}
We then define 
\[
 \rho_\varepsilon (x,t) = \frac{1}{\varepsilon^{d+1}} \rho_1 \left(\frac{t}{\varepsilon}\right) \rho_2 \left(\frac{x}{\varepsilon}\right),
\]
and a family of corresponding smooth approximation operators 
\[
 \mathbf{P}_\varepsilon u = \rho_\varepsilon \ast u.
\]
Let $u$ be a local solution to \eqref{nonlinear equation pre} defined in an open time interval $J$. We consider the approximated solution $u_\varepsilon = \mathbf{P}_\varepsilon u$, which solves the equation 
\[
 \partial_t^2 u_\varepsilon - \Delta u_\varepsilon = f_\varepsilon; 
\]
with 
\[
 f_\varepsilon = \mathbf{P}_\varepsilon \left(-\frac{u}{|x|} - |u|^{p-1} u\right).
\]
It is not difficult to see ($I \subset J$ is a bounded closed interval)
\begin{align*}
 (u_\varepsilon, \partial_t u_\varepsilon) &\rightarrow (u,u_t) & & {\rm in}\quad \mathcal{C}(I; (\dot{H}^1\cap L^{\frac{2d}{d-2}}) \times L^2);\\
 |x|^{-1} u_\varepsilon & \rightarrow |x|^{-1} u & & {\rm in}\quad \mathcal{C}(I; L^2(\Rm^d));\\
 u_\varepsilon & \rightarrow u & & {\rm in}\quad L^{\frac{2p}{(d-2)p-d}}(I; L^{2p}(\Rm^d)); \\
 f_\varepsilon & \rightarrow -\frac{u}{|x|} - |u|^{p-1} u & & {\rm in}\quad L^1(I; L^2(\Rm^d)); 
\end{align*}
Next we let 
\[
 e_\varepsilon(x,t) = \frac{1}{2} |\nabla u_\varepsilon(x,t)|^2 + \frac{1}{2} |\partial_t u_\varepsilon(x,t)|^2 + \frac{|u_\varepsilon(x,t)|^2}{2|x|} + \frac{|u_\varepsilon(x,t)|^{p+1}}{p+1}. 
\]
A direct calculation shows that ($t_1<t_2$, $R>0$)
\begin{align*}
 & \int_{|x|<R+t_2}e_\varepsilon(x,t_2) {\rm d} x  -  \int_{|x|<R+t_1}e_\varepsilon(x,t_1) {\rm d} x\\
 & \quad = \int_{t_1}^{t_2} \int_{|x|<R+t} \left(\nabla u_\varepsilon \cdot \nabla (\partial_t u_\varepsilon) + \partial_t u_\varepsilon \partial_{tt} u_\varepsilon + \frac{u_\varepsilon}{|x|} \partial_t u_\varepsilon + |u_\varepsilon|^{p-1} u_\varepsilon \partial_t u_\varepsilon \right) {\rm d} x {\rm d} t\\
 & \qquad + \int_{t_1}^{t_2} \int_{|x|=R+t} e_\varepsilon (x,t) {\rm d} S {\rm d} t\\
 & \quad = \int_{t_1}^{t_2} \int_{|x|<R+t} \partial_t u_\varepsilon \left(\partial_{tt} u_\varepsilon - \Delta u_\varepsilon + \frac{u_\varepsilon}{|x|} + |u_\varepsilon|^{p-1} u_\varepsilon\right) {\rm d} x {\rm d} t\\
 & \qquad +  \int_{t_1}^{t_2} \int_{|x|=R+t} \left(e_\varepsilon (x,t) + \partial_t u_\varepsilon \partial_r u_\varepsilon \right) {\rm d} S {\rm d} t\\
 & \quad = \int_{t_1}^{t_2} \int_{|x|<R+t} \partial_t u_\varepsilon \left( f_\varepsilon + \frac{u_\varepsilon}{|x|} + |u_\varepsilon|^{p-1} u_\varepsilon\right) {\rm d} x {\rm d} t\\
 & \qquad + \int_{t_1}^{t_2} \int_{|x|=R+t} \left(e_\varepsilon (x,t) + \partial_t u_\varepsilon \partial_r u_\varepsilon \right) {\rm d} S {\rm d} t. 
\end{align*}
The convergence above implies that 
\[
  f_\varepsilon + \frac{u_\varepsilon}{|x|} + |u_\varepsilon|^{p-1} u_\varepsilon \rightarrow 0 \quad {\rm in}\quad L^1 L^2([t_1,t_2]\times \Rm^d); 
\]
and that $\|\partial_t u_\varepsilon\|_{L^\infty L^2([t_1,t_2]\times \Rm^d)}$ is uniformly bounded as $\varepsilon \rightarrow 0$. In addition, we always have 
\[
 e_\varepsilon (x,t) + \partial_t u_\varepsilon \partial_r u_\varepsilon \geq 0
\]
Thus we may let $\varepsilon \rightarrow 0$ in the identity above to deduce that 
\begin{equation} \label{finite speed of energy}
 \int_{|x|<R+t_2} e(x,t_2) {\rm d} x  \geq  \int_{|x|<R+t_1} e(x,t_1) {\rm d} x, \qquad \forall R> 0. 
\end{equation}
Here we use the various convergence of $u_\eps$ to $u$. We then let $R\rightarrow +\infty$ and obtain $E(t_1) \leq E(t_2)$. If we substitute the forward light cone $|x|= R+t$ by the backward light cone $|x|=R-t$ in the argument above, we also have $E(t_2) \leq E(t_1)$. This immediately gives the energy conservation law. 

\paragraph{Global well-posedness} The global well-posedness in the energy sub-critical case immediately follows the energy conservation law and Lemma \ref{local well-posedness sub-critical}. 
Let $t$ be an arbitray time in the maximal lifespan of $u$. Then the energy conservation law implies that 
\[
 \|(u(t),u_t(t))\|_{\dot{H}^1\times L^2} \leq (2E)^{1/2} < +\infty. 
\]
Therefore we are able to choose $M = (2E)^{1/2}$, which does not depend on $t$, then apply Lemma \ref{local well-posedness sub-critical} to solve the equation starting from time $t$ and conclude that the maximal lifespan of the solution also contains at least the time interval $[t, t+T(d,p,M)]$. Since $T(d,p,M)$ is independent of $t$, it follows that the solution can never break down before $t=+\infty$. The negative time is similar because the wave equation is time-reversible. 

\paragraph{Energy critical case} The energy critical case $p = 1 + \frac{4}{d-2}$ is more subtle. In this case the argument given above does not work any more, since the minimal lifespan given in the fixed-point argument depends on not only the norm but also the profile of initial data. In order to prove the global well-posedness of solutions, we view \eqref{nonlinear equation pre} as a small perturbation of the classic defocusing wave equation 
\[
 \partial_t^2 u - \Delta u + |u|^\frac{4}{d-2} u = 0
\]
in a small time interval. More precisely we have

\begin{lemma} \label{perturbation lemma}
 Let $M, T$ be positive constants. Then there exist two positive numbers $\eps_0 = \eps_0(M,T)$ and $C(M,T)$ so that if a solution $v$ to the approximated equation
 \[
  \left\{\begin{array}{ll}\displaystyle \partial_t^2 v - \Delta v + \frac{v}{|x|} + |v|^\frac{4}{d-2} v = e(x,t), & (x,t)\in \Rm^d \times [0,T'];\\ (v(0),v_t(0)) = (v_0,v_1);& \end{array}\right.
 \]
 and another pair of initial data $(u_0,u_1)$ satisfy 
 \[
  \varepsilon \doteq \|e\|_{L^1 L^2([0,T']\times \Rm^d)} + \|(u_0,u_1)-(v_0,v_1)\|_{\dot{H}^1\times L^2} < \varepsilon_0,
 \]
 and 
 \begin{align*}
   &\|v\|_{L^\frac{d+2}{d-2} L^\frac{2(d+2)}{d-2}([0,T']\times \Rm^d)} < M,& &T'\leq T;&
 \end{align*}
 then the corresponding solution to the equation 
  \begin{equation} \label{perturbation equation 1} 
  \left\{\begin{array}{l}\displaystyle \partial_t^2 u - \Delta u + \frac{u}{|x|} + |u|^\frac{4}{d-2} u = 0;\\ (u(0),u_t(0)) = (u_0,u_1);  \end{array}\right.
 \end{equation} 
 is well-defined in the time interval $[0,T']$ with 
 \[
  \sup_{t\in [0,T']} \left\|((u(t)),u_t(t))-(v(t),v_t(t))\right\|_{\dot{H}^1\times L^2} + \|u-v\|_{L^\frac{d+2}{d-2} L^\frac{2(d+2)}{d-2}([0,T']\times \Rm^d)} \leq C(M,T)\eps. 
 \]
\end{lemma}

\begin{proof}
 Let $u$ be the solution to \eqref{perturbation equation 1} and $[0,T'']$ be a time interval contained in its maximal lifespan with $T''\leq T'$. For convenience we introduce the notation 
 \[
  \|w\|_{X(T'')} = \|w\|_{L^\frac{d+2}{d-2} L^\frac{2(d+2)}{d-2}([0,T'']\times \Rm^d)} + \max_{t\in [0,T'']} \|(w,w_t)\|_{\dot{H}^1\times L^2}.
 \]
 By the Strichartz estimates we have 
 \begin{align}
  \|u-v\|_{X(T'')} & \leq C_1 \|(u_0,u_1)-(v_0,v_1)\|_{\dot{H}^1\times L^2} \nonumber \\
   & \qquad + C_1 \left\|e - \frac{v}{|x|}-|v|^\frac{4}{d-2} + \frac{u}{|x|} + |u|^\frac{4}{d-2} u\right\|_{L^1 L^2 ([0,T'']\times \Rm^d)} \label{combined Strichartz}\\
   & \leq C_1 \varepsilon + C_2 T'' \|u-v\|_{X(T'')} + C_3 \|u-v\|_{X(T'')} \left(M^\frac{4}{d-2} + \|u-v\|_{X(T'')}^\frac{4}{d-2}\right) \nonumber
 \end{align}
 Here $C_1, C_2, C_3>1$ are all positive constants solely determined by the dimension $d$. We may choose small constants $M_0, T_0, \varepsilon_0$ so that if  $\varepsilon \leq \varepsilon_0(d)$ is a sufficiently small positive constant, then the inequality 
 \begin{equation} \label{contradiction inequality}
  2 C_1 \varepsilon > C_1 \varepsilon + C_2 T_0 (2C_1 \varepsilon) + C_3 (2C_1 \varepsilon)\left(M_0^\frac{4}{d-2} + (2C_1\varepsilon)^\frac{4}{d-2}\right)
 \end{equation}
 holds. We claim that if the approximation solution $v$ and the initial data $(u_0,u_1)$ satisfy 
 \begin{align*}
  \varepsilon \doteq \|e\|_{L^1 L^2([0,T']\times \Rm^d)} + \|(u_0,u_1)-(v(0),v_t(0))\|_{\dot{H}^1\times L^2} < \varepsilon_0
 \end{align*}
 \begin{align*}
  &\|v\|_{L^\frac{d+2}{d-2} L^\frac{2(d+2)}{d-2}([0,T']\times \Rm^d)} < M_0;& &T' \leq T_0;&
 \end{align*}
 then the corresponding solution $u$ is well-defined in $[0,T']$ and satisfies 
 \[
  \|u-v\|_{X(T')} < 2C_1 \varepsilon. 
 \]
 In fact, if this were false, i.e. the solution $u$ blows up before time $T'$ or $\|u-v\|_{X(T')} \geq 2 C_1 \varepsilon$, then we would find a time $T'' \in (0,T']$ so that $u$ still exists at the time $T''$ and $\|u-v\|_{X(T'')} = 2C_1 \varepsilon$, thanks to the continuity of the norm $X(T)$ with respect to $T$ and the finite time blow-up criterion. We then utilize the Strichartz estimates \eqref{combined Strichartz} and obtain 
 \[
  2C_1\varepsilon \leq  C_1 \varepsilon + C_2 T'' (2C_1 \varepsilon) + C_3 (2C_1 \varepsilon) \left(M_0^\frac{4}{d-2} + (2C_1 \varepsilon)^\frac{4}{d-2}\right). 
 \]
 This contradicts with \eqref{contradiction inequality}. This proves the lemma when $M, T$ are both small. If $M$ and/or $T$ is large, we may first split the time interval $[0,T']$ into a few sub-intervals $J_1, J_2, \cdots, J_n$ so that 
 \begin{align*}
  &\|v\|_{L^\frac{d+2}{d-2} L^\frac{2(d+2)}{d-2} (J_k \times \Rm^d)} \leq M_0,& &|J_k| \leq T_0, & & \forall k \in \{1,2,\cdots,n\}
 \end{align*}
 with 
 \[
  n \leq n_0(M,T) \doteq \left(\frac{M}{M_0}\right)^\frac{d+2}{d-2} + \frac{T}{T_0} + 2;
 \]
 and apply a bootstrap argument. The final upper bounds $\varepsilon_0 (M,T)$ and $C(M,T)$ can be chosen as 
 \begin{align}
  &\varepsilon_0 (M,T) = (2C_1)^{1-n_0(M,T)} \varepsilon_0(d);& &C(M,T) = 2(2C_1)^{n_0(M,T)}.
 \end{align}
 Therefore they only depend on $M,T$ and are independent of the specific approximated solution $v$ and initial data $(u_0,u_1)$. 
\end{proof}

Now we are at the position to prove the global well-posedness of the solutions in the energy critical case. Let $u$ be such a solution with an energy $E$ and $T$ is an arbitrary time in its maximal lifespan. By the energy conservation law, the norm
\[
 \|(u(T), u_t(T))\|_{\dot{H}^1\times L^2} \leq (2E)^{1/2} 
\]
is uniformly bounded. According to the global existence and scattering theory of the classic energy-critical wave equation (see Shatah-Struwe \cite{ss2} and Nakanishi \cite{enscatter1}, for instance)
\[
 \partial_t^2 v - \Delta v + |v|^\frac{4}{d-2} v = 0,
\] 
its solution $v$ with initial data $(v(T), v_t(T)) = (u(T), u_t(T))$ exists for all time $t \in \Rm$ and satisfies the inequality 
\[ 
 \|v\|_{L^\frac{d+2}{d-2} L^\frac{2(d+2)}{d-2} (\Rm \times \Rm^d)} \leq C_0.
\]
Here $C_0$ solely depends on the upper bound of $\|(u(T), u_t(T))\|_{\dot{H}^1\times L^2}$, thus can be chosen independent of $T$. In fact, it can be determined by the energy $E$ only. In addition, the energy conservation law implies that 
\begin{align*}
 \int_{\Rm^d} & \left(\frac{1}{2}|\nabla v(x,t)|^2 + \frac{1}{2}|v_t(x,t)|^2 + \frac{d-2}{2d}|v(x,t)|^{\frac{2d}{d-2}} \right) {\rm d} x \\
 & = \int_{\Rm^d} \left(\frac{1}{2}|\nabla u(x,T)|^2 + \frac{1}{2}|u_t(x,T)|^2 + \frac{d-2}{2d}|u(x,T)|^{\frac{2d}{d-2}} \right) {\rm d} x \leq E. 
\end{align*}
It follows Hardy's inequality that 
\[
 \left\|\frac{v(x,t)}{|x|}\right\|_{L^2(\Rm^d)} \leq C_1, \qquad \forall t\in \Rm. 
\]
Here $C_1$ is again a constant solely determined by the dimension $d$ and the energy $E$. Next we let $\varepsilon_0 = \varepsilon_0 (C_0,1)$ be the small constant given by Lemma \ref{perturbation lemma} and let $T_0 < 1$ be a small number so that 
\[
 C_1 T_0 < \varepsilon_0. 
\]
All these constants depend on only the dimension and the energy $E$. Now we observe that $v$ solves the approximated equation 
\[
 \partial_t^2 v - \Delta v + \frac{v}{|x|} + |v|^\frac{4}{d-2} v = \frac{v}{|x|}
\]
with 
\begin{align*}
 & \|v\|_{L^\frac{d+2}{d-2} L^\frac{2(d+2)}{d-2} ([T,T+T_0] \times \Rm^d)} \leq C_0; & &\left\|\frac{v(x,t)}{|x|}\right\|_{L^1 L^2([T, T+T_0] \times \Rm^d)} \leq C_1 T_0 < \varepsilon_0. 
\end{align*}
Thus we may apply Lemma \ref{perturbation lemma} on the approximated solution $v$ and the initial data 
\[
 (u(T), u_t(T))= (v(T), v_t(T))
\]
 to conclude that the solution $u$ is still well defined in the time interval $[T,T+T_0]$. Because $T$ is an arbitrary time in the maximal lifespan of $u$ but $T_0$ is a constant independent of $T$, we conclude that the solution $u$ never blows up in finite time. 
 
 \paragraph{Continuity in $\mathcal{H}^1\times L^2$} The remaining work to prove Proposition \ref{prop global well} is to show that $u \in \mathcal{C}(\Rm; \mathcal{H}^1)$. Let us recall $u\in \mathcal{C}(\Rm; \dot{H}^1)$. By the definition of $\mathcal{H}^1$ it suffices to show that 
 \begin{equation} \label{continuity of u weight L2}
  u \in \mathcal{C}(\Rm; L^2(\Rm^d; |x|^{-1} {\mathrm d} x)). 
 \end{equation}
 A combination of the continuity of $u$ in $\dot{H}^1$ and the Hardy inequality immediately gives that 
 \[
  u \in \mathcal{C}(\Rm; L^2(\Rm^d; |x|^{-2} {\mathrm d} x)).
 \]
 As a result, given any $R>0$, we have 
 \begin{equation} \label{continuity of u com 1} 
  u \in \mathcal{C}(\Rm; L^2(\{x: |x|<R\}; |x|^{-1} {\mathrm d} x)).
 \end{equation} 
 On the other hand, we may combine the energy conservation law and \eqref{finite speed of energy}, as well as the time symmetric property of wave equation, to deduce the energy inequality 
 \[
  \int_{|x|>|t|+R_1} e(x,t) {\mathrm d} x \leq \int_{|x|>R_1} e(x,0) {\mathrm d} x, \quad R_1>0. 
 \]
 This implies that for any given bounded time interval $J$, we have 
 \[
  \lim_{R\rightarrow \infty} \sup_{t\in J} \int_{|x|>R} e(x,t) {\mathrm d} x = 0. 
 \]
 Combining this tail estimate with \eqref{continuity of u com 1}, we verify \eqref{continuity of u weight L2} and finish the proof of Proposition \ref{prop global well}. 
 
\begin{remark}
 We view $\partial_t^2 u - \Delta u$ as the major linear part and all the other terms $+u/|x|+|u|^{p-1} u$ as the perturbation part in the local theory given above.  Of course, one may also view the linear part $\partial_t^2 u - \Delta u + u/|x|$ as the major part and the nonlinear term $|u|^{p-1}u$ as the perturbation part, as long as we might prove the corresponding Strichartz-type estimates. The reasons of our choice in this work include the following 
 \begin{itemize}
  \item We are only able to prove the Strichartz estimates with Coulomb potential in the radial case, which means that we can only give a local theory for radial solutions if we choose the second way. But our local theory given above may deal with non-radial solutions as well. 
  \item Our proof of radial Strichartz estimates depends on the inward/outward energy theory, or at least the Morawetz estimates for linear waves with Coulomb potential. In order to simplify our argument we hope to give an inward/outward energy theory for both linear and nonlinear solutions in a unified way. This means that we need to prove the global existence of nonlinear solutions before we may introduce the inward/outward energy theory and then prove the radial Strichartz estimates. 
 \end{itemize} 
\end{remark}

\section{Regularity of radial solutions}

We also consider the regularity of radial solutions to the homogeneous linear wave equation with the Coulomb potential. This guarantees that we have sufficient derivatives to work with in the argument, at least for sufficiently good initial data. 

\begin{lemma} \label{continuity of radial solutions} 
 Let $u$ be a radial free Coulomb wave with a finite energy. Then $u(x,t)$ is continuous in $(\Rm^d\setminus \{0\})\times \Rm$. 
\end{lemma}
\begin{proof}
 We recall that radial $\dot{H}^1(\Rm^d)$ functions $u$ are continuous in $\Rm^d \setminus \{0\}$ and satisfies the point-wise estimate (see Kenig-Merle \cite{km} and Shen \cite{shenhighradial})
 \[
  |u(x)| \lesssim_d \frac{|u|_{\dot{H}^1(\Rm^d)}}{|x|^{\frac{d-2}{2}}}.
 \]
 The continuity of a solution $u$ to linear homogeneous wave equation with the Coulomb potential immediately follows the observation above and the fact 
 \[ 
  u(t) \in \mathcal{C}(\Rm; \mathcal{H}^1(\Rm^d)) \hookrightarrow \mathcal{C}(\Rm; \dot{H}^1(\Rm^d)). 
 \]
\end{proof}

\begin{proposition} \label{C2 continuity} 
 Let $u$ be a radial free Coulomb wave with initial data $u_0,u_1\in \mathcal{C}_0^\infty (\Rm^d \setminus \{0\})$. Then $u$ satisfies $u \in \mathcal{C}^2((\Rm^d \setminus \{0\})\times \Rm)$. 
\end{proposition}

In order to prove this proposition, we start by proving an embedding property for radial $\mathcal{H}^s$ functions. 

\begin{lemma}
 Given $0<r_1<r_2<+\infty$, we let $K = \{x\in \Rm^d: r_1<|x|<r_2\}$. Let 
 \begin{align*}
  & \mathcal{H}_{\rm rad}^k (\Rm^d)= \left\{f: f \; {\rm radial}, \;  f\in \mathcal{H}^j(\Rm^d), \; j=1,2,\cdots, k \right\};& &\|u\|_{\mathcal{H}_{\rm rad}^k}= \sum_{j=1}^k \|u\|_{\mathcal{H}^j}. 
 \end{align*}
 Then we have the following embedding of spaces 
 \begin{align*}
  \mathcal{H}_{\rm rad}^k \hookrightarrow \mathcal{C}_{b}^{k-1}(K), \qquad k=1,2,3.  
 \end{align*}
 Here $\mathcal{C}_b^{k-1}(K)$ is the space of $\mathcal{C}^{k-1}(K)$ functions with a finite norm 
 \[
  \|f\|_{\mathcal{C}_b^{k-1}(K)} = \sup_{x\in K, |\alpha|\leq k-1} \left|\partial^\alpha f(x)\right|.
 \] 
\end{lemma}
\begin{proof}
 Let us start by $k=1$. Any function $u \in \mathcal{H}^1$ satisfies $u \in H_{\rm loc}^1(\Rm^d)$. Combining this with the radial assumption, we obtain that if $u$ is viewed as a function of the radius, then $u(r)$ is absolutely continuous in $[r_1,r_2]$ with $u_r (r) \in L^2([r_1,r_2])$. Therefore $u \in C(K)$. In addition, the point-wise estimate of $\dot{H}^1$ function gives 
 \[
  \|u\|_{\mathcal{C}_b (K)} = \max_{x\in K} |u(x)| \lesssim_d r_1^{-\frac{d-2}{2}} \|u\|_{\dot{H}^1(\Rm^d)} \lesssim_d r_1^{-\frac{d-2}{2}} \|u\|_{\mathcal{H}_{\rm rad}^1(\Rm^d)}. 
 \]
 Next we consider $k=2$. We recall that any function $u \in \mathcal{H}^2$ satisfies $u \in H_{\rm loc}^2(\Rm^d)$. Again a combination of this with the radial assumption implies that $u$ and $u_r$ are both absolutely continuous in $[r_1,r_2]$ with $u_{rr} (r) \in L^2([r_1,r_2])$. Thus we have $u \in \mathcal{C}^1(K)$. In addition, we utilize the formula of $\mathcal{H}^2$ norm given in  Section \ref{sec: notation} and deduce
 \[
  \int_{0}^\infty \left(\left|u_{rr} + \frac{d-1}{r} u_r\right|^2 + \frac{2|u_r|^2}{r} + \frac{|u|^2}{r^2} \right) r^{d-1} {\rm d} r \lesssim_d \|u\|_{\mathcal{H}^2}^2. 
 \]
 It follows that 
 \[
  \int_{r_1}^{r_2} \left(|w_{rr}|^2 + |w_r|^2\right) {\rm d} r \lesssim_{d,r_1,r_2} \|u\|_{\mathcal{H}^2}^2. 
 \]
 The Sobolev embedding $W^{1,2} (J)\hookrightarrow \mathcal{C}_b (J)$ with $J=(r_1,r_2)$ then gives 
 \[
  \sup_{r\in (r_1,r_2)} u_r (r) \lesssim_K \|u\|_{\mathcal{H}^2}. 
 \]
 Combining this with the upper bound on $\|u\|_{\mathcal{C}_b(K)}$, we obtain 
 \[
  \|u\|_{\mathcal{C}_b^1(K)} \lesssim_K \|u\|_{\mathcal{H}_{\rm rad}^2(\Rm^d)}. 
 \]
 Finally we consider $k=3$. Let $u \in \mathcal{H}_{\rm rad}^3 \hookrightarrow \mathcal{H}_{\rm rad}^2$. Then we have 
 \[
  -u_{rr} - \frac{d-1}{r} u_r + \frac{u}{r} = -\Delta u + \frac{u}{|x|} \in \mathcal{H}_{\rm rad}^1 (\Rm^d) \hookrightarrow \mathcal{C}(K). 
 \]
 We also have $u_r, u\in \mathcal{C}(K)$ by the fact $u \in \mathcal{H}_{\rm rad}^2 \hookrightarrow \mathcal{C}^1(K)$. Thus we have $u_{rr} \in \mathcal{C}(\mathcal{K})$. It follows that $u \in \mathcal{C}^2(K)$. In addition, we have 
 \[
  \sup_{r_1<r<r_2} \left|-u_{rr} - \frac{d-1}{r} u_r + \frac{u}{r}\right| = \left\|-\Delta u + \frac{u}{|x|}\right\|_{\mathcal{C}_b(K)} \lesssim_K \left\|-\Delta u + \frac{u}{|x|}\right\|_{\mathcal{H}_{\rm rad}^1} \leq \|u\|_{\mathcal{H}_{\rm rad}^3 (\Rm^d)}. 
 \]
 We have already obtained that 
 \[
  \sup_{r_1<r<r_2} (|u(r)|+|u_r(r)|) \lesssim_K \|u\|_{\mathcal{H}_{\rm rad}^2 (\Rm^d)}\leq  \|u\|_{\mathcal{H}_{\rm rad}^3 (\Rm^d)}.
 \]
 It follows that 
 \[
  \sup_{r_1<r<r_2} |u_{rr}| \lesssim_K \|u\|_{\mathcal{H}_{\rm rad}^3 (\Rm^d)}. 
 \]
 Thus 
 \[
  \|u\|_{\mathcal{C}_b^2(K)}\lesssim_K \|u\|_{\mathcal{H}_{\rm rad}^3 (\Rm^d)}. 
 \]
\end{proof}

\begin{proof}[Proof of Proposition \ref{C2 continuity}]
 First of all, we observe that the initial data satisfy $(u_0,u_1) \in \mathcal{H}^k \times \mathcal{H}^{k-1}$ for all integers $k \geq 1$. The conservation law immediately gives 
 \[
  (u,u_t) \in \mathcal{C}(\Rm; \mathcal{H}^k \times \mathcal{H}^{k-1}), \qquad \forall k\geq 1. 
 \]
 Therefore we have 
 \[
  (u,u_t) \in \mathcal{C}(\Rm; \mathcal{H}_{\rm rad}^3 \times \mathcal{H}_{\rm rad}^2). 
 \]
By the embedding given above, we have for any $K= \{x: r_1<|x|<r_2\}$ that
\begin{equation} \label{C2C1K}
 (u,u_t) \in \mathcal{C}(\Rm; \mathcal{C}^2(K) \times \mathcal{C}^1(K)). 
\end{equation}
Furthermore we have 
\begin{equation} \label{CK}
 u_{tt} = -\left(-\Delta + \frac{1}{|x|}\right) u \in \mathcal{C}(\Rm; \mathcal{H}_{\rm rad}^1 (\Rm^d)) \hookrightarrow \mathcal{C}(\Rm; \mathcal{C}(K)).
\end{equation}
Combining \eqref{C2C1K} and \eqref{CK}, we conclude that $u \in \mathcal{C}^2 (\Rm \times K)$ and finish the proof. 
\end{proof}

\begin{remark} \label{C3 continuity}
  The solutions $u$ are not likely to be $\mathcal{C}^2$ at the origin for all $t$ even if the initial data are as smooth as possible. Since this means that $u(x,t)/|x|$ is also continuous at the origin, or $u(0,t) = 0$, which is not true in general. 
\end{remark} 

\begin{remark}
 A similar argument as above gives that $\mathcal{H}_{\rm rad}^4 \hookrightarrow \mathcal{C}_{b}^{3}(K)$ for $K=\{x: a<|x|<b\}$ with $0<a<b<+\infty$. Thus $u\in \mathcal{C}^3(\Rm \times (\Rm^d\setminus \{0\}))$ if initial data $u_0,u_1\in \mathcal{H}_{\rm rad}^4 \times (\mathcal{H}_{\rm rad}^3\cap L^2)$.  In fact this assumption guarantees that $(u,u_t) \in \mathcal{C}(\Rm; \mathcal{H}_{\rm rad}^4 \times \mathcal{H}_{\rm rad}^3)$.
\end{remark}

\chapter{Inward/outward energy theory}

In this chapter we introduce the inward/outward energy theory of solutions to \eqref{unified wave equation}. Generally speaking, the inward/outward energy theory generalizes the classic Morawetz estimates of wave equations and gives plentiful information about energy distribution of solutions. We shall split the total energy into two parts, i.e. the inward and outward energy. Roughly speaking, the inward energy moves toward and the outward energy moves away from the origin. We then discuss the space-time distribution and transformation of the inward/outward energies, especially their asymptotic behaviours as the time tends to infinity. The results in this chapter will be frequently used in the subsequent chapters. We first introduce our assumptions on the Coulomb equation and a few notations for convenience. 

\paragraph{Assumptions} In this chapter we consider two kinds of Coulomb wave equations 
\begin{itemize}
 \item The homogeneous linear Coulomb wave equation $\partial_t^2 u - \Delta u + \frac{u}{|x|} = 0$ in all dimensions $d\geq 3$;  
 \item The defocusing Coulomb wave equation with immediate defocusing power-like nonlinearity
\[
 \partial_t^2 u - \Delta u + \frac{u}{|x|} + |u|^{p-1} u = 0, \qquad (x,t) \in \Rm^d \times \Rm
\]
with $3\leq d\leq 5$ and $1+\frac{4}{d-1} \leq p \leq 1+\frac{4}{d-2}$. 
\end{itemize} 
We may write them in a unified way 
\begin{equation} \label{unified equation inout}
  \partial_t^2 u - \Delta u + \frac{u}{|x|} + \zeta |u|^{p-1} u = 0, \qquad \zeta \in \{0,1\}.  
\end{equation} 

\begin{remark}
 The inward/outward energy theory also applies to wave equation $\partial_t^2 - \Delta u + f(x,u) = 0$ with more general nonlinear terms/potentials and/or in higher dimensions. The essential condition is that the following inequality holds for a constant $\mu > 0$ and all $(x,u)\in \Rm^d \times \Rm$
 \[
  H(x,u) \doteq \frac{n-1}{4} u f(x,u) - \frac{n-1}{2} F(x,u) - \frac{1}{2} x \cdot \nabla_x F (x, u) \geq \mu F(x,u) \geq 0,
 \]
 with 
 \[
  F(x,u) = \int_0^u f(x,v) {\rm d}v. 
 \]
 Our assumptions on $d$ and $p$ in the defocusing case above are technical conditions for the convenience of our discussion. 
\end{remark}

\paragraph{Energy density functions} We define $e(x,t)$ to be the regular energy density function 
\[  
 e(x,t)  = \frac{1}{2} |\nabla u|^2 + \frac{1}{2} |u_t|^2 + \frac{1}{2}\frac{|u|^2}{|x|} + \frac{\zeta}{p+1}|u|^{p+1}. 
\]
To define the inward/outward energy, we first define 
\begin{align*}
 \mathbf{L} u  & = \frac{x}{|x|}\cdot \nabla u + \frac{d-1}{2}\cdot \frac{u}{|x|}; \\
 \mathbf{L}_\pm u & = \frac{x}{|x|}\cdot \nabla u + \frac{d-1}{2}\cdot \frac{u}{|x|} \pm u_t.
\end{align*}
For initial data $(u_0,u_1)$ we also define 
\[
 \mathbf{L}_\pm (u_0,u_1) = \frac{x}{|x|}\cdot \nabla u_0 + \frac{d-1}{2}\cdot \frac{u_0}{|x|} \pm u_1.
\]
By Hardy's inequality, we have $\mathbf{L} u, \mathbf{L}_\pm u \in \mathcal{C}(\Rm; L^2(\Rm^d))$ as long as the solution $u$ satisfies $(u,u_t)\in \mathcal{C}(\Rm; \dot{H}^1(\Rm^d)\times L^2(\Rm^d))$. About the norm of $\mathbf{L} u$, we have 
\begin{lemma}[see \cite{shen3dnonradial}] \label{lemma L}
Assume $d\geq 3$ and $u \in \dot{H}^1(\Rm^d)$. Then 
\begin{align*}
 &\int_{\mathbb R^d} \left(\left|\mathbf{L} u\right|^2+ \lambda \cdot \frac{|u|^2}{|x|^2}\right) {\rm d} x = \int_{\Bbb R^n} |u_r|^2 {\rm} x; & & \lambda = \frac{(d-1)(d-3)}{4}.
\end{align*}
More generally, we have ($0<a<b<+\infty$)
 \begin{align*}
  \int_{a<|x|<b}  & \left(\left|\mathbf{L} u\right|^2 + \lambda \frac{|u|^2}{|x|^2}\right) {\rm d} x \\
   & = \int_{a<|x|<b} |u_r|^2 {\rm d} x + \frac{n-1}{2b} \int_{|x|=b} |u|^2 {\rm d}\sigma (x) - \frac{n-1}{2a} \int_{|x|=a} |u|^2 {\rm d}\sigma (x). 
\end{align*}
The surface integral over the sphere $|x|=a$ or $|x|=b$ can be ignored if $a=0$ or $b=+\infty$. 
\end{lemma}
Now we are able to define the inward/outward energy density function 
\[
 e_\mp (x,t) = \frac{1}{4}|\mathbf{L}_\pm u(x,t)|^2 + \frac{1}{2} e'(x,t).
\]
Here $e'(x,t)$ is the non-directional energy density defined by
\[
 e'(x,t)  = \frac{1}{2} |\slashed{\nabla} u|^2 + \frac{\lambda}{2} \frac{|u|^2}{|x|^2} + \frac{1}{2}\frac{|u|^2}{|x|} + \frac{\zeta}{p+1}|u|^{p+1};
\]
and $\slashed{\nabla} u$ is non-radial derivative defined by 
\[
 \slashed{\nabla} u = \nabla u - \left(\frac{x}{|x|}\cdot \nabla u\right)\frac{x}{|x|}.
\]
Throughout this chapter, the letter $\lambda$ represents the same constant given in Lemma \ref{lemma L}. 

\paragraph{Inward/outward energy} Given a time $t$,  we define the outward/inward energy at time $t$ by 
\[
 E_\pm (t) = \int_{\Rm^d} e_\pm (x,t) {\rm d}x. 
\]
If $\Sigma\subset \Rm^d$ is a radially symmetric region, then we may also define the outward/inward energy contained in $\Sigma$ at time $t$ to be 
\[
 E_\pm (\Sigma;t) = \int_\Sigma e_\pm (x,t) {\rm d}x. 
\]
In particular we use the notations 
\begin{align*}
 E_\pm (r_1,r_2; t) = E_\pm (\{x: r_1<|x|<r_2\};t) = \int_{r_1<|x|<r_2} e_\pm(x,t) {\rm d} x. 
\end{align*}
Please note that 
\[
 e_+(x,t) + e_-(x,t) = \frac{1}{2}|\mathbf{L} u|^2 + \frac{1}{2} |u_t|^2 + \frac{1}{2} |\slashed{\nabla} u|^2 + \frac{\lambda}{2} \frac{|u|^2}{|x|^2} + \frac{1}{2}\frac{|u|^2}{|x|} + \frac{\zeta}{p+1}|u|^{p+1} \neq e(x,t). 
\]
Thus in general, we have 
\[
 E_-(\Sigma;t) + E_+(\Sigma; t) \neq \int_{\Sigma} e(x,t) {\rm d} x.
\]
But in view of Lemma \ref{lemma L} we still have 
\begin{align}
 E_-(0,R;t) + E_+(0,R;t) & \geq \int_{|x|<R} e(x,t) {\rm d} x; \label{inside inward outward total}\\
 E_-(R,+\infty;t) + E_+(R,\infty;t) & \leq \int_{|x|>R} e(x,t) {\rm d}x; \label{outside inward outward total}\\
 E_-(t) + E_+(t) & = E. 
\end{align}
Finally we define the Morawetz density function 
\[
 M(x,t) = \frac{1}{2} \cdot \frac{|\slashed{\nabla} u(x,t)|^2}{|x|} + \frac{\lambda}{2} \cdot \frac{|u(x,t)|^2}{|x|^3} + \frac{\zeta (d-1)(p-1)}{4(p+1)}\frac{|u|^{p+1}}{|x|} + \frac{|u|^2}{4|x|^2}.
\]
\begin{remark}
 The inward/outward energy theory for the defocusing wave equation with power nonlinearity in dimension $3$ 
 \[
  \partial_t^2 u - \Delta u = -|u|^{p-1} u, \qquad (x,t)\in \Rm^3\times \Rm
 \]
 has been discussed by third author in \cite{shenenergy, shen3dnonradial}. A generalized version of this theory for the wave equation 
 \[ 
   \partial_t^2 u - \Delta u + f(x,u) =0
 \]
 for all dimensions $d\geq 3$ has been discussed in Chapter 7 of \cite{wavebook} if the nonlinear term $f(x,u)$ satisfies suitable conditions. Roughly speaking, the function $f(x,u)$ there is assumed to satisfy a repulsive condition and grow at a rate between those of $|u|^{1+\frac{4}{d-1}}$ and $|u|^{1+\frac{4}{d-2}}$, which does not cover the situation $f(x,u)=u/|x|$. It turns out that the same general idea and a similar argument works for the Coulomb wave equation as well, with insignificant modifications. We still give details about this theory here not only for completeness of our theory but also for convenience of use in subsequent chapters. 
\end{remark}
 
\begin{remark}
 High-order weighted Morawetz estimates for free Coulomb waves given in this chapter are new results and not available for nonlinear Coulomb or classic wave equations. The $L^2$-level conservation law plays an important role in this argument. 
\end{remark} 

Before we start our discussion on the inward/outward energy theory, we give a few preliminary results. 

\begin{remark} \label{finite speed of propagation}
By applying divergence theorem on the energy identity 
\[
  \partial_t e(x,t) + \nabla \cdot (-u_t \nabla u) = u_t {\rm Eq} (u)
 \]
in the region $\{(x,t): |x|<t+R; 0<t <t_0\}$, we obtain 
\[
 \int_{|x|<R+t_0} e(x,t_0) {\rm d} x = \int_{|x|<R} e(x,0) {\rm d} x + \frac{1}{\sqrt{2}} \int_{C} \left(\left|\nabla u + \frac{x}{|x|} u_t\right|^2+ \frac{|u|^2}{2|x|} + \frac{\zeta |u|^{p+1}}{p+1}\right) {\rm d} x. 
\]
Here $C=\{(x,t): |x|=t+R, 0\leq t\leq t_0\}$ is a cone surface. It immediately follows 
\[
  \int_{|x|<R+t_0} e(x,t_0) {\rm d} x \geq \int_{|x|<R} e(x,0) {\rm d} x. 
\]
By the energy conservation law, we also have
\[
  \int_{|x|>R+t_0} e(x,t_0) {\rm d} x \leq \int_{|x|>R} e(x,0) {\rm d} x. 
\]
This property is usually called the finite speed of energy propagation. 
\end{remark}

\begin{remark} \label{weighted energy plus negative}
 Given any increasing absolutely continuous function $a: [0,+\infty) \rightarrow [0,+\infty)$, we always have 
 \[ 
  \int_{\Rm^d} a(|x|) [e_+(x,0)+e_-(x,0)] {\rm d} x \leq \int_{\Rm^d} a(|x|) e(x,0) {\rm d} x. 
 \]
 Without loss of generality we assume $a(0)=0$. Indeed, we may split $a(r)$ into $a_1(r) = a(r)-a(0)$ and $a_2(r) = a(0)$, then prove that the inequality holds for each $a_j$ with $j=1,2$.  The first function $a_1(r)$ satisfies the same assumption given above and $a_1(0)=0$. The function $a_2(r)$ is a constant, the inequality above clearly holds as both sides are $a(0)$ times of the energy. If $a(0)=0$, then we have
 \begin{align*}
   \int_{\Rm^d} a(|x|) [e_-(x,0)+e_+(x,0)] {\rm d} x  &= \int_0^\infty \left(a'(r) \int_{|x|>r} [e_-(x,0)+e_+(x,0)] {\rm d} x\right) {\rm d} r \\
   & \leq \int_0^\infty \left(a'(r) \int_{|x|>r} e(x,0) {\rm d} x\right) {\rm d} r \\
   & \leq  \int_{\Rm^d} a(|x|) e(x,0) {\rm d} x.
  \end{align*}
  Here we use the inequality \eqref{outside inward outward total}. 
\end{remark}

\section{Energy flux formula}

In this section we introduce the energy flux formula for inward/outward energies. This formula is the major tool to develop the inward/outward energy theory. 

\paragraph{Regions} We first make a few assumptions on the regions involved for convenience. Let $\Omega \subset \mathbb R^d\times \mathbb R$ be a radially symmetric region so that it can be expressed by spherical coordinates $(r,\theta)$ in the following way:
\[
 \Omega = \{(r\theta, t)\in \Bbb R^n \times \Bbb R: (r,t)\in \Phi, \theta \in \mathbb{S}^{n-1}\}.
\]
Here $\Phi \subset [0,\infty) \times {\mathbb R}_t$ is a bounded, closed region, whose boundary is a simple curve consisting of finite line segments paralleled to either $t$-axis, $r$-axis or the characteristic lines $t \pm r=0$. As a result, the boundary surface $\partial \Omega$ of $\Omega$ consists of finite pieces of annulus, circular cylinders or light cones. If the boundary $\partial \Phi$ contains a line segment of the $t$-axis, then the boundary surface $\partial \Omega$ contains a degenerate part, i.e. the same line segment of $t$-axis. 

\begin{proposition}[Inward/outward energy flux] \label{energy flux formula}
 Let $d\geq 3$ and define two vector fields in $\mathbb R^n \times \mathbb R$:
 \begin{align*}
  \mathbf{V}_- = & \left[ -\frac{1}{4} \left|\mathbf{L}_+ u\right|^2 + \frac{1}{2} e'(x,t)\right]\frac{\vec{x}}{|x|} + \left[\frac{1}{4}\left|\mathbf{L}_+ u\right|^2 + \frac{1}{2}e'(x,t)\right]\vec{e};& &{\rm (Inward)}\\
  \mathbf{V}_+ = & \left[ +\frac{1}{4} \left|\mathbf{L}_- u\right|^2 - \frac{1}{2} e'(x,t)\right]\frac{\vec{x}}{|x|} + \left[\frac{1}{4}\left|\mathbf{L}_- u\right|^2 + \frac{1}{2}e'(x,t)\right]\vec{e}.& &{\rm (Outward)}
 \end{align*}
Here $\vec{x}$ and $\vec{e}$ are both vectors in $\mathbb R^d \times \mathbb R$.
 \begin{align*}
  &\vec{x}=(x,0);& &\vec{e}=(0,0,\cdots,0,1).
 \end{align*}
If $\Omega = \{(r\theta, t)\in \Bbb R^d \times \Bbb R: (r,t) \in \Phi, \theta \in \mathbb{S}^{d-1}\}$ is a region as described above and $u$ is a solution to \eqref{unified equation inout} with a finite energy, then we have
 \[
  \int_{\partial \Omega} \mathbf{V}_{\pm} \cdot {\rm d}\mathbf{S} = \pm  \iint_\Omega M(x,t) {\rm d} x {\rm d} t.
 \]
 If $d=3$ and a line segment $[t_1,t_2]$ of $t$-axis is part of the boundary $\partial \Phi$, i.e. the boundary surface $\partial \Omega$ contains a line segment of $t$-axis, we understand the surface integral over this degenerate part of $\partial \Omega$ as
 \[
  \mp \pi \int_{t_1}^{t_2} |u(0,t)|^2 {\rm d} t \quad \hbox{(plus for inward, minus for outward)}
 \]
 In higher dimensions $d\geq 4$, the surface integral over degenerate part of $\partial \Omega$ is simply ignored.
\end{proposition}

\begin{figure}[h]
 \centering
 \includegraphics[scale=1.1]{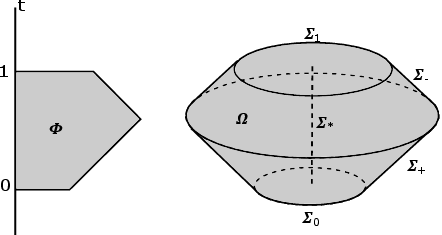}
 \caption{energy flux example} \label{figure flux}
\end{figure}

Before we give the proof of Proposition \ref{energy flux formula}, we first demonstrate why this is called energy flux formula. We consider the 3-dimensional case. Let $\Phi \subset \Rm_r^+ \times \Rm_t$ be a pentagon as illustrated in left part of figure \ref{figure flux}.  Then the boundary surface of the corresponding region $\Omega \subset \Rm^d \times \Rm$ consists of five parts: two disks $\Sigma_1$ and $\Sigma_0$, two truncated cones $\Sigma_+$ and $\Sigma_-$, as well as a line segment on the $t$-axis. We apply Proposition \ref{energy flux formula} for the inward energy on this region and obtain 
\begin{align*}
 \int_{\Sigma_1} e_-(x,t) {\rm d} S + \pi \int_0^1 |u(0,t)|^2 {\rm d} t - \int_{\Sigma_0} e_-(x,t) {\rm d} S - \frac{1}{2\sqrt{2}} \int_{\Sigma_+} |\mathbf{L}_+ u|^2 {\rm d} S & +  \frac{1}{\sqrt{2}} \int_{\Sigma_-} e'(x,t) {\rm d} S \\
 & = -\iint_{\Omega} M(x,t) {\rm d} x {\rm d}t. 
\end{align*}
We may rewrite it in the form of 
\begin{align*}
 \int_{\Sigma_1} & e_-(x,t) {\rm d} S  - \int_{\Sigma_0} e_-(x,t) {\rm d} S \\
& =  -\pi \int_0^1 |u(0,t)|^2 {\rm d} t + \frac{1}{2\sqrt{2}} \int_{\Sigma_+} |\mathbf{L}_+ u|^2 {\rm d} S  -  \frac{1}{\sqrt{2}} \int_{\Sigma_-} e'(x,t) {\rm d} S  -\iint_{\Omega} M(x,t) {\rm d} x {\rm d}t. 
\end{align*}
The left hand side is the difference of the inward energy contained in the region $\Sigma_1$ and $\Sigma_0$. This difference is a sum of four terms. The first term is a loss, which is the amount of energy carried by inward waves which move through the origin and transform to outward waves. The second term is a gain, which is the amount of energy carried by inward waves which move through the surface $\Sigma_+$. The third term is negative, which is the amount of inward energy loss along the surface $\Sigma_-$ due to both linear and nonlinear effect. Finally the last term is also a loss, which is exactly the amount of inward energy transforming to outward energy in the region $\Omega$. 

Now we give the proof of Proposition \ref{energy flux formula} by an application of the divergence theorem. We first prove 
\begin{lemma} \label{lemma energy flux}
 If we introduce the notation 
 \begin{align*}
  {\rm Eq}(u) = \partial_t^2 u - \Delta u + \frac{u}{|x|} + \zeta |u|^{p-1} u,
 \end{align*}
 then 
 \begin{align*}
  \partial_t \left[\frac{1}{2}|\mathbf{L}_- u|^2 + e' \right] + \nabla \cdot \left[+\slashed{\nabla} u \mathbf{L}_- u + \frac{x}{|x|}\left(+\frac{1}{2}|\mathbf{L}_- u|^2 -e' \right)\right] & = -  {\rm Eq} (u) \mathbf{L}_- u + 2M(x,t);\\
  \partial_t \left[\frac{1}{2}|\mathbf{L}_+ u|^2 + e' \right] + \nabla \cdot \left[-\slashed{\nabla} u \mathbf{L}_+ u + \frac{x}{|x|}\left(-\frac{1}{2}|\mathbf{L}_+ u|^2 +e' \right)\right] & = +{\rm Eq} (u) \mathbf{L}_+ u  -2M(x,t). 
 \end{align*}
\end{lemma}
\begin{proof}
 Taking the sum and difference of the two identities gives above, we obtain 
 \begin{align}
 & \partial_t \left[2e+\frac{(d-1)(d-2)|u|^2}{2|x|^2}+(d-1)\frac{x\cdot u\nabla u}{|x|^2}\right] + \nabla  \cdot \left[-2u_t \slashed{\nabla} u -2 \frac{x}{|x|} u_t \mathbf{L} u\right] = 2 u_t {\rm Eq} (u)  \label{energy divergence}\\
  &\partial_t \left[-2 u_t \mathbf{L} u\right] + \nabla\cdot \left[2\slashed{\nabla} u \mathbf{L} u + \frac{x}{|x|} \left(|\mathbf{L} u|^2 + |u_t|^2 - 2 e'\right)\right]  = -2{\rm Eq}(u) \mathbf{L} u + 4M. \label{L divergence}
 \end{align}
 It suffices to prove these two identities. In view of the identity 
 \[
  -2u_t \slashed{\nabla} u -2 \frac{x}{|x|} u_t \mathbf{L} u = -2 u_t \nabla u - \frac{(d-1)x}{|x|^2} u u_t, 
 \]
 the first identity \eqref{energy divergence} is simply the combination of the energy identity 
 \[
  \partial_t e(x,t) + \nabla \cdot (-u_t \nabla u) = u_t {\rm Eq} (u)
 \]
 and the identity 
 \[
  \partial_t \left[\frac{(d-1)(d-2)|u|^2}{2|x|^2}+(d-1)\frac{x\cdot u\nabla u}{|x|^2}\right] + \nabla  \cdot \left[- \frac{(d-1)x}{|x|^2} u u_t\right] = 0.
 \]
 Now we prove the second identity \eqref{L divergence}. We split the left hand side of \eqref{L divergence} into three parts $I_1+I_2+I_3$ with 
 \begin{align*}
  I_1 & = \partial_t [-2u_t \mathbf{L} u] + \nabla \cdot \frac{x}{|x|} |u_t|^2;\\
  I_2 & = \nabla \cdot \left[2\slashed{\nabla} u \mathbf{L} u + \frac{x}{|x|} \left(|\mathbf{L} u|^2 - |\slashed{\nabla} u|^2 - \lambda \frac{|u|^2}{|x|^2}\right)\right];\\
  I_3 & = \nabla \cdot \frac{x}{|x|} \left( -\frac{|u|^2}{|x|} - \frac{2\zeta}{p+1} |u|^{p+1} \right);
 \end{align*}
 and then calculate each term individually:
 \begin{align}
  I_1 & = - 2 u_{tt} \mathbf{L} u -2 u_t \partial_t \left(\frac{x}{|x|}\cdot \nabla u\right) -(d-1) \frac{|u_t|^2}{|x|} + \nabla \cdot \frac{x|u_t|^2}{|x|}\nonumber \\
  & = - 2 u_{tt} \mathbf{L} u. \label{L divergence part 1}
 \end{align}
 \begin{align*}
  I_2 & =  \nabla \cdot \left[2\nabla  u \mathbf{L} u + \frac{x}{|x|}\left( \left(\mathbf{L} u- 2u_r \right)\mathbf{L} u - |\slashed{\nabla} u|^2 - \lambda \frac{|u|^2}{|x|^2}\right)\right] \\
  & = 2 \Delta u \mathbf{L} u + 2 \nabla u \cdot \nabla \mathbf{L} u + \nabla \cdot \frac{x}{|x|} \left(\frac{(d-1) |u|^2}{2|x|^2} - |\nabla u|^2\right)\\
  & = 2 \Delta u \mathbf{L} u + 2\sum_{j} u_j \partial_j \left(\sum_{k} \frac{x_k}{|x|} u_k + \frac{(d-1) u}{2|x|}\right) + \nabla \cdot \frac{(d-1)x}{2|x|^3} |u|^2 - \sum_{j,k} \partial_j \left(\frac{x_j}{|x|} u_k^2\right)\\
  & = 2 \Delta u \mathbf{L} u + 2 \sum_{j,k} \left(\frac{\delta_{jk}}{|x|} u_j u_k - \frac{x_j x_k}{|x|^3} u_j u_k + \frac{x_k}{|x|}u_j u_{jk}\right) + (d-1)\frac{|\nabla u|^2}{|x|} - \nabla u \cdot \frac{(d-1)x}{|x|^3} u\\
  & \qquad + \frac{(d-1)x}{|x|^3} \cdot (u\nabla u) + 2\lambda \frac{|u|^2}{|x|^3} -  (d-1)\frac{|\nabla u|^2}{|x|} - 2\sum_{j,k} \frac{x_j}{|x|} u_k u_{jk}\\
  & = 2 \Delta u \mathbf{L} u + \frac{2|\slashed{\nabla} u|^2}{|x|} + 2\lambda \frac{|u|^2}{|x|^3}.
 \end{align*}
 \begin{align*}
  I_3 & = - \nabla \cdot \left(\frac{x}{|x|^2}|u|^2\right) - \frac{2\zeta}{p+1} \nabla \cdot \left(\frac{x}{|x|} |u|^{p+1}\right)\\
  & = -\frac{d-2}{|x|^2}|u|^2 - 2 \frac{u}{|x|} \left(\frac{x}{|x|}\cdot \nabla u\right) - \frac{2\zeta(d-1)}{p+1} \frac{|u|^{p+1}}{|x|} - 2\zeta |u|^{p-1} u \left(\frac{x}{|x|}\cdot \nabla u\right)\\
  & = -2 \left(\frac{u}{|x|} + \zeta |u|^{p-1} u\right)\left(\frac{x}{|x|}\cdot \nabla u + \frac{(d-1)u}{2|x|}\right) + \frac{|u|^2}{|x|^2} + \frac{\zeta(d-1)(p-1)}{p+1} \frac{|u|^{p+1}}{|x|}\\
  & = -2 \left(\frac{u}{|x|} + \zeta |u|^{p-1} u\right) \mathbf{L} u  + \frac{|u|^2}{|x|^2} + \frac{\zeta(d-1)(p-1)}{p+1} \frac{|u|^{p+1}}{|x|}.
 \end{align*}
 We collect all these terms and conclude 
 \begin{align*}
  I_1 + I_2 + I_3 & = -2\left(u_{tt}-\Delta u + \frac{u}{|x|} + \zeta |u|^{p-1} u\right) \mathbf{L} u +  \frac{2|\slashed{\nabla} u|^2}{|x|} + 2\lambda \frac{|u|^2}{|x|^3}\\& \qquad \quad + \frac{|u|^2}{|x|^2} + \frac{\zeta(d-1)(p-1)}{p+1} \frac{|u|^{p+1}}{|x|}\\
  & = -2 {\rm Eq}(u) \mathbf{L} u + 4 M(x,t). 
 \end{align*}
\end{proof}

\begin{proof}[Proof of Proposition \ref{energy flux formula}]
 We consider the energy flux of the outward energy as an example. The proof of the inward energy is similar. We calculate as though $u$ is sufficiently smooth. Strictly speaking, we need to apply smooth approximation techniques here, as we did for the proof of energy conservation law. More details about these techniques can be found in Miao-Shen \cite{wavebook}. Let $u$ be a finite-energy solution to \eqref{unified wave equation} and $\Omega$ be a radially symmetric region as described above. We also temporarily assume that $\Omega$ is away from the $t$-axis.  We apply the divergence theorem on the first divergence identity given in Lemma \ref{lemma energy flux} in the region $\Omega$, insert the equation ${\rm Eq}(u)=0$ and obtain 
 \begin{align*}
  2  \int_{\partial \Omega} \mathbf{V}_{+} \cdot {\rm d}\mathbf{S} + \int_{\partial \Omega} \mathbf{V}_\ast  \cdot {\rm d}\mathbf{S}  = 2 \int_{\Omega} M(x,t) {\rm d} x {\rm d} t. 
 \end{align*}
 Here the vector field $\mathbf{V}_\ast$ is defined by 
 \[
  \mathbf{V}_\ast = (\slashed{\nabla} u \mathbf{L}_- u, 0).
 \]
 We observe that $\mathbf{V}_\ast$ is orthogonal to the outer normal vector of the surface $\partial \Omega$ as long as $\Omega$ is a radially symmetric region, since we have 
 \[
  (\slashed{\nabla} u \mathbf{L}_- u ) \cdot \frac{x}{|x|} = 0. 
 \]
 As a result, the surface integral of $\mathbf{V}_\ast$ is zero. This verifies the energy flux formula when $\Omega$ is away from the $t$-axis. If part of the boundary $\partial \Omega$ is a line segment $[t_1,t_2]$ on the $t$-axis, we may first apply the energy flux formula on the region 
 \[
  \Omega_r = \Omega \cap \{(x,t): |x|\geq r\},
 \]
 and then let $r\rightarrow 0^+$. In this limiting procedure we need to consider the surface integral over the cylinder 
 \[
  \Sigma_r = \{(x,t): |x|=r, t_1 + \kappa_1 r \leq t \leq t_2 + \kappa_2 r\}, \qquad \kappa_1, \kappa_2 \in \{-1,0,1\}.
 \]
 We have 
 \[ 
  \int_{\Sigma_r} \mathbf{V}_+ \cdot {\rm d} \mathbf{S} = -\int_{t_1 + \kappa_1 r}^{t_2 + \kappa_2 r} \int_{|x|=r} \left(\frac{1}{4} \left|\mathbf{L}_- u\right|^2 - \frac{1}{2} e'(x,t)\right) {\rm d} \sigma {\rm d} t.
 \]
 We recall that 
 \begin{align*}
  &\mathbf{L}_- u = \frac{x}{|x|}\cdot \nabla u  + \frac{(d-1)u}{2|x|} - u_t;& &e' = \frac{1}{2} |\slashed{\nabla} u|^2 + \frac{\lambda}{2} \frac{|u|^2}{|x|^2} + \frac{1}{2}\frac{|u|^2}{|x|} + \frac{\zeta}{p+1}|u|^{p+1}.
 \end{align*}
 It is easily seen that if $u$ is sufficiently smooth, then 
 \begin{align*}
  &\lim_{r\rightarrow 0^+}\int_{\Sigma_r} \mathbf{V}_+ \cdot {\rm d} \mathbf{S} = - \pi \int_{t_1}^{t_2} |u(0,t)|^2 {\rm d} t,& & d=3; \\
  &\lim_{r\rightarrow 0^+}\int_{\Sigma_r} \mathbf{V}_+ \cdot {\rm d} \mathbf{S} = 0,& & d\geq 4.
 \end{align*}
 This verifies the energy flux formula when the boundary of $\Omega$ contains a degenerate part on the $t$-axis. 
\end{proof}

\begin{remark}
 The following Morawetz-Strichartz estimates are very useful. Let $u$ solve the linear classic wave equation 
 \[
  \left\{\begin{array}{l} \partial_t^2 u - \Delta u = g; \\ (u(0),u_t(0)) = (u_0,u_1). \end{array} \right. 
 \]
 Then we have
  \begin{align*}
  \big\||x|^{-1/2}\slashed{\nabla} u\big\|_{L^2(\Bbb R^{1+d})}+\big\||x|^{-3/2} u\big\|_{L^2(\Bbb R^{1+d})} &\leq C\|(u_0,u_1)\|_{\dot{H}^1\times L^2} + C\|g\|_{L^1 L^2(\Bbb R\times \Bbb R^d)}, & & d\geq 4;\\
  \big\||x|^{-1/2}\slashed{\nabla} u\big\|_{L^2(\Bbb R^{1+3})}+ \|u(0,t)\|_{L^2(\Bbb R)} &\leq C\|(u_0,u_1)\|_{\dot{H}^1\times L^2} + C\|g\|_{L^1 L^2(\Bbb R\times \Bbb R^3)}, & & d= 3.
 \end{align*}
 In addition, we also have ($2^*=2d/(d-2)$)
 \[
  \big\|u\big\|_{L^{2^*}(\Bbb R^{1+d};|x|^{-1} {\rm d}x {\rm d} t)} \leq C\|(u_0,u_1)\|_{\dot{H}^1\times L^2} + C\|g\|_{L^1 L^2(\Bbb R\times \Bbb R^d)}.
 \]
 These imply that $M(x,t)$ are always locally integrable since
 \begin{itemize}
  \item The local theory guarantees that $|u|^{p-1}u \in L^1 L^2 (I\times \Rm^d)$ for any finite time interval;
  \item The Hardy inequality means that $u/|x| \in L^2(\Rm^d)$ for any time $t$.
 \end{itemize} 
 The Morawetz-Strichartz estimates can be seen in the following way. A similar calculation as in Lemma \ref{lemma energy flux} shows that any solution to the free classic wave equation satisfies 
 \begin{align*}
  \partial_t \left[-2 u_t \mathbf{L} u\right] + \nabla\cdot \left[2\slashed{\nabla} u \mathbf{L} u + \frac{x}{|x|} \left(|\mathbf{L} u|^2 + |u_t|^2 - |\slashed{\nabla} u|^2 - \lambda \frac{|u|^2}{|x|^2}\right)\right]  = \frac{2|\slashed{\nabla} u|^2}{|x|} + 2\lambda \frac{|u|^2}{|x|^3}. 
 \end{align*}
We may apply the divergence theorem on the divergence identity above in a region $[t_1,t_2]\times \{x: |x|>r\}$ and then let $r\rightarrow 0^+$, $t_1\rightarrow -\infty$, $t_2\rightarrow +\infty$ to deduce the Morawetz-Strichartz estimates for the free wave equation with initial data $(u_0,u_1) \in C_0^\infty (\Rm^d)$. Since $C_0^\infty$ is dense in $\dot{H}^1\times L^2$, these inequalities also hold for general free waves. We observe that the solution to the inhomogeneous wave equation $\partial_t^2 u - \Delta u = g$ satisfies 
\begin{equation} \label{Duhamel}
 u(t) = \mathbf{S}_{\mathcal{W}}(t) (u_0,u_1) + \int_0^\infty \chi(s,t) \mathbf{S}_{\mathcal{W}}(t-s) (0, g(s)) {\rm d} s 
\end{equation}
Here $\chi(s,t)$ is the characteristic function of the region $\{(s,t): t>s\}$. The Strichartz estimate of the inhomogeneous part then follows \eqref{Duhamel} and the corresponding estimate for the homogeneous part. Finally we may combine the Hardy inequality and the Strichartz estimates to deduce ($3\leq d \leq 6$)
\begin{align*}
 \iint \frac{|u(x,t)|^{2^*}}{|x|} {\rm d} x {\rm d} t 
 & \lesssim \|u\|_{L^\infty(\Bbb R; \dot{H}^1(\Bbb R^n))} \|u\|_{L^\frac{n+2}{n-2}L^\frac{2(n+2)}{n-2}(\Bbb R\times \Bbb R^n)}^{2^*-1}\\
 & \lesssim \left(C\|(u_0,u_1)\|_{\dot{H}^1\times L^2} + C\|g\|_{L^1 L^2(\Bbb R\times \Bbb R^n)}\right)^{2^*}.
\end{align*}
\end{remark}

\paragraph{Notations} For convenience of future use, we introduce a few notations. Given $s,\tau \in \Bbb R$, we use the following notations for backward($-$) and forward($+$) light cones
\begin{align*}
 &C^-(s) = \{(x,t): |x|+t = s\};& &C^+(\tau) = \{(x,t): t-|x|=\tau\}.
\end{align*}
We also need to consider their truncated version, i.e. the cones between two given times $t_1, t_2$.
\begin{align*}
 &C^-(s;t_1,t_2) = \{(x,t): |x|+t = s, t_1\leq t\leq t_2\};\\
  &C^+(\tau; t_1,t_2) = \{(x,t): t-|x|=\tau, t_1\leq t\leq t_2\}.
\end{align*}
For convenience we introduce the following notations for the region inside the cones $C^-(s;t_1,t_2)$ or $C^+(\tau,t_1,t_2)$.
\begin{align*}
 &K^-(s;t_1,t_2) = \{(x,t): |x|+t<s, t_1\leq t\leq t_2\};\\
  &K^+(\tau; t_1,t_2) = \{(x,t): t-|x|>\tau, t_1\leq t\leq t_2\}.
\end{align*}
We also consider the energy flux through the truncated cones:
\begin{align*}
 Q_-^- (s; t_1,t_2) & =  \frac{1}{\sqrt{2}} \int_{C^-(s;t_1,t_2)} e'(x,t) {\rm d} S,& &t_1<t_2\leq s;\\
 Q_+^-(s; t_1,t_2) & = \frac{1}{2\sqrt{2}} \int_{C^-(s;t_1,t_2)} |\mathbf{L}_- u|^2 {\rm d} S, & & t_1<t_2\leq s;\\
 Q_-^+(\tau; t_1,t_2) & = \frac{1}{2\sqrt{2}} \int_{C^+(\tau;t_1,t_2)} |\mathbf{L}_+ u|^2 {\rm d} S,& &\tau\leq t_1<t_2;\\
 Q_+^+(\tau; t_1,t_2) & = \frac{1}{\sqrt{2}} \int_{C^+(\tau;t_1,t_2)} e'(x,t) {\rm d} S,& &\tau\leq t_1<t_2.
\end{align*}
The upper index tells us whether this is energy flux through backward light cone($-$), or forward light cone($+$). The lower index indicates whether this is energy flux of inward energy ($-$) or outward energy ($+$). 

\paragraph{The cone law} Before we introduce the basic theory of inward/outward energy, we give a simple corollary of the energy flux formula. 

\begin{lemma}[The cone law]
 Let $u$ be a global solution to \eqref{unified equation inout} with a finite energy. Assume that $t_0 \leq s$. Then 
 \begin{align*}
  E_-(0,s-t_0;t_0) &= \iint_{K^-(s;t_0,s)} M(x,t) {\rm d} x {\rm d} t + Q_-^-(s;t_0,s), & &d\geq 4;\\
  E_-(0,s-t_0;t_0) &= \iint_{K^-(s;t_0,s)} M(x,t) {\rm d} x {\rm d} t + Q_-^-(s;t_0,s) + \pi \int_{t_0}^s |u(0,t)|^2 {\rm d} t, & &d=3.
 \end{align*}
\end{lemma}
\begin{proof}
 We simply apply the energy flux formula of the inward energy in the cone region $K^-(s;t_0,s)$. 
\end{proof}

\begin{remark}
 The cone law also holds for outward energy and forward light cones. A similar argument shows that if $\tau\leq t_0$, then 
  \begin{align*}
  E_+(0,t_0-\tau;t_0) &= \iint_{K^+(\tau;\tau,t_0)} M(x,t) {\rm d} x {\rm d} t + Q_+^+(\tau;\tau,t_0), & &d\geq 4;\\
  E_+(0,t_0-\tau;t_0) &= \iint_{K^+(\tau;\tau,t_0)} M(x,t) {\rm d} x {\rm d} t + Q_+^+(\tau;\tau,t_0) + \pi \int_{\tau}^{t_0} |u(0,t)|^2 {\rm d} t, & &d=3.
 \end{align*}
\end{remark}

\section{Theory of inward/outward energy}

In this section we always assume that $u$ is a global solution to the Coulomb wave equation \eqref{unified equation inout} with a finite energy and introduce the basic theory of inward/outward energy. 

\begin{lemma}[Morawetz estimates] \label{global Morawetz inequality} 
Let $u$ be a global solution to \eqref{unified wave equation} with a finite energy. Then the following global integral estimates hold:
\begin{align*}
 \iint_{\mathbb R^{1+d}} M(x,t) {\rm d} x {\rm d} t & \leq E, & &d\geq 4; \\
 \iint_{\mathbb R^{1+d}} M(x,t) {\rm d} x {\rm d} t + \pi \int_{-\infty}^{\infty} |u(0,t)|^2 {\rm d} t &\leq E, & &d=3.
\end{align*}
\end{lemma}
\begin{proof}
 We observe that the left hand of each identity in the cone law is bounded by the energy $E$ and the right hand is a decreasing function of $t_0 \in (-\infty,t_0]$. We let $t_0 \rightarrow -\infty$ and obtain 
  \begin{align*}
  \iint_{|x|+t\leq s} M(x,t){\rm d} x {\rm d} t + Q_-^-(s; -\infty, s)&\leq E, & &d\geq 4;\\
  \iint_{|x|+t\leq s} M(x,t){\rm d} x {\rm d} t + Q_-^-(s; -\infty, s) + \pi \int_{-\infty}^s |u(0,t)|^2 {\rm d} t &\leq E, & &d=3.
 \end{align*}
 We then let $s\rightarrow +\infty$ to finish the proof. 
\end{proof}

Next we observe that the following inequality holds for almost everywhere $(x,t) \in \Rm^d \times \Rm$. 
\begin{equation} \label{key estimate for in out energy theory}
 M(x,t)\geq \mu \frac{e'(x,t)}{|x|}
\end{equation} 
where 
\[
 \mu = \min\left\{\frac{1}{2}, \frac{(d-1)(p-1)}{4}\right\}
\]
Please note that $\mu = 1/2$ for immediate Coulomb wave equation $1+\frac{4}{d-1} \leq p \leq 1+ \frac{4}{d-2}$. This estimate plays an essential role in the inward/outward energy theory. It immediately follows that 
\begin{lemma} \label{average cone}
For all $t_0\in \Bbb R$ and radii $0\leq r_1<r_2$, the following inequality holds
 \[
  \int_{t_0+r_1}^{t_0+r_2} Q_-^-(s; t_0, s) {\rm d} s \leq \frac{r_2}{\mu} \iint_{\Omega^-(t_0;r_1,r_2)} M(x,t) {\rm d} x {\rm d} t,
 \]
 where $\mu$ is the constant in \eqref{key estimate for in out energy theory} and $\Omega^-(t_0;t_1,t_2)$ is a cone shell defined by
 \begin{equation} 
 \Omega^-(t_0;r_1,r_2) = \{(x,t): t_0+r_1\leq |x|+t\leq t_0+r_2, t\geq t_0\} \subset B(0,r_2)\times \Bbb R.
 \end{equation}
 %as illustrated in figure \ref{figure inf3D1}.
\end{lemma}
\begin{proof}
 This inequality follows a straight forward calculation by Fubini's theorem
\begin{align*}
 \int_{t_0+r_1}^{t_0+r_2} Q_-^-(s;t_0,s) {\rm d} s & = \frac{1}{\sqrt{2}} \int_{t_0+r_1}^{t_0+r_2} \int_{C^-(s;t_0,s)} e'(x,t) {\rm d} S {\rm d} s \\
 & = \iint_{\Omega^-(t_0;r_1,r_2)} e'(x,t) {\rm d} x {\rm d} t \\
 & \leq r_2\iint_{\Omega^-(t_0;r_1,r_2)} \frac{e'(x,t)}{|x|} {\rm d} x {\rm d} t\\
 & \leq \frac{r_2}{\mu} \iint_{\Omega^-(t_0;r_1,r_2)} M(x,t) {\rm d} x {\rm d} t.
\end{align*} 
Here we use the inequality \eqref{key estimate for in out energy theory}. 
\end{proof}

\begin{corollary} \label{bigger cone law}
 Let $\kappa>1$ be a constant. Then the following identities hold for all $r>0$ and $t_0 \in \Rm$. 
 \begin{align*}
  E_-(0,r; t_0) &\leq C(\kappa,\mu) \int_{K^-(t_0+\kappa r; t_0, t_0+\kappa r)} M(x,t) {\rm d} x {\rm d} t, & & d\geq 4; \\
  E_-(0,r; t_0) &\leq C(\kappa,\mu) \int_{K^-(t_0+\kappa r; t_0, t_0+\kappa r)} M(x,t) {\rm d} x {\rm d} t + \pi \int_{t_0}^{t_0+\kappa r} |u(0,t)|^2 {\rm d} t, & & d=3. 
 \end{align*}
 Here the constant $C(\kappa,\mu)$ is given by 
 \[
  C(\kappa,\mu) = 1 + \frac{\kappa}{(\kappa-1)\mu}. 
 \]
\end{corollary}
\begin{proof}
 Let us prove the three-dimensional case. Higher dimensional case is similar. We apply the cone law on the cone 
 \[
  K^-(t_0+r'; t_0, t_0+r') = \{(x,t): |x|+t < t_0+r', t_0<t < t_0+r'\}, \qquad r' \in [r,\kappa r],
 \]
 and obtain 
 \begin{align*}
  E_-(0,r';t_0) & = \iint_{K^-(t_0+r';t_0,t_0+r')} M(x,t) {\rm d} x {\rm d} t + \pi \int_{t_0}^{t_0+r'} |u(0,t)|^2 {\rm d} t + Q_-^-(t_0+r';t_0,t_0+r')\\
  & \leq \iint_{K^-(t_0+\kappa r;t_0,t_0+\kappa r)} M(x,t) {\rm d} x {\rm d} t + \pi \int_{t_0}^{t_0+\kappa r} |u(0,t)|^2 {\rm d} t + Q_-^-(t_0+r';t_0,t_0+r')\\
  & \doteq J(t_0,\kappa r) + Q_-^-(t_0+r';t_0,t_0+r').
 \end{align*}
 We then integrate this inequality from $r' = r$ to $r'=\kappa r$, divide both sides by $(\kappa-1)r$ and apply Lemma \ref{average cone} to deduce
 \begin{align*}
  E_-(0,r;t_0) & \leq \frac{1}{(\kappa-1)r} \int_{r}^{\kappa r} E_-(0,r';t_0){\rm d} r' \\
  & \leq J(t_0,\kappa r) + \frac{1}{(\kappa-1)r} \int_{r}^{\kappa r} Q_-^-(t_0+r';t_0,t_0+r') {\rm d} r'\\
  & \leq J(t_0,\kappa r) + \frac{1}{(\kappa-1)r} \int_{t_0 + r}^{t_0 +\kappa r} Q_-^-(s;t_0,s) {\rm d} s\\
  & \leq J(t_0,\kappa r) + \frac{\kappa}{(\kappa-1)\mu} \int_{\Omega^-(t_0;r,\kappa r)} M(x,t) {\rm d} x {\rm d} t. 
 \end{align*}
 Finally we observe that $\Omega^-(t_0;r,\kappa r) \subset K^-(t_0+\kappa r; t_0, t_0+\kappa r)$ and finish the proof. 
\end{proof}

\begin{remark} \label{bigger cone law 2} 
 By time symmetry a similar argument to Corollary \ref{bigger cone law} gives 
  \begin{align*}
  E_+(0,r; t_0) &\leq C(\kappa,\mu) \int_{K^+(t_0-\kappa r; t_0-\kappa r, t_0)} M(x,t) {\rm d} x {\rm d} t, & & d\geq 4; \\
  E_+(0,r; t_0) &\leq C(\kappa,\mu) \int_{K^+(t_0-\kappa r; t_0-\kappa r, t_0)} M(x,t) {\rm d} x {\rm d} t + \pi \int_{t_0-\kappa r}^{t_0} |u(0,t)|^2 {\rm d} t, & & d=3. 
 \end{align*}
\end{remark} 
A direct consequence of Corollary \ref{bigger cone law} is 

\begin{proposition} \label{summation of energy}
Let $u$ be a global solution to \eqref{unified equation inout} with a finite energy $E$. Assume that $R\geq T>0$ are both constants. Then   
 \[
  \sum_{k=-\infty}^\infty \int_{|x|\leq R} e(x,kT) {\rm d}x  \lesssim \frac{RE}{T}.
 \]
\end{proposition}
\begin{proof}
Again we only consider the case $d=3$. Higher dimensions are similar. We apply Corollary \ref{bigger cone law} with $t_0 = kT$, $r = R$, $\kappa =2$ and obtain 
\begin{align*}
 E_-(0,R; kT) &\leq \left(1+ \frac{2}{\mu}\right) \int_{K^-(kT+2R; kT, kT+2R)} M(x,t) {\rm d} x {\rm d} t + \pi \int_{kT}^{kT+2R} |u(0,t)|^2 {\rm d} t\\
 & \leq \left(1+ \frac{2}{\mu}\right) \int_{kT}^{kT+2R} \int_{\Rm^d} M(x,t) {\rm d} x {\rm d} t + \pi \int_{kT}^{kT+2R} |u(0,t)|^2 {\rm d} t. 
\end{align*}
By Remark \ref{bigger cone law 2} we also have 
\[
 E_+(0,R;kT) \leq \left(1+ \frac{2}{\mu}\right) \int_{kT-2R}^{kT} \int_{\Rm^d} M(x,t) {\rm d} x {\rm d} t + \pi \int_{kT-2R}^{kT} |u(0,t)|^2 {\rm d} t.
\]
Combining this together we have 
\begin{align*}
 \int_{|x|\leq R} e(x,kT) {\rm d} x & \leq E_-(0,R;kT) + E_+(0,R;kT)\\
 & \leq \left(1+ \frac{2}{\mu}\right) \int_{kT-2R}^{kT+2R} \int_{\Rm^d} M(x,t) {\rm d} x {\rm d} t + \pi \int_{kT-2R}^{kT+2R} |u(0,t)|^2 {\rm d} t. 
\end{align*}
We may take a sum and obtain 
\begin{align*}
 \sum_{k=-\infty}^\infty \int_{|x|\leq R} e(x,kT) {\rm d} x & \leq N(4R/T)  \left[\left(1+ \frac{2}{\mu}\right) \iint_{\Rm^d\times \Rm} M(x,t) {\rm d} x {\rm d} t + \pi \int_{\Rm} |u(0,t)|^2 {\rm d} t\right]\\
  & \leq N(4R/T) \left(1+ \frac{2}{\mu}\right) E. 
\end{align*}
Here $N(4R/T) \lesssim R/T$ is the minimal integer greater than or equal to $4R/T$. We use the global Morawetz estimate given in Lemma \ref{global Morawetz inequality} in the final step. 
\end{proof}

\begin{lemma} \label{integral expression of energy} 
Given any time $t_0$, the following identities hold:
\begin{align*}
 E_-(t_0) & = \int_{t_0}^\infty \int_{\Bbb R^d} M(x,t) {\rm d} x {\rm d} t, & & d\geq 4;\\
 E_-(t_0) & = \int_{t_0}^\infty \int_{\Bbb R^d} M(x,t) {\rm d} x {\rm d} t + \pi \int_{t_0}^\infty |u(0,t)|^2 {\rm d} t, & & d=3.
\end{align*}
\end{lemma}
\begin{proof}
 We consider the higher dimensional case $d \geq 4$. The proof for three-dimensional case is similar. We apply the cone law on the cone 
 \[
  K^-(s; t_0, s) = \{(x,t): |x|+t < s, t_0<t < s\}, \qquad s>t_0
 \]
 and obtain
 \begin{equation} \label{cone law identity for expression}
 E_-(0,s-t_0;t_0) = \iint_{K^-(s;t_0,s)} M(x,t) {\rm d} x {\rm d} t + Q_-^-(s;t_0,s);
 \end{equation}
 then let $s\rightarrow +\infty$. By the global integrability of $e_-(x,t_0)$ in $\Rm^d$ and the global integrability of $M(x,t)$ in $\Rm^d \times \Rm$, it is clear that 
 \begin{align*}
  \lim_{s \rightarrow +\infty} E_-(0,s-t_0;t_0) & = E_-(t_0); \\
  \lim_{s \rightarrow +\infty}\iint_{K^-(s;t_0,s)} M(x,t) {\rm d} x {\rm d} t & = \int_{t_0}^\infty \int_{\Bbb R^n} M(x,t) {\rm d} x {\rm d} t. 
 \end{align*} 
 In order to complete the proof, it suffice to show 
 \begin{equation} \label{cone surface limit}
  \lim_{s\rightarrow +\infty} Q_-^-(s;t_0,s) = 0. 
 \end{equation}
 Since the other two terms in \eqref{cone law identity for expression} both converge to a finite number, $Q_-^-(s;t_0,s)$ has to converge to a finite number as $s\rightarrow +\infty$. To figure out the limit of $Q_-^-(s;t_0,s)$, we consider its limit in the average sense. Indeed, by Lemma \ref{average cone}, we have 
 \[
 \frac{1}{r}\int_{t_0+r}^{t_0+2r} Q_-^-(s;t_0,s) {\rm d} s \leq \frac{2}{\mu}\iint_{\Omega^-(t_0;r,2r)} M(x,t) {\rm d} x {\rm d} t \rightarrow 0, \quad {\rm as}\; r\rightarrow +\infty.
\]
This implies that the limit of $Q_-^-(s;t_0,s)$ as $s\rightarrow 0$ must be zero and finishes the proof. 
\end{proof}

\begin{corollary} \label{Morawetz estimate cor}
 The inward energy $E_-(t)$ is a decreasing function of time $t$ satisfying
 \begin{align*}
  &\lim_{t\rightarrow +\infty} E_-(t) = 0;& &\lim_{t\rightarrow -\infty} E_-(t) = E.
 \end{align*}
 In addition, we have the Morawetz identity 
 \begin{align*}
 \iint_{\mathbb R^{1+d}} M(x,t) {\rm d} x {\rm d} t & = E, & &d\geq 4; \\
 \iint_{\mathbb R^{1+d}} M(x,t) {\rm d} x {\rm d} t + \pi \int_{-\infty}^{\infty} |u(0,t)|^2 {\rm d} t &= E, & &d=3.
\end{align*}
\end{corollary} 

\begin{proof}
The monotonicity of $E_-(t)$ and the limit of $E_-(t)$ at positive infinity directly follows the integral expression of $E_-(t)$ given in Lemma \ref{integral expression of energy} and the global integrability of $M(x,t)$ (as well as the integrability of $|u(0,t)|^2$ in dimension 3). By time symmetry we also have 
 \[
  \lim_{t\rightarrow -\infty} E_+(t) = 0. 
 \]
 The limit of $E_-(t)$ at the negative infinity then follows the identity $E_-(t)+E_+(t) = E$. Finally we let $t_0\rightarrow -\infty$ in the integral expression of inward energy $E_-(t)$ and deduce the Morawetz identity. 
\end{proof}

Finally we give a finer property of energy distribution. 
\begin{lemma}
 Given any $c\in (0,1)$, we have 
 \begin{align*}
  &\lim_{t\rightarrow +\infty} E_+(0,ct; t) = 0;& &\lim_{t\rightarrow -\infty} E_-(0,c|t|;t) = 0. 
 \end{align*}
\end{lemma}
\begin{proof}
 By time symmetry it suffices to prove $E_-(0,-ct,t)\rightarrow 0$ as $t\rightarrow -\infty$. We apply Lemma \ref{bigger cone law} with $\kappa = \frac{c+1}{2c}$ and obtain 
 \begin{align*}
  E_-(0,-ct;t) &\lesssim_{c,\mu} \int_{K^-(\frac{1-c}{2}t; t, \frac{1-c}{2}t)} M(x,t') {\rm d} x {\rm d} t' + \pi \int_t^{\frac{1-c}{2}t} |u(0,t')|^2 {\rm d} t'\\
  & \lesssim_{c,\mu} \int_{-\infty}^{\frac{1-c}{2}t} \int_{\Rm^d} M(x,t') {\rm d} x {\rm d} t' + \pi \int_{-\infty}^{\frac{1-c}{2}t} |u(0,t')|^2 {\rm d} t'. 
 \end{align*}
 Since we have $\frac{1-c}{2} t \rightarrow -\infty$, the desired result follows. 
\end{proof}

Now we summarize all the energy distribution properties given above. 

\begin{proposition} \label{energy distribution summary}
 Let $u$ be a solution to \eqref{unified equation inout} with a finite energy $E$. Then it satisfies all the energy distribution properties below. 
 \begin{itemize}
  \item The inward energy $E_-(t)$ is a decreasing function of $t$ satisfying 
  \begin{align*}
   & \lim_{t\rightarrow -\infty} E_-(t) = E; & & \lim_{t\rightarrow +\infty} E_-(t) = 0.
  \end{align*}
  \item The outward energy $E_+(t)$ is an increasing function of $t$ satisfying
  \begin{align*}
   & \lim_{t\rightarrow -\infty} E_+(t) = 0; & & \lim_{t\rightarrow +\infty} E_+(t) = E.
  \end{align*}
  \item The non-directional energy converges to zero in both two time directions
  \[ 
   \lim_{t \rightarrow \pm \infty} \int_{\Rm^d} e'(x,t) {\rm d} x = 0.
  \]
  \item Given any $c\in (0,1)$, we have
  \[
   \lim_{t\rightarrow \pm \infty} \int_{|x|<c|t|} e(x,t) {\rm d} x = 0. 
  \]
 \end{itemize}
\end{proposition}

\begin{remark}
 The proposition above implies that all the energy eventually transforms from the inward energy to outward energy as time moves from $-\infty$ to $\infty$. This transformation may happen at any time and any place, whose amount is given by the integral of the Morawetz density function $M(x,t)$. In dimension 3 this transformation also concentrates near the origin. Intuitively this represents the amount of energy carried by the inward wave which moves through the origin and transforms to outward wave. 
\end{remark}

As an application of the inward/outward energy theory, we may  give some basic asymptotic behaviours of linear Coulomb waves. 

\begin{proposition} \label{half energy}
 Let $u$ be a Coulomb free wave with a finite energy $E$. Then we have 
 \[
  \lim_{t\rightarrow +\infty} \|u\|_{\mathcal{H}^1}^2 = \lim_{t\rightarrow +\infty} \|u\|_{\dot{H}^1}^2 =  \lim_{t\rightarrow +\infty} \|u_t\|_{L^2}^2 = E. 
 \]
\end{proposition}
\begin{proof}
 We start by writting the energy conservation law with more details
 \[
  \int_{\Rm^d} \left(|\mathbf{L}u|^2 + |u_t|^2 + |\slashed{\nabla} u|^2 + \lambda \frac{|u|^2}{|x|^2} + \frac{|u|^2}{|x|}\right) {\rm d} x = 2E.
 \]
 By the inward/outward energy theory, we also have 
 \begin{equation} \label{limit of part of energy}
  \lim_{t\rightarrow +\infty} \int_{\Rm^d} \left(\frac{1}{4}|\mathbf{L}u+u_t|^2 + \frac{1}{4}|\slashed{\nabla} u|^2 + \frac{\lambda}{4} \frac{|u|^2}{|x|^2} + \frac{|u|^2}{4|x|}\right) {\rm d} x = \lim_{t\rightarrow +\infty} E_-(t) = 0.
 \end{equation}
 Combining these two identities we have 
 \begin{align*} 
   \lim_{t\rightarrow +\infty} \left|\|\mathbf{L} u\|_{L^2} - \|u_t\|_{L^2}\right| \leq \lim_{t\rightarrow +\infty} \|\mathbf{L}u +u_t\|_{L^2(\Rm^d)} = 0
 \end{align*}
 and 
 \[
  \lim_{t\rightarrow +\infty} \left(\|\mathbf{L} u\|_{L^2}^2 + \|u_t\|_{L^2}^2\right) = 2E. 
 \]
 This immediately gives 
 \[ 
  \lim_{t\rightarrow +\infty} \|u_t\|_{L^2}^2 = E. 
 \]
 The limit of $\|u\|_{\mathcal{H}^1}^2$ then follows from the energy conservation law
 \[
  \|u\|_{\mathcal{H}^1}^2 + \|u_t\|_{L^2}^2 = 2 E, \qquad \forall t \in \Rm.
 \]
 Finally we obtain the limit of $\|u\|_{\dot{H}^1}^2$ by a combination of \eqref{limit of part of energy} and the identity 
 \[ 
  \|u\|_{\mathcal{H}^1}^2 = \|u\|_{\dot{H}^1}^2 + \int_{\Rm^d} \frac{|u|^2}{|x|} {\rm d} x. 
 \]
\end{proof}

\begin{corollary} \label{balance for L2}
 Let $u$ be a free Coulomb wave with initial data $(u_0,u_1)\in L^2 \times \mathcal{H}^{-1}$. Then we have 
 \[
  \lim_{t\rightarrow +\infty} \|u\|_{L^2}^2 = \frac{1}{2}\|u_0\|_{L^2}^2 + \frac{1}{2}\|u_1\|_{\mathcal{H}^{-1}}^2. 
 \]
\end{corollary}
\begin{proof}
 Clearly $\mathbf{H}^{-1/2} u$ is also a free Coulomb wave with initial data $\mathbf{H}^{-1/2} (u_0,u_1)\in \mathcal{H}^1 \times L^2$. Thus we may apply Proposition \ref{half energy} on $\mathbf{H}^{-1/2} u$ and obtain  
 \begin{align*}
  \lim_{t\rightarrow +\infty} \|u\|_{L^2}^2  = \lim_{t\rightarrow +\infty} \|\mathbf{H}^{-1/2} u\|_{\mathcal{H}^1}^2 & = \frac{1}{2} \|\mathbf{H}^{-1/2} u_0\|_{\mathcal{H}^1}^2 + \frac{1}{2} \|\mathbf{H}^{-1/2} u_1\|_{L^2}^2\\
  & = \frac{1}{2}\|u_0\|_{L^2}^2 + \frac{1}{2}\|u_1\|_{\mathcal{H}^{-1}}^2.
 \end{align*}
\end{proof}

\section{Weighted Morawetz estimates}

In this section we give a few weighted Morawetz estimates, which are perhaps the most powerful applications of the inward/outward energy theory. We start by giving a basic weighted Morawetz estimate, which works in both the linear and nonlinear case, as long as the initial data decay at certain rate as the spatial variable tends to infinity. 

\begin{proposition} [Classic weighted Morawetz inequality] \label{prop weighted Morawetz}
Assume $\kappa \in (0,1/2)$ and let $u$ be a solution to the Coulomb wave equation \eqref{unified equation inout} with a finite energy satisfying
 \[
  \int_{\Rm^d} |x|^\kappa e_-(x,0) {\rm d} x < +\infty.
 \]
 Then $u$ satisfies the global weighted Morawetz estimate 
 \begin{align*}
  \iint_{\Rm^d \times \Rm^+} (t+|x|)^\kappa M(x,t) {\rm d} x {\rm d} t & \lesssim_{\kappa} \int_{\Rm^d} |x|^\kappa e_-(x,0) {\rm d} x, & & d\geq 4;\\
  \int_0^\infty t^\kappa |u(0,t)|^2 {\rm d} t + \iint_{\Rm^3 \times \Rm^+} (t+|x|)^\kappa M(x,t) {\rm d} x {\rm d} t & \lesssim_{\kappa} \int_{\Rm^d} |x|^\kappa e_-(x,0) {\rm d} x, & & d=3.
 \end{align*}
\end{proposition}
\begin{proof}
 We consider the 3-dimensional case. The proof in higher dimensional case is similar. We start by applying the energy flux formula of the inward energy on the region 
 \[
  \Omega(s, s') = \left\{(x,t)\in \mathbb R^3 \times [0,+\infty): s<|x|+t<s' \right\}
 \]
 and obtain 
 \begin{align*}
    -E_-(s,s';0)  + Q_-^-(s';0,s')+\pi \int_s^{s'} |u(0,t)|^2 {\rm d} t - Q_-^-(s; 0,s)  = - \iint_{\Omega(s,s')} M(x,t) {\rm d} x {\rm d} t.
 \end{align*}
 Next we let $s'\rightarrow +\infty$. By \eqref{cone surface limit} we may ignore the term $Q_-^-(s; 0, s)$ and obtain
  \[
  -E_-(0;s,+\infty)  +\pi \int_s^{\infty} |u(0,t)|^2 {\rm d} t -Q_-^-(s; 0,s) + \iint_{\Omega(s)} M(x,t) {\rm d} x {\rm d} t= 0.
 \]
 This is actually the energy flux formula for the unbounded region 
 \[
  \Omega(s) = \{(x,t)\in \Rm^3 \times [0,+\infty) : |x|+t>s\}.
 \]   
   We move the terms with a negative sign to the other side of the equality and insert the details of $Q_-^-$ and $E_-$ to deduce 
 \begin{align*} 
  \pi \int_s^{\infty} |u(0,t)|^2 {\rm d} t + \iint_{\Omega(s)} M(x,t) {\rm d} x {\rm d} t  =  \int_{|x|>s} e_-(x,0) {\rm d} x + \frac{1}{\sqrt{2}} \int_{C^-(s;0,s)} e'(x,t) {\rm d} S.
 \end{align*}
 Multiplying both sides by $\kappa s^{\kappa-1}$ and integrating for $s$ from $0$ to $R\gg 1$, we obtain
 \begin{align}
  & \pi  \int_0^\infty \min\{t,R\}^\kappa |u(0,t)|^2 {\rm d} t + \iint_{\Bbb R^3\times \Bbb R^+} \min\{|x|+t,R\}^\kappa M(x,t) {\rm d} x {\rm d} t \nonumber\\
  &\quad =  \int_{\Bbb R^3} \min\{|x|,R\}^\kappa e_-(x,0) {\rm d} x + \iint_{\Omega(0,R)} \kappa (|x|+t)^{\kappa-1} e'(x,t)  {\rm d} x {\rm d} t. \label{m weighted M}
 \end{align}
We combine the key inequality \eqref{key estimate for in out energy theory} and the fact $\mu = 1/2$ to deduce 
\[
 e'(x,t) \leq 2|x| M(x,t) \; \Rightarrow \; \kappa (|x|+t)^{\kappa-1} e'(x,t) \leq 2\kappa |x| (|x|+t)^{\kappa-1} M(x,t) \leq 2\kappa (|x|+t)^\kappa M(x,t).
\]
Inserting this into \eqref{m weighted M}, we obtain 
 \begin{align*}
   \pi  \int_0^\infty & \min\{t,R\}^\kappa |u(0,t)|^2 {\rm d} t + \iint_{\Bbb R^3\times \Bbb R^+} \min\{|x|+t,R\}^\kappa M(x,t) {\rm d} x {\rm d} t \\
  & \leq   \int_{\Bbb R^3} \min\{|x|,R\}^\kappa e_-(x,0) {\rm d} x + \iint_{\Omega(0,R)} 2\kappa (|x|+t)^\kappa M(x,t)  {\rm d} x {\rm d} t. 
 \end{align*}
 Since $2\kappa < 1$, the last term in the right hand side can be absorbed by the second term in the left hand side. As a result, we have 
 \begin{align*}
   \pi  \int_0^\infty  \min\{t,R\}^\kappa |u(0,t)|^2 {\rm d} t + (1-2\kappa)\iint_{\Bbb R^3\times \Bbb R^+} & \min\{|x|+t,R\}^\kappa M(x,t) {\rm d} x {\rm d} t \\
  & \leq   \int_{\Bbb R^3} \min\{|x|,R\}^\kappa e_-(x,0) {\rm d} x. 
 \end{align*}
 Finally we may let $R\rightarrow +\infty$ and conclude 
 \begin{align*}
  \pi  \int_0^\infty  t^\kappa |u(0,t)|^2 {\rm d} t + (1-2\kappa)\iint_{\Bbb R^3\times \Bbb R^+} & (|x|+t)^\kappa M(x,t) {\rm d} x {\rm d} t \leq   \int_{\Bbb R^3} |x|^\kappa e_-(x,0) {\rm d} x. 
 \end{align*}
 This finishes the proof.
\end{proof}

\begin{corollary} \label{coro weighted Morawetz} 
 Let $u$ be a solution as in Proposition \ref{prop weighted Morawetz}. Then we have 
 \begin{align*}
  &E_-(t) \lesssim_{d,\kappa} t^{-\kappa};& &E_-(t) \in L^{1/\kappa} (\Rm).
 \end{align*}
\end{corollary}
\begin{proof}
 Again we only give the proof for $d=3$. We recall the expression of the inward energy by Morawetz integral (see Lemma \ref{integral expression of energy}) and obtain
\begin{align*}
 E_-(t_0) & = \int_{t_0}^\infty \int_{\Bbb R^n} M(x,t) {\rm d} x {\rm d} t + \pi \int_{t_0}^\infty |u(0,t)|^2 {\rm d} t \nonumber \\
 & \leq \int_{t_0}^\infty t^{-\kappa}\left(\int_{\Bbb R^n} (t+|x|)^\kappa M(x,t) {\rm d} x\right) {\rm d} t + \pi  \int_{t_0}^\infty t^{-\kappa} \cdot t^\kappa |u(0,t)|^2 {\rm} t. 
\end{align*} 
The conclusion then follows an application of the following lemma with the measure $\mu'$ given by ${\rm d} \mu' = \left(\int_{\Rm^d} (t+|x|)^\kappa M(x,t) {\rm d} x\right) {\rm d} t$ and ${\rm d}\mu' = t^\kappa |u(0,t)|^2 {\rm d} t $. The finiteness of this measure is exactly the consequence of weighted Morawetz estimates given in the proposition above. 
\end{proof}

\begin{lemma}
 Let $\mu'$ be a continuous, nonnegative, finite measure on $\Rm^+$ and $\kappa \in (0,1)$ be a constant. Then the function $g: \Bbb R^+ \rightarrow \Bbb R$ defined by $g(x) = \int_x^\infty y^{-\kappa} {\rm d}\mu'(y)$ satisfies $g(x)\leq x^{-\kappa} \mu'(\Rm^+)$ and $g \in L^{1/\kappa} (\Bbb R^+)$. In fact we have $\|g\|_{L^{1/\kappa}(\Bbb R^+)} \leq \mu'(\Bbb R^+)$.
\end{lemma}
\begin{proof}
 It is clear that 
 \[
  g(x) = \int_x^\infty y^{-\kappa} {\rm d}\mu'(y) \leq \int_x^\infty x^{-\kappa} {\rm d}\mu'(y)  \leq x^{-\kappa} \mu'([x,+\infty))\leq x^{-\kappa} \mu'(\Rm^+). 
 \]
 To verify $g\in L^{1/\kappa} (\Rm^+)$, we start by finding an upper bound of $g(x)$ via H\"{o}lder inequality
 \begin{align*}
  g(x)  =  \int_x^\infty (1\cdot y^{-\kappa}) {\rm d}\mu'(y) & \leq  \left(\int_x^\infty 1^{1/(1-\kappa)} {\rm d}\mu'(y)\right)^{1-\kappa} \left(\int_x^\infty \left(y^{-\kappa}\right)^{1/\kappa} {\rm d}\mu'(y)\right)^{\kappa} \\
  & \leq \mu'(\Bbb R^+)^{1-\kappa} \left(\int_x^\infty y^{-1} {\rm d}\mu'(y)\right)^{\kappa}
 \end{align*}
 Therefore we have
 \begin{align*}
  \int_0^\infty |g(x)|^{1/\kappa} {\rm d} x & \leq \mu'(\Bbb R^+)^{(1-\kappa)/\kappa} \int_0^\infty \left(\int_x^\infty y^{-1} {\rm d}\mu'(y)\right) {\rm d} x \\
  & =  \mu'(\Bbb R^+)^{(1-\kappa)/\kappa} \int_0^\infty \left(\int_0^y y^{-1} {\rm d} x \right) {\rm d}\mu'(y)\\
  & = \mu'(\Bbb R^+)^{1/\kappa}.
 \end{align*}
\end{proof}

Next we consider stronger weighted Morawetz estimates if decay rate of initial data is higher. We first prove a technical lemma. 

\begin{lemma} \label{total energy flux}
 Let $u$ be a global solution to \eqref{unified equation inout} with a finite energy. Then given $R>0$, we have 
 \[
  \frac{1}{\sqrt{2}} \int_{C^+(-R; 0,+\infty)} \left(\frac{1}{2} |\mathbf{L}_+ u|^2 + e'(x,t)\right) {\rm d} S \leq \int_{|x|>R} e(x,0) {\rm d} x. 
 \]
\end{lemma}
\begin{proof}
 We apply the energy flux formula on the cone $K^+(-R, 0, t_0)$ for inward and outward energy respectively, and then add them together to deduce 
 \begin{align}
  E_+(0,R+t_0; t_0) + E_-(0,R+t_0;t_0) = & E_+(0,R;0) + E_-(0,R;0) \nonumber \\
  & + \frac{1}{\sqrt{2}} \int_{C^+(-R; 0,t_0)} \left(\frac{1}{2} |\mathbf{L}_+ u|^2+ e'(x,t)\right) {\rm d} S. \label{in out sum lemma 1}
 \end{align}
 For any $t_0>0$, the left hand side satisfies 
 \[
  {\rm LHS} \leq E = E_+(0) + E_-(0) = E_+(R,+\infty; 0) + E_-(R,+\infty;0) + E_+(0,R;0) + E_-(0,R;0)
 \]
 Inserting this into \eqref{in out sum lemma 1} and cancelling $E_\pm (0,R;0)$, we obtain 
 \[
  \frac{1}{\sqrt{2}} \int_{C^+(-R; 0,t_0)} \left(\frac{1}{2} |\mathbf{L}_+ u|^2 + e'(x,t)\right) {\rm d} S \leq E_+(R,+\infty;0) + E_-(R,\infty;0). 
 \]
 We then finish the proof by making $t_0\rightarrow +\infty$ and using the inequality \eqref{outside inward outward total}. 
\end{proof}

\begin{proposition} \label{weighted Morawetz kappa 1}
 Let $(u_0,u_1)\in (\mathcal{H}^1 \cap L^2) \times (L^2 \cap \mathcal{H}^{-1})$ be initial data satisfying the weighted energy estimate 
 \[
  K \doteq \int_{\Rm^d} |x| \left(\frac{1}{2}|\nabla u_0|^2 + \frac{|u_0|^2}{2|x|} + \frac{1}{2} |u_1|^2\right) {\rm d} x < +\infty. 
 \]
 Then the corresponding free Coulomb wave $u$ satisfies 
 \begin{itemize}
  \item[(i)] The following global integral estimate holds
  \[
   \iint_{\Rm^d \times \Rm^+} t M_0(x,t) {\rm d} x {\rm d} t + \int_{0}^\infty \int_{|x|<t} \frac{(t-|x|)|u|^2}{4|x|^2} {\rm d} x {\rm d} t \leq \frac{3}{2} K,
  \]
  where
   \[ 
    M_0(x,t) = \frac{1}{2} \cdot \frac{|\slashed{\nabla} u(x,t)|^2}{|x|} + \frac{\lambda}{2} \cdot \frac{|u(x,t)|^2}{|x|^3} 
   \]
   In the three-dimensional case, we may also add the following term to the left hand side of the inequality above 
   \[
    \pi \int_0^\infty t|u(0,t)|^2 {\rm d} t.
   \]
  \item[(ii)] The inward energy satisfies the decay estimate $E_-(t) \lesssim t^{-1}$. In addition, we have ($E$ is the energy of the solution $u$)
  \[
   \|u_t\|_{L^2(\Rm^d)}^2 = E + O(t^{-1/2}). 
  \]
  \item[(iii)] We also have the integral estimate 
  \[
   \int_{\Rm^d} \left(\frac{1}{2}|\mathbf{L} u(x,t)|^2 + \frac{1}{2}|u_t(x,t)|^2 + e'(x,t)\right)\left(|t-|x||+1\right) {\rm d} x \lesssim \ln t, \qquad t\gg 1.  
  \]
 % Given any $\kappa \in [1/2,1)$ and $c>0$, the energy in the centre part satisfies the decay estimate
 % \[
 %  \int_{|x|<t-ct^{\kappa}} \left(|\nabla u|^2 + |u_t|^2 + e'(x,t)\right) {\rm d} x \lesssim t^{-\kappa} \ln t, \qquad t\gg 1. 
%  \]                  
 \end{itemize}
\end{proposition}

\begin{proof}
 The proof is similar to the regular weighted Morawetz estimate. We apply the energy flux formula of inward energy in the unbounded region 
 \[
  \Omega(s) = \{(x,t): t>0, |x|+t \geq s\}
 \]
 and obtain 
 \[
  \pi \int_s^\infty |u(0,t)|^2 {\rm d} t + \iint_{\Omega(s)} M(x,t) {\rm d} x {\rm d} t = \int_{|x|>s} e_-(x,0) {\rm d} x + \frac{1}{\sqrt{2}}\int_{C^-(s;0,s)} e'(x,t) {\rm d} S. 
 \]
Please note that the term regarding $|u(0,t)|^2$ is ignored in the higher dimensional case $d\geq 4$. We integrate from $s=0$ to $R\gg 1$ to deduce 
 \begin{align*}
  \pi \int_0^\infty \min \{t, R\} |u(0,t)|^2 {\rm d} t & + \iint_{\Rm^d \times \Rm^+} \min\{|x|+t, R\} M(x,t) {\rm d} x {\rm d} t\\
  & = \int_{\Rm^d} \min\{|x|,R\} e_-(x,0) {\rm d} x + \iint_{\Omega(0,R)} e'(x,t) {\rm d} x {\rm d} t. 
 \end{align*}
We observe 
\begin{align*}
 &M(x,t) = M_0(x,t) + \frac{|u|^2}{4|x|^2};& &e'(x,t) = |x| M_0(x,t) + \frac{|u|^2}{2|x|};& 
\end{align*}
and deduce
\begin{align}
  &\pi \int_0^\infty \min \{t, R\} |u(0,t)|^2 {\rm d} t + \iint_{\Omega(s)} R M(x,t) {\rm d} x {\rm d} t + \iint_{\Omega(0,R)} t M_0(x,t){\rm d} x {\rm d}t\nonumber \\
  & \qquad + \iint_{\Omega_1(R)} \frac{(t-|x|)|u|^2}{4|x|^2} {\rm d} x {\rm d} t= \int_{\Rm^d} \min\{|x|,R\} e_-(x,0) {\rm d} x + \iint_{\Omega_2(R)} \frac{(|x|-t)|u|^2}{4|x|^2} {\rm d} x {\rm d} t. \label{weighted Morawetz kappa 1 e1}
 \end{align}
 Here 
 \begin{align*}
  \Omega_1(R) &= \left\{(x,t)\in \Rm^d \times \Rm^+: |x|+t<R, |x|<t\right\};\\
  \Omega_2(R) &= \left\{(x,t)\in \Rm^d \times \Rm^+: |x|+t<R, |x|>t\right\}.
 \end{align*}
 Next we give an upper bound of the integral over $\Omega_2(R)$ by making use of Lemma \ref{total energy flux}
 \begin{align*} 
  \iint_{\Omega_2(R)} \frac{(|x|-t)|u|^2}{4|x|^2} {\rm d} x {\rm d} t & \leq \iint_{\Omega_2(R)} \frac{|u|^2}{4|x|} {\rm d} x {\rm d} t \leq \frac{1}{2} \iint_{\Omega_2(R)} e'(x,t) {\rm d} x {\rm d} t\\
  & \leq \frac{1}{2} \int_0^R \left(\frac{1}{\sqrt{2}} \int_{C^+(-s; 0, +\infty)} e'(x,t) {\rm d} S \right) {\rm d} s\\
  & \leq \frac{1}{2} \int_0^R \int_{|x|>s} e(x,0) {\rm d} x {\rm d} s\\
  & \leq \frac{1}{2} \int_{\Rm^d} \min\{|x|,R\} e(x,0) {\rm d} x. 
 \end{align*}
 Inserting this into \eqref{weighted Morawetz kappa 1 e1} and discarding the integral over $\Omega(s)$, we have 
 \begin{align}
  \pi \int_0^\infty \min \{t, R\} |u(0,t)|^2 {\rm d} t + \iint_{\Omega(0,R)} t M_0(x,t){\rm d} x {\rm d}t & + \iint_{\Omega_1(R)} \frac{(t-|x|)|u|^2}{4|x|^2} {\rm d} x {\rm d} t \nonumber\\
  & \leq \frac{3}{2}\int_{\Rm^d} \min\{|x|,R\} e(x,0) {\rm d} x. \label{kappa 1 to take limit}
  \end{align}
  Here we use the inequality given in Remark \ref{weighted energy plus negative}.  Letting $R\rightarrow +\infty$ in \eqref{kappa 1 to take limit}, we immediately obtain (i). In order to prove (ii), we start by observing
  \begin{align*}
  \int_{t_0}^{4t_0} \int_{0<|x|<3t_0} \frac{t|u|^2}{|x|^2} {\rm d} x {\rm d} t & \leq \int_{t_0}^{4t_0} \int_{|x|<t/2} \frac{2(t-|x|)|u|^2}{|x|^2} {\rm d} x {\rm d} t + \int_{t_0}^{4t_0} \int_{t/2<|x|<3t_0} \frac{t|u|^2}{|x|^2} {\rm d} x {\rm d} t\\
  & \lesssim_1 \int_{\Rm^d} |x|e(x,0) {\rm d} x + \frac{1}{t_0} \int_{t_0}^{4t_0} \int_{\Rm^d} |u|^2 {\rm d} x {\rm d} t\\
  & \lesssim_1 \int_{\Rm^d} |x|e(x,0) {\rm d} x + \|u_0\|_{L^2}^2 + \|u_1\|_{\mathcal{H}^{-1}}^2.
  \end{align*}
  Here we use (i) and the $L^2$-level conservation law, respectively. Combining this with the estimate of $M_0(x,t)$ in (i), we obtain that 
  \begin{equation} \label{kappa 1 equation 1}
   \int_{t_0}^{4t_0} \int_{0<|x|<3t_0} t M(x,t) {\rm d} x {\rm d} t \lesssim_1  \int_{\Rm^d} |x|e(x,0) {\rm d} x + \|u_0\|_{L^2}^2 + \|u_1\|_{\mathcal{H}^{-1}}^2
  \end{equation}
  is uniformly bounded for all $t\gg 1$. Now we are at the position to prove the decay estimate in (ii). On one hand, we apply Corollary \ref{bigger cone law} with $r=2t_0$,  and $\kappa = 3/2$ to deduce 
  \begin{align*}
   E_-(0,2t_0;t_0) &\lesssim_1 \int_{K^-(4t_0; t_0, 4t_0)} M(x,t) {\rm d} x {\rm d} t + \pi \int_{t_0}^{4t_0} |u(0,t)|^2 {\rm d} t\\
   & \lesssim_1 t_0^{-1} \int_{t_0}^{4t_0} \int_{0<|x|<3t_0} t M(x,t) {\rm d} x {\rm d} t + t_0^{-1} \int_{t_0}^{4t_0} t |u(0,t)|^2 {\rm d} t\\
   & \lesssim_1 t_0^{-1}.
  \end{align*}
  Here we use \eqref{kappa 1 equation 1} and (i). On the other hand, by \eqref{outside inward outward total} and finite speed of energy propagation (see Remark \ref{finite speed of propagation}), we have 
  \[
   E_-(2t,\infty;t) \leq \int_{|x|>2t} e(x,t) {\rm d} x \leq \int_{|x|>t} e(x,0) {\rm d} x \leq t^{-1} K. 
  \]
  A combination of these two inequalities proves the decay estimate in (ii). It immediately follows that 
  \begin{equation} \label{estimate on the difference Lu ut}
   \left|\|\mathbf{L} u(\cdot,t)\|_{L^2(\Rm^d)}-\|u_t(\cdot,t)\|_{L^2(\Rm^d)}\right| \leq \left\|\mathbf{L} u (\cdot,t) + u_t(\cdot,t)\right\|_{L^2(\Rm^d)}\lesssim_1 E_-(t)^{1/2} \lesssim t^{-1/2}. 
  \end{equation} 
  In addition, we may combine the energy conservation law with the estimate $E_-(t) \lesssim t^{-1}$ to deduce
  \begin{equation} 
   \|\mathbf{L} u(\cdot,t)\|_{L^2(\Rm^d)}^2 + \|u_t(\cdot,t)\|_{L^2(\Rm^d)}^2 = 2 E - 2\int_{\Rm^d} e'(x,t) {\rm d} x = 2 E + O(t^{-1}). 
  \end{equation}
  Therefore we have 
  \begin{align*}
   \|\mathbf{L} u(\cdot,t)\|_{L^2}+\|u_t(\cdot,t)\|_{L^2} &= \left(2\|\mathbf{L} u(\cdot,t)\|_{L^2}^2 + 2\|u_t(\cdot,t)\|_{L^2}^2 - \big|\|\mathbf{L} u(\cdot,t)\|_{L^2}-\|u_t(\cdot,t)\|_{L^2}\big|^2\right)^{1/2}\\
   & = 2 E^{1/2} + O(t^{-1}). 
  \end{align*}
  Combining this with \eqref{estimate on the difference Lu ut}, we obtain 
  \[
   \|u_t\|_{L^2(\Rm^d)} = E^{1/2} + O(t^{-1/2}).
  \]
  This finishes the proof of part (ii). To verify (iii), we observe that (iii) is equivalent to 
  \[
   \int_{\Rm^d} (e_+(x,t)+e_-(x,t))\left(\left||x|-t\right|+1\right) {\rm d} x \lesssim \ln t, \qquad t\gg 1. 
  \]
  By the energy conservation law, it suffices to prove 
  \[
   \int_{\Rm^d} (e_+(x,t)+e_-(x,t)) \left||x|-t\right| {\rm d} x \lesssim \ln t, \qquad t\gg 1. 
  \]
  Let us consider the exterior part $|x|>t$ and interior part $|x|<t$ separately. For the exterior part, by the finite speed of energy propagation 
  \[
   \int_{|x|>r} (e_+(x,t)+e_-(x,t)) {\rm d} x \leq \int_{|x|>r} e(x,t) {\rm d} x \leq \int_{|x|>r-t} e(x,0) {\rm d} x, \quad r>t. 
  \]
  An integration of this inequality for $r \in (t,+\infty)$ immediately gives 
  \begin{equation} \label{weighted Morawetz kappa 1 exterior} 
   \int_{|x|>t} (e_+(x,t)+e_-(x,t)) (|x|-t) {\rm d} x \leq \int_{\Rm^d} |x| e(x,0) {\rm d} x = K< +\infty. 
  \end{equation}
  For the interior part, by (ii) we have
  \begin{equation} \label{weighted Morawetz kappa 1 interior inward}
   \int_{|x|<t} e_-(x,t) (t-|x|) {\rm d} x \leq t \int_{|x|<t} e_-(x,t) {\rm d} x \lesssim 1. 
  \end{equation}
  Next we deal with the outward energy part. The cone law gives the identity 
  \begin{align*}
   \int_{|x|<r} e_+(x,t) {\rm d} x = \pi  &\int_{t-r}^t |u(0,t')|^2 {\rm d} t' + \frac{1}{\sqrt{2}} \int_{C^+(t-r;t-r,t)} e'(x,t') {\rm d} S \\
    & \qquad + \int_{t-r}^t \int_{|x|<t'-t+r} M(x,t') {\rm d} x {\rm d} t'.
  \end{align*}
  We integrate this identity for $r \in (0,t)$ and obtain
  \begin{align*}
   \int_{|x|<t} (t-|x|) e_+(x,t) {\rm d} x = & \pi \int_{0}^t t' |u(0,t')|^2 {\rm d} t' + \int_0^t \int_{|x|<t'} e'(x,t') {\rm d x} {\rm d}t' \\
   & \qquad + \int_0^t \int_{|x|<t'} (t'-|x|) M(x,t') {\rm d} x {\rm d} t'.
  \end{align*} 
  Inserting 
  \begin{align*}
   &e'(x,t') = |x| M_0(x,t') + \frac{|u(x,t')|^2}{2|x|};& &M(x,t') = M_0(x,t')+\frac{|u(x,t')|^2}{4|x|^2},&
  \end{align*}
  we have 
  \begin{align}
   \int_{|x|<t} & (t-|x|) e_+(x,t) {\rm d} x \nonumber\\
   & = \pi \int_{0}^t t' |u(0,t')|^2 {\rm d} t' + \int_0^t \int_{|x|<t'} t' M_0(x,t') {\rm d x} {\rm d}t' + \int_0^t \int_{|x|<t'} \frac{(t'+|x|)|u|^2}{4|x|^2} {\rm d} x {\rm d} t' \nonumber\\
    &\leq  \pi \int_{0}^t t' |u(0,t')|^2 {\rm d} t' + \int_0^t \int_{|x|<t'} t' M_0(x,t') {\rm d x} {\rm d}t'  + \int_0^t \int_{|x|<t'/3} \frac{(t'-|x|)|u|^2}{2|x|^2} {\rm d} x {\rm d} t' \nonumber \\
    & \qquad  + \int_0^t \int_{t'/3<|x|<t'} \frac{(t'+|x|)|u|^2}{4|x|^2} {\rm d} x {\rm d} t'. \label{estimate on outward interior kappa 1}
  \end{align} 
  The first three terms in the right hand are all uniformly bounded for all $t\gg 1$ by the conclusion of part (i). The upper bound of the last term depends on the estimate $E_-(t)\lesssim t^{-1}$:
  \begin{align*}
   \int_0^t \int_{t'/3<|x|<t'} \frac{(t'+|x|)|u|^2}{4|x|^2} {\rm d} x {\rm d} t' &\lesssim \int_0^t \int_{t'/3<|x|<t'} \frac{|u|^2}{|x|} {\rm d} x {\rm d} t'\\
   & \lesssim \int_0^t E_-(t') {\rm d} t' \lesssim \ln t. 
  \end{align*}
  Inserting these upper bounds into \eqref{estimate on outward interior kappa 1} gives us 
  \[
   \int_{|x|<t} (t-|x|) e_+(x,t) {\rm d} x \lesssim \ln t, \qquad t\gg 1. 
  \]
 Finally we combine this with \eqref{weighted Morawetz kappa 1 exterior} and \eqref{weighted Morawetz kappa 1 interior inward} to complete the proof of part (iii). 
\end{proof}

\section{Higher order decay estimates}

In this subsection we prove weighted Morawetz estimates of higher order for sufficiently good radial initial data. We start by proving a few technical results. 

\begin{lemma} \label{limit zero infinity case}
 There exist two radial solutions $\phi$, $\psi$ to the Laplace equation 
 \[
  -\Delta u + \frac{u}{|x|} = 0, \qquad x\in \Rm^3\setminus\{0\}
 \]
 satisfying $\phi, \psi \in C^\infty(\Rm^3\setminus \{0\})$ and 
 \begin{itemize}
  \item The solution $\phi(r)$ can extend to a smooth function of $r \in \Rm$; in particular, $\phi$ is a continuous function in $\Rm^3$ and $\phi \in H^2(B(0,1))$;
  \item The solution $\phi(x)$ tends to infinity as $x \rightarrow +\infty$; 
  \item $r \psi(r)$ converges to a nonzero number as $r\rightarrow 0^+$; 
  \item $\psi, \psi_r \lesssim r^{-N}$ for any fixed constant $N>0$ and any sufficiently large radius $r$.
 \end{itemize}
\end{lemma}

\begin{proof}
 The solutions $\phi$ and $\psi$ can be given explicitly. For simplicity we consider the solutions $\Phi (r) = r\phi(r)$ and $\Psi (r) = r \psi(r)$ to the equation
 \[
  -w_{rr} + \frac{w}{r} = 0.
 \]
 We define 
 \begin{align*}
  &\Phi (r) = \sum_{k = 1}^\infty \frac{r^k}{k[(k-1)!]^2};& &\phi(x) = \sum_{k = 1}^\infty \frac{|x|^{k-1}}{k[(k-1)!]^2}.&
 \end{align*}
 These series converges for all $r \in [0,+\infty)$ and $x\in \Rm^3$. A simple calculation verifies that they solve the corresponding Laplace equations. The properties of $\phi$ then follows straight-forward observations. The other linear independent solutions can be given by the formula 
 \begin{align*}
  &\Psi (r) = \Phi (r) \int_r^\infty \Phi^{-2} (s) {\rm d} s;& &\psi (x) = \frac{1}{|x|} \Phi (|x|) \int_{|x|}^\infty \Phi^{-2} (s) {\rm d} s.&
 \end{align*}
 Now we consider the asymptotic behaviour of $\Psi$. Since $\Phi$ is an increasing function of $r$ and satisfies $\Phi(r)\gtrsim r^{N+1}$, we have 
 \begin{align*}
  \Psi(r) = \int_r^\infty \frac{\Phi(r)}{\Phi(s)} \cdot \Phi^{-1}(s) {\rm d} s \leq \int_r^\infty  \Phi^{-1}(s) {\rm d} s \leq \int_r^\infty r^{-N-1} \lesssim r^{-N}, \quad r\gg 1.
 \end{align*}
 Furthermore, we have
 \begin{align*}
  \Phi'(r) = \sum_{k = 1}^\infty \frac{r^{k-1}}{[(k-1)!]^2} = 1 + \sum_{k=1}^\infty \frac{r^k}{(k!)^2} \leq \Phi(r) + 1.
 \end{align*}
 A similar argument as above gives 
 \begin{align*}
  |\Psi'(r)| & = \left|\Phi'(r) \int_r^\infty \Phi^{-2} (s) {\rm d}s - \Phi^{-1}(r) \right|\\
  & \leq \Phi(r) \int_r^\infty \Phi^{-2} (s) {\rm d}s + \int_r^\infty \Phi^{-2} (s) {\rm d}s + \Phi^{-1}(r) \\
  & \lesssim r^{-N}. 
 \end{align*}
 Finally we consider the asymptotic behaviour $\Psi$ near zero. In fact we have 
 \begin{align*}
  \Phi(r) \simeq r\quad \Rightarrow \quad \int_r^\infty \Phi^{-2}(s) {\rm d} s \simeq r^{-1} \quad \Rightarrow \quad \Psi(r) \simeq 1. 
 \end{align*}
 This finishes the proof. 
\end{proof}

\paragraph{Higher dimensional case} A similar result to Lemma \ref{limit zero infinity case} holds for higher dimensions $d\geq 4$ as well. In fact we let $w(r,t)= r^\frac{d-1}{2} u(r,t)$ and rewrite the Laplace equation in the form of 
\[
 - w_{rr} + \lambda \frac{w}{r^2} + \frac{w}{r} = 0. 
\]
It is clear that one of the solution is 
\[
 \Phi(r) = \sum_{k=0}^\infty \frac{(d-2)! }{k! (k+d-2)!}r^{\frac{d-1}{2}+k}
\]
Again another solution can be given by 
\[
 \Psi (r) = \Phi (r) \int_r^\infty \Phi^{-2} (s) {\rm d} s. 
\]
Following the same argument as above, we obtain 
\begin{lemma} \label{limit zero infinity case d4}
 There exist two radial solutions $\phi$, $\psi$ to the Laplace equation 
 \[
  -\Delta u + \frac{u}{|x|} = 0, \qquad x\in \Rm^d\setminus\{0\}
 \]
 satisfying $\phi, \psi \in C^\infty(\Rm^d\setminus \{0\})$ and 
 \begin{itemize}
  \item The solution $\phi(r)$ can extend to a smooth function of $r \in \Rm$; in particular, $\phi$ is a continuous function in $\Rm^d$ and $\phi \in H^2(B(0,1))$;
  \item The solution $\phi(x)$ tends to infinity as $x \rightarrow +\infty$; 
  \item $r^{d-2} \psi(r)$ converges to a nonzero number as $r\rightarrow 0^+$; 
  \item $\psi, \psi_r \lesssim r^{-N}$ for any fixed constant $N>0$ and any sufficiently large radius $r$.
 \end{itemize}
\end{lemma}

\begin{lemma} \label{estimate on H minus 2}
 Let $f \in C_0^\infty(\Rm^d\setminus\{0\})$ be radial. Then there exists $u \in L^2(\Rm^d)$ such that $\mathbf{H}u = f$, thus $u\in \mathcal{H}^{-2}$. In fact there exists a radial solution $u \in C^\infty(\Rm^d \setminus \{0\}) \cap H^1(\Rm^d)$ to the Laplace equation 
 \[
  -\Delta u + \frac{u}{|x|} = f
 \]
 with two constants $c_1,c_2$ so that $u(x) = c_1 \phi(x)$ in a neighbourhood of the origin and $u(x) = c_2 \psi(x)$ near the infinity. 
\end{lemma}
\begin{proof}
 Without loss of generality we assume ${\rm Supp} f \in \{x: r_1\leq |x|\leq R_1\}$. We first let $w = r^\frac{d-1}{2} u(r)$ and solve the equivalent equation 
 \[
  -w_{tt} + \lambda\frac{w}{r^2} + \frac{w}{r} = r^{\frac{d-1}{2}} f(r),\qquad r\in (0,+\infty),
 \]
 with an pair of arbitrary initial data at $r=1$. Our assumption that $f$ vanishes for $r<r_1$ and $r>R_1$ guarantees that there exist constants $a_1, b_1, a_2, b_2$ so that 
 \begin{align*}
  &w = a_1 \Phi+ b_1 \Psi, \; r<r_1;& &w = a_2 \Phi + b_2 \Psi,\; r>R_1.
 \end{align*}
 Clearly $\tilde{w} = w-b_1\Psi - a_2 \Phi$ and $u = r^{-\frac{d-1}{2}} \tilde{w}$ solve the corresponding equations and satisfy
  \begin{align*}
  &u = (a_1-a_2) \phi , \; |x|<r_1;& &u = (b_2-b_1) \psi,\; |x|>R_1.
 \end{align*}
\end{proof}

\begin{proposition} \label{prop weighted Morawetz kappa large}
 Assume that $d\geq 3$ and $\kappa \in (1,3/2)$. Let $u$ be a radial free Coulomb wave and $t_0\in \Rm$. Assume that the data of $u$ at time $t_0$ satisfy 
 \begin{align*}
  &u(t_0)\in C_0^\infty(\Rm^d);& & u_t(t_0)\in C_0^\infty(\Rm^d\setminus \{0\}).
 \end{align*}
Then $u$ satisfies 
 \begin{equation} \label{weighted Morawetz kappa large d4} 
   \int_{0}^\infty \int_{|x|<t/3} t^\kappa M(x,t) {\rm d} x {\rm d} t + \int_0^\infty \int_{|x|>t/3} \frac{t^{\kappa-\frac{1}{2}}\left(\left|t-|x|\right|+1\right) |u|^2}{|x|^2} {\rm d}x {\rm d} t< +\infty.
 \end{equation}
 In particular, if $d=3$, then we have
 \begin{equation} \label{weighted Morawetz kappa large}
  \int_0^\infty t^\kappa |u(0,t)|^2 {\rm d} t + \int_{0}^\infty \int_{|x|<t/3} \frac{t^\kappa |u|^2}{|x|^2} {\rm d} x {\rm d} t + \int_0^\infty \int_{|x|>t/3} \frac{t^{\kappa-\frac{1}{2}}\left(\left|t-|x|\right|+1\right) |u|^2}{|x|^2} {\rm d}x {\rm d} t< +\infty.
 \end{equation}
\end{proposition}
\begin{proof}
 By finite speed of propagation the initial data $(u_0,u_1) = (u(0), u_t(0))$ are also compactly supported. Let us assume that the support of initial data is contained in a ball $B(0,R_0)$. By finite of propagation we have  
 \begin{equation} \label{support of u}
  u(x,t) = 0,\qquad |x|>|t|+R_0.
 \end{equation}
 The general idea is similar to the proof of usual weighted Morawetz estimate. The extra ingredient here is a very precise estimate on the norm $\|u(t)\|_{L^2}$ as well as the spatial distribution of $|u|^2$ as time tends to infinity. We consider the auxiliary wave equation 
 \[
  \left\{\begin{array}{l} \displaystyle \partial_t^2 v + \mathbf{H} v  = 0; \\ v(t_0) = - \mathbf{H}^{-1} u_t(t_0); \\ v_t(t_0) = u(t_0). \end{array}\right.
 \]
One one hand, according to Lemma \ref{estimate on H minus 2}, the initial data $v(t_0) \in \mathcal{H}^1 \cap L^2$ and decay very fast at the infinity. Thus $v$ satisfies the assumptions of Proposition \ref{weighted Morawetz kappa 1}, up to a time translation. On the other hand, $v_t$ is also a Coulomb free wave with the initial data 
\begin{align*}
 v_t (t_0) & = u(t_0); \\
 v_{tt} (t_0) & = -\mathbf{H} v(t_0) = u_t(t_0).
\end{align*}
This implies that $u=v_t$. Combining this fact with Proposition \ref{weighted Morawetz kappa 1}, we obtain that 
\begin{itemize} 
 \item There exists a constant $E_1$, so that $\|u(\cdot,t)\|_{L^2(\Rm^d)}^2 = E_1 + O(|t-t_0|^{-1/2})$ for $t\gg 1$; 
 \item The following integral estimate holds for large time $t \gg 1$
 \begin{equation*}
  \int_{\Rm^d} \left(\left|t-t_0-|x|\right|+1\right) |u(x,t)|^2 {\rm d} x \lesssim \ln (t-t_0). 
 \end{equation*}
 \end{itemize}
 Combining these facts with the uniform boundedness of $\|u(\cdot,t)\|_{L^2(\Rm^d)}$, we obtain 
 \begin{equation} \label{details 3 in distribution of u2}
  \|u(\cdot,t)\|_{L^2(\Rm^d)}^2 = E_1 + O(t^{-1/2}), \qquad \forall t > 0;
 \end{equation} 
 and 
 \begin{equation} \label{details in distribution of u2}
  \int_{\Rm^d} \left(\left|t-|x|\right|+1\right) |u(x,t)|^2 {\rm d} x \lesssim \ln t, \qquad t\gg 1.  
 \end{equation}
 In particular, we have 
 \begin{equation} \label{details 2 in distribution of u2}
  \int_{|x|<t-t^{1/2}} |u(x,t)|^2 {\rm d} x \lesssim t^{-1/2} \ln t, \qquad t\gg 1. 
 \end{equation}
Now we prove the integral estimate. We first consider the case $d=3$. We apply the energy flux formula of inward energy in the unbounded region 
 \[
  \Omega(s) = \{(x,t): t>0, |x|+t \geq s\}
 \]
 and obtain 
 \[
  \pi \int_s^\infty |u(0,t)|^2 {\rm d} t + \iint_{\Omega(s)} M(x,t) {\rm d} x {\rm d} t = \int_{|x|>s} e_-(x,0) {\rm d} x + \frac{1}{\sqrt{2}}\int_{C^-(s;0,s)} e'(x,t) {\rm d} S. 
 \]
 Here 
 \[
  C^-(s;0,s) = \{(x,t): |x|+t=s, 0\leq t \leq s\}. 
 \]
We multiply both sides by $\kappa s^{\kappa-1}$ and integrate from $s=0$ to $R\gg 1$ to deduce 
 \begin{align*}
  &\pi \int_0^\infty \min \{t, R\}^\kappa |u(0,t)|^2 {\rm d} t + \iint_{\Rm^3 \times \Rm^+} \min\{|x|+t, R\}^\kappa M(x,t) {\rm d} x {\rm d} t\\
  & \qquad = \int_{\Rm^3} \min\{|x|,R\}^\kappa e_-(x,0) {\rm d} x + \iint_{\Omega(0,R)} \kappa (|x|+t)^{\kappa-1} e'(x,t) {\rm d} x {\rm d} t. 
 \end{align*}
 Here 
 \[
  \Omega(0,R) = \{(x,t): |x|+t < R, t>0\}. 
 \]
We recall in 3-dimensional case with radial data that 
 \begin{align*}
  & M(x,t) = \frac{|u|^2}{4|x|^2}; & &e'(x,t) =  \frac{1}{2}\frac{|u|^2}{|x|};
 \end{align*}
 and split the region $\Omega(0,R)$ into two parts 
 \begin{align*}
  &\Omega_1 = \left\{(x,t)\in \Omega(0,R): t>(2\kappa-1)|x|\right\}; & &\Omega_2 =  \left\{(x,t)\in \Omega(0,R): t<(2\kappa-1)|x|\right\}.&
 \end{align*}
 to deduce 
 \begin{align}
  \pi \int_0^\infty \min \{t, R\}^\kappa |u(0,t)|^2 {\rm d} t + J_1 & +  R^\kappa \iint_{\Omega(R)} \frac{|u|^2}{4|x|^2} {\rm d} x {\rm d} t \nonumber \\
  & = \int_{\Rm^3} \min\{|x|,R\}^\kappa e_-(x,0) {\rm d} x + J_2. \label{weighted Morawetz cancelled}
 \end{align}
 Here 
 \begin{align*}
  J_1 & =  \iint_{\Omega_1} \frac{(|x|+t)^{\kappa-1}(t+|x|-2\kappa |x|)|u|^2}{4|x|^2}{\rm d} x {\rm d} t\\
  J_2 & = \iint_{\Omega_2} \frac{(|x|+t)^{\kappa-1}(2\kappa |x|-|x|-t)|u|^2}{4|x|^2}{\rm d} x {\rm d} t.
  \end{align*}
  \begin{figure}[h]
 \centering
 \includegraphics[scale=1.1]{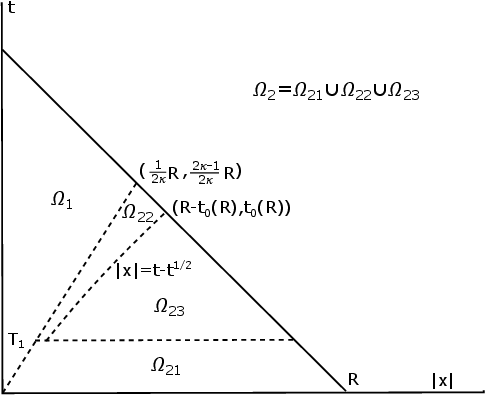}
 \caption{Illustration for regions} \label{figure regions}
\end{figure}
 We consider the upper bound of $J_2$.  We further divide the region $\Omega_2$ into three parts 
 \[
  \Omega_2 = \Omega_{21} \cup \Omega_{22} \cup \Omega_{23}
 \]
 with 
 \begin{align*}
  \Omega_{21} & = \left\{(x,t) \in \Omega_2: 0<t<T_1\right\}; \\
  \Omega_{22} & = \left\{(x,t) \in \Omega_2: t>T_1, |x|<t-t^{1/2} \right\}; \\
  \Omega_{23} & = \left\{(x,t) \in \Omega_2: t>T_1, t-t^{1/2} < |x|<R-t\right\}.
 \end{align*}
 and consider the integrals $J_{21}, J_{22}, J_{23}$ over these regions respectively. Here $T_1=T_1(\kappa)$ is a fixed but sufficiently large number. Please refer to figure \ref{figure regions} for a visual illustration of these regions. We start by observing that 
 \begin{align*}
  J_{21} =  \int_0^{T_1} \int_{\frac{t}{2\kappa-1}<|x|<t+R_0} \frac{(|x|+t)^{\kappa-1}(2\kappa |x|-|x|-t)|u|^2}{4|x|^2}{\rm d} x {\rm d} t 
 \end{align*}
 is a finite constant independent of $R\gg 1$. In addition, we use the spatial distribution property of $|u|^2$ given in \eqref{details 2 in distribution of u2} and obtain 
 \begin{align*}
  J_{22} & \leq \int_{T_1}^{\frac{2\kappa-1}{2\kappa}R} \int_{\frac{t}{2\kappa-1}<|x|<t-t^{1/2}} \frac{(|x|+t)^{\kappa-1}(2\kappa |x|-|x|-t)|u|^2}{4|x|^2}{\rm d} x {\rm d} t \\
  & \lesssim \int_{T_1}^{\frac{2\kappa-1}{2\kappa}R} \int_{\frac{t}{2\kappa-1}<|x|<t-t^{1/2}} \frac{|u|^2}{t^{2-\kappa}} {\rm d} x {\rm d} t\\
  & \lesssim \int_{T_1}^{\frac{2\kappa-1}{2\kappa}R} t^{\kappa- \frac{5}{2}} \ln t \; {\rm d} t \lesssim 1. 
 \end{align*}
 Finally we have 
 \begin{align*}
  J_{23} = \int_{T_1}^{t_0(R)} \int_{t-t^{1/2}<|x|<R-t} \frac{(|x|+t)^{\kappa-1}(2\kappa |x|-|x|-t)|u|^2}{4|x|^2}{\rm d} x {\rm d} t 
 \end{align*}
 Here $t_0(R)$ solves the equation 
 \[
  t-t^{1/2} = R-t,
 \]
 thus 
 \[
  \frac{R}{2} < t_0(R) < \frac{1}{2} (R + R^{1/2}). 
 \]
 By the support of $u$, we have 
 \begin{align*}
  J_{23} \leq \int_{T_1}^{t_0(R)} \int_{t-t^{1/2}<|x|<t+R_0} \frac{(|x|+t)^{\kappa-1}(2\kappa |x|-|x|-t)|u|^2}{4|x|^2}{\rm d} x {\rm d} t 
 \end{align*}
 Here 
 \begin{align*}
 \frac{(|x|+t)^{\kappa-1}(2\kappa |x|-|x|-t)}{4|x|^2} & = \frac{(2t)^{\kappa-1} (2\kappa t - 2t)}{4t^2} + O\left(t^{1/2} \cdot t^{\kappa-3}\right)\\
 & = 2^{\kappa-2} (\kappa-1) t^{\kappa-2} + O(t^{\kappa-5/2}).
 \end{align*}
 Therefore we have 
 \begin{align*}
  J_{23} &\leq  \int_{T_1}^{t_0(R)} \left[\left(2^{\kappa-2} (\kappa-1) t^{\kappa-2} + O(t^{\kappa-5/2})\right) \int_{\Rm^3} |u|^2 {\rm d} x\right] {\rm d} t \\
  & \leq \int_{T_1}^{t_0(R)} \left[\left(2^{\kappa-2} (\kappa-1) t^{\kappa-2} + O(t^{\kappa-5/2})\right) \left(E_1 + O(t^{-1/2})\right)\right] {\rm d} t\\
  & \leq \int_{T_1}^{t_0(R)} \left[2^{\kappa-2} (\kappa-1) t^{\kappa-2} E_1 + O(t^{\kappa-5/2})\right] {\rm d} t \\
  & \leq 2^{\kappa-2} (t_0(R))^{\kappa-1} E_1 + O(1)\\
  & \leq \frac{1}{2} R^{\kappa-1} E_1 + O(1).
 \end{align*}
 In summary we have 
 \[
  J_2 \leq \frac{1}{2} R^{\kappa-1} E_1 + O(1). 
 \]
 Next we recall \eqref{details 3 in distribution of u2}, \eqref{details 2 in distribution of u2} and consider the lower bound
 \begin{align*}
  \iint_{\Omega(R)} \frac{|u|^2}{4|x|^2} {\rm d} x {\rm d} t &\geq \int_{t_0(R)}^\infty \int_{t-t^{1/2}<|x|<t+R_0} \frac{|u|^2}{4|x|^2} {\rm d} x {\rm d} t\\
  & \geq \int_{t_0(R)}^\infty \frac{1}{4(t+R_0)^2} \left(\int_{t-t^{1/2}<|x|<t+R_0} |u|^2 {\rm d} x\right) {\rm d} t\\
  & \geq \int_{t_0(R)}^\infty \left(\frac{1}{4t^2} - O(t^{-3})\right) \left(E_1 - O\left(t^{-1/2}\ln t\right)\right) {\rm d} t\\
  & \geq \int_{t_0(R)}^\infty \left(\frac{E_1}{4t^2} - O\left(t^{-5/2}\ln t\right)\right) {\rm d} t\\
  & \geq \frac{E_1}{4 t_0(R)} - O\left(t_0(R)^{-3/2} \ln t_0(R)\right)\\
  & \geq \frac{E_1}{2R} - O\left(R^{-3/2} \ln R\right).
 \end{align*}
 Inserting these lower/upper bounds into \eqref{weighted Morawetz cancelled} and cancelling the term $E_1 R^{\kappa-1}/2$, we obtain 
 \begin{align*} 
 \pi \int_0^\infty \min \{t, R\}^\kappa |u(0,t)|^2 {\rm d} t + J_1  \leq \int_{\Rm^d} \min\{|x|,R\}^\kappa e_-(x,0) {\rm d} x + O(1).
 \end{align*}
 We then let $R\rightarrow +\infty$ and obtain 
 \[
  \pi \int_0^\infty t^\kappa |u(0,t)|^2 {\rm d} t + \int_0^\infty \int_{|x|<\frac{t}{2\kappa-1}} \frac{(|x|+t)^{\kappa-1}(t+|x|-2\kappa |x|)|u|^2}{4|x|^2}{\rm d} x {\rm d} t < +\infty.
 \]
 It immediately follows that 
 \begin{equation} \label{weighted Morawetz kappa large 3}
  \pi \int_0^\infty t^\kappa |u(0,t)|^2 {\rm d} t + \int_0^\infty \int_{|x|<t/3} \frac{t^{\kappa}|u|^2}{|x|^2}{\rm d} x {\rm d} t < +\infty.
 \end{equation}
 In addition, we use \eqref{details in distribution of u2} to deduce 
 \begin{align*}
  \int_{0}^\infty \int_{|x|>t/3} \frac{t^{\kappa-\frac{1}{2}}(|t-|x||+1) |u|^2}{|x|^2} {\rm d} x {\rm d} t & \lesssim 1+ \int_{2}^\infty \int_{|x|>t/3} t^{\kappa - \frac{5}{2}} \left(\left|t-|x|\right|+1\right) |u|^2 {\rm d} x {\rm d} t\\
  & \lesssim 1 + \int_2^\infty t^{\kappa-\frac{5}{2}}\ln t {\rm d} t < +\infty. 
 \end{align*}
 A combination of this with \eqref{weighted Morawetz kappa large 3} immediately gives the integral estimate \eqref{weighted Morawetz kappa large}.  
 
 Next we consider the higher dimensional case $d\geq 4$. The proof is similar to the 3-dimensional case. Thus we only sketch the proof and emphasize the differences between these two cases. Following the same argument as in dimension 3, we obtain 
  \begin{align*}
   \iint_{\Rm^d \times \Rm^+} & \min\{|x|+t, R\}^\kappa M(x,t) {\rm d} x {\rm d} t\\
  & = \int_{\Rm^d} \min\{|x|,R\}^\kappa e_-(x,0) {\rm d} x + \iint_{\Omega(0,R)} \kappa (|x|+t)^{\kappa-1} e'(x,t) {\rm d} x {\rm d} t. 
 \end{align*}
 In the higher dimensional case we have
 \begin{align} \label{def of M d4}
  & M(x,t) = \lambda \frac{|u|^2}{2|x|^3} + \frac{|u|^2}{4|x|^2}; & &e'(x,t) =  \lambda \frac{|u|^2}{2|x|^2} + \frac{1}{2}\frac{|u|^2}{|x|};& 
 \end{align}
 Thus we may write the identity above in the form of 
 \begin{align}
 J_1  + J_3 +  R^\kappa \iint_{\Omega(R)} M(x,t) {\rm d} x {\rm d} t = \int_{\Rm^d} \min\{|x|,R\}^\kappa e_-(x,0) {\rm d} x + J_2 + J_4. \label{weighted Morawetz cancelled d4}
 \end{align}
 Here $J_1$, $J_2$ are defined in exactly the same manner as in dimension 3. While $J_3$ and $J_4$ are defined by
  \begin{align*}
  J_3 & =  \iint_{\Omega_3} \frac{\lambda (|x|+t)^{\kappa-1}(t+|x|-\kappa |x|)|u|^2}{2|x|^3}{\rm d} x {\rm d} t,& &\Omega_3=\{(x,t)\in \Omega(0,R): t>(\kappa-1)|x|\}; &\\
  J_4 & = \iint_{\Omega_4} \frac{\lambda (|x|+t)^{\kappa-1}(\kappa |x|-|x|-t)|u|^2}{2|x|^3}{\rm d} x {\rm d} t.& &\Omega_4=\{(x,t)\in \Omega(0,R): t<(\kappa-1)|x|\}.&
  \end{align*}
  The same argument as in 3-dimensional case gives the upper bound 
  \begin{equation} \label{upper bound J2 d4}
   J_2 \leq \frac{1}{2} R^{\kappa-1} E_1 + O(1); 
  \end{equation}
  and the lower bound
  \begin{equation} \label{lower bound outside d4}
   \iint_{\Omega(R)} M(x,t) {\rm d} x {\rm d} t \geq \iint_{\Omega(R)} \frac{|u|^2}{4|x|^2} {\rm d} x {\rm d} t \geq \frac{E_1}{2R} - O\left(R^{-3/2} \ln R\right).
  \end{equation}
  We also need to consider the upper bound of $J_4$. In fact we have 
  \begin{align*}
   J_4 & \leq \int_0^{\frac{\kappa-1}{\kappa} R} \int_{|x|>\frac{t}{\kappa-1}} \frac{\lambda (|x|+t)^{\kappa-1}(\kappa |x|-|x|-t)|u|^2}{2|x|^3}{\rm d} x {\rm d} t\\
   & = \int_0^{T_1} \int_{\frac{t}{\kappa-1}<|x|<t+R_0} \frac{\lambda (|x|+t)^{\kappa-1}(\kappa |x|-|x|-t)|u|^2}{2|x|^3}{\rm d} x {\rm d} t\\
    & \qquad + \int_{T_1}^{\frac{\kappa-1}{\kappa} R} \int_{|x|>\frac{t}{\kappa-1}} \frac{\lambda (|x|+t)^{\kappa-1}(\kappa |x|-|x|-t)|u|^2}{2|x|^3}{\rm d} x {\rm d} t\\
    & \leq O(1) + C(d,\kappa)  \int_{T_1}^{\frac{\kappa-1}{\kappa} R} \int_{|x|>\frac{t}{\kappa-1}} \frac{|u(x,t)|^2}{t^{3-\kappa}} {\rm d} x {\rm d} t \leq O(1).
   \end{align*}
  Inserting these upper or lower bounds into \eqref{weighted Morawetz cancelled d4}, we obtain 
 \[
  J_1  + J_3  \leq \int_{\Rm^d} \min\{|x|,R\}^\kappa e_-(x,0) {\rm d} x + O(1).
 \]
 A limit of $R\rightarrow +\infty$ gives 
 \begin{align*}
   \int_0^\infty \int_{|x|<\frac{t}{2\kappa-1}} \frac{(|x|+t)^{\kappa-1}(t+|x|-2\kappa |x|)|u|^2}{4|x|^2}{\rm d} x {\rm d} t &< +\infty;\\
   \int_0^\infty \int_{|x|<\frac{t}{\kappa-1}}  \frac{\lambda (|x|+t)^{\kappa-1}(t+|x|-\kappa |x|)|u|^2}{2|x|^3}{\rm d} x {\rm d} t &< +\infty. 
 \end{align*}
 We recall \eqref{def of M d4} and obtain 
 \[
  \int_0^\infty \int_{|x|<t/3} t^\kappa M(x,t) {\rm d} x {\rm d} t < +\infty. 
 \]
 This give the finiteness of the first term in \eqref{weighted Morawetz kappa large d4}. The finiteness of the second term can be dealt with in the same manner as $d=3$. 
 \end{proof}

\begin{remark} \label{faster center decay}
 Let $u$ be a free Coulomb wave as in Proposition \ref{prop weighted Morawetz kappa large} and $\kappa \in (1,3/2)$. Then we have 
 \begin{align*}
  \int_1^\infty \int_{|x|>t/3} t^{\kappa-\frac{1}{2}}\left(\left|t-|x|\right|+1\right) M(x,t) {\rm d}x {\rm d} t \lesssim_1  \int_1^\infty \int_{|x|>t/3} \frac{t^{\kappa-\frac{1}{2}}\left(\left|t-|x|\right|+1\right)|u|^2}{|x|^2} {\rm d}x {\rm d} t<\infty.
 \end{align*}
 It immediately follows that 
 \begin{equation} \label{weighted Morawetz kappa large d4 M}
  \int_0^\infty \int_{|x|>t/3} t^{\kappa-\frac{1}{2}}\left(\left|t-|x|\right|+1\right) M(x,t) {\rm d}x {\rm d} t < +\infty. 
 \end{equation}
 A combination of this with \eqref{weighted Morawetz kappa large d4} and \eqref{weighted Morawetz kappa large} implies for any $c\in (0,1)$ that
 \begin{align*}
  \int_{0}^\infty t^\kappa |u(0,t)|^2 {\rm d} t + \int_{0}^\infty \int_{|x|<ct} \frac{t^\kappa |u|^2}{|x|^2} {\rm d} x {\rm d} t &< +\infty, & &d=3; \\
  \int_{0}^\infty \int_{|x|<ct} t^\kappa M(x,t) {\rm d} x {\rm d} t &< +\infty, & &d\geq 4.
 \end{align*}
 It immediately follows Corollary \ref{bigger cone law} and Remark \ref{bigger cone law 2} that 
 \[
  \int_{|x|<c t} (e_+(x,t)+e_-(x,t)) {\rm d} x \lesssim t^{-\kappa}, \qquad \forall c\in (0,1).
 \]
\end{remark}

\begin{corollary} \label{corollary in out pointwise}
 Let $u$ be a free Coulomb wave as in Proposition \ref{prop weighted Morawetz kappa large}. Assume that $\kappa \in (1,3/2)$ is a constant. Then we have the following decay estimate on the energy contained in the centre part of space. Here $t\gg 1$ and $1<r<t/2$. The implicit constants in the inequalities do not depend on $t,r$. 
 \begin{align}
  \int_{|x|<t-r} e_-(x,t) {\rm d} x + \int_{C^+(r; t,+\infty)} |\mathbf{L}_+ u|^2 {\rm d} S &\lesssim \max\left\{t^{-\kappa}, r^{-1} t^{-\kappa+\frac{1}{2}}\right\}; \label{decay of E minus tr}\\
  \int_{|x|<t-r} e_+(x,t) {\rm d} x & \lesssim \max\left\{r^{-\kappa}, r^{-2} \ln^2 t\right\}. \nonumber
 \end{align}
 In addition, we have the point-wise estimate: If $t\gg 1$ and $\ln t < r < t/2$, then 
 \[
  \sup_{|x|<t-r} |x|^\frac{d-1}{2} |u(x,t)| \lesssim \left(r^{-1}\ln t\right)^{1/4} \left(\max\left\{r^{-\kappa}, r^{-2} \ln^2 t\right\}\right)^{1/4}.
 \]
\end{corollary}
\begin{proof}
 By the energy flux in the unbounded region 
 \[
  \Omega = \left\{(x,t'): |x|<t'-r, t'>t\right\}, 
 \]
 we have 
 \begin{equation} \label{center estimate e minus}
   \int_{|x|<t-r} e_-(x,t) {\rm d} x + \frac{1}{2\sqrt{2}}\int_{C^+(r;t,+\infty)} |\mathbf{L}_+ u|^2 {\rm d} S = \iint_{\Omega} M(x,t') {\rm d} x {\rm d} t' + \pi \int_{t}^\infty |u(0,t')|^2 {\rm d} t'.
 \end{equation} 
 Here the integral of $|u(0,t')|^2$ is ignored if $d\geq 4$. We split $\Omega$ into two parts $\Omega_1, \Omega_2$ with 
 \begin{align*}
  &\Omega_1 =\{(x,t')\in \Omega: |x|<t'/3\};& &\Omega_2 = \{(x,t')\in \Omega: |x|>t'/3\};&
 \end{align*}
 and observe 
 \begin{align*}
 \pi \int_{t}^\infty |u(0,t')|^2 {\rm d} t' & \lesssim t^{-\kappa} \int_{t}^\infty (t')^\kappa |u(0,t')|^2 {\rm d} t' \lesssim t^{-\kappa}; \\
  \iint_{\Omega_1} M(x,t') {\rm d} x {\rm d} t' &\lesssim t^{-\kappa} \iint_{\Omega_1} (t')^\kappa M(x,t') {\rm d} x {\rm d} t' \lesssim t^{-\kappa};\\
  \iint_{\Omega_2}  M(x,t') {\rm d} x {\rm d} t' &\lesssim t^{-\kappa+\frac{1}{2}} r^{-1} \iint_{\Omega_2} \frac{(t')^{\kappa-\frac{1}{2}} \left(|t'-|x||+1\right) |u|^2}{|x|^2} {\rm d} x {\rm d} t'\\
  & \lesssim t^{-\kappa+\frac{1}{2}} r^{-1}.
 \end{align*}
 Here we use the integral estimates given in Proposition \ref{prop weighted Morawetz kappa large}. Inserting these upper bounds into \eqref{center estimate e minus} verifies the decay estimate of inward energy. Next we consider the outward energy.  For any $s\in (r/2,r)$, the energy flux for the region 
 \[
  \left\{(x,t'): |x|<t'-s, 3r/2<t'<t\right\}
 \]
 gives 
 \[
  \int_{|x|<t-s} e_+(x,t) {\rm d} x  = J_1 + J_2 + J_3 + J_4,
 \]
 with
 \begin{align*}
  J_1 = &\int_{|x|<3r/2-s} e_+(x,3r/2) {\rm d} x; & J_2 = &\pi  \int_{3r/2}^t |u(0,t')|^2 {\rm d} t'; & \\
  J_3 = &\int_{3r/2}^t \int_{|x|<t'-s} M(x,t'){\rm d} x {\rm d} t'; &  J_4 = &\frac{1}{\sqrt{2}} \int_{C^+(s;3r/2,t)} e'(x,t') {\rm d} S.&
 \end{align*}
 Again the term $J_2$ is ignored if $d\geq 4$. It follows Remark \ref{faster center decay} that 
 \[
  J_1 \leq \int_{|x|<r} e_+(x,3r/2) {\rm d} x \lesssim r^{-\kappa}; 
 \]
  By a comparison with \eqref{weighted Morawetz kappa large d4 M} and the integral estimate in Proposition \ref{prop weighted Morawetz kappa large} , we have
 \begin{align*}
  J_2 &\lesssim r^{-\kappa}; \\
  J_3 &\leq \int_{3r/2}^t \int_{|x|<t'-r/2} M(x,t') {\rm d} x {\rm d} t' \lesssim r^{-\kappa}.
 \end{align*}   
 Therefore we have 
  \begin{align*}
   \int_{|x|<t-r} e_+(x,t) {\rm d} x & \leq \frac{2}{r} \int_{r/2}^r \left(\int_{|x|<t-s} e_+(x,t) {\rm d} x\right) {\rm d} s\\
  & \lesssim \frac{2}{r}\int_{r/2}^r \left(\frac{1}{\sqrt{2}} \int_{C^+(s;3r/2,t)} e'(x,t') {\rm d} S\right) {\rm d} s + r^{-\kappa}\\
  & \lesssim \frac{1}{r} \int_{3r/2}^t \int_{t'-r<|x|<t'-r/2} e'(x,t') {\rm d} x {\rm d} t' + r^{-\kappa}\\
  & \lesssim \frac{1}{r} \int_{3r/2}^t \left(\frac{1}{rt'} \int_{t'-r<|x|<t'-r/2} (t'-|x|)|u(x,t')|^2 {\rm d} x\right) {\rm d} t' + r^{-\kappa}\\
  & \lesssim \frac{1}{r} \int_{3r/2}^t \frac{\ln t'}{rt'} {\rm d} t' + r^{-\kappa} \lesssim r^{-2} (\ln t)^2 + r^{-\kappa}. 
  \end{align*}
 Here we use \eqref{details in distribution of u2}.  In order to verify the point-wise estimate, we observe that the function $w (r',t) = (r')^\frac{d-1}{2} u(r',t)$ satisfies 
  \begin{align*}
   \int_0^{t-r} |w(r',t)|^2 {\rm d} r'  & = \frac{1}{\sigma_{d-1}} \int_{|x|<t-r} |u(x,t)|^2 {\rm d} x \lesssim r^{-1} \ln t;\\
   \int_0^{t-r} |\partial_{r'} w(r',t)|^2 {\rm d} r'  & = \frac{1}{\sigma_{d-1}} \int_{|x|<t-r} |\mathbf{L} u(x,t)|^2 {\rm d} x \lesssim \int_{|x|<t-r} [e_+ (x,t) + e_-(x,t)] {\rm d} x\\
   & \lesssim \max\{r^{-\kappa}, r^{-2} \ln^2 t\}. 
   \end{align*}
   In addition, we have $w(0,t) = 0$. Thus we may utilize the basic identity 
 \begin{align*}
  |w(R,t)|^2 & =2 \int_0^R w(r',t) \partial_{r'} w(r,'t) {\rm d} r' \\
  & \lesssim_1 \left(\int_0^{R} |w(r',t)|^2 {\rm d} r'\right)^{1/2} \left(\int_0^{R} |\partial_{r'} w(r',t)|^2 {\rm d} r'\right)^{1/2}, \quad \forall R>0; 
 \end{align*}
 to deduce 
   \begin{align*}
    \sup_{|x|<t-r} |x|^\frac{d-1}{2} |u(x,t)| & = \sup_{r' \in [0,t-r]} |w(r',t)|\\
     &\lesssim_1 \left(\int_0^{t-r} |w(r',t)|^2 {\rm d} r'\right)^{1/4} \left(\int_0^{t-r} |\partial_{r'} w(r',t)|^2 {\rm d} r'\right)^{1/4} \\
    & \lesssim \left(r^{-1} \ln t\right)^{1/4} \left(\max\{r^{-\kappa}, r^{-2} \ln^2 t\}\right)^{1/4}. 
   \end{align*}
\end{proof}

\begin{remark} \label{uniform boundedness of w} 
 Let $u$ be a radial free Coulomb wave with initial data $(u_0,u_1)\in (\mathcal{H}^1\cap L^2)\times (L^2 \cap \mathcal{H}^{-1})$. Then the argument above actually gives that 
 \[
  \sup_{(x,t)\in \Rm^d\times \Rm} |x|^\frac{d-1}{2} |u(x,t)| \lesssim_d \|(u_0,u_1)\|_{\mathcal{H}^1\times L^2}^{1/2} \|(u_0,u_1)\|_{L^2\times \mathcal{H}^{-1}}^{1/2}. 
 \]
\end{remark}

\begin{corollary}
 Let $u$ be a solution as in Proposition \ref{prop weighted Morawetz kappa large}. Assume that $\kappa \in (1,3/2)$ is a constant. Then the following estimate holds for $t\gg 1$. 
 \[
  \int_{\Rm^d} \left(|\mathbf{L}_+ u (x,t)|^2 + \lambda \frac{|u(x,t)|^2}{|x|^2}\right){\rm d} x \lesssim t^{-\kappa}. 
 \]
\end{corollary}
\begin{proof}
 The integral expression of inward energy gives 
 \begin{align*}
  \int_{\Rm^d} \left(\frac{1}{4}|\mathbf{L}_+ u (x,t)|^2+\frac{|u|^2}{4|x|} + \frac{\lambda |u|^2}{4|x|^2}\right) {\rm d} x  = \pi \int_t^\infty |u(0,t')|^2 {\rm d} t'+  \int_{t}^{\infty} \int_{\Rm^d} \left(\frac{|u|^2}{4|x|^2} + \frac{\lambda |u|^2}{2|x|^3}\right) {\rm d} x {\rm d} t'. 
 \end{align*}
 Again the integral of $|u(0,t')|^2$ is ignored if $d \geq 4$. We recall that fact that $u$ is supported in $\{(x,t'): |x|\leq |t'|+R_0\}$, the $L^2$ estimate $\|u(\cdot,t)\|_{L^2(\Rm^d)}^2 = E_1 + O(t^{-1/2})$, the decay estimates \eqref{details in distribution of u2}, \eqref{details 2 in distribution of u2}, \eqref{decay of E minus tr}, the global integral estimates \eqref{weighted Morawetz kappa large d4}, \eqref{weighted Morawetz kappa large} and obtain 
 \begin{align*}
  \frac{1}{4} \int_{\Rm^d} \frac{|u(x,t)|^2}{|x|} {\rm d} x &= \frac{1}{4} \int_{|x|<t/2} \frac{|u(x,t)|^2}{|x|} {\rm d} x + \frac{1}{4} \int_{|x|>t/2} \frac{|u(x,t)|^2}{|x|} {\rm d} x\\
  & = O(t^{-\kappa}) + \frac{1}{4} \int_{|x|>t/2} \frac{|u(x,t)|^2}{t} {\rm d} x + \frac{1}{4} \int_{|x|>t/2} \frac{(t-|x|)|u(x,t)|^2}{|x|t} {\rm d} x\\
  & = \frac{E_1}{4t} + O(t^{-\kappa});
 \end{align*}
 \begin{align*}
  \pi \int_t^\infty |u(0,t')|^2 {\rm d} t' = O(t^{-\kappa});
 \end{align*}
 \begin{align*}
  \int_{t}^{\infty} \int_{\Rm^d} \frac{|u|^2}{4|x|^2} {\rm d} x {\rm d} t' & = \int_{t}^{\infty} \int_{|x|<t'/3} \frac{|u|^2}{4|x|^2} {\rm d} x {\rm d} t' + \int_{t}^{\infty} \int_{|x|>t'/3} \frac{|u|^2}{4|x|^2} {\rm d} x {\rm d} t' \\
  & = O(t^{-\kappa}) + \int_{t}^{\infty} \int_{|x|>t'/3} \frac{|u|^2}{4t'^2} {\rm d} x {\rm d} t' + \int_{t}^{\infty} \int_{|x|>t'/3} \frac{(t'^2-|x|^2)|u|^2}{4|x|^2 t'^2} {\rm d} x {\rm d} t'\\
  & = O(t^{-\kappa}) + \int_{t}^{\infty} \frac{1}{4t'^2} (E_1+O(t'^{-1/2})) {\rm d} t' + O(t^{-\kappa-1/2})\\
  & = \frac{E_1}{4t} + O(t^{-\kappa});
 \end{align*}
 and 
 \begin{align*}
  \int_t^\infty \int_{\Rm^d} \frac{\lambda |u|^2}{2|x|^3} {\rm d} x {\rm d} t' & = \int_t^\infty \int_{|x|<t'/3} \frac{\lambda |u|^2}{2|x|^3} {\rm d} x {\rm d} t' + \int_t^\infty \int_{|x|>t'/3} \frac{\lambda |u|^2}{2|x|^3} {\rm d} x {\rm d} t'\\
  & = O(t^{-\kappa}) + O(t^{-2}) = O(t^{-\kappa}).
 \end{align*}
 A combination of these estimates with the energy flux identity finishes the proof. 
\end{proof}

\chapter{Asymptotic behaviour of linear waves}

In this chapter we shall prove a basic asymptotic behaviour of free Coulomb waves, namely, energy retraction. More precisely we show that given any forward light cone, almost all the energy is located inside the light cone as the time tends to infinity. We temporarily assume that the initial data are radial. The non-radial case will be dealt with in subsequent chapters by a decomposition of linear waves by spherical harmonic functions. 
\begin{lemma} \label{inner decay asymptotic}
 Let $u$ be a radial free Coulomb wave with a finite energy $E$ and $\eps$ be a positive number. Then there exists a constant $c_2$ such that 
 \[
  \limsup_{t\rightarrow +\infty} \int_{|x|<t-c_2 \ln t} e(x,t) {\rm d} x < \eps.
 \]
\end{lemma}
\begin{proof} 
 This is a direct consequence of the inward/outward energy theory. By linearity it suffices to consider smooth and compactly supported initial data since these initial data are dense in the energy space $\mathcal{H}^1 \times L^2$. Part (iii) of Proposition \ref{weighted Morawetz kappa 1} then gives 
\[
   \int_{|x|<t-c_2 \ln t} \left(\frac{1}{2}|\mathbf{L} u(x,t)|^2 + \frac{1}{2}|u_t(x,t)|^2 + e'(x,t)\right) {\rm d} x \lesssim c_2^{-1}, \qquad t\gg 1.  
\]
Namely 
\[
  \int_{|x|<t-c_2 \ln t} \left(e_-(x,t) + e_+ (x,t)\right) {\rm d} x \lesssim c_2^{-1}, \qquad t\gg 1.  
\]
A combination of this inequality with the energy inequality \eqref{inside inward outward total} yields
\[
  \int_{|x|<t-c_2 \ln t} e(x,t) {\rm d} x \lesssim c_2^{-1}, \qquad t\gg 1.  
\]
This immediately finishes the proof. 
\end{proof}
\begin{proposition} \label{prop energy retraction}
 Let $u$ be a radial free Coulomb wave with a finite energy $E$. Then given any function $\ell: \Rm^+ \rightarrow \Rm^+$ satisfying 
 \[
   \lim_{t\rightarrow +\infty} \frac{\ell (t)}{\ln t} = +\infty, 
 \]
 and $\eta \in \Rm$, we have 
 \[
  \lim_{t\rightarrow \pm\infty} \int_{|t|-\ell(|t|)<|x|<|t|-\eta} e(x,t) {\rm d}x = E.
 \]
\end{proposition}
\begin{proof}
In order to prove this proposition, by energy conservation law and the time symmetry it suffices to show that 
\begin{align}
 &\lim_{t\rightarrow \infty} \int_{|x|>t-\eta} e(x,t) {\rm d}x = 0; & & \forall \eta \in \Rm; \label{out vanishes}\\
 &\lim_{t\rightarrow \infty} \int_{|x|<t- \ell(t)} e(x,t) {\rm d}x = 0. & & \label{in vanishes}
\end{align}
The second limit \eqref{in vanishes} is a direct consequence of Lemma \ref{inner decay asymptotic}. 
The rest of this chapter is devoted to the proof of \eqref{out vanishes}. By a time translation it suffices to consider the case $\eta =0$. By smooth approximation and center/tail cut-off techniques, we assume that the initial data are smooth and supported in a annulus $B(0,R_0)\setminus B(0,r_0)$, without loss of generality. It immediately follows that 
\begin{align*}
 \|u(\cdot,t)\|_{\mathcal{H}^2}^2 + \|u_t(\cdot,t)\|_{\mathcal{H}^1}^2 &= C_1 \doteq \|u_0\|_{\mathcal{H}^2}^2 + \|u_1\|_{\mathcal{H}^1}^2 < +\infty;\\
 \|u(\cdot,t)\|_{\mathcal{H}^1}^2 + \|u_t(\cdot,t)\|_{L^2}^2 &= 2E < +\infty. 
\end{align*}
Thus 
\[
 \int_{\Rm^d}\left(|\Delta u|^2 + \frac{2|\nabla u|^2}{|x|} + \frac{|u|^2}{|x|^2}\right) {\rm d} x \leq C_1, \qquad \forall t\in \Rm. 
\]
Please note our assumption on the initial data also implies that  $u \in \mathcal{C}^2(\Rm \times ({\Rm^d \setminus \{0\}}))$, thanks to Proposition \ref{C2 continuity}. Since $u$ is a radial solution, we rewrite the integral of $|\Delta u|^2$ in the form of polar coordinates and obtain
\[
 \int_0^{\infty} \left|\frac{\partial^2}{\partial r^2} (r^{\frac{d-1}{2}} u) -\lambda r^{\frac{d-5}{2}} u\right|^2 {\rm d} r \lesssim_1 C_1. 
\]
Combining this with the Hardy inequality
\[
 \int_0^\infty r^{d-3} |u(r)|^2 {\rm d} r \lesssim_d \|u\|_{\dot{H}^1}^2 \leq \|u\|_{\mathcal{H}^1}^2 \lesssim_1 E; 
\]
we obtain the uniform $L^2$ boundedness of the second derivative of $w(r,t) = r^{\frac{d-1}{2}} u(r,t)$
\begin{equation} \label{bound of wrr} 
  \sup_{t\in \Rm} \int_1^\infty |w_{rr}(r,t)|^2 {\rm d} r < +\infty. 
\end{equation}
In addition, the energy conservation law gives that 
\begin{equation} \label{bound of wr} 
  \sup_{t\in \Rm} \int_0^\infty |w_{r}(r,t)|^2 {\rm d} r \lesssim_1 E. 
\end{equation}
Furthermore, our assumption on the support of initial data and the finite speed of propagation yields 
\[
 w(r,t) = 0, \quad \forall r\geq |t| +R_0.
\]
Now we consider two families of functions defined on the interval $s\in [0,R_0]$:
\begin{align*}
 &\{w(s+t,t): t>1\};& &\{w_r(s+t,t): t>1\}.
\end{align*}
A combination of \eqref{bound of wrr}, \eqref{bound of wr} and the support of function $w$ implies that both families of functions are equicontinuous and uniformly bounded. Thus we may extract a sequence of time $t_k \rightarrow +\infty$, so that the following uniform convergence holds 
\begin{align} \label{uniform convergence of f}
 &w(s+t_k,t_k) \rightrightarrows f(s);& &w_r(s+t_k,t_k)\rightrightarrows f'(s).
\end{align}
Here $f(s)$ is a continuously differentiable function on $[0,R_0]$. We claim that the following limit holds (without extracting a sequence of time)
\begin{equation} \label{continuous limit} 
 \lim_{t\rightarrow \infty} w_r(s+t,t) = f'(s), \qquad \forall s\in [0,R_0].
\end{equation} 
If this were false, then we might find an $s_0\in [0,R_0]$ so that 
\[
 \limsup_{t\rightarrow \infty} w_r(s+t,t) > \liminf_{t\rightarrow \infty} w_r(s+t,t). 
\]
For any number $A$ between these two limits, by continuity of $w_r$ we may extract a sequence of time $t'_k$, so that 
\[
 \lim_{k\rightarrow \infty} w_r(s_0+t'_k, t'_k) = A. 
\]
Without loss of generality, we choose a real number $A$ such that $A\neq f'(s_0), -f'(s_0)$. By extracting a subsequence of $\{t'_k\}$ if necessary, we may assume that the following uniform convergence holds as well by the equicontinuous and uniformly bounded properties of $w, w_r$. 
\begin{align} \label{uniform convergence of g}
 &w(s+t'_k,t'_k) \rightrightarrows g(s);& &w_r(s+t'_k,t'_k)\rightrightarrows g'(s).
\end{align}
Please note that $f(s), f'(s), g(s), g'(s)$ are still continuous and $g'(s_0) = A \neq f'(s_0), -f'(s_0)$. Next we recall the energy flux formula of inward/outward energies and obtain that the following limit exists as $t\rightarrow +\infty$: 
\begin{align*}
 \lim_{t\rightarrow +\infty} \int_{r_1+t<|x|<R_1+t} e_+ (x,t) {\rm d} t, \qquad  \forall r_1<R_1. 
\end{align*}
Rewriting this in term of $w$, we obtain the existence of 
\[
 \lim_{t\rightarrow +\infty} \int_{r_1+t}^{R_1+t} \left(\frac{1}{4}|(w_r-w_t)(r,t)|^2+\frac{\lambda}{4}\frac{|w(r,t)|^2}{r^2} + \frac{1}{4} \frac{|w(r,t)|^2}{r}\right) {\rm d} r, \quad \forall r_1<R_1. 
\]
Similarly we have the convergence of inward energy 
\begin{equation} \label{limit of inward equi}
 \lim_{t\rightarrow +\infty} \int_{r_1+t}^{R_1+t} \left(\frac{1}{4}|(w_r+w_t)(r,t)|^2+\frac{\lambda}{4}\frac{|w(r,t)|^2}{r^2} + \frac{1}{4} \frac{|w(r,t)|^2}{r}\right) {\rm d} r = 0, \quad \forall r_1<R_1. 
\end{equation} 
A combination of these two limits gives the existence of the following limit (without extracting a sequence)
\[
  \lim_{t\rightarrow +\infty} \int_{r_1+t}^{R_1+t} |w_r(r,t)|^2 {\rm d} r. 
\]
By the uniform convergence \eqref{uniform convergence of f} and \eqref{uniform convergence of g}, this immediately gives 
\[
 \int_{r_1}^{R_1} |f'(s)|^2 {\rm d} s = \int_{r_1}^{R_1} |g'(s)|^2 {\rm d} s, \qquad 0\leq r_1 < R_1 \leq R_0. 
\]
By the continuity of $f'$ and $g'$, we obtain 
\[
 |f'(s)| = |g'(s)|, \qquad s\in [0,R_0].
\]
This immediately gives a contraction with $g'(s_0) \neq f'(s_0), -f'(s_0)$ thus verifies \eqref{continuous limit}. Since $w_r$ is equicontinuous and $w(R_0+t,t)\equiv 0$, we actually have the uniform convergence for $s \in [0,R_0]$
\begin{align*}
  &w(s+t,t) \rightrightarrows f(s);& &w_r(s+t,t)\rightrightarrows f'(s).
\end{align*}
Combining this with \eqref{limit of inward equi}, we have
\begin{equation} \label{limit of L2 equi}
 w_t(s+t,t) - w_r(s+t,t) \rightarrow -2f'(s) \quad {\rm in}\; L^2([0,R_0]). 
\end{equation}
Now we claim that $f'(s) = f(s) = 0$ for all $s\in [0,R_0]$. If this were false, then there would exist $s_0\in [0,R_0]$ so that $f(s_0)\neq 0$. Without loss of generality, we assume $f(s_0)>0$. It follows the continuity that there exists an interval $J \subset [0,R_0]$ and $C>0$, so that 
\[
 f(s) > 2C, \qquad s\in J. 
\] 
By the uniform convergence there exists a large time $t_0>1$ so that 
\[
 w(s+t,t) > C, \qquad s\in J,\; t>t_0. 
\]
We observe that 
\[
 \frac{\partial}{\partial t} \left[w_t(s+t,t) - w_r(s+t,t)\right] = (w_{tt} - w_{rr})(s+t,t) = - \frac{w(s+t,t)}{s+t} - \lambda \frac{w(s+t,t)}{(s+t)^2}. 
\]
Therefore 
\begin{align*}
 &\left[w_t(s+t_1,t_1) - w_r(s+t_1,t_1)\right] - \left[w_t(s+t_0,t_0) - w_r(s+t_0,t_0)\right]\\
 & \qquad = - \int_{t_0}^{t_1} \left(\frac{w(s+t,t)}{s+t} + \lambda \frac{w(s+t,t)}{(s+t)^2}\right){\rm d} t \leq - C \ln \frac{s+t_1}{s+t_0}, \quad s\in J.
\end{align*}
This gives a contradiction with \eqref{limit of L2 equi} thus we must have $f(s)=f'(s) = 0$. Thus we have 
\[
 \lim_{t\rightarrow +\infty} \int_t^{t+R_0} |w_r(r,t)|^2 {\rm d} r = 0. 
\]
Combining this with the decay of inward energy and the inequality
\[
 \sup_{r\in [t,t+R_0]} |w(r,t)| \leq R_0^{1/2} \left(\int_t^{t+R_0} |w_r(r,t)|^2 {\rm d} r \right)^{1/2}, 
\]
we have 
\[
  \lim_{t\rightarrow +\infty} \int_{t}^{t+R_0} \left(|w_r(r,t)|^2+|w_t(r,t)|^2 + |w(r,t)|^2\right) {\rm d} r = 0.
\]
Finally we observe that 
\[
 r^{\frac{d-1}{2}} u_r = w_r - \frac{d-1}{2} \frac{w}{r}
\]
and obtain 
\[
 \lim_{t\rightarrow +\infty} \int_{t<|x|<t+R_0} \left(|\nabla u(x,t)|^2 + |u_t(x,t)|^2 + |u(x,t)|^2 \right) {\rm d} x = 0.
\]
This completes the proof in the radial case. 
\end{proof} 

\chapter{An isometry between the scattering profiles}

In this section we give an isometry between the scattering profiles of one-dimensional Klein-Gordon equation and those of Coulomb wave equation in the radial case. In the major part of this chapter we consider the Coulomb wave equation in 3-dimensional space for simplicity. In the final section we explain why the same argument works in higher dimensions as well. We consider the Hilbert space $\mathcal{K}$ of all finite-energy solutions to the linear homogeneous Klein-Gordon equation  
\begin{equation} \label{klein-gordon}
 v_{\tau \tau}-v_{yy}+2v = 0
\end{equation}
equipped with the norm 
\[
 \|v\|_{\mathcal{K}}^2 = \int_{\mathbb{R}} \left(|v_\tau|^2 + |v_y|^2 + 2|v|^2\right) {\rm d}y. 
\]
Here $\tau$ represents the time. Similarly we consider the space $\mathcal{C}$ of all radial solutions $u$ to the wave equation with a Coulomb potential 
\[
 \partial_t^2 u - \Delta u + \frac{u}{|x|} = 0, \qquad (x,t)\in \Rm^3 \times \Rm,
\]
with finite norm 
\begin{equation} \label{norm of C}
 \|u\|_{\mathcal{C}}^2 = \|\mathbf{H}^{1/2} u\|_{L^2}^2 + \frac{1}{2}\|u\|_{L^2}^2 + \|u_t\|_{L^2}^2 + \frac{1}{2}\|\mathbf{H}^{-1/2} u_t\|_{L^2}^2.
\end{equation}
Please note that the right hand of \eqref{norm of C} is a conserved quantity independent of the time $t$.  

\section{The idea}

The basic idea is to apply the following geometric transformation 
\[
 (y,\tau) = \left(\frac{r-t+\ln (t+r)}{2}, \frac{t-r+\ln(t+r)}{2}\right)
\]
and its inverse 
\[ 
 (r,t) = \left(\frac{e^{y+\tau}+y-\tau}{2}, \frac{e^{y+\tau}-y+\tau}{2}\right). 
\]
Let $u(x,t)$ be a radial free Coulomb wave in 3-dimensional space. It follows that $w(r,t) = r u(r,t)$ satisfies the Coulomb wave equation 
\[
 w_{tt} - w_{rr} + \frac{w}{r} = 0. 
\]
A direct calculation shows that $v(y,\tau) =  w(r,t) = w\left(\frac{e^{y+\tau}+y-\tau}{2}, \frac{e^{y+\tau}-y+\tau}{2}\right)$ solves the equation 
\[
 v_{\tau \tau} - v_{yy} + 2 v = \frac{r-t}{r} w(r,t). 
\]
This is NOT an exact solution to the linear free Klein-Gordon equation but the error term can be ignored if we consider the asymptotic behaviour of solutions. In fact, the inward/outward energy theory implies that almost all energy of $u$ or $w$ concentrates around the light cone $r=t$ as $t$ tends to infinity, where the factor $\frac{r-t}{t}$ is small. Of course, when $y$ is close to $-\tau$, the corresponding value of $r$ is negative, which is not reasonable. Thus we need to cut the insignificant part off by applying suitable cut-off techniques. We shall give the details of our argument in the subsequent sections. 

\section{Construction of operator} \label{sec: construction of operator}

We start by defining a linear bounded operator $\mathbf{T}$ from $\mathcal{K}$ to $\mathcal{C}$. For simplicity we first define it on a sub-space of $\mathcal{K}$ in this section and then extend it to the whole space $\mathcal{K}$ by continuity in subsequent sections. We consider the sub-space $\mathcal{K}_0$ consisting of $v\in \mathcal{K}$ satisfying  
\begin{itemize}
 \item The solution $v(y,\tau)$ is smooth;  
 \item There exists a time $\tau_0>0$ so that the data $(v(\cdot,\tau_0), v_{\tau}(\cdot,\tau_0))$ is supported in the interval $[-\tau_0,\tau_0]$. 
\end{itemize}
It is not difficult to check that $\mathcal{K}_0$ is a dense sub-space of $\mathcal{K}$, as shown in Corollary \ref{dense subset K0} in the appendix. 
Now we define $\mathbf{T} v$ for $v\in \mathcal{K}_0$. We first fix a smooth, decreasing cut-off function $\rho: \Rm \rightarrow [0,1]$ such that 
\[
 \rho(s) = \left\{\begin{array}{ll} 0, & s\in [2,+\infty); \\ 1, & s\in (-\infty,1/2];\end{array}\right.
\]
and define a smooth function 
\begin{equation} \label{def of w t}
 w(r,t) = \rho\left(\frac{t-r}{t^{1/2}}\right) v\left(\frac{r-t+\ln (t+r)}{2}, \frac{t-r+\ln(t+r)}{2}\right), \quad (r,t)\in \Rm^+ \times [T_0,+\infty).
\end{equation}
Here $T_0$ is a large constant. This function $w$ always vanishes unless $r\in (t-2t^{1/2},t)$. A direct calculation shows that 
\begin{align*}
 w_r & = -\frac{1}{t^{1/2}}\rho' v + \rho v_y \left(\frac{1}{2}+\frac{1}{2(t+r)}\right) +\rho v_\tau \left(-\frac{1}{2}+\frac{1}{2(t+r)}\right);\\
 w_t & = \frac{t+r}{2t^{3/2}}\rho' v + \rho v_y \left(-\frac{1}{2}+\frac{1}{2(t+r)}\right) + \rho v_\tau \left(\frac{1}{2}+\frac{1}{2(t+r)}\right);
\end{align*}
and 
\begin{align*}
 w_{rr} & = \frac{1}{t}\rho'' v - \frac{2}{t^{1/2}} \rho' \left[v_y \left(\frac{1}{2}+\frac{1}{2(t+r)}\right) + v_\tau \left(-\frac{1}{2}+\frac{1}{2(t+r)}\right)\right]\\
 & \qquad + \rho\left[v_{yy}\left(\frac{1}{2}+\frac{1}{2(t+r)}\right)^2 + 2v_{\tau y}\left(\frac{1}{4(t+r)^2}-\frac{1}{4}\right) + v_{\tau \tau} \left(-\frac{1}{2}+\frac{1}{2(t+r)}\right)^2 \right]\\
 & \qquad \qquad + \rho (v_y+v_\tau)\cdot \frac{-1}{2(t+r)^2};\\
 w_{tt} & =  -\frac{t+3r}{4t^{5/2}}\rho' v +\frac{(t+r)^2}{4t^3}\rho'' v  +\frac{t+r}{t^{3/2}} \rho' \left[v_y \left(-\frac{1}{2}+\frac{1}{2(t+r)}\right) + v_\tau \left(\frac{1}{2}+\frac{1}{2(t+r)}\right)\right]\\
 & \qquad + \rho\left[v_{yy}\left(-\frac{1}{2}+\frac{1}{2(t+r)}\right)^2 + 2v_{\tau y}\left(\frac{1}{4(t+r)^2}-\frac{1}{4}\right) + v_{\tau \tau} \left(\frac{1}{2}+\frac{1}{2(t+r)}\right)^2 \right]\\
 & \qquad \qquad + \rho (v_y+v_\tau)\cdot \frac{-1}{2(t+r)^2};
\end{align*}
Thus we have 
\begin{align*}
 w_{tt} - w_{rr} &=  -\frac{t+3r}{4t^{5/2}}\rho' v + \frac{(r-t)(r+3t)}{4t^3}\rho'' v + \frac{t-r}{t^{3/2}} \rho' \left(\frac{1}{2}v_y - \frac{1}{2} v_\tau\right)\\
 & \qquad + \frac{3t+r}{t^{3/2}}\rho' \cdot \frac{1}{2(t+r)}(v_y + v_\tau)  + \rho \cdot \frac{v_{\tau \tau} - v_{yy}}{t+r}.
\end{align*}
Since $v$ satisfies the equation \eqref{klein-gordon}, $w$ satisfies the approximated linear wave equation with a Coulomb potential 
\begin{equation} \label{approximated equation w}
 w_{tt} - w_{rr} + \frac{w}{r} = f(r,t). 
\end{equation}
Here $f(r,t)$ is defined by 
\begin{align*}
 f(r,t) &=  -\frac{t+3r}{4t^{5/2}}\rho' v + \frac{(r-t)(r+3t)}{4t^3}\rho'' v + \frac{t-r}{2t^{3/2}} \rho' \left(v_y -  v_\tau\right) \\
& \qquad + \frac{(3t+r)\rho' (v_y + v_\tau) }{2t^{3/2}(t+r)} +\rho \frac{(t-r)v}{r(t+r)}\\
 &= J_1 + J_2 + J_3 + J_4 + J_5.
\end{align*}
If we define $u(x,t) =|x|^{-1} w(|x|,t)$ be a radial function, then $u$ satisfies the approximated linear wave equation with a Coulomb potential in $\Rm^3$. 
\begin{equation} \label{approximated equation u}
 u_{tt} - \Delta u + \frac{u}{|x|} = \frac{f(|x|,t)}{|x|}. 
\end{equation}
We claim that 
\begin{equation}  \label{estimate of f}
 \left\|(r^{1/2}+1)f(r,t)\right\|_{L_t^1([T_0,+\infty); L_r^2(\Rm^+))}<+\infty,
\end{equation} 
It immediately follows that 
\[
 |x|^{-1} f(|x|,t) \in L_t^1([T_0,+\infty);L^2 \cap \mathcal{H}^{-1}(\Rm^3)).
\]
Here we use the embedding $\|g\|_{\mathcal{H}^{-1}(\Rm^3)} \leq \||x|^{1/2} g\|_{L^2(\Rm^3)}$, which is a direct consequence of 
\[
 \left|\int_{\Rm^3} g(x) h(x) {\rm d} x\right| \leq \left\||x|^{1/2} g(x)\right\|_{L^2} \left\||x|^{-1/2} h(x)\right\|_{L^2} \leq \left\||x|^{1/2} g(x)\right\|_{L^2} \|h\|_{\mathcal{H}^1(\Rm^3)}. 
\]
Thus $u$ scatters in the positive time direction. More precisely, if we define 
\[ 
 \vec{u}_+ = \vec{\mathbf S}_{\mathcal C} (-T_0)\vec{u}(\cdot,T_0)+\int_{T_0}^\infty \vec{\mathbf S}_{\mathcal C} (-t) (0, |x|^{-1}f(|x|,t)) {\rm d} t  \in \vec{\mathcal{H}}^1 \cap \vec{\mathcal{H}}^0,
\]
then a straight-forward calculation shows that 
\[
 \lim_{t\rightarrow +\infty} \left\|\vec{u}(\cdot,t) - \vec{\mathbf{S}}_{\mathcal C} (t) \vec{u}_+\right\|_{\vec{\mathcal{H}}^1 \cap \vec{\mathcal{H}}^0}  = 0.
\]
We then define $\mathbf{T} v$ to be the scattering target of $u$, i.e.
\begin{equation}
 \mathbf{T} v = \mathbf{S}_{\mathcal C} (t) \vec{u}_+ \in \mathcal{C}. 
\end{equation}
Now we verify our claim \eqref{estimate of f}. Since $f$ is smooth and vanishes unless $t - 2t^{1/2}<r<t$, it suffices to show 
\begin{equation} \label{estimate of f 1} 
 \|f(\cdot,t)\|_{L^2(\Rm^+)} \lesssim t^{-3/2-\varepsilon} 
\end{equation}
for some small constant $\varepsilon > 0$ and all sufficiently large time $t$. 

%\begin{lemma} \label{decay of u}
% Assume that $v \in {\mathcal K}_0$ and let $w$, $u$ be approximated solutions to the wave equation with a Coulomb potential as defined above. If $(r,t)$ satisfies 
% \begin{align*}
%  &t-2t^{1/2} < r < t-t^{1/4};& &t\gg 1;&
% \end{align*}
% then we have 
% \begin{itemize} 
%  \item The following decay estimate holds for any fixed $N>0$
%   \[
%    |v(y,\tau)|+ |v_y(y,\tau)| + |v_\tau(y,\tau)| \lesssim t^{-N};
%   \]
%   where 
%   \[
%    (y,\tau) = \left(\frac{r-t+\ln (t+r)}{2}, \frac{t-r+\ln(t+r)}{2}\right). 
%   \]
%  \item We also have the following decay estimate for any fixed constant $N>0$;
%   \begin{align*}
%     &|w(r,t)| + |w_r(r,t)| + |w_t(r,t)| \lesssim t^{-N};& &|u(r,t)|+|u_r(r,t)|+|u_t(r,t)|\lesssim t^{-N}.&
%   \end{align*}
% \end{itemize}
%\end{lemma}

Recall that $f(r,t) = J_1 + J_2 + \cdots + J_5$. We first consider $J_k$ for $1\leq k\leq 4$. The presence of $\rho'$ or $\rho''$ implies that 
\[
 \hbox{Supp} J_k \subset \left\{(r,t): t\geq T_0, t-2t^{1/2}\leq r\leq t-(1/2)t^{1/2} \right\}. 
\]
Thus if $(r,t)\in \hbox{Supp} J_k$, then the corresponding $(y,\tau)$ satisfies  
\begin{align*}
& \tau = \frac{t-r+\ln(t+r)}{2} > \frac{1}{4} t^{1/2}; & & |y+\tau| = \ln (t+r) < \ln(2t).&
\end{align*}
Therefore we have 
\[
 |y| > |\tau| - \ln(2t) > |\tau| - |\tau|^{1/3}. 
\]
It immediately follows Lemma \ref{fast decay outside parabola} that if $(r,t)\in \hbox{Supp} J_k$ and $t$ is sufficiently large, then the inequality 
\[
 |v|+ |v_y| + |v_\tau| \lesssim |\tau|^{-N} \lesssim t^{-N/2} 
\]
holds for any given $N>0$. The arbitrarity of $N$ then immediately gives the following estimate for large time $t$:
\[
 \|J_k(\cdot,t)\|_{L^2(\Rm^+)} \lesssim t^{-2}, \qquad k\in \{1,2,3,4\}. 
\]
Next we consider the major term $J_5$.  We split the interval $(t-2t^{1/2},t)$ into two sub-intervals. If $r\in (t-2t^{1/2},t-t^{1/4})$, then we have 
\begin{align*}
& \tau = \frac{t-r+\ln(t+r)}{2} > \frac{1}{2} t^{1/4}; & & |y+\tau| = \ln (t+r) < \ln(2t).&
\end{align*}
Again we have
\[
 |y| > |\tau| - \ln(2t) > |\tau| - |\tau|^{1/3}. 
\]
We apply Lemma \ref{fast decay outside parabola} for another time and deduce for any $N>0$ that 
\[
 |v|  \lesssim |\tau|^{-N} \lesssim t^{-N/4}, \qquad t\gg 1. 
\]
Thus we have 
\[
 \left(\int_{0}^{t-t^{1/4}} |J_5(r,t)|^2 {\rm d} r \right)^{1/2} \lesssim t^{-2}. 
\]
On the other hand, if $r\in (t-t^{1/4},t)$, then we have 
\[
 |J_5(r,t)| \leq \frac{t^{1/4}}{t^2}, \qquad \Rightarrow \qquad \left(\int_{t-t^{1/4}}^t |J_5(r,t)|^2 {\rm d} r \right)^{1/2} \lesssim t^{-13/8}. 
\]
In summary we have 
\begin{equation} \label{estimate of f detail}
 \|J_5(\cdot,t)\|_{L^2(\Rm^+)} \lesssim t^{-13/8} \qquad \Rightarrow \qquad \|f(\cdot,t)\|_{L^2(\Rm^+)} \lesssim t^{-3/2-1/8}.  
\end{equation}
This verifies \eqref{estimate of f 1} thus \eqref{estimate of f}. 
\begin{remark} \label{decay of w}
 In the argument above, we actually show that if $t$ is sufficiently large and $0<r<t-(1/2) t^{1/2}$, then 
 \[
  |w(r,t)|+|w_r(r,t)|+|w_t(r,t)| \lesssim t^{-N}, \qquad \forall N>0. 
 \]
\end{remark}

\section{Preservation of norm} \label{sec: preservation of norm}

In this section we show that the operator $\mathbf{T}$ defined on $\mathcal{K}_0$ preserve the norms, i.e. 
\[
 \|\mathbf{T} v\|_{\mathcal{C}} = 4\pi \|v\|_{\mathcal K}, \qquad v\in \mathcal{K}_0,
\]
thus it extends to a norm-preserving operator from $\mathcal{K}$ to $\mathcal{C}$. First of all, by Corollary \ref{balance for L2} we have 
\begin{align}
 \|\mathbf{T} v\|_{\mathcal C}^2 &= \lim_{t\rightarrow +\infty} \left(\|\mathbf{T} v(\cdot,t)\|_{\mathcal{H}^1}^2 + \|\partial_t \mathbf{T} v(\cdot,t)\|_{L^2}^2 + \|\mathbf{T} v(\cdot,t)\|_{L^2}^2 \right)\nonumber \\
 & = \lim_{t\rightarrow +\infty} \left(\|u(\cdot,t)\|_{\mathcal{H}^1}^2 + \|u_t (\cdot,t)\|_{L^2}^2 + \|u(\cdot,t)\|_{L^2}^2 \right). \label{norm of Tv}
\end{align}
A basic calculation shows 
\begin{align*}
 J(t) \doteq \|u(\cdot,t)\|_{\mathcal{H}^1}^2 + \|u_t (\cdot,t)\|_{L^2}^2 + \|u(\cdot,t)\|_{L^2}^2  &= \int_{\Rm^3} \left(|\nabla u|^2 + \frac{|u|^2}{|x|}+|u_t|^2 + |u|^2\right) {\rm d} x.
\end{align*}
By the fact $\|u(\cdot, t)-(\mathbf{T} v) (\cdot,t)\|_{L^2(\Rm^3)} \rightarrow 0$ and Corollary \ref{balance for L2}, we have
\begin{equation} \label{limit of u2 R3}
 \lim_{t\rightarrow +\infty} \|u(\cdot,t)\|_{L^2} = \lim_{t\rightarrow +\infty} \|(\mathbf{T} v)(\cdot,t)\|_{L^2} = \frac{1}{2}\|(\mathbf{T}v, \partial_t \mathbf{T}v)\|_{L^2 \times \mathcal{H}^{-1}} \lesssim_1 \|\mathbf{T} v\|_{\mathcal{C}} < +\infty.
\end{equation}
A combination of this limit  with the support of $u$ gives 
\begin{equation} \label{limit of Vu R3}
 \lim_{t\rightarrow +\infty} \int_{\Rm^3} \frac{|u(x,t)|^2}{|x|} {\rm d} x = 0.
\end{equation}
It immediately follows that 
\begin{align*}
 \lim_{t\rightarrow +\infty}  J(t)
 & = \lim_{t\rightarrow +\infty} \int_{\Rm^3} \left(|\nabla u|^2 +|u_t|^2 + |u|^2\right) {\rm d} x \\
 & = 4\pi \lim_{t\rightarrow +\infty} \int_0^\infty \left(|w_r|^2 + |w_t|^2 + |w|^2\right) {\rm d} r\\
 & = 2\pi \lim_{t\rightarrow +\infty} \int_0^\infty \left(|w_r+w_t|^2 + |w_r-w_t|^2 + 2|w|^2\right) {\rm d} r\\
 & = 2\pi \lim_{t\rightarrow +\infty} \int_0^t \left(|w_r+w_t|^2 + |w_r-w_t|^2 + 2|w|^2\right) {\rm d} r.
\end{align*}
Here we use the support of $w$. Now we recall $w = \rho v$ and observe 
\begin{align*}
 w_r + w_t & = \frac{r-t}{2t^{3/2}} \rho' v + \frac{1}{t+r}\rho (v_y + v_\tau);\\
 w_r - w_t & = - \frac{3t+r}{2t^{3/2}}\rho' v + \rho(v_y - v_\tau).
\end{align*}
In view of Remark \ref{decay of w}, we may ignore the integral over the interval $(0,t-(1/2)t^{1/2})$ and obtain 
\begin{align}
  \lim_{t\rightarrow +\infty} J(t) &= 2\pi \lim_{t\rightarrow +\infty} \int_{t-\frac{1}{2}t^{1/2}}^t \left(|v_y-v_\tau|^2 + \frac{|v_y+v_\tau|^2}{(t+r)^2} + 2 |v|^2\right) {\rm d} r. \label{limit of Jt}
\end{align}

 \begin{figure}[h]
 \centering
 \includegraphics[scale=1.1]{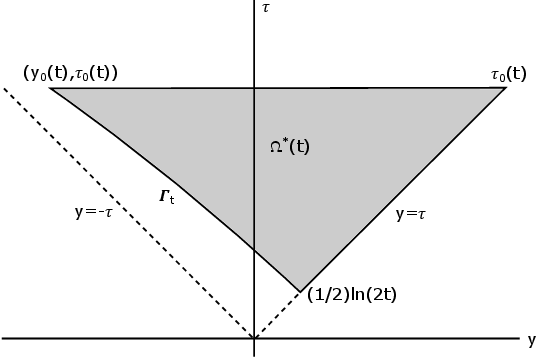}
 \caption{Illustration of region $\Omega^*(t)$} \label{figure gammat}
\end{figure}

In order to evaluate this limit, we switch to the $y$-$\tau$ space. We start by observing the divergence identity for energy of the Klein-Gordon equation
\[
 \frac{\partial}{\partial \tau} \left(|v_\tau|^2 + |v_y|^2 + 2|v|^2\right) + \frac{\partial}{\partial y} \left(-2 v_y v_\tau\right) = 2v_\tau(v_{\tau \tau}-v_{yy}+2v) = 0,
\]
We then apply Green's formula on the line integral 
\[
 \int 2 v_y v_\tau {\rm d} \tau + \left(|v_\tau|^2 + |v_y|^2 + 2|v|^2\right) {\rm d} y
\]
 over the boundary of the region $\Omega^\ast (t)$ enclosed by the (oriented) curve 
\[
 \Gamma_t = \left\{\left(\frac{r-t+\ln (t+r)}{2}, \frac{t-r+\ln(t+r)}{2}\right): t-\frac{1}{2}t^{1/2} < r < t\right\} \subset \Rm_y \times \Rm_\tau
\]
and the straight lines
\begin{align*}
 &\tau = \tau_0(t) = \frac{1}{4} t^{1/2} + \frac{1}{2}\ln\left(2t - \frac{1}{2} t^{1/2}\right);& &\tau = y;&
\end{align*}
as illustrated in figure \ref{figure gammat} to deduce the energy flux formula 
\begin{align}
 \int_{y_0(t)}^{\tau_0(t)} & \left(|v_\tau(y,\tau_0(t))|^2 + |v_y(y,\tau_0(t))|^2 + 2|v(y,\tau_0(t))|^2\right) {\rm d} y\nonumber\\
  & = \int_{\Gamma_t} 2 v_y v_\tau {\rm d} \tau + \left(|v_\tau|^2 + |v_y|^2 + 2|v|^2\right) {\rm d} y. \label{energy flux formula v}
\end{align}
where 
\[
 y_0(t) = -\frac{1}{4} t^{1/2} + \frac{1}{2}\ln\left(2t - \frac{1}{2} t^{1/2}\right).
\]
By Lemma \ref{fast decay outside parabola}, we obtain 
\begin{equation} \label{limit of energy tau}
 \lim_{t\rightarrow +\infty} \int_{y_0(t)}^{\tau_0(t)}  \left(|v_\tau(y,\tau_0(t))|^2 + |v_y(y,\tau_0(t))|^2 + 2|v(y,\tau_0(t))|^2\right) {\rm d} y = \|v\|_{\mathcal K}^2.
\end{equation}
By the parameter representation of $\Gamma_t$ given above,  the energy flux 
\[
 {\rm Flux} (v, \Gamma_t) \doteq \int_{\Gamma_t}  2 v_y v_\tau {\rm d} \tau + \left(|v_\tau|^2 + |v_y|^2 + 2|v|^2\right) {\rm d} y
\]
can be given by 
\begin{align}
 {\rm Flux} (v, \Gamma_t) & = \int_{t-\frac{1}{2}t^{1/2}}^t \left[ 2v_y v_\tau \cdot \frac{1}{2}\left(-1+\frac{1}{t+r}\right) + \left(|v_\tau|^2 + |v_y|^2 + 2|v|^2\right) \cdot \frac{1}{2}\left(1+\frac{1}{t+r}\right)\right] {\rm d} r\nonumber \\
  & = \frac{1}{2}\int_{t-\frac{1}{2}t^{1/2}}^t \left[|v_y-v_\tau|^2+\frac{1}{t+r}|v_y+v_\tau|^2+2\left(1+\frac{1}{t+r}\right)|v|^2\right] {\rm d} r. \label{expression of flux}
\end{align}
It immediately follows \eqref{energy flux formula v} and \eqref{limit of energy tau} that ${\rm Flux}(v,\Gamma_t)$ satisfies 
\begin{align} \label{limit of energy flux v} 
 &{\rm Flux}(v,\Gamma_t) \leq \|v\|_{\mathcal K}^2; & &\lim_{t\rightarrow +\infty} {\rm Flux}(v, \Gamma_t) = \|v\|_{\mathcal K}^2. &
\end{align}
Now we compare the expressions of $J(t)$ and ${\rm Flux}(v,\Gamma_t)$ given in \eqref{limit of Jt} and \eqref{expression of flux}. It turns out that the coefficients of $|v_y-v_\tau|^2$ and $|v|^2$ match nicely as $t\rightarrow +\infty$. However, the coefficient of $v_y+v_\tau$ is much different. We claim that this term can be ignored because 
\begin{equation} \label{ignore the difference}
 \int_{t-\frac{1}{2}t^{1/2}}^t \frac{1}{t+r}|v_y+v_\tau|^2 {\rm d} r =  \int_{t-\frac{1}{2}t^{1/2}}^t (t+r)|w_r+w_t|^2 {\rm d} r
\end{equation}
actually converges to zero as $t$ tends to infinity.  In order to verify this we start by integrating the divergence formula 
\[
 \frac{\partial}{\partial t} \left(|w_t + w_r|^2 + \frac{|w|^2}{r}\right) + \frac{\partial}{\partial r}\left(\frac{|w|^2}{r}-|w_t+w_r|^2\right) = 2(w_t+w_r)\left(w_{tt}-w_{rr}+\frac{w}{r}\right)- \frac{|w|^2}{r^2},
\]
in the region $\Rm^+ \times [t_0,T_0]$ and obtain 
\begin{align*}
 \left.\int_0^\infty \left(|w_t + w_r|^2 + \frac{|w|^2}{r}\right) {\rm d} r\right|_{t=t_0}^{T_0}  = 2 \int_{t_0}^{T_0} \int_0^\infty (w_t+w_r) f(r,t){\rm d}r {\rm d}t - \int_{t_0}^{T_0}  \int_0^\infty \frac{|w|^2}{r^2} {\rm d}r {\rm d}t.
\end{align*}
We first combine \eqref{expression of flux}, \eqref{limit of energy flux v} and \eqref{ignore the difference} to deduce 
\[
  \int_{t-\frac{1}{2}t^{1/2}}^t |w_r+w_t|^2 {\rm d} r \lesssim t^{-1}. 
\]
We combine this with the decay estimates given in Remark \ref{decay of w} to deduce  
\begin{equation} \label{decay of wrt whole}
 \|(w_r+w_t)(\cdot,t)\|_{L^2(\Rm^+)} \lesssim t^{-1/2}, \qquad t\gg 1. 
\end{equation}
We then recall \eqref{estimate of f detail} and conclude 
\[
 \left|\int_{t_0}^{T_0} \int_0^\infty (w_t+w_r) f(r,t){\rm d}r {\rm d}t\right| \lesssim t_0^{-9/8}, \qquad t_0\gg 1.
\]
Therefore we have 
\begin{align*}
 \left.\int_0^\infty \left(|w_t + w_r|^2 + \frac{|w|^2}{r}\right) {\rm d} r\right|_{t=t_0}^{T_0}  = O(t_0^{-9/8}) - \int_{t_0}^{T_0}  \int_0^\infty \frac{|w|^2}{r^2} {\rm d}r {\rm d}t.
\end{align*}
Letting $T_0\rightarrow +\infty$, by \eqref{limit of Vu R3} and \eqref{decay of wrt whole} we have 
\begin{align*}
 \int_0^\infty \left(|(w_t + w_r)(r,t_0)|^2 + \frac{|w(r,t_0)|^2}{r}\right) {\rm d} r = O(t_0^{-9/8}) + \int_{t_0}^{\infty}  \int_0^\infty \frac{|w|^2}{r^2} {\rm d}r {\rm d}t.
\end{align*}
Thus 
\begin{align} \label{limit of t wt plus wr} 
 \lim_{t\rightarrow +\infty} t\int_{0}^\infty |w_r+w_t|^2 {\rm d} r & = \lim_{t\rightarrow +\infty} \left(t\int_{t}^{\infty}  \int_0^\infty \frac{|w|^2}{r^2} {\rm d}r {\rm d} t'- t\int_0^\infty \frac{|w|^2}{r}{\rm d}r\right)
\end{align} 
It follows \eqref{limit of u2 R3} that $\|u\|_{L^2(\Rm^3)}$ and $\|w\|_{L^2(\Rm^+)}$ converges. We let 
\begin{equation} \label{limit of w L2}
 \lim_{t\rightarrow +\infty} \int_0^\infty |w(r,t)|^2 {\rm d} r = E_1. 
\end{equation} 
Next we observe
\begin{align*}
 t\int_{t}^{\infty}  \int_0^\infty \frac{|w|^2}{r^2} {\rm d}r {\rm d} t' - E_1 &=  t\int_{t}^{\infty} \int_0^\infty \left(\frac{|w|^2}{t'^2} + \frac{(t'^2-r^2) |w|^2}{r^2 t'^2}\right){\rm d}r {\rm d} t' - E_1\\
 &=  t \int_t^\infty \frac{1}{t'^2} \left(\int_0^\infty |w|^2 {\rm d} r - E_1 \right) {\rm d} t' + t \int_{t}^{\infty} \int_0^\infty \frac{(t'^2-r^2) |w|^2}{r^2 t'^2} {\rm d} r {\rm d}t'.
\end{align*}
The first term in the right hand side converges to zero as $t\rightarrow +\infty$ by the convergence of $\|w\|_{L^2}$. We may also gives an upper bound of the second term by the support of $w$ and \eqref{limit of w L2}:
\begin{align*}
 t \left|\int_{t}^{\infty} \int_0^\infty \frac{(t'^2-r^2) |w|^2}{r^2 t'^2} {\rm d} r {\rm d}t' \right| & \lesssim t \int_{t}^{\infty} \int_0^\infty \frac{|w|^2}{t'^{5/2}} {\rm d} r {\rm d}t' \lesssim t^{-1/2}. 
\end{align*}
Therefore we have 
\[
 \lim_{t\rightarrow +\infty} t\int_{t}^{\infty}  \int_0^\infty \frac{|w|^2}{r^2} {\rm d}r {\rm d} t' = E_1. 
\]
Similarly we have 
\[
 \lim_{t\rightarrow +\infty} t\int_0^\infty \frac{|w|^2}{r}{\rm d}r = E_1. 
\]
Inserting these two limits into \eqref{limit of t wt plus wr} yields
\begin{equation} \label{decay of w inward}
 \lim_{t\rightarrow +\infty} t\int_{0}^\infty |w_r+w_t|^2 {\rm d} r  = 0.
\end{equation}  
By \eqref{ignore the difference} we also have 
\begin{equation} 
 \lim_{t\rightarrow +\infty} \int_{t-\frac{1}{2}t^{1/2}}^t \frac{1}{t+r}|v_y+v_\tau|^2 {\rm d} r = 0.
\end{equation}
We combine this limit with \eqref{norm of Tv}, \eqref{limit of Jt}, \eqref{expression of flux}, \eqref{limit of energy flux v} and obtain  
\begin{align*}
 \|\mathbf{T} v\|_{\mathcal{C}}^2 & = \lim_{t\rightarrow \infty} J(t) = 2\pi \lim_{t\rightarrow +\infty} \int_{t-\frac{1}{2}t^{1/2}}^t \left(|v_y-v_\tau|^2 + 2 |v|^2\right) {\rm d} r\\
 \|v\|_{\mathcal{K}}^2 & =  \lim_{t\rightarrow +\infty} \frac{1}{2}\int_{t-\frac{1}{2}t^{1/2}}^t \left[|v_t-v_\tau|^2+2\left(1+\frac{1}{t+r}\right)|v|^2\right] {\rm d} r.
\end{align*}
A simple comparison implies that 
\[
 \|\mathbf{T} v\|_{\mathcal{C}}^2 = 4\pi \|v\|_{\mathcal{K}}^2,\qquad v\in \mathcal{K}_0. 
\]
It follows that $\mathbf{T}$ can naturally extend to a bounded linear operator from $\mathcal{K}$ to $\mathcal{C}$ with 
\begin{equation} \label{norm identity}
 \|\mathbf{T} v\|_{\mathcal{C}}^2 = 4\pi \|v\|_{\mathcal{K}}^2,\qquad  v\in \mathcal{K}.
\end{equation}

\section{Homomorphism between $\mathcal{K}$ and $\mathcal{C}$}

Next we show that the linear transformation $\mathbf{T}: \mathcal{K}\rightarrow \mathcal{C}$ defined above is actually a homomorphism. In view of the norm identity \eqref{norm identity}, it suffices to show the range of $\mathbf{T}$ is dense in the space $\mathcal{C}$. We consider the following subset $\mathcal{C}_0$ of the space $\mathcal{C}$:

\begin{lemma}
 Let $\mathcal{C}_0$ be the subset of $\mathcal{C}$ consisting of the radial free Coulomb waves satisfying the following conditions:  
 \begin{itemize} 
  \item There exists a time $t_0>0$, so that the data $(u(t_0), u_t(t_0))$ are smooth and supported in the closed ball $\bar{B}(0,t_0)$;
  \item The data $(u(t_0), u_t(t_0))$ vanishes in a neighbourhood of zero.
 \end{itemize}
Then $\mathcal{C}_0$ is dense in the space $\mathcal{C}$. 
\end{lemma}
\begin{proof}
 Given $u \in \mathcal{C}$ and $\varepsilon >0$, we need to find a solution $\tilde{u} \in \mathcal{C}_0$ so that $\|u-\tilde{u}\|_{\mathcal{C}} < \varepsilon$.  Since radial data $(u_0,u_1)\in \mathcal{C}_0^\infty(\Rm^3\setminus \{0\})$ are dense in the radial subspace of $(\mathcal{H}^1 \cap L^2)\times (L^2\cap \mathcal{H}^{-1})$, it suffices to consider solutions with these data. Let us assume the initial data are supported in the ball $B(0,R_0)$. Thanks to Proposition \ref{prop energy retraction}, we have 
 \[
   J_1(t) = \int_{|x|>t-2} \left(|\nabla u(x,t)|^2 + |u_t(x,t)|^2  \right) {\rm d} x  \rightarrow 0.
 \]
 Since $u(\cdot,t)$ is supported in the ball $B(0,t+R_0)$, we have ($t\gg 1$)
\begin{align*}
 \sup_{|x|>t-2} |u(x)| \lesssim_1 (R_0+2)^{1/2} \left(\int_{t-2}^{t+R} |u_r(r',t)|^2 {\rm d} r' \right)^{1/2} \lesssim_1 (R_0+2)^{1/2} t^{-1} J_1(t)^{1/2}.
\end{align*}
Thus we have 
\begin{equation} \label{decay of u cone}
 \lim_{t\rightarrow +\infty} \int_{|x|>t-2} \left(|\nabla u(x,t)|^2 + |u_t(x,t)|^2 + |u(x,t)|^2 \right) {\rm d} x = 0.
\end{equation}
We observe that $v = \mathbf{H}^{-1} u$ solves the same equation with initial data 
 \[
  \left\{\begin{array}{l} v(0) = \mathbf{H}^{-1} u_0 = (-\Delta +\frac{1}{|x|})^{-1} u_0 \in  \mathcal{H}^1; \\ 
   v_t(0) = \mathbf{H}^{-1} u_1 =  (-\Delta +\frac{1}{|x|})^{-1} u_1 \in  L^2.  \end{array}\right.
 \]
 Here we apply Lemma \ref{estimate on H minus 2}. Thus by Proposition \ref{prop energy retraction} we have 
 \begin{equation} \label{estimate of J2}
  J_2(t) = \int_{|x|>t-2} \left(|\nabla v(x,t)|^2 + |v_t(x,t)|^2  \right) {\rm d} x \rightarrow 0. 
 \end{equation}
 Next we recall that Lemma \ref{estimate on H minus 2} also gives
 \[
  (v(0),v_t(0)) = (c_1 \psi(x), c_2 \psi(x)), \qquad |x|>R_0. 
 \]
 Here $\psi$ is a fast-decaying solution to the elliptic equation 
 \[
  \left(-\Delta + \frac{1}{|x|}\right) \psi = 0. 
 \]
  By finite speed of propagation we have 
 \[
  v(x,t) = (c_1 + c_2 t)\psi(x), \qquad |x|>R_0+|t|. 
 \]
 Thus for any positive integer $N$ we have
 \[
  |v_t(x,t)| + |\nabla v(x,t)| + |v(x,t)| \lesssim  |x|^{-N}, \qquad |x|>R_0+|t|.
 \]
 Combining this with \eqref{estimate of J2}, we also have 
 \begin{equation} \label{decay of v1}
  \lim_{t\rightarrow \infty} \int_{|x|>t-2} \left(|\nabla v(x,t)|^2 + |v_t(x,t)|^2 +|v(x,t)|^2\right) {\rm d} x = 0.
 \end{equation}
 Furthermore, $v_t$ also solve the Coulomb wave equation $\partial_t^2 \tilde{v} + \mathbf{H} \tilde{v} = 0$ with the initial data 
 \[
  \left\{\begin{array}{l} \tilde{v} (0) = v_t(0) = \mathbf{H}^{-1} u_1 = (-\Delta +\frac{1}{|x|})^{-1} u_1 \in  \mathcal{H}^1; \\ 
   \tilde{v}_t(0) = v_{tt}(0) = -\mathbf{H} v(0) = - u_0 \in C_0^\infty(\Rm^3).  \end{array}\right.
 \]
 It follows from Proposition \ref{prop energy retraction} that 
 \begin{equation} \label{decay of v2}
  \lim_{t\rightarrow +\infty} \int_{|x|>t-2} \left(|\nabla v_t(x,t)|^2 + |v_{tt}(x,t)|^2 \right) {\rm d} x = 0.
 \end{equation} 
 We choose 
 \[
  \bar{v}_1(t) = \rho (|x|-t) v_t(\cdot,t)
 \]
 Here $\rho: \Rm \rightarrow [0,1]$ is a smooth cut-off function so that 
 \begin{align*}
  &\rho(s) = 0, \; s\leq -2;& &\rho(s) = 1, \; s\geq -1.&
 \end{align*}
 In view of \eqref{decay of v1}, it is clear that
 \begin{equation} \label{decay of bar v}
  \lim_{t\rightarrow +\infty} \|\bar{v}_1(t)\|_{L^2} = 0. 
 \end{equation}
Next we choose 
\[
 (\bar{u}_0(t), \bar{u}_1(t)) = \left(\rho(|x|-t) u(t), \mathbf{H} \bar{v}_1(t)\right) = \left(\rho u(t), \left(-\Delta + \frac{1}{|x|}\right) \bar{v}_1(t)\right). 
\]
A basic calculation shows that 
\[
 \|\bar{u}_0\|_{\mathcal{H}^1\cap L^2}^2 \lesssim \int_{|x|>t-2} \left(|\nabla u(x,t)|^2 + |u(x,t)|^2 \right){\rm d} x \rightarrow 0, \quad {\rm as}\; t\rightarrow \infty.
\]
In addition, inserting the definition of $\bar{v}_1$ yields
\begin{align*}
 \bar{u}_1(t) &= \left(-\frac{\partial^2}{\partial r^2} - \frac{2}{r}\cdot\frac{\partial}{\partial r} + \frac{1}{r}\right) \left[\rho(|x|-t) v_{t} (t)\right]\\
 & = \rho u_t (t) - (\rho'' + 2\rho'/r) v_t(t) - 2\rho' v_{tr}(t).
\end{align*}
Please note that we may take the third derivative of $v$ above because $v\in \mathcal{C}^3(\Rm\times (\Rm^3 \setminus \{0\}))$ is sufficiently smooth by Remark \ref{C3 continuity}. Here the initial data $(v(0),v_t(0)) = \mathbf{H}^{-1} (u_0,u_1) \in \mathcal{H}_{\rm rad}^4 \times (\mathcal{H}_{\rm rad}^3 \cap L^2)$. As a result of the expression of $\bar{u}_1(t)$, we have
\begin{align*}
 \|\bar{u}_1(t)\|_{L^2(\Rm^3)}^2 \lesssim \|\rho u_t\|_{L^2}^2 + \int_{t-2<|x|<t-1} \left(|\nabla v_t(x)|^2 + |v_t(x)|^2 \right) {\rm d} x. 
\end{align*}
It follows from \eqref{decay of u cone}, \eqref{decay of v1} and \eqref{decay of v2} that 
\[
 \lim_{t\rightarrow +\infty} \|\bar{u}_1(t)\|_{L^2(\Rm^3)} = 0. 
\]
By \eqref{decay of bar v} we also have 
\[
 \lim_{t\rightarrow +\infty} \|\mathbf{H}^{-1} \bar{u}_1(t)\|_{L^2} = \lim_{t\rightarrow +\infty} \|\bar{v}_1(t)\|_{L^2} = 0; \quad \Rightarrow \quad\lim_{t\rightarrow \infty}\|\bar{u}_1(t)\|_{\mathcal{H}^{-2}} = 0. 
\]
An interpolation then gives that 
\[
 \lim_{t\rightarrow +\infty} \|\bar{u}_1(t)\|_{L^2 \cap \mathcal{H}^{-1}} = 0. 
\]
Thus the following inequality holds for a sufficiently large time $t_0$
\[
 \left\|(\bar{u}_0(t_0), \bar{u}_1(t_0))\right\|_{(\mathcal{H}^1 \cap L^2)\times (L^2 \cap \mathcal{H}^{-1})} < \varepsilon/2. 
\]
In addition, the expression of $\bar{u}_1$ given above implies that the data 
\begin{align*}
 (u'_0, u'_1) & = (u(t_0), u_t(t_0)) - (\bar{u}_0(t_0), \bar{u}_1(t_0)) \\
 & = \left((1-\rho)u(t_0), (1-\rho)u_t(t_0)+(\rho'' + 2\rho'/r) v_t(t_0) + 2\rho' v_{tr}(t_0)\right)
\end{align*}
are supported in the ball $B(0,t_0-1)$. Obviously we have 
\[
  \left\|(u(t_0),u_t(t_0))-(u'_0,u'_1)\right\|_{(\mathcal{H}^1 \cap L^2)\times (L^2 \cap \mathcal{H}^{-1})} < \varepsilon/2.
\]
By a standard smooth approximation of $(u'_0,u'_1)$ we may find a pair of smooth initial data $(\tilde{u}_0, \tilde{u}_1)$, so that they are supported in $B(0,t_0)$, vanish near the origin and satisfy 
\[
 \left\|(u(t_0), u_t(t_0))-(\tilde{u}_0,\tilde{u}_1)\right\|_{(\mathcal{H}^1 \cap L^2)\times (L^2 \cap \mathcal{H}^{-1})} < \varepsilon. 
\]
Let $\tilde{u} = \mathbf{S}_{\mathcal{C}}(t-t_0) (\tilde{u}_0,\tilde{u}_1)$. Clearly we have $\tilde{u} \in \mathcal{C}_0$ and 
\[
 \|\tilde{u}-u\|_{\mathcal{C}} < \varepsilon,
\]
thus finish the proof. 
\end{proof}

Next we show that $\mathcal{C}_0 \subset {\rm Ran} \mathbf{T}$. The following argument actually gives the inverse $\mathbf{T}^{-1} u$ for $u \in \mathcal{C}_0$. As usual we define $w(r,t) = r u(r,t)$. We first observe that the inverse of the geometric transform  
\[
 (y,\tau) = \left(\frac{r-t+\ln(t+r)}{2}, \frac{t-r+\ln(t+r)}{2}\right)
\]
is 
\[ 
 (r,t) = \left(\frac{e^{y+\tau}+y-\tau}{2}, \frac{e^{y+\tau}-y+\tau}{2}\right). 
\]
We define 
\begin{equation} \label{def of v tau}
 v(y,\tau) = \chi(y+\tau-2\ln \tau) w\left(\frac{e^{y+\tau}+y-\tau}{2}, \frac{e^{y+\tau}-y+\tau}{2}\right), \qquad \tau \gg 1. 
\end{equation}
Here $\chi: \Rm \rightarrow [0,1]$ is smooth cut-off function so that 
\begin{align*}
 &\chi(s) = 0, \; s\leq 1;& &\chi(s) = 1, \; s\geq 2.&
\end{align*}
It is easy to see the function $\tilde{v}$ vanishes unless $-\tau + 2\ln \tau + 1 \leq y \leq \tau$. 
A direct calculation shows that 
\begin{align*}
 v_\tau & = \chi' \left(1-\frac{2}{\tau}\right) w + \chi w_r \cdot \frac{e^{y+\tau}-1}{2} + \chi w_t \cdot \frac{e^{y+\tau}+1}{2};\\
 v_y & = \chi' w + \chi w_r \cdot \frac{e^{y+\tau}+1}{2} + \chi w_t \cdot \frac{e^{y+\tau}-1}{2};
\end{align*}
and 
\begin{align*}
 v_{\tau \tau} &= \left[\left(1-\frac{2}{\tau}\right)^2\chi'' + \frac{2}{\tau^2}\chi' \right]w + 2\chi' \left(1-\frac{2}{\tau}\right)\left[\frac{e^{y+\tau}-1}{2} w_r + \frac{e^{y+\tau}+1}{2}w_t\right]\\
 & \quad + \chi \left[w_{rr} \left(\frac{e^{y+\tau}-1}{2}\right)^2 + 2 w_{rt}  \left(\frac{e^{y+\tau}-1}{2}\right)\left(\frac{e^{y+\tau}+1}{2}\right) +w_{tt} \left(\frac{e^{y+\tau}+1}{2}\right)^2\right]\\
 & \qquad + \chi (w_r + w_t)\frac{e^{y+\tau}}{2};
\end{align*}
\begin{align*}
 v_{yy} &= \chi'' w + 2\chi' \left[\frac{e^{y+\tau}+1}{2} w_r + \frac{e^{y+\tau}-1}{2} w_t\right]\\
 & \quad + \chi \left[w_{rr} \left(\frac{e^{y+\tau}+1}{2}\right)^2 + 2 w_{rt}  \left(\frac{e^{y+\tau}-1}{2}\right)\left(\frac{e^{y+\tau}+1}{2}\right) +w_{tt} \left(\frac{e^{y+\tau}-1}{2}\right)^2\right]\\
 & \qquad + \chi (w_r + w_t)\frac{e^{y+\tau}}{2}.
\end{align*}
Thus $v$ solves the approximated equation 
\begin{align}
 v_{\tau \tau} - v_{yy} + 2 v& = \left[\frac{4-4\tau}{\tau^2} \chi'' + \frac{2}{\tau^2}\chi' \right] w -\frac{2}{\tau} \chi' e^{y+\tau} (w_r + w_t) \nonumber\\
 & \qquad  + \left(2-\frac{2}{\tau}\right) \chi' (w_t - w_r) + \chi e^{y+\tau} (w_{tt} - w_{rr}) + 2 \chi w. \label{equation of tilde v}
\end{align}
Since $\tilde{w}$ solve the equation 
\[
 w_{tt} - w_{rr} = -\frac{w}{r} = - \frac{2w}{e^{y+\tau} + y -\tau},
\]
we may rewrite the equation above in the form 
\begin{align}
 v_{\tau \tau} - v_{yy} + 2 v& = \left[\frac{4-4\tau}{\tau^2} \chi'' + \frac{2}{\tau^2}\chi' \right]w -\frac{2}{\tau} \chi' e^{y+\tau} (w_r + w_t)\nonumber \\
 & \qquad  + \left(2-\frac{2}{\tau}\right) \chi' (w_t - w_r) + \chi \frac{2y-2\tau}{e^{y+\tau}+y -\tau} w\nonumber \\
 & = J_1 + J_2 + J_3 + J_4. \label{approximated equation tilde}
\end{align} 
We recall that $w(r,t) \in \mathcal{C}^2 (\Rm^+\times \Rm)$ by Proposition \ref{C2 continuity}. Thus $v(y,\tau) \in \mathcal{C}^2$ for large time $\tau$. Combining this with the fact 
\[
  v (y,\tau) = 0, \qquad |y|\geq \tau, \; \tau \gg 1;
\]
 we can see that the data $(v(\tau), v_\tau(\tau)) \in H^1 \times L^2$ for large time $\tau\geq \tau_0$. Next we claim that 
 \[
  J_k \in L^1([\tau_0,+\infty); L^2(\Rm))
 \]
 Let us first consider the major contribution $J_4$.  We split the interval $y \in (-\tau + 2\ln \tau +1, \tau)$ into two pieces: If $y \geq -\tau + 3\ln \tau$, then 
 \[
  \left|\frac{2y-2\tau}{e^{y+\tau}+y -\tau} \right| \lesssim_1 \frac{1}{\tau^2}.
 \]
 Combining this with the uniform boundedness of $w$ given in Remark \ref{uniform boundedness of w}, we have 
 \[
  \int_{-\tau + 3\ln \tau}^\tau |J_4(y,\tau)|^2 {\rm d} y \lesssim \frac{1}{\tau^3}. 
 \]
 On the other hand, if $-\tau+2\ln \tau + 1 \leq y \leq -\tau + 3\ln \tau$, then 
 \begin{align}
  &e\tau^2 \leq e^{y+\tau} \leq \tau^3;& &t = \frac{e^{y+\tau}-y+\tau}{2} \in (\tau^2, \tau^3);& &t-r = \tau - y \in (\tau,2\tau).&  \label{distribution of tr}
 \end{align}
 Thus we may apply the point-wise estimate of $w$ given in Corollary \ref{corollary in out pointwise} and obtain 
 \begin{align*}
  |J_4(y,\tau)| \leq \left|\frac{2y-2\tau}{e^{y+\tau}+y -\tau} w\right| &\lesssim \frac{1}{\tau} \left(\tau^{-1}\ln \tau\right)^{1/4}\left(\max\{\tau^{-\kappa}, \tau^{-2} \ln^2 \tau\}\right)^{1/4}\\
  & \lesssim \tau^{-1-\frac{\kappa+1}{4}} \ln^{1/4} \tau.
 \end{align*}
 Here $\kappa$ is a constant slightly smaller than $3/2$. Thus 
 \[
  \int_{-\tau+2\ln \tau + 1}^{-\tau + 3\ln \tau} |J_4(y,\tau)|^2 {\rm d} y \lesssim \tau^{-2-\frac{\kappa+1}{2}} \ln^{3/2} \tau.
 \]
 In summary we have 
 \[
  \|J_4(\tau)\|_{L^2(\Rm)} \lesssim \tau^{-3/2}\qquad \Rightarrow \qquad J_4 \in L^1([\tau_0,+\infty); L^2(\Rm)). 
 \]
 Next we deal with $J_1$. Clearly $J_1$ vanishes unless $-\tau + 2\ln \tau + 1 \leq y \leq -\tau + 2\ln \tau +2$. This implies 
 \begin{align*}
  &e\tau^2 \leq e^{y+\tau} \leq e^2 \tau^2;& &t = \frac{e^{y+\tau}-y+\tau}{2} \in (\tau^2, 5\tau^2);& &t-r = \tau - y \in (\tau,2\tau).&
 \end{align*}
 Again by the point-wise estimate of $w$ we have 
 \[
  |J_1(y,\tau)| \lesssim \tau^{-1-\frac{\kappa+1}{4}} \ln^{1/4} \tau \qquad \Rightarrow \qquad J_1 \in L^1([\tau_0,+\infty); L^2(\Rm)). 
 \]
 In order to deal with $J_2$ and $J_3$, we integrate 
 \begin{align*}
 \int_{\tau_1}^{2\tau_1} \int_{-\infty}^\infty \left(|J_2|^2 + |J_3|^2\right) {\rm d} y {\rm d}\tau & =  \int_{\tau_1}^{2\tau_1} \int_{-\tau+2\ln \tau +1}^{-\tau+2\ln \tau +2} \left(|J_2|^2 + |J_3|^2\right) {\rm d} y {\rm d}\tau\\
 & \lesssim \int_{\tau_1}^{2\tau_1} \int_{-\tau+2\ln \tau +1}^{-\tau+2\ln \tau +2} \left(\tau^{-2} e^{2y+2\tau} |w_t + w_r|^2 + |w_t - w_r|^2\right) {\rm d} y {\rm d}\tau\\
 & \lesssim \iint_{\Sigma(\tau_1)} \left(\tau_1^{-2} (t+r)^2 |w_t + w_r|^2 + |w_t - w_r|^2\right) \frac{{\rm d} r {\rm d} t}{t+r}.
 \end{align*}
 Here we use the change of variables formula 
 \[
  {\rm d} y {\rm d}\tau = \frac{{\rm d} r {\rm d} t}{t+r}.
 \]
 A basic calculation shows that the integral region 
 \[ 
  \Sigma(\tau_1) = \left\{\left(\frac{e^{y+\tau}+y-\tau}{2}, \frac{e^{y+\tau}-y+\tau}{2}\right): \tau_1 < \tau< 2\tau_1, 1<y+\tau-2\ln \tau<2\right\}
  \] 
  satisfies 
 \[
  \Sigma(\tau_1) \subset \Sigma^+(\tau_1) \doteq \{(r,t): \tau_1^2 \leq t \leq 20 \tau_1^2; \tau_1 \leq t-r \leq 4\tau_1\}. 
 \]
 In this region $\Sigma^+(\tau_1)$ we have 
 \[
  t+r \simeq_1 \tau_1^2. 
 \]
 Therefore the integral 
 \[
  I \doteq \int_{\tau_1}^{2\tau_1} \int_{-\infty}^\infty \left(|J_2|^2 + |J_3|^2\right) {\rm d} y {\rm d}\tau
 \]
  is dominated by 
 \begin{align*}
  I &\lesssim \iint_{\Sigma^+(\tau_1)} \left( |w_t + w_r|^2 + \tau_1^{-2}|w_t - w_r|^2\right) {\rm d} r {\rm d} t\\
  & \lesssim \int_{\tau_1^2}^{20\tau_1^2} \int_{t-4\tau_1<|x|<t-\tau_1} \left(|\mathbf{L}_+ u|^2 + \tau_1^{-2} |\mathbf{L}_- u|^2 \right) {\rm d} x {\rm d} t\\
  & \lesssim \int_{\tau_1}^{4\tau_1} \left(\int_{C^+(s; \tau_1^2; +\infty)} |\mathbf{L}_+ u|^2 {\rm d} S\right) {\rm d} s + \frac{1}{\tau_1^2} \int_{\tau_1^2}^{20 \tau_1^2} \int_{|x|<t-\tau_1} e_+(x,t) {\rm d} x {\rm d} t\\
  & \lesssim \int_{\tau_1}^{4\tau_1} \max\{\tau_1^{-2\kappa}, s^{-1} \tau_1^{-2\kappa+1}\} {\rm d} s + \frac{1}{\tau_1^2} \int_{\tau_1^2}^{20 \tau_1^2}  \max\{\tau_1^{-\kappa}, \tau_1^{-2} \ln^2 t\}{\rm d} t\\
  & \lesssim \tau_1^{-2\kappa+1} + \tau_1^{-\kappa} \lesssim \tau_1^{-\kappa}.
 \end{align*}
 Here we utilize the conclusion of Corollary \ref{corollary in out pointwise}; $\kappa$ is a constant slightly smaller than $3/2$. Thus we have ($k=2,3$)
 \[
  \|J_k\|_{L^2([\tau_1,2\tau_1]; L^2(\Rm))} \lesssim \tau_1^{-\kappa/2}; \quad \Rightarrow \quad \|J_k\|_{L^1([\tau_1,2\tau_1]; L^2(\Rm))} \lesssim \tau_1^{-\frac{\kappa-1}{2}}. 
 \]
 Thus we have 
 \[
  \|J_2\|_{L^1([\tau_0,+\infty); L^2(\Rm))} +  \|J_3\|_{L^1([\tau_0,+\infty); L^2(\Rm))} < +\infty. 
 \]
 Collecting all the upper bounds of $J_k$, we conclude that $v$ scatters in the energy space. More precisely, there exists a solution $v^\ast \in \mathcal{K}$, so that 
 \[
  \lim_{\tau \rightarrow +\infty} \int_{\Rm} \left(|v_y^\ast - v_y|^2 + 2|v^\ast-v|^2 + |v_t^\ast - v_t|^2 \right) {\rm d} y = 0.
 \]
 In fact we may choose $v^\ast = \mathbf{S}_{\mathcal{K}}(t) \vec{v}_+$ with 
 \[
  \vec{v}_+  = \mathbf{S}_{\mathcal{K}}(-\tau_0) \vec{v} (\tau_0) + \int_{\tau_0}^\infty \mathbf{S}_{\mathcal{K}}(-t) (0, J_1+J_2+J_3+J_4) {\rm d} t\in H^1\times L^2. 
 \]
 We claim that $\mathbf{T} v^\ast = u$ thus $u \in {\rm Ran} \mathbf{T}$. In fact we may choose $v^k \in \mathcal{K}_0$ so that 
 \[
  v^k \rightarrow v^\ast, \quad {in}\; \mathcal{K}. 
 \]
 Clearly we have 
 \[
  \lim_{\tau \rightarrow +\infty} \int_{\Rm} \left(|v_y^k - v_y|^2 + 2|v^k -v|^2 + |v_t^k - v_t|^2 \right) {\rm d} y = \|v^\ast - v^k\|_{\mathcal{K}}^2.
 \]
 For any large time $t$, we let $\Omega^\ast(t)$, $\Gamma_t$, $\tau_0(t)$, $y_0(t)$ be the region, curve and coordinates introduced in the previous section. Then by the fact $\tau_0(t) \rightarrow +\infty$ as $t\rightarrow +\infty$ we have 
 \[
  \limsup_{t\rightarrow +\infty} E_k(t) \leq  \|v^\ast-v^k\|_{\mathcal{K}}^2. 
 \]
 where 
 \[
  E_k(t) = \int_{y_0(t)}^{\tau_0(t)} \left(|v_y^k - v_y|^2 + 2|v^k -v|^2 + |v_t^k - v_t|^2 \right)|_{\tau =\tau_0(t)} {\rm d} y.
 \]
 We recall the divergence formula of energy
 \[
  \int_{\partial \Omega^\ast(t)}  2 \varphi_y \varphi_\tau {\rm d} \tau + \left(|\varphi_\tau|^2 + |\varphi_y|^2 + 2|\varphi|^2\right) {\rm d} y = \iint_{\Omega^\ast(t)} 2v_\tau(\varphi_{\tau \tau}-\varphi_{yy}+2\varphi) {\rm d} y {\rm d} \tau,
 \]
 and apply it on $\varphi = v - v^k$ to obtain the energy formula as we did in the previous section 
 \begin{align} \label{energy flux difference}
  E_k(t) + \iint_{\Omega^\ast(t)} 2(v_\tau - v_\tau^k)(v_{\tau \tau}-v_{yy}+2v) {\rm d} y {\rm d} \tau  = {\rm Flux}(v-v^k, \Gamma_t).
 \end{align}
 Here 
 \begin{equation} \label{def of energy flux rebuild}
  {\rm Flux}(\varphi, \Gamma_t) = \frac{1}{2}\int_{t-\frac{1}{2}t^{1/2}}^t \left[|\varphi_y-\varphi_\tau|^2+\frac{1}{t+r}|\varphi_y+\varphi_\tau|^2+2\left(1+\frac{1}{t+r}\right)|\varphi|^2\right] {\rm d} r
 \end{equation}
 We recall that $v_{\tau \tau}-v_{yy}+2v \in L_\tau^1([\tau_0,+\infty); L_y^2(\Rm))$, and $\|v_\tau - v_\tau^k\|_{L_y^2(\Rm)}$ are uniformly bounded for all large time $\tau$ and large $k$. It follows that 
 \[
  \lim_{t\rightarrow +\infty} \iint_{\Omega^\ast(t)} 2(v_\tau - v_\tau^k)(v_{\tau \tau}-v_{yy}+2v) {\rm d} y {\rm d} \tau = 0. 
 \]
 Thus we may let $t\rightarrow +\infty$ in \eqref{energy flux difference} to deduce 
 \begin{equation} \label{limit of energy flux rebuild} 
  \limsup_{t\rightarrow +\infty} {\rm Flux}(v-v^k, \Gamma_t) \leq  \|v^\ast-v^k\|_{\mathcal{K}}^2. 
 \end{equation}
 Next we let $w^k(r,t)$ and $u^k(x,t)$ be approximated solutions associated to $v^k$ as we defined in Section \ref{sec: construction of operator}. We also observe that when $t$ is large and $t-\frac{1}{2} t^{1/2} < r < t$, then the corresponding $\tau$ and $y$ satisfies $\tau \gg 1$ and 
 \begin{align*}
  y + \tau - 2\ln \tau & = \ln (t+r) - 2 \ln \frac{t-r+\ln(t+r)}{2}  = \ln \frac{4(t+r)}{(t-r+\ln(t+r))^2}\\
  & \geq \ln \frac{8t - 2 t^{1/2} }{(\frac{1}{2}t^{1/2}+\ln (2t))^2} \geq 3. 
 \end{align*}
 It follows the definitions \eqref{def of w t} and \eqref{def of v tau} that for these $(t,r)$, we have
 \begin{align*}
  w^k (r,t) & = v^k \left(\frac{r-t+\ln(t+r)}{2}, \frac{t-r+\ln(t+r)}{2}\right);\\
  w (r,t) & = v \left(\frac{r-t+\ln(t+r)}{2}, \frac{t-r+\ln(t+r)}{2}\right). 
 \end{align*}
 Thus we may follow a similar argument as in Section \ref{sec: preservation of norm} to conduct a direct calculation and deduce 
 \begin{align*}
  & \limsup_{t\rightarrow +\infty} \int_{t-\frac{1}{2}t^{1/2} }^\infty \left(|w_r - w_r^k|^2 + |w_t-w_t^k|^2 + |w-w^k|^2\right) {\rm d} r \\
  & =  \frac{1}{2} \limsup_{t\rightarrow +\infty} \int_{t-\frac{1}{2}t^{1/2}}^t \left(|\varphi_y-\varphi_\tau|^2 + \frac{|\varphi_y+\varphi_\tau|^2}{(t+r)^2} + 2 |\varphi|^2\right) {\rm d} r, \quad \varphi= v - v^k. 
 \end{align*}
 A comparison with the energy flux \eqref{def of energy flux rebuild} immediately gives that 
 \begin{align*}
   \limsup_{t\rightarrow +\infty} \int_{t-\frac{1}{2}t^{1/2} }^\infty & \left(|w_r - w_r^k|^2 + |w_t-w_t^k|^2 + |w-w^k|^2\right) {\rm d} r \\
   &\qquad \quad  \leq \limsup_{t\rightarrow +\infty} {\rm Flux} (v-v^k, \Gamma_t) \leq \|v^\ast - v^k\|_{\mathcal{K}}^2. 
 \end{align*}
 Here we recall the limit given in \eqref{limit of energy flux rebuild}. We then combine Remark \ref{decay of w}, \eqref{details in distribution of u2} and Proposition \ref{weighted Morawetz kappa 1}, part (iii) to deduce 
 \[
  \limsup_{t\rightarrow +\infty} \int_0^{t-\frac{1}{2}t^{1/2}}\left(|w_r - w_r^k|^2 + |w_t-w_t^k|^2 + |w-w^k|^2\right) {\rm d} r = 0.
 \]
 Therefore we have 
  \[
  \limsup_{t\rightarrow +\infty} \int_0^{\infty}\left(|w_r - w_r^k|^2 + |w_t-w_t^k|^2 + |w-w^k|^2\right) {\rm d} r \leq \|v^\ast - v^k\|_{\mathcal{K}}^2.
 \]
 This gives 
  \[
  \limsup_{t\rightarrow +\infty} \int_{\Rm^3}\left(|\nabla u - \nabla u^k|^2 + |u_t-u_t^k|^2 + |u-u^k|^2\right) {\rm d} x \leq 4\pi \|v^\ast - v^k\|_{\mathcal{K}}^2.
 \]
 By the scattering property we also have 
 \[
  \limsup_{t\rightarrow +\infty} \int_{\Rm^3}\left(|\nabla (\mathbf{T} v^k)  - \nabla u^k|^2 + |\partial_t (\mathbf{T} v^k)-u_t^k|^2 + |\mathbf{T} v^k -u^k|^2\right) {\rm d} x = 0.
 \]
Thus 
\[
  \limsup_{t\rightarrow +\infty} \int_{\Rm^3}\left(|\nabla u - \nabla (\mathbf{T} v^k)|^2 + |u_t-\partial_t (\mathbf{T} v^k)|^2 + |u-\mathbf{T} v^k|^2\right) {\rm d} r \leq 4\pi \|v^\ast - v^k\|_{\mathcal{K}}^2.
 \]
Since $u - \mathbf{T} v^k \in \mathcal{C}$, by Proposition \ref{half energy} and Corollary \ref{balance for L2} we have 
\[
 \|u-\mathbf{T} v^k\|_{\mathcal{C}}^2 \leq 4\pi \|v^\ast - v^k\|_{\mathcal{K}}^2, \quad \forall k\gg 1. 
\]
Therefore we must have $\mathbf{T} v^\ast = u$ and finish the proof. This implies that the map $\mathbf{T}$ is a homeomorphism from $\mathcal{K}$ to $\mathcal{C}$ in dimension $3$.

\section{Higher dimensional case}

In this section we explain how to define a homeomorphism from $\mathcal{K}$ to $\mathcal{C}$ in the higher dimensional case. In fact it can be done in almost the same way. We only give main statements and skip the details. Let $v\in \mathcal{K}_0$. We can define the same smooth function as in the 3-dimensional case 
\[
 w(r,t) = \rho\left(\frac{t-r}{t^{1/2}}\right) v\left(\frac{r-t+\ln (t+r)}{2}, \frac{t-r+\ln(t+r)}{2}\right), \qquad (r,t)\in \Rm^+ \times [T_0,+\infty),
\]
which satisfies the approximated equation 
\begin{align*}
 w_{tt} - w_{rr} + \frac{w}{r} = f(r,t)
\end{align*}
with 
\begin{align} \label{upper bound of f higher}
 \|f(\cdot,t)\|_{L^2(\Rm^+)} \lesssim t^{-3/2-1/8}.  
\end{align}
We then define a radial function $u(x,t) = |x|^{-\frac{d-1}{2}} w(|x|,t)$ for $(x,t)\in \Rm^d \times [T_0,+\infty)$, which solves the approximated Coulomb wave equation 
\[
 \partial_t^2 u - \Delta u + \frac{u}{|x|} = |x|^{-\frac{d-1}{2}} f(|x|,t) + \lambda |x|^{-\frac{d+3}{2}} w(|x|,t).
\]
Here the right hand side $g(x,t) = |x|^{-\frac{d-1}{2}} f(|x|,t) + \lambda |x|^{-\frac{d+1}{2}} w(|x|,t)$ satisfies 
\[
 |g(x,t)| \lesssim |x|^{-\frac{d-1}{2}} f(|x|,t) + |x|^{-\frac{d+3}{2}} 
\]
We observe that $g(x,t)$ is always zero if $|x|<t-2t^{1/2}$ or $|x|>t$. Thus we may recall \eqref{upper bound of f higher} and deduce  
\begin{align*}
 \|g\|_{L^1([T_0,+\infty);L^2(\Rm^d)))} & < +\infty; \\
 \|g\|_{L^1([T_0,+\infty); \mathcal{H}^{-1}(\Rm^d))} \leq \||x|^{1/2} g\|_{L^1([T_0,+\infty); L^2(\Rm^d))} & < +\infty.
\end{align*}
As in the 3-dimensional case, we may define 
\[ 
 \vec{u}_+ = \mathbf{S}_{\mathcal{C}}(-T_0) \vec{u}(\cdot,T_0)+\int_{T_0}^\infty \vec{\mathbf S}_{\mathcal C} (-t) (0, g(\cdot,t)) {\rm d} t  \in \vec{\mathcal{H}}^1 \cap \vec{\mathcal{H}}^0,
\]
It follows that 
\[
 \lim_{t\rightarrow +\infty} \left\|\vec{u}(\cdot,t) - \vec{\mathbf{S}}_{\mathcal C} (t) \vec{u}_+\right\|_{\vec{\mathcal{H}}^1 \cap \vec{\mathcal{H}}^0}  = 0.
\]
We then define $\mathbf{T} v$ to be the scattering target of $u$, i.e.
\begin{equation}
 \mathbf{T} v = \mathbf{S}_{\mathcal C} (t) \vec{u}_+ \in \mathcal{C}. 
\end{equation}
Next we show for $v\in \mathcal{K}_0$ that 
\begin{equation} \label{energy identity dimension higher T}
 \|\mathbf{T} v\|_{\mathcal{C}}^2 = \sigma_{d-1} \|v\|_{\mathcal{K}}^2. 
\end{equation}
Here $\sigma_{d-1}$ is the area of the sphere $\mathbb{S}^{d-1}$. In fact, we may recall Proposition \ref{half energy} and Corollary \ref{balance for L2} to deduce
\begin{align*}
 \|\mathbf{T} v\|_{\mathcal C}^2 &= \lim_{t\rightarrow +\infty} \left(\|\mathbf{T} v(\cdot,t)\|_{\mathcal{H}^1(\Rm^d)}^2 + \|\partial_t \mathbf{T} v(\cdot,t)\|_{L^2(\Rm^d)}^2 + \|\mathbf{T} v(\cdot,t)\|_{L^2(\Rm^d)}^2 \right)\\
 & = \lim_{t\rightarrow +\infty} \left(\|u(\cdot,t)\|_{\mathcal{H}^1(\Rm^d)}^2 + \|u_t (\cdot,t)\|_{L^2(\Rm^d)}^2 + \|u(\cdot,t)\|_{L^2(\Rm^d)}^2 \right) \\
 & = \lim_{t\rightarrow +\infty} \int_{\Rm^d} \left(|\nabla u(x,t)|^2 + |u_t(x,t)|^2 + \frac{|u(x,t)|^2}{|x|} + |u(x,t)|^2 \right){\rm d} x\\
 & = \lim_{t\rightarrow +\infty} \sigma_{d-1} \int_0^\infty \left(|w_r|^2 + \lambda \frac{|w|^2}{r^2} + |w_t|^2 + \frac{|w|^2}{r} + |w|^2\right) {\rm d} r\\
 & = \sigma_{d-1} \lim_{t\rightarrow +\infty}  \int_0^\infty \left(|w_r|^2 +  |w_t|^2 + |w|^2\right) {\rm d} r.
\end{align*}
We observe that the limit in the last line does NOT depend on the dimension $d$ because $w$ is independent of $d$. Thus we may use the result in 3-dimensional case and obtain 
\[
 \lim_{t\rightarrow +\infty} \int_0^\infty \left(|w_r|^2 +  |w_t|^2 + |w|^2\right) {\rm d} r = \|v\|_{\mathcal{K}}^2. 
\]
The energy identity \eqref{energy identity dimension higher T} then immediately follows. 

Finally we show that the range of $\mathbf{T}$ is the whole space $\mathcal{C}$.  As in the 3-dimensional case, since $\mathcal{C}_0$ is dense in $\mathcal{C}$, it suffices to show that any $u \in \mathcal{C}_0$ is contained in the range of $\mathbf{T}$. Let $w(r,t) = r^\frac{d-1}{2} u(r,t)$ and 
\[
 v(y,\tau) = \chi(y+\tau-2\ln \tau) w\left(\frac{e^{y+\tau}+y-\tau}{2}, \frac{e^{y+\tau}-y+\tau}{2}\right), \qquad \tau \geq \tau_0 \gg 1. 
\]
It follows \eqref{equation of tilde v} and the identity $w_{tt} - w_{rr} = - w/r - \lambda w/r^2$  that 
\begin{align*}
 v_{\tau \tau} - v_{yy} + 2 v& = \left[\frac{4-4\tau}{\tau^2} \chi'' + \frac{2}{\tau^2}\chi' \right]w -\frac{2}{\tau} \chi' e^{y+\tau} (w_r + w_t)\\
 & \qquad  + \left(2-\frac{2}{\tau}\right) \chi' (w_t - w_r) + \chi e^{y+\tau} (w_{tt} - w_{rr}) + 2 \chi w\\
 & = \left[\frac{4-4\tau}{\tau^2} \chi'' + \frac{2}{\tau^2}\chi' \right]w -\frac{2}{\tau} \chi' e^{y+\tau} (w_r + w_t)\\
 & \qquad  +  \left(2-\frac{2}{\tau}\right) \chi' (w_t - w_r) + \chi \frac{2y-2\tau}{e^{y+\tau}+y -\tau} w - \lambda \chi \frac{e^{y+\tau}}{r^2} w\\
 & = J_1 + J_2 + J_3 + J_4 + J_5. 
\end{align*}
Next we show that $J_k \in L^1 ([\tau_0,+\infty); L^2(\Rm))$. The first four terms are defined in the same manner as in dimension 3 thus can be dealt with in the same way. The last term $J_5$ satisfies 
\[
 |J_5| \lesssim  \frac{1}{\tau^2},
\]
Here we use the uniform boundedness of $w$ given in Remark \ref{uniform boundedness of w} and the fact that $r \simeq e^{y+\tau} \gtrsim \tau^2$ if $\chi$ is nonzero. Therefore we have 
\[
 \|J_5(\cdot,t)\|_{L^2(\Rm)} \lesssim \tau^{-3/2} \quad \Rightarrow \quad \|J_5\|_{L^1([\tau_0,+\infty); L^2(\Rm))} < +\infty. 
\]
It follows that $v$ scatters in the positive time direction. More precisely there exists a finite-energy solution $v^*$ to the free Klein-Gordon equation so that 
\[
 \lim_{t\rightarrow +\infty} \int_\Rm \left(|v_x - v_x^*|^2 + |v_t - v_t^*|^2 + 2|v-v^*|^2 \right) {\rm d} x = 0.
\]
We then choose $v^k \in \mathcal{K}_0$ so that $v^k \rightarrow v^*$ in $\mathcal{K}$, then let $w^k(r,t)$ and $u^k(x,t)$ be approximated solutions as we defined in the previous sections of this chapter. Following the same argument as in the three dimensional case, we have 
 \begin{align*}
   \limsup_{t\rightarrow +\infty} \int_{t-\frac{1}{2}t^{1/2} }^\infty  \left(|w_r - w_r^k|^2 + |w_t-w_t^k|^2 + |w-w^k|^2\right) {\rm d} r \leq \|v^\ast - v^k\|_{\mathcal{K}}^2. 
 \end{align*}
In fact this part does not depend on the dimension. A similar argument as in the 3-dimensional case gives
\[
  \limsup_{t\rightarrow +\infty} \int_0^{t-\frac{1}{2}t^{1/2}}\left(|w_r - w_r^k|^2 + |w_t-w_t^k|^2 + |w-w^k|^2\right) {\rm d} r = 0.
 \]
In summary we have 
\begin{align*}
   \limsup_{t\rightarrow +\infty} \int_{0}^\infty  \left(|w_r - w_r^k|^2 + |w_t-w_t^k|^2 + |w-w^k|^2\right) {\rm d} r \leq \|v^\ast - v^k\|_{\mathcal{K}}^2. 
 \end{align*}
By the support and uniform boundedness of $w^k$, as well as the inward/outward energy theory for $u$, we also have 
\begin{align*}
 \lim_{t\rightarrow +\infty} \int_0^\infty \frac{|w^k|^2}{r^2} {\rm d} r & \lesssim \lim_{t\rightarrow +\infty} t^{-3/2} = 0; \\
 \lim_{t\rightarrow +\infty} \int_0^\infty \frac{|w|^2}{r^2} {\rm d} r & \lesssim \lim_{t\rightarrow +\infty} \int_{\Rm^d} \frac{|u|^2}{|x|^2} {\rm d} x= 0.
\end{align*}
We combine these two limits by making use of Lemma \ref{lemma L} and deduce that 
\begin{align*}
   \limsup_{t\rightarrow +\infty} \int_{\Rm^d} \left(|\nabla u - \nabla u^k|^2 + |u_t-u_t^k|^2 + |u-u^k|^2\right) {\rm d} x \leq \sigma_{d-1} \|v^\ast - v^k\|_{\mathcal{K}}^2. 
 \end{align*}
 This implies that 
 \begin{align*}
   \limsup_{t\rightarrow +\infty} \int_{\Rm^d} \left(|\nabla u - \nabla (\mathbf{T} v^k)|^2 + |u_t- \partial_t (\mathbf{T} v^k)|^2 + |u-\mathbf{T} v^k|^2\right) {\rm d} x \leq \sigma_{d-1} \|v^\ast - v^k\|_{\mathcal{K}}^2.
 \end{align*}
 Finally we follow the same argument as in dimension 3 to conclude that $\mathbf{T} v^* = u$ and finish the proof. 
 
 \begin{remark}
 To complete the proof of Theorem \ref{main transformation}, we still need to show that if the initial data $v$ are smooth and compact supported in $B(0,R)$, then  
  \[
  \lim_{t\rightarrow +\infty} \left\|(\mathbf{T} v - u, \partial_t(\mathbf{T} v - u))\right\|_{(\mathcal{H}^1 \cap L^2) \times (L^2 \cap \mathcal{H}^{-1})} = 0,
 \]
 where $u$ is defined in the same manner as in the case of $v\in \mathcal{K}_0$:
  \[
 u(x,t) = \rho\left(\frac{t-|x|}{t^{1/2}}\right) |x|^{-\frac{d-1}{2}} v\left(\frac{|x|-t+\ln (t+|x|)}{2}, \frac{t-|x|+\ln(t+|x|)}{2}\right). 
\]
 In fact, we may follow a similar argument to the case $v\in \mathcal{K}_0$ and conclude that $u$ solves an approximated equation 
 \begin{align*}
  &\partial_t^2 u - \Delta u + \frac{u}{|x|} = f, & &f \in L^1([T,+\infty); L^2\cap \mathcal{H}^{-1}).
 \end{align*}
 Therefore there exists a solution $\tilde{u} \in \mathcal{C}$ so that 
  \[
  \lim_{t\rightarrow +\infty} \left\|(\tilde{u} - u, \tilde{u}_t - u_t)\right\|_{(\mathcal{H}^1 \cap L^2) \times (L^2 \cap \mathcal{H}^{-1})} = 0,
 \]
 It suffices to show that $\mathbf{T} v = \tilde{u}$. Let $v_k \in \mathcal{K}_0$ so that $v_k \rightarrow v$ in $\mathcal{K}$. We then combine the energy flux formula of $v-v_k$ in an enlarged region enclosed by the curve 
\[
 \Gamma_t = \left\{\left(\frac{r-t+\ln (t+r)}{2}, \frac{t-r+\ln(t+r)}{2}\right): t-\frac{1}{2}t^{1/2} < r < t+R\right\}
\]
and the straight lines
\begin{align*}
 &\tau = \tau_0(t) = \frac{1}{4} t^{1/2} + \frac{1}{2}\ln\left(2t - \frac{1}{2} t^{1/2}\right);& & y = \tau + R;&
\end{align*}
with the basic asymptotic property of free Coulomb waves to deduce that 
 \[
  \|\tilde{u} - \mathbf{T} v_k\|_{\mathcal{C}} \lesssim_d \|v-v^k\|_{\mathcal{K}} \rightarrow 0,
 \]
 as we did in Section \ref{sec: preservation of norm}. This immediately verifies $\tilde{u} = \mathbf{T} v$. 
\end{remark}

\section{Energy dispersive rate in the radial direction}

Now we give an application of the isomorphism $\mathbf{T}$ of scattering profiles given above. In fact we may show that the dispersive rate of the energy in the radial direction is roughly $(\ln t)^{-1}$. More precisely we may prove the following radial version of our first main theorem. 

\begin{proposition} \label{radial non concentration}
 Let $u$ be a radial free Coulomb wave with a finite energy and $e(x,t)$ be its corresponding energy density function. Then 
 \begin{itemize}
  \item[(i)] Given any $\eps>0$, there exists two constants $c_1, c_2 > 0$ so that 
  \begin{align*}
   &\limsup_{t\rightarrow +\infty} \int_{|x|>t-c_1 \ln t} e(x,t) {\rm d} x < \eps;& &  \limsup_{t\rightarrow +\infty} \int_{|x|<t-c_2 \ln t} e(x,t) {\rm d} x < \eps. 
  \end{align*}
  \item[(ii)] Assume that $\ell(t)$ is a positive function satisfying the growth condition 
 \[
  \lim_{t\rightarrow +\infty} \frac{\ell(t)}{\ln t} = 0. 
 \]
 Then we also have
 \[
  \lim_{t\rightarrow +\infty} \left(\sup_{r\geq 0} \int_{r<|x|<r+\ell(t)} e(x,t) {\rm d} x \right) = 0. 
 \]
 \end{itemize}
\end{proposition}
\begin{proof}
The existence of $c_2$ has been proved in Lemma \ref{inner decay asymptotic}. We only prove the existence of $c_1$ and part (ii). Since the map $\mathbf{T}$ is a bijective isometry (up to a constant) from $\mathcal{K}$ to $\mathcal{C}$ and $\mathcal{K}_0$ is a dense subspace of $\mathcal{K}$, $\mathbf{T} \mathcal{K}_0$ must be a dense subspace of $\mathcal{C}$. Let $\mathcal{C}_e$ be the space of all finite-energy radial Coulomb free waves equipped with the natural norm 
 \[
  \|u\|_{\mathcal{C}_e}^2 = \|(u,u_t)\|_{\mathcal{H}^1 \times L^2}^2 = \int_{\Rm^d} \left(|\nabla u(x,t)|^2 + |u_t(x,t)|^2 + \frac{|u(x,t)|^2}{|x|}\right) {\rm d} x.
 \]
 Clearly the natural embedding $i: \mathcal{C}\hookrightarrow \mathcal{C}_e$ is a bounded linear map. In addition, $i \mathcal{C}$ is a dense subspace of $\mathcal{C}_e$. It follows that $\mathbf{T} \mathcal{K}_0$ must be a dense subspace of $\mathcal{C}_e$. By linearity it suffices to prove the corresponding asymptotic behaviour for radial free Coulomb waves  $\tilde{u} \in \mathbf{T} \mathcal{K}_0$. Assume that $v \in \mathcal{K}_0$ with $\tilde{u} = \mathbf{T} v$. Let $w$, $u$ be functions associated to $v$ as in the definition of $\mathbf{T}$. We have already shown that $u, \tilde{u}$ share the same asymptotic behaviour in the energy space as $t\rightarrow +\infty$. More precisely we have
 \[
  \lim_{t\rightarrow +\infty} \int_{\Rm^d} \left(|\nabla (u-\tilde{u})(x,t)|^2 + |(u_t-\tilde{u}_t)(x,t)|^2 + \frac{|(u-\tilde{u})(x,t)|^2}{|x|}\right) {\rm d} x = 0.
 \]
 Thus it suffices to prove the following properties of the associated approximated solution $u$ for each $v\in \mathcal{K}_0$ in order to finish the proof of the proposition:
 \begin{itemize} 
   \item[(a)] The following limit holds
 \begin{equation*} 
  \lim_{t\rightarrow +\infty} \left(\sup_{r\geq 0} \int_{r<|x|<r+\ell(t)} \left(|\nabla u(x,t)|^2 + |u_t(x,t)|^2 + \frac{|u(x,t)|^2}{|x|}\right) {\rm d} x \right) = 0.
 \end{equation*} 
  \item[(b)] Given any $\eps>0$, there exists a constant $c_1$ such that 
  \[
   \limsup_{t\rightarrow +\infty} \int_{|x|>t-c_1 \ln t} \left(|\nabla u(x,t)|^2 + |u_t(x,t)|^2 + \frac{|u(x,t)|^2}{|x|}\right) {\rm d} x < \eps. 
  \]
\end{itemize} 
Let us first consider (a). Because $u(x,t)=0$ unless $t-2t^{1/2} < |x| < t$, we have the following identity for large time $t$:
 \[
  \int_{r<|x|<r+\ell(t)} \left(|\nabla u(x,t)|^2 + |u_t(x,t)|^2 + \frac{|u(x,t)|^2}{|x|}\right) {\rm d} x = 0, \qquad r<t/2. 
 \]
 Therefore we only need to show 
  \begin{equation} \label{small ln decay to prove}
  \lim_{t\rightarrow +\infty} \left(\sup_{r\geq t/2} \int_{r<|x|<r+\ell(t)} \left(|\nabla u(x,t)|^2 + |u_t(x,t)|^2 + \frac{|u(x,t)|^2}{|x|}\right) {\rm d} x \right) = 0.
 \end{equation} 
 Next we rewrite the integral above in term of $w(r,t) = r u(r,t)$ by making use of Lemma \ref{lemma L}
 \begin{align}
  \frac{1}{\sigma_{d-1}} & \int_{r<|x|<r+\ell(t)} \left(|\nabla u(x,t)|^2 + |u_t(x,t)|^2 + \frac{|u(x,t)|^2}{|x|}\right) {\rm d} x \nonumber\\
  & = \int_r^{r+\ell(t)} \left(|w_r(r',t)|^2 + |w_t(r',t)|^2 + \frac{|w(r',t)|^2}{r'} + \lambda \frac{|w(r',t)|^2}{r'^2}\right) {\rm d} r' \label{relation u w app}\\
   & \qquad + \frac{d-1}{2}\cdot \frac{|w(r,t)|^2}{r} - \frac{d-1}{2}\cdot \frac{|w(r+\ell(t),t)|^2}{r+\ell(t)}. \nonumber
 \end{align}
 Here the function $w(r',t)$ always vanishes unless $r'\in (t-2t^{1/2},t)$. We also have the following detailed expression of $w$, $w_r$ and $w_t$
\begin{align*}
 w(r,t) & = \rho\left(\frac{t-r}{t^{1/2}}\right) v\left(\frac{r-t+\ln (t+r)}{2}, \frac{t-r+\ln(t+r)}{2}\right);\\
 w_r(r,t)& = -\frac{1}{t^{1/2}}\rho' v + \rho v_y \left(\frac{1}{2}+\frac{1}{2(t+r)}\right) +\rho v_\tau \left(-\frac{1}{2}+\frac{1}{2(t+r)}\right);\\
 w_t(r,t) & = \frac{t+r}{2t^{3/2}}\rho' v + \rho v_y \left(-\frac{1}{2}+\frac{1}{2(t+r)}\right) + \rho v_\tau \left(\frac{1}{2}+\frac{1}{2(t+r)}\right);
\end{align*}
Here we have 
\[
 \frac{t-r+\ln(t+r)}{2} \geq \frac{\ln t}{2}, \qquad r\in [t-2t^{1/2},t]. 
\]
It follows Lemma \ref{general decay of KG} that 
\begin{equation} \label{upper bound of w a half}
 |w(r,t)| + |w_r(r,t)| + |w_t(r,t)| \lesssim (\ln t)^{-1/2}, \qquad t\gg 1. 
\end{equation}
Inserting these upper bounds into \eqref{relation u w app}, we obtain 
\[
 \int_{r<|x|<r+\ell(t)} \left(|\nabla u(x,t)|^2 + |u_t(x,t)|^2 + \frac{|u(x,t)|^2}{|x|}\right) {\rm d} x \lesssim \frac{\ell(t)}{\ln t}+\frac{1}{t\ln t}, \qquad \forall r\geq t/2. 
\]
This immediately verifies \eqref{small ln decay to prove} thus completes the proof of (a). Next we consider (b). Similarly we may rewrite the integral of $u$ in the form of 
 \begin{align*}
  \frac{1}{\sigma_{d-1}} & \int_{t-c_1 \ln t <|x|<t} \left(|\nabla u(x,t)|^2 + |u_t(x,t)|^2 + \frac{|u(x,t)|^2}{|x|}\right) {\rm d} x \\
  & = \int_{t-c_1 \ln t}^{t} \left(|w_r(r',t)|^2 + |w_t(r',t)|^2 + \frac{|w(r',t)|^2}{r'} + \lambda \frac{|w(r',t)|^2}{r'^2}\right) {\rm d} r' \\
   & \qquad + \frac{d-1}{2}\cdot \frac{|w(t-c_1 \ln t,t)|^2}{t-c_1 \ln t}. 
 \end{align*}
 We then make use of \eqref{upper bound of w a half} to deduce 
 \[
  \limsup_{t\rightarrow +\infty} \int_{t-c_1 \ln t <|x|<t} \left(|\nabla u(x,t)|^2 + |u_t(x,t)|^2 + \frac{|u(x,t)|^2}{|x|}\right) {\rm d} x \lesssim c_1.
 \]
 This immediately verifies (b). 
\end{proof}

\chapter{Asymptotic behaviour in the non-radial case} 

In this chapter we prove the first main theorem for non-radial free Coulomb waves. The main tool is the spherically harmonic function decomposition. In this way we may write a free wave as a sum of orthogonal parts, each of them corresponds to a radial free Coulomb wave in a possibly higher dimensional space. This situation is similar to that of the free wave equation. 

\paragraph{Harmonic polynomials} Let us first introduce the harmonic polynomials. We recall that the eigenfunctions of the Laplace-Beltrami operator on $\mathbb{S}^{d-1}$ are the homogeneous harmonic polynomials of the variables $x_1, x_2, \cdots, x_d$. Such a polynomial $\Phi$ of degree $\nu$ satisfies 
\[
 -\Delta_{\mathbb{S}^{d-1}} \Phi = \nu(\nu+d-2) \Phi. 
\]
We choose a Hilbert basis $\{\Phi_k(\theta)\}_{k\geq 0}$ of the operator $-\Delta_{\mathbb{S}^{d-1}}$ on the sphere $\mathbb{S}^{d-1}$. Here we assume that the harmonic polynomial $\Phi_k$ is of degree $\nu_k$. In particular we assume $\nu_0 = 0$ and $\nu_k > 0$ if $k \geq 1$. For more details, please refer to M\"{u}ller \cite{sharmonics}. 

\paragraph{Orthogonal decomposition} Now we apply an orthogonal decomposition of a Coulomb wave $u$ on the sphere $\{x: |x|=r\}$ for each given $r>0$ and time $t>0$. More precisely we define 
\begin{equation} \label{definition of uk}
 u_k(r,t) = r^{-\nu_k} \int_{\mathbb{S}^{d-1}} u(r\theta, t) \Phi_k(\theta) d\theta.
\end{equation}
Thus we have 
\begin{align} \label{reconstruction of u}  
 &u(r\theta,t) = \sum_{k=0}^\infty r^{\nu_k} u_k(r,t) \Phi_k(\theta);& &u_t(r\theta,t) = \sum_{k=0}^\infty r^{\nu_k} \partial_t u_k(r,t) \Phi_k(\theta).
\end{align} 
Furthermore, if we let $ \Box = \partial_t^2 - \partial_r^2 - \frac{d+2\nu_k-1}{r} \partial_r$. A basic calculation shows (we calculate as though $u$ is sufficiently smooth, for general solutions one may apply smooth approximation techniques)
\begin{align*}
 \Box u_k & = (\Box r^{-\nu_k}) \int_{\mathbb{S}^{d-1}} u(r\theta, t) \Phi_k(\theta) d\theta + r^{-\nu_k} \int_{\mathbb{S}^{d-1}} \Box u(r\theta, t) \Phi_k(\theta) d\theta\\
 & \qquad - 2 \partial_r(r^{-\nu_k}) \int_{\mathbb{S}^{d-1}} \partial_r u(r\theta,t) \Phi_k(\theta) d\theta\\
 & = r^{-\nu_k} \int_{\mathbb{S}^{d-1}} \left(\partial_t^2 - \partial_r^2 - \frac{d-1}{r}\partial_r \right) u (r\theta, t) \Phi_k(\theta) d\theta + \nu_k(d-2+\nu_k) r^{-2} u_k\\
 & = r^{-\nu_k} \int_{\mathbb{S}^{d-1}} \left[\left(r^{-2} \Delta_{\mathbb{S}^{d-1}} -r^{-1}\right)u (r\theta, t)\right] \Phi_k(\theta) d\theta + \nu_k(d-2+\nu_k) r^{-2} u_k\\
 & = r^{-\nu_k-2} \int_{\mathbb{S}^{d-1}}  u (r\theta, t) \Delta_{\mathbb{S}^{d-1}} \Phi_k(\theta) d\theta - r^{-1} u_k + \nu_k(d-2+\nu_k) r^{-2} u_k\\
 & = - r^{-1} u_k. 
\end{align*}
Thus if we view $u_k(r,t)$ as a radial function defined in $\Rm^{d+2\nu_k}$, then it solves the Coulomb wave equation 
\[
 \partial_t^2 u_k - \Delta u_k + \frac{u_k}{|x|} = 0, \qquad (x,t) \in \Rm^{d+2\nu_k} \times \Rm.
\]
In addition, we may deduce from \eqref{reconstruction of u} and the orthogonality that ($\sigma_{\ell}$ is the area of the unit sphere $\mathbb S^{\ell}$)
\begin{align*}
 \|u\|_{\dot{H}^1}^2 &=   \int_0^\infty r^{d-1} \int_{\mathbb S^{d-1}} \left|\left(\sum_{k=0}^\infty \partial_r (r^{\nu_k} u_k) \Phi_k \right) \vec{n} + r^{-1} \sum_{k=0}^\infty r^{\nu_k} u_k \nabla_\theta \Phi_k\right|^2 {\rm d}\theta {\rm d} r\\
 & =   \int_0^\infty r^{d-1} \sum_{k=0}^\infty \left[(\nu_k r^{\nu_k-1} u_k + r^{\nu_k} \partial_r u_k)^2 + \nu_k (d-2+\nu_k) r^{2\nu_k-2} u_k^2 \right] {\rm d} r\\
 & =  \sum_{k=0}^\infty \int_0^\infty \left[r^{d+2\nu_k-1} |\partial_r u_k|^2 + 2 \nu_k r^{d+2\nu_k-2} u_k \partial_r u_k + \nu_k(d+2\nu_k-2)r^{d+2\nu_k-3} |u_k|^2\right] {\rm d} r\\
 & = \sum_{k=0}^\infty \int_0^\infty r^{d+2\nu_k-1} |\partial_r u_k|^2 {\rm d} r\\
 & = \sum_{k=0}^\infty \sigma_{d+2\nu_k-1}^{-1} \|u_k\|_{\dot{H}^1(\Rm^{d+2\nu_k})}^2; 
\end{align*}
and 
\begin{align*}
\int_{\Rm^d} \frac{|u(x,t)|^2}{|x|} {\rm d} x& = \int_0^\infty r^{d-2} \int_{\mathbb S^{d-1}}\left|\sum_{k=0}^\infty r^{\nu_k} u_k(r,t) \Phi_k(\theta) \right|^2 {\rm d \theta} {\rm d} r\\
 & = \int_0^\infty r^{d-2} \sum_{k=0}^\infty r^{2\nu_k} |u_k(r,t)|^2  {\rm d} r\\
 & = \sum_{k=0}^\infty \sigma_{d+2\nu_k-1}^{-1} \int_{\Rm^{d+2\nu_k}} \frac{|u_k(x,t)|^2}{|x|} {\rm d} x;
\end{align*}
and 
\begin{align*}
\|u_t\|_{L^2(\Rm^d)}^2 & = \int_0^\infty r^{d-1} \int_{\mathbb S^{d-1}}\left|\sum_{k=0}^\infty r^{\nu_k} \partial_t u_k(r,t) \Phi_k(\theta) \right|^2 {\rm d \theta} {\rm d} r\\
 & = \int_0^\infty r^{d-1} \sum_{k=0}^\infty r^{2\nu_k} |\partial_t u_k(r,t)|^2  {\rm d} r\\
 & = \sum_{k=0}^\infty \sigma_{d+2\nu_k-1}^{-1} \|\partial_t u_k\|_{L^2(\Rm^{d+2\nu_k})}^2. 
\end{align*}
In summary each $u_k(r,t)$ is a radial Coulomb free wave with a finite energy in $\Rm^{d+2\nu_k}$. We also have the energy identity 
\begin{equation} \label{full energy identity}
  E(u) = \sum_{k=0}^\infty \sigma_{d+2\nu_k-1}^{-1} E(u_k)
\end{equation} 
\begin{remark}
 If $0\leq a < b \leq \infty$, then we may follow the same argument as above to deduce
 \[
  \int_{a<|x|<b} |\nabla u|^2 {\rm d} x  = \sum_{k=0}^\infty \int_a^b r^{d-1} \left[(\nu_k r^{\nu_k-1} u_k + r^{\nu_k} \partial_r u_k)^2 + \nu_k (d-2+\nu_k) r^{2\nu_k-2} u_k^2 \right] {\rm d} r.
 \]
 Therefore for any given $N \geq 1$, we have 
 \begin{align*}
   \int_{a<|x|<b} |\nabla u|^2 {\rm d} x  &\leq  \sum_{k=0}^N \int_a^b r^{d-1} \left[(\nu_k r^{\nu_k-1} u_k + r^{\nu_k} \partial_r u_k)^2 + \nu_k (d-2+\nu_k) r^{2\nu_k-2} u_k^2 \right] {\rm d} r\\
   & \; +  \sum_{k=N+1}^\infty \int_0^\infty r^{d-1} \left[(\nu_k r^{\nu_k-1} u_k + r^{\nu_k} \partial_r u_k)^2 + \nu_k (d-2+\nu_k) r^{2\nu_k-2} u_k^2 \right] {\rm d} r.
 \end{align*}
 We then integrate by parts and obtain the following inequality for any $t\in \Rm$:
 \begin{align*}
  \int_{a<|x|<b} |\nabla u|^2 {\rm d} x & \leq  \sum_{k=0}^N \left(\int_a^b r^{d+2\nu_k-1} |\partial_r u_k|^2 {\rm d} r + \nu_k r^{d+2\nu_k-2} |u_k(r)|^2\bigg|_{r=a}^b \right)\\
  & \qquad +  \sum_{k=N+1}^\infty \int_0^\infty r^{d+2\nu_k-1} |\partial_r u_k|^2 {\rm d} r. 
 \end{align*}
Here if $a=0$ or $b=+\infty$, then the term $\nu_k r^{d+2\nu_k-2} |u_k(r)|^2$ at the corresponding endpoint can be ignored. Similarly we may write the integrals $\displaystyle \int_a^b |u_t|^2 {\rm d} x$ and $\displaystyle \int_a^b \frac{|u|^2}{|x|} {\rm d} x$ in term of $u_k$ and deduce the estimate
 \begin{align}
   \int_{a<|x|<b} e(u,x,t) {\rm d} x &\leq \sum_{k=0}^N \left(\sigma_{d+2\nu_k-1}^{-1} \int_{a<|x|<b} e(u_k,x,t) {\rm d} x + \nu_k r^{d+2\nu_k-2} |u_k(r,t)|^2\bigg|_{r=a}^b \right)\nonumber \\
  & \qquad +  \sum_{k=N+1}^\infty \sigma_{d+2\nu_k-1}^{-1} E(u_k).  \label{partial energy identity N} 
 \end{align}
Here $e(u,x,t)$ and $e(u_k, x,t)$ are the energy density function of $u, u_k$ respectively. 
\begin{align*}
 &e(u,x) = \frac{1}{2} |\nabla u|^2 + \frac{1}{2} |u_t|^2 + \frac{|u|^2}{2|x|}; & &e(u_k,x) = \frac{1}{2} |\nabla u_k|^2 + \frac{1}{2} |\partial_t u_k|^2 + \frac{|u_k|^2}{2|x|}.
\end{align*}
\end{remark}

\proof[Proof of Theorem \ref{thm general energy distribution}] Now we prove theorem \ref{thm general energy distribution} without the radial assumption. The proof is a combination of the corresponding radial estimate and the spherical harmonic function decomposition given above. Let us first show that given $u$ and $\varepsilon > 0$, there exists two constant $c_1, c_2 > 0$, so that 
\begin{align}
 \limsup_{t\rightarrow +\infty} \int_{|x|<t-c_2 \ln t} e(x,t) {\rm d} x & < \varepsilon; \label{c2 decay}\\
 \limsup_{t\rightarrow +\infty} \int_{|x|>t-c_1 \ln t} e(x,t) {\rm d} x & < \varepsilon. \label{c1 decay}
\end{align}
We first recall the energy identity \eqref{full energy identity} and choose a large integer $N$, so that 
\begin{equation} \label{tail energy estimate N}
 \sum_{k=N+1}^\infty  \sigma_{d+2\nu_k-1}^{-1} E(u_k) < \varepsilon/2. 
\end{equation}
Next we recall that $u_k$ is a radial finite-energy free Coulomb wave in $\Rm^{d+2\nu_k}$. By the decay estimate of radial $\mathcal{H}^1$ functions (see Lemma \ref{radial pointwise estimate}) 
\[ 
 |u_k(r,t)| \lesssim_{d+2\nu_k} r^{-\frac{2d+4\nu_k-3}{4}}\|u_k(\cdot,t)\|_{\mathcal{H}^1(\Rm^{d+2\nu_k})} \lesssim_{d+2\nu_k} r^{-\frac{2d+4\nu_k-3}{4}} E(u_k)^{1/2},
\]
given any $k$, we have the following uniform limit for all $t$
\[
   r^{d+2\nu_k-2} |u_k(r,t)|^2 \rightrightarrows 0, \qquad {\rm as}\; r\rightarrow +\infty. 
\]
Inserting these into the energy inequality \eqref{partial energy identity N}, we obtain for any given $c>0$ that
\begin{align*}
 \limsup_{t\rightarrow +\infty}  \int_{0<|x|<t-c\ln t} e(u,x) {\rm d} x &\leq \sum_{k=0}^N \left(\frac{1}{\sigma_{d+2\nu_k-1}} \limsup_{t\rightarrow +\infty}\int_{0<|x|<t-c\ln t} e(u_k,x,t) {\rm d} x  \right) + \frac{\varepsilon}{2}. 
\end{align*}
By the corresponding result concerning energy distribution of radial Coulomb free waves, i.e. Proposition \ref{radial non concentration}, there exists $c'_k$ such that 
\[
 \limsup_{t\rightarrow +\infty}\int_{0<|x|<t-c'_k \ln t} e(u_k,x,t) {\rm d} x < \frac{\sigma_{d+\nu_k-1} \varepsilon}{2(N+1)}, \qquad k=0,1,\cdots, N. 
\]
Choosing $c_2 = \max\{c'_0, c'_1, \cdots, c'_N\}$, we immediately have 
\begin{equation} \label{c2 decay proved}
 \limsup_{t\rightarrow +\infty}  \int_{0<|x|<t-c_2\ln t} e(u,x,t) {\rm d} x  < \varepsilon. 
\end{equation}
In the same manner we can find a constant $c_1$ so that \eqref{c1 decay} holds. A direct consequence of \eqref{c2 decay proved} is that 
\begin{equation} \label{inner decay proved}
 \lim_{t\rightarrow +\infty} \int_{|x|<2t/3} e(u,x,t) {\rm d} x = 0.  
\end{equation}
Next we assume that $\ell(t)$ satisfies the growth condition 
\[
 \lim_{t\rightarrow +\infty} \frac{\ell(t)}{\ln t} = 0. 
\]
Thanks to \eqref{inner decay proved}, it suffices to show that 
\begin{equation} \label{to prove ell}
 \lim_{t\rightarrow +\infty} \left(\sup_{r>t/2} \int_{r<|x|<r+\ell(t)} e(u,x,t) {\rm d} x\right) = 0. 
\end{equation}
Given $\varepsilon > 0$, we may follow a similar argument as above to deduce 
\begin{align*}
 \limsup_{t\rightarrow +\infty} &\left(\sup_{r>t/2} \int_{r<|x|<r+\ell(t)} e(u,x,t) {\rm d} x\right)\\
  &\leq \limsup_{t\rightarrow +\infty} \left(\sup_{r>t/2} \sum_{k=0}^N \frac{1}{\sigma_{d+2\nu_k-1}}\int_{r<|x|<r+\ell(t)} e(u_k,x,t) {\rm d} x\right) + \frac{\varepsilon}{2}\\
 &\leq \limsup_{t\rightarrow +\infty} \left(\sum_{k=0}^N \frac{1}{\sigma_{d+2\nu_k-1}}\sup_{r>0}\int_{r<|x|<r+\ell(t)} e(u_k,x,t) {\rm d} x\right) + \frac{\varepsilon}{2}\\
 &\leq \frac{\varepsilon}{2}. 
\end{align*}
Here we use the corresponding result in the radial case given in Proposition \ref{radial non concentration}. Since $\varepsilon > 0$ is arbitrary, we obtain \eqref{to prove ell} and finish the proof. 

\chapter{Strichartz estimates in the radial case}

In this chapter we give a family of Strichartz estimates in the radial case, which play an important role in the scattering theory. We focus on the Strichartz estimates of linear homogeneous equation whose initial data come with a finite energy. First of all, we recall the generalized Strichartz estimate of free wave equation (see Ginibre-Velo \cite{strichartz}). 

\begin{proposition}[Strichartz estimates of free wave equations]\label{Strichartz estimates} 
 Let $2\leq p_1,p_2 \leq \infty$, $2\leq q_1,q_2 < \infty$ and $\rho_1,\rho_2,s\in \Rm$ be constants with
 \begin{align}
  &\frac{2}{p_i} + \frac{d-1}{q_i} \leq \frac{d-1}{2},& &(p_i,q_i)\neq \left(2,\frac{2(d-1)}{d-3}\right),& &i=1,2; \label{admissible condition St} \\
  &\frac{1}{p_1} + \frac{d}{q_1} = \frac{d}{2} + \rho_1 - s;& &\frac{1}{p_2} + \frac{d}{q_2} = \frac{d-2}{2} + \rho_2 +s.& \label{scaling condition St}
 \end{align}
 Assume that $u$ is the solution to the linear wave equation
\[
 \left\{\begin{array}{ll} \partial_t^2 u - \Delta u = F(x,t), & (x,t) \in \Rm^d \times [0,T];\\
 u|_{t=0} = u_0 \in \dot{H}^s; & \\
 \partial_t u|_{t=0} = u_1 \in \dot{H}^{s-1}. &
 \end{array}\right.
\]
Then we have
\begin{align*}
 \left\|\left(u(\cdot,T), \partial_t u(\cdot,T)\right)\right\|_{\dot{H}^s \times \dot{H}^{s-1}} & +\|D_x^{\rho_1} u\|_{L^{p_1} L^{q_1}([0,T]\times \Rm^d)} \\
 & \leq C\left(\left\|(u_0,u_1)\right\|_{\dot{H}^s \times \dot{H}^{s-1}} + \left\|D_x^{-\rho_2} F(x,t) \right\|_{L^{\bar{p}_2} L^{\bar{q}_2} ([0,T]\times \Rm^d)}\right).
\end{align*}
Here the coefficients $\bar{p}_2$ and $\bar{q}_2$ satisfy $1/p_2 + 1/\bar{p}_2 = 1$, $1/q_2 + 1/\bar{q}_2 = 1$. The constant $C$ does not depend on $T$ or $u$. 
\end{proposition}

In particular, we call $(p,q)$ a free wave admissible pair, if and only if it satisfies the conditions \eqref{admissible condition St} in Proposition \ref{Strichartz estimates}. 

\paragraph{Local Strichartz estimates} If $(p,q)$ is a free wave admissible pair with 
\begin{equation} \label{energy level identity}
 \frac{1}{p} + \frac{d}{q} = \frac{d-2}{2}, 
\end{equation}
then the solution to the Coulomb wave equation
\[
 \partial_t^2 u - \Delta u + \frac{u}{|x|} = 0
 \]
with initial data $(u_0,u_1) \in \mathcal{H}^1 \times L^2$ satisfies 
\begin{align}
 \|u\|_{L^p L^q([0,T]\times \Rm^d)} &\leq C \left(\|(u_0,u_1)\|_{\dot{H}^1\times L^2} + \|u/|x|\|_{L^1 L^2([0,T]\times \Rm^d)}\right)\nonumber \\
 & \leq C\left(\|(u_0,u_1)\|_{\mathcal{H}^1\times L^2} + T^{1/2} \|u/|x|\|_{L^2 L^2([0,T]\times \Rm^d)}\right)\nonumber \\
 & \leq C(1+T^{1/2}) \|(u_0,u_1)\|_{\mathcal{H}^1\times L^2}. \label{local strichartz estimate}
\end{align}
Here we use the following Morawetz estimate given in Corollary \ref{Morawetz estimate cor}
\begin{equation} \label{classic Morawetz Coulomb}
 \int_{\Rm} \int_{\Rm^d} \frac{|u(x,t)|^2}{|x|^2} {\rm d} x {\rm d} t \lesssim_d \|(u_0,u_1)\|_{\mathcal{H}^1 \times L^2}. 
\end{equation}
Thus the local Strichartz estimate also holds for finite-energy free Coulomb waves as long as $(p,q)$ is a free wave admissible pair with \eqref{energy level identity}. This may work fairly well in the local theory of nonlinear Coulomb wave equation but probably won't help us with the scattering theory, as the upper bound blows up as $T$ tends to infinity. 

In this chapter we prove the following Strichartz estimates for radial solutions to the Coulomb wave equation

\begin{proposition}[Strichartz estimates of radial solutions] \label{strichartz estimate}
 Assume that $d\geq 3$ and $(p,q)$ satisfies either 
 \begin{align*}
  &2<p \leq \infty;&  &\frac{1}{p} + \frac{d}{q} \geq \frac{d}{2} - 1;& &\frac{2}{p} + \frac{2d-1}{q}< d - \frac{3}{2};&
 \end{align*}
 or 
 \begin{align*}
  &\frac{2(2d+1)}{2d-3} \leq p \leq \infty;& &\frac{2}{p} + \frac{2d-1}{q}= d - \frac{3}{2}.&
 \end{align*}
 Then any radial solution to the linear Coulomb wave equation
 \[
  \left\{\begin{array}{l} \displaystyle \partial_t^2 - \Delta u + \frac{u}{|x|} = f, \quad (x,t)\in \Rm^d \times \Rm; \\ u(0)=u_0 \in \mathcal{H}^1;\\ u_t(0) = u_1 \in L^2\end{array}\right.
 \] 
 satisfies
 \[
  \|u\|_{L^p L^q (\Rm \times \Rm^d)} \lesssim_{d,p,q} \|(u_0,u_1)\|_{\mathcal{H}^1\times L^2} + \|f\|_{L^1 L^2(\Rm \times \Rm^d)}. 
 \]
\end{proposition} 

\begin{figure}[h]
 \centering
 \includegraphics[scale=1.5]{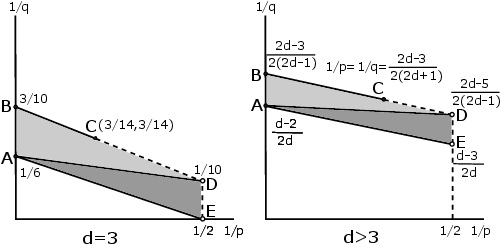}
 \caption{Illustration for allowed pairs} \label{figure admpair}
\end{figure}

\begin{remark}
 The allowed pairs $(1/p,1/q)$ in Proposition \ref{strichartz estimate} are illustrated in figure \ref{figure admpair}, with dimension $d=3$ and $d>3$ respectively. The region containing all possible pairs $(1/p,1/q)$ is exactly the interior of a trapezoid $ABDE$ plus part of its boundary, as shown in the figure \ref{figure admpair}. The point $C$ is on the line segment connecting $BD$. Please note that the open line segments $EA$, $AB$ and $BC$, as well as the points $A, B, C$ are included, while the other part of boundary is not included. We call this region the Coulomb allowed region.
\end{remark}

\begin{remark} \label{nonradial energy level}
 Clearly allowed pairs of the radial Strichartz estimates given above contain all admissible pairs of classic wave equation at the energy level, i.e. those pairs $(1/p,1/q)$ satisfying 
 \begin{equation} \label{pair energy level} 
  \frac{1}{p} + \frac{d}{q} = \frac{d}{2} - 1,
 \end{equation}
 except for the boundary case 
 \[
  \left(\frac{1}{p}, \frac{1}{q}\right) = \left(\frac{1}{2}, \frac{d-3}{2d}\right).
 \]
 In addition, it also contains much more pairs, which are actually at the energy sub-critical level, i.e. $1/p+d/q>d/2-1$. However, the non-radial case is much different. In fact, we claim that if the Strichartz estimates 
 \[
  \|u\|_{L^p L^q(\Rm \times \Rm^d)} \lesssim \|(u_0,u_1)\|_{\mathcal{H}^1 \times L^2} 
 \]
 holds for all free Coulomb waves, without the radial assumption, then $(p,q)$ must satisfies \eqref{pair energy level}. We sketch a proof of this claim at the final part of this chapter. 
 \end{remark}

The majority of this chapter is devoted to the proof of radial Strichartz estimates. 

\section{Preliminary results}

We start by giving a few decay estimates of radial $\mathcal{H}^1$ functions.
\paragraph{Improved decay estimate} It is well known that any radial $\dot{H}^1(\Rm^d)$ function $u$ satisfies a point-wise decay estimate 
\[
 |u(r)| \lesssim_d r^{-\frac{d-2}{2}} \|u\|_{\dot{H}^1(\Rm^d)}. 
\]
Since $\mathcal{H}^1$ norm is stronger than $\dot{H}^1(\Rm^d)$, any radial $\mathcal{H}^1(\Rm^d)$ function $u$ satisfies the same point-wise decay estimate. In fact, the presence of the potential term $|u|^2/|x|$ implies a stronger decay near infinity.

\begin{lemma} \label{radial pointwise estimate}
Let $u\in \mathcal{H}^1(\Rm^d)$ be radial. Then 
\[
 |u(r)| \lesssim_d r^{-\frac{2d-3}{4}} \|u\|_{\mathcal{H}^1 (\Rm^d)}.
\]
\end{lemma}
\begin{proof}
 Clearly it suffices to consider the case $r>1$, since the case $r\leq 1$ has already been covered by the usual $\dot{H}^1$ estimate. We compare 
 \[
  m \doteq \inf_{r<|x|<r+r^{1/2}} |u(x)|
 \]
 with $|u(r)|/2$. There are two cases. 
 \begin{itemize}
  \item If $m \geq |u(r)|/2$, then we have 
  \begin{align*}
   r^{d-\frac{3}{2}} |u(r)|^2\lesssim_1 r^{d-\frac{3}{2}} m^2 \lesssim_d \int_{r<|x|<r+r^{1/2}} \frac{|u(x)|^2}{|x|} {\rm d} x \leq \|u\|_{\mathcal{H}^1}^2. 
  \end{align*}
  This gives the desired result.
  \item If $m\leq |u(r)|/2$, then
  \begin{align*}
   \frac{|u(r)|}{2} & \leq \sup_{r<|x|<r+r^{1/2}} |u(x)-u(r)| = \int_{r}^{r+r^{1/2}} |u_r(r')| {\rm d} r'\\
    & \leq r^{1/4} \left(\int_r^{r+r^{1/2}} |u_{r}(r')|^2 {\rm d} r'\right)^{1/2}\\
    & \lesssim_d r^{-\frac{2d-3}{4}}\left(\int_r^{r+r^{1/2}} (r')^{d-1}|u_{r}(r')|^2 {\rm d} r'\right)^{1/2}\\
    & \lesssim_d r^{-\frac{2d-3}{4}} \|u\|_{\mathcal{H}^1}.
  \end{align*}
 \end{itemize}
 In summary, we always have $|u(r)| \lesssim_d r^{-\frac{2d-3}{4}} \|u\|_{\mathcal{H}^1 (\Rm^d)}$. 
\end{proof}

\begin{corollary} \label{radial q0 infinity estimate}
 If $u \in \mathcal{H}^1(\Rm^d)$ is radial, then we have 
 \[
  \|u\|_{L^{2+\frac{4}{2d-3}}(\Rm^d)} \lesssim_d \|u\|_{\mathcal{H}^1}.
 \]
 In addition, for $q\geq 2+\frac{4}{2d-3}$ and $R>0$ we have 
 \[
  \|u\|_{L^q(\{x:|x|>R\})} \lesssim_d R^{-\frac{2d-3}{4}\cdot\frac{q-2}{q} + \frac{1}{q}} \|u\|_{\mathcal{H}^1}. 
 \]
\end{corollary}
\begin{proof}
 The first estimate is a direct consequence of the definition of $\mathcal{H}^1$ norm and the point-wise decay 
 \[
   |u(x)| \leq c |x|^{-\frac{2d-3}{4}} \|u\|_{\mathcal{H}^1 (\Rm^d)},
 \]
 as given in Lemma \ref{radial pointwise estimate}. Here $c=c(d)$ is a constant. We have
 \begin{align*}
  \int_{\Rm^d} |u(x)|^{2+\frac{4}{2d-3}} {\rm d} x &\lesssim_d \int_{\Rm^d} |u(x)|^2 \left(|x|^{-\frac{2d-3}{4}}\|u\|_{\mathcal{H}^1}\right)^\frac{4}{2d-3} {\rm d} x \\
  & \lesssim_d \|u\|_{\mathcal{H}^1}^\frac{4}{2d-3} \int_{\Rm^d} \frac{|u(x)|^2}{|x|} {\rm d} x \lesssim_d \|u\|_{\mathcal{H}^1}^{2+\frac{4}{2d-3}}.
 \end{align*}
 The proof of the second inequality is similar. 
 \begin{align*}
  \int_{|x|>R} |u(x)|^{q} {\rm d} x &\leq  \int_{|x|>R} |u(x)|^2 \left(c |x|^{-\frac{2d-3}{4}}\|u\|_{\mathcal{H}^1}\right)^{q-2} {\rm d} x \\
  & \leq c^{q-2} \|u\|_{\mathcal{H}^1}^{q-2} \int_{|x|>R} \frac{|u(x)|^2}{|x|} \cdot |x|^{-\frac{(2d-3)}{4}(q-2-\frac{4}{2d-3})} {\rm d} x\\
  & \leq c^{q-2} R^{-\frac{(2d-3)}{4}(q-2-\frac{4}{2d-3})}\|u\|_{\mathcal{H}^1}^{q-2} \int_{|x|>R} \frac{|u(x)|^2}{|x|} {\rm d} x\\
  & \leq c^{q-2} R^{-\frac{(2d-3)}{4}(q-2-\frac{4}{2d-3})}\|u\|_{\mathcal{H}^1}^{q}
 \end{align*}
 Therefore we have
 \[
  \|u\|_{L^q(\{x:|x|>R\})} \leq c^{1-\frac{2}{q}} R^{-\frac{2d-3}{4}\cdot\frac{q-2}{q} + \frac{1}{q}} \|u\|_{\mathcal{H}^1} \leq \max\{c,1\} R^{-\frac{2d-3}{4}\cdot\frac{q-2}{q} + \frac{1}{q}} \|u\|_{\mathcal{H}^1}.
 \]
\end{proof}

\begin{remark} \label{radial Sobolev embedding} 
 A combination of Corollary \ref{radial q0 infinity estimate} and the embedding $\mathcal{H}^1 \hookrightarrow \dot{H}^1 \hookrightarrow L^\frac{2d}{d-2}$ immediately give the radial Sobolev embedding 
 \[
  \|u\|_{L^q (\Rm^d)} \lesssim \|u\|_{\mathcal{H}^1(\Rm^d)}, \qquad 2+\frac{4}{2d-3} \leq q \leq 2 + \frac{4}{d-2}. 
 \]
 From the proof of Corollary \ref{radial q0 infinity estimate}, we may also deduce the inequality for radial $\mathcal{H}^1$ functions 
 \begin{equation} \label{GN inequality}
   \int_{\Rm^d} |u(x)|^{2+\frac{4}{2d-3}} {\rm d} x \lesssim_d \|u\|_{\mathcal{H}^1}^\frac{4}{2d-3} \int_{\Rm^d} \frac{|u(x)|^2}{|x|} {\rm d} x. 
 \end{equation} 
\end{remark}

\begin{proposition}[Special Strichartz estimates] \label{special strichartz estimates}
 Let $u$ be a radial solution to the free Coulomb wave equation 
 \[
  \left\{\begin{array}{l} \displaystyle \partial_t^2 u - \Delta u + \frac{u}{|x|} = 0, \quad (x,t)\in \Rm^d \times \Rm; \\ u(0)=u_0 \in \mathcal{H}^1;\\ u_t(0) = u_1 \in L^2. \end{array}\right.
 \] 
 Then the following global space-time estimates hold. 
 \begin{align*}
  \|u\|_{L^\infty L^{2^*}(\Rm \times \Rm^d)} + \|u\|_{L^\infty L^{2+\frac{4}{2d-3}}(\Rm \times \Rm^d)} &+ \|u\|_{L^{2+\frac{8}{2d-3}} L^{2+\frac{8}{2d-3}}(\Rm \times \Rm^d)} \\
  & \lesssim_d \|(u_0,u_1)\|_{\mathcal{H}^1\times L^2(\Rm^d)}.  
 \end{align*}
 Here $2^*$ satisfies $1/2^* = 1/2-1/d$. 
\end{proposition}
\begin{proof}
 The $L^\infty L^{2^*}$ estimate immediately follows the Sobolev embedding 
 \[
  \|u\|_{L^\infty L^{2^*}} \lesssim_d \|u\|_{L^\infty(\Rm; \dot{H}^1(\Rm^d))} \leq \|u\|_{L^\infty(\Rm; \mathcal{H}^1)} \leq \|(u_0,u_1)\|_{\mathcal{H}^1\times L^2(\Rm^d)}.
 \]
 Meanwhile the $L^\infty L^{2+\frac{4}{2d-3}}$ estimate depends on Corollary \ref{radial q0 infinity estimate}. 
 \[
  \|u\|_{L^\infty L^{2+\frac{4}{2d-3}}} \leq \|u\|_{L^\infty(\Rm; \mathcal{H}^1)} \leq \|(u_0,u_1)\|_{\mathcal{H}^1\times L^2(\Rm^d)}.
 \]
 The final special Strichartz estimate is a consequence of the Morawetz estimate \eqref{classic Morawetz Coulomb} and the point-wise estimate in Lemma \ref{radial pointwise estimate}.
 \begin{align*}
  \int_{\Rm^d \times \Rm} |u(x,t)|^{2+\frac{8}{2d-3}} {\rm d} x {\rm d} t & \lesssim_d \iint_{\Rm^d \times \Rm} |u(x,t)|^2 \left(|x|^{-\frac{2d-3}{4}} \|u(\cdot,t)\|_{\mathcal{H}^1}\right)^{\frac{8}{2d-3}} {\rm d} x {\rm d} t\\
  & \lesssim_d \left(\sup_{t} \|u(\cdot,t)\|_{\mathcal{H}^1}\right)^{\frac{8}{2d-3}} \iint_{\Rm^d \times \Rm} \frac{|u|^2}{|x|^2} {\rm d} x {\rm d} t\\
  & \lesssim_d \|(u_0,u_1)\|_{\mathcal{H}^1\times L^2}^{2+\frac{8}{2d-3}}.
 \end{align*}
\end{proof}

\section{Slowly-growing norms} 

As we showed above, local Strichartz estimates hold for all free Coulomb waves if the corresponding Strichartz estimates hold for the free waves. In particular, the space-time norm in the time interval $[0,T]$ grows at a rate lower than $T^{1/2}$. Next we show that for some pairs $(2,q)$, the norm grows at very low rate, if it does grow, as $T$ tends to the infinity. 

\begin{lemma} \label{slow growing lemma}
 Let $u$ be a radial solution to the free Coulomb wave equation
 \[
  \left\{\begin{array}{l} \displaystyle \partial_t^2 u - \Delta u + \frac{u}{|x|} = 0, \quad (x,t)\in \Rm^d \times \Rm; \\ u(0)=u_0 \in \mathcal{H}^1;\\ u_t(0) = u_1 \in L^2. \end{array}\right.
 \] 
Assume that  
\[ 
  \frac{2(2d-1)}{2d-5}\leq  q \leq \frac{2d}{d-3}.
\]
Then given any positive constant $\alpha > 0$, the following local Strichartz estimate holds 
\[
 \|u\|_{L^2 L^q ([0,T]\times \Rm^d)} \lesssim_{d,q,\alpha} T^\alpha \|(u_0,u_1)\|_{\mathcal{H}^1\times L^2}, \qquad T\geq 1. 
\]
\end{lemma}
\begin{proof}
 First of all, we recall that the pair $(2,\frac{2d}{d-3})$ is a radial free wave admissible pair, i.e. the Strichartz estiamte 
 \begin{equation} \label{radial Strichartz 2d3}
  \|u\|_{L^2 L^\frac{2d}{d-3}(\Rm \times \Rm^d)} \lesssim_d \|(u_0,u_1)\|_{\dot{H}^1 \times L^2} + \|F\|_{L^1 L^2(\Rm \times \Rm^d)} 
 \end{equation}
 holds for all radial solutions to the linear wave equation $\partial_t^2 u - \Delta u = F$ with initial data $(u_0,u_1)$, because
 \begin{itemize} 
  \item The case $d\geq 4$ is covered by the regular Strichartz estimates given in Proposition \ref{Strichartz estimates}.
  \item In the 3-dimensional case, the pair $(2,+\infty)$ is NOT free wave admissible. But the Strichartz estimate \eqref{radial Strichartz 2d3} still holds for radial solutions. In fact it suffices to consider the homogeneous case $F=0$ by the Duhamel's formula. This can be proved by the theory of radiation fields and classic maximal function. Please refer to Section 8.3 of Miao-Shen \cite{wavebook}. 
 \end{itemize}
 The local Strichartz estimate of radial free Coulomb waves immediately follows: 
 \[
  \|u\|_{L^2 L^\frac{2d}{d-3} ([0,T]\times \Rm^d)} \lesssim_d T^{1/2} \|(u_0,u_1)\|_{\mathcal{H}^1\times L^2}, \qquad T\geq 1. 
 \] 
 This gives the Strichartz estimates of the centre part
 \begin{equation} \label{center step 0} 
  \|u\|_{L^2([0,T]; L^q (\{x:|x|\leq 1\}))} \lesssim_d T^{1/2} \|(u_0,u_1)\|_{\mathcal{H}^1\times L^2}. 
 \end{equation}
  The exterior part estimate depends on Corollary \ref{radial q0 infinity estimate}, since we always have $\frac{2(2d-1)}{2d-5}>2+\frac{4}{2d-3}$.
 \begin{equation} \label{exterior step 0} 
  \|u\|_{L^2([0,T]; L^q (\{x:|x|> 1\}))} \lesssim_{d,q} \|u\|_{L^2([0,T]; \mathcal{H}^1)} \lesssim_{d,q} T^{1/2}\|(u_0,u_1)\|_{\mathcal{H}^1\times L^2}. 
 \end{equation}
 A combination of \eqref{center step 0} and \eqref{exterior step 0} gives the local Strichartz estimate with $\alpha  = 1/2$. 
 \begin{equation*}
   \|u\|_{L^2([0,T]; L^q (\Rm^d))} \lesssim_{d,q} T^{1/2} \|(u_0,u_1)\|_{\mathcal{H}^1\times L^2}, \qquad T\geq 1. 
 \end{equation*}
 Next we apply an induction argument to show that the local Strichartz estimate holds for any $\alpha > 0$. We assume that the following local Strichartz estimate holds for some $\alpha \in (0,1/2]$. 
 \begin{equation} \label{step 1 assumption}
   \|u\|_{L^2([0,T]; L^q (\Rm^d))} \lesssim_{d,q,\alpha} T^{\alpha} \|(u_0,u_1)\|_{\mathcal{H}^1\times L^2},\qquad \forall T\geq 1.
 \end{equation}
Given a time $T\geq 1$, we split the space into the interior part $\{x: |x|<T^\beta\}$ and the exterior part $\{x: |x|>T^\beta\}$, then give the upper bound of the norms separately. Here $\beta\in (0,1)$ is a constant to be determined later. For the interior part we recall the energy estimate
 \[
 \sum_{k=-\infty}^\infty \int_{|x|\leq R} e(x,kT') {\rm d}x  \lesssim_d \frac{RE}{T'}
 \]
 for $R\geq T'>0$ given in Proposition \ref{summation of energy}. We let $R = 4T^\beta$, $T' = T^\beta$ and obtain 
 \[
  \sum_{k=-\infty}^\infty \int_{|x|\leq 4T^\beta} e(x,kT^\beta) {\rm d} x  \lesssim_d \|(u_0,u_1)\|_{\mathcal{H}^1\times L^2}^2. 
 \]
 Next we define 
 \[
  (u_{0,k},u_{1,k}) = \rho(T^{-\beta}x) (u(x, kT^\beta), u_t(x,kT^\beta)). 
 \]
 Here $\rho$ is a fixed smooth cut-off function satisfying 
 \begin{align*}
  &\rho(s) = 1, \; s\leq 2;& &\rho(s) = 0, \; s\geq 4.&
 \end{align*}
 A basic calculation shows that 
 \begin{align*}
  \|(u_{0,k}, u_{0,k})\|_{\mathcal{H}^1\times L^2}^2 & \lesssim_d \int_{|x|<4T^\beta} \left(|\nabla u(x,kT^\beta)|^2 + \frac{|u(x,kT^\beta)|^2}{|x|} + |u_t(x,kT^\beta)|^2\right) {\rm d} x\\
  & \lesssim_d \int_{|x|<4T^\beta} e(x,kT^\beta) {\rm d} x.
 \end{align*}
 Thus we have 
 \begin{equation} \label{interior step 1 sum energy}
  \sum_{k=-\infty}^\infty \|(u_{0,k}, u_{0,k})\|_{\mathcal{H}^1\times L^2}^2 \lesssim_d \|(u_0,u_1)\|_{\mathcal{H}^1\times L^2}^2. 
 \end{equation} 
 It immediately follows \eqref{step 1 assumption} that $u^{(k)} = \mathbf{S}_{\mathcal{C}} (u_{0,k}, u_{1,k})$ satisfies  
 \[
  \|u^{(k)}\|_{L^2 L^q ([0,T^\beta]\times \Rm^d)} \lesssim_{d,p,\alpha} T^{\alpha \beta} \|(u_{0,k}, u_{1,k})\|_{\mathcal{H}^1\times L^2}. 
 \]
 By finite speed of propagation and the fact that 
 \[
  (u_{0,k}(x),u_{1,k}(x)) = (u(x, kT^\beta), u_t(x,kT^\beta)), \qquad |x|<2 T^\beta.
 \]
 we have 
 \[
  u(x,t+kT^\beta) = u^{(k)}(x,t), \qquad |x|+|t|<2T^\beta. 
 \]
 Thus we have 
 \[
   \|u\|_{L^2 L^q ([kT^\beta, (k+1)T^\beta] \times \{x: |x|<T^\beta\})} \lesssim_{d,p,\alpha} T^{\alpha \beta} \|(u_{0,k}, u_{1,k})\|_{\mathcal{H}^1\times L^2}.
 \]
 A combination of this with \eqref{interior step 1 sum energy} yields 
 \begin{align*}
  \|u\|_{L^2 L^q(\Rm \times \{x: |x|<T^\beta\})} & = \left(\sum_{k=-\infty}^\infty \|u\|_{L^2 L^q ([kT^\beta, (k+1)T^\beta] \times \{x: |x|<T^\beta\})}^2 \right)^{1/2}\\
  & \lesssim_{d,q,\alpha} T^{\alpha \beta} \left(\sum_{k=-\infty}^\infty \|(u_{0,k}, u_{1,k})\|_{\mathcal{H}^1\times L^2}^2 \right)^{1/2}\\
  & \lesssim_{d,q,\alpha} T^{\alpha \beta} \|(u_0,u_1)\|_{\mathcal{H}^1\times L^2}. 
 \end{align*}
 Next we consider the exterior part. We apply Corollary \ref{radial q0 infinity estimate} and obtain 
 \begin{align*}
  \|u\|_{L^2([0,T];L^q(\{x:|x|>T^\beta\}))} & \leq T^{1/2} \sup_{t\in [0,T]} \|u(\cdot,t)\|_{L^q(\{x: |x|>T^\beta\})} \\
  & \lesssim_d T^{1/2} T^{\beta(-\frac{2d-3}{4}\cdot\frac{q-2}{q} + \frac{1}{q})} \sup_{t\in [0,T]} \|u(\cdot,t)\|_{\mathcal{H}^1} \\
  & \lesssim_d T^{\frac{1}{2} + \beta(-\frac{2d-3}{4}\cdot\frac{q-2}{q} + \frac{1}{q})} \|(u_0,u_1)\|_{\mathcal{H}^1\times L^2}. 
 \end{align*}
 Please note that our assumption on $q$ implies that 
 \[
  q > 2 + \frac{4}{2d-3}. 
 \]
 In summary we obtain 
 \[
   \|u\|_{L^2([0,T];L^q(\Rm^d))} \lesssim_{d,q,\alpha} T^{\max\left\{\alpha \beta, \frac{1}{2} + \beta(-\frac{2d-3}{4}\cdot\frac{q-2}{q} + \frac{1}{q})\right\}} \|(u_0,u_1)\|_{\mathcal{H}^1\times L^2}. 
 \]
 Namely the estimate \eqref{step 1 assumption} still holds if we substitute $\alpha$ by 
 \[
  \max\left\{\alpha \beta, \frac{1}{2} + \beta\left(-\frac{2d-3}{4}\cdot\frac{q-2}{q} + \frac{1}{q}\right)\right\}. 
 \]
 Here $\beta\in (0,1)$ is an arbitrary constant. Please note that our assumption on $q$ guarantees that 
 \[
  -\frac{2d-3}{4}\cdot\frac{q-2}{q} + \frac{1}{q} \leq -\frac{1}{2}. 
 \]
 Thus given any small constant $\varepsilon > 0$, we may choose a constant $\beta\in (0,1)$ such that 
 \[
  \varepsilon =  \frac{1}{2} + \beta\left(-\frac{2d-3}{4}\cdot\frac{q-2}{q} + \frac{1}{q}\right). 
 \]
 Now we may start with $\alpha = 1/2$ and iterate the argument above to deduce that the estimate \eqref{step 1 assumption} holds for all parameters 
 \[
  \alpha = \max\left\{\frac{1}{2}\beta^k, \varepsilon\right\}, \qquad k\geq 0. 
 \]
 Finally we let $k\rightarrow +\infty$ to conclude that 
 \[
   \|u\|_{L^2([0,T]; L^q (\Rm^d))} \lesssim_{d,q,\varepsilon} T^{\varepsilon} \|(u_0,u_1)\|_{\mathcal{H}^1\times L^2}
 \]
 thus finishes the proof.
\end{proof}

\section{Proof of Strichartz estimates}

A direct consequence of the slow growing $L^2 L^q$ norm is the following partial result of the Strichartz estimates: 

\begin{corollary} \label{strichartz low triangle}
 Assume that $(p,q)$ satisfies $2<p \leq \infty$ and 
 \begin{align*} 
  &\frac{1}{p} + \frac{d}{q} \geq \frac{d}{2} - 1;& &\frac{2}{p} + \frac{d(2d-1)}{q} \leq \frac{(2d-1)(d-2)}{2}.
 \end{align*}
 Then any radial solution to 
 \[
  \left\{\begin{array}{l} \displaystyle \partial_t^2 u - \Delta u + \frac{u}{|x|} = 0, \quad (x,t)\in \Rm^d \times \Rm; \\ u(0)=u_0 \in \mathcal{H}^1;\\ u_t(0) = u_1 \in L^2. \end{array}\right.
 \] 
 satisfies the Strichartz estimate
 \[
  \|u\|_{L^p L^q (\Rm \times \Rm^d)} \lesssim_{d,p,q} \|(u_0,u_1)\|_{\mathcal{H}^1\times L^2}. 
 \] 
\end{corollary}
\begin{proof}
 First of all, we observe that these pairs are exactly those in the darker triangle region of figure \ref{figure admpair}. We start by following the same argument as in the proof of Lemma \ref{slow growing lemma} and obtain 
 \[
  \sum_{k=-\infty}^\infty \int_{|x|\leq 4R} e(x,kR) {\rm d}x  \lesssim_d \|(u_0,u_1)\|_{\mathcal{H}^1\times L^2}^2, \qquad \forall R>0. 
 \]
 Thus if we define 
 \[
  (u_{0,k},u_{1,k}) = \rho(x/R) (u(x, kR), u_t(x,kR))
 \]
with a fixed smooth cut-off function $\rho$ satisfying 
\begin{align*}
  &\rho(s) = 1, \; s\leq 2;& &\rho(s) = 0, \; s\geq 4;&
 \end{align*}
then we have 
\begin{equation} \label{interior R sum energy}
  \sum_{k=-\infty}^\infty \|(u_{0,k}, u_{0,k})\|_{\mathcal{H}^1\times L^2}^2 \lesssim_d \|(u_0,u_1)\|_{\mathcal{H}^1\times L^2}^2, \qquad \forall R\geq 1. 
 \end{equation} 
 By the finite speed of propagation and the conclusion of Lemma \ref{slow growing lemma}, we have 
 \[
  \|u\|_{L^2([kR, (k+1)R]; L^{\tilde{q}}(\{x: |x|\leq R\}))} \lesssim_{d, \tilde{q}, \alpha}  R^\alpha \|(u_{0,k}, u_{1,k})\|_{\mathcal{H}^1\times L^2}, \qquad \forall \alpha > 0,\; R\geq 1. 
 \]
 Here $\tilde{q}$ is an arbitrary parameter satisfying 
 \[
   \frac{2(2d-1)}{2d-5}\leq  \tilde{q} \leq \frac{2d}{d-3}.
 \]
 It immediately follows this upper bound and \eqref{interior R sum energy} that 
 \begin{equation} \label{2q estimate}
  \|u\|_{L^2(\Rm; L^{\tilde{q}}(\{x: |x|\leq R\}))} \lesssim_{d, \tilde{q}, \alpha}  R^\alpha \|(u_0, u_1)\|_{\mathcal{H}^1\times L^2}, \qquad \forall \alpha > 0, \; R\geq 1. 
 \end{equation} 
 Next we recall the universal $L^\infty L^{2^*}$ estimate given by the Sobolev embedding:
 \[
  \|u\|_{L^\infty(\Rm; L^{2^*}(\Rm^d))} \lesssim_d \|u\|_{L^\infty(\Rm; \dot{H}^1(\Rm^d))} \leq \|(u_0, u_1)\|_{\mathcal{H}^1\times L^2}. 
 \]
 Please note that the pair $(0, 1/2^*)$ corresponding to the point $A$ in the figure \ref{figure admpair} and that the pairs $(1/2, 1/\tilde{q})$ are exactly those pairs on the line segment $DE$. An interpolation between this $L^\infty L^{2^*}$ estimate and \eqref{2q estimate} with $R=1$ then gives 
 \begin{equation} \label{low triangle centre part}
   \|u\|_{L^p (\Rm; L^{q}(\{x: |x|\leq 1\}))} \lesssim_{d,p,q} \|(u_0, u_1)\|_{\mathcal{H}^1\times L^2}
 \end{equation}
 for any pair $(1/p,1/q)$ in the closed triangle $ADE$. Next we consider the cylinder region $\{(x,t): 2^k \leq |x| < 2^{k+1}, t\in \Rm\}$. On one hand, we utilize \eqref{2q estimate} and write 
 \begin{equation} \label{2q estimate cylinder}
  \|u\|_{L^2(\Rm; L^{\tilde{q}}(\{x: 2^k \leq |x|<2^{k+1}\}))} \lesssim_{d, \tilde{q}, \alpha}  2^{\alpha k} \|(u_0, u_1)\|_{\mathcal{H}^1\times L^2}, \qquad \forall \alpha > 0, \; k\geq 0. 
 \end{equation} 
 On the other hand, we recall the decay estimate given in Corollary \ref{radial q0 infinity estimate} 
 \[
  \|u\|_{L^\infty L^{2^*} (\Rm \times \{x: |x|\geq 2^k\})} \lesssim_d 2^{-\frac{d-1}{2d}k} \|(u_0, u_1)\|_{\mathcal{H}^1\times L^2}, \qquad k\geq 0. 
 \] 
 Again an interpolation of this inequality with \eqref{2q estimate cylinder} shows that for any $(1/p,1/q)$ satisfying the assumption, there exists $\beta = \beta(d,p,q)>0$ such that 
 \begin{equation} \label{pq estimate cylinder}
  \|u\|_{L^p(\Rm; L^{q}(\{x: 2^k \leq |x|<2^{k+1}\}))} \lesssim_{d, p, q}  2^{-\beta k} \|(u_0, u_1)\|_{\mathcal{H}^1\times L^2}, \qquad  k\geq 0. 
 \end{equation} 
 Taking a sum for all $k \geq 0$, we obtain 
 \begin{equation} \label{pq estimate exterior}
  \|u\|_{L^p(\Rm; L^{q}(\{x: |x|\geq 1\}))} \lesssim_{d, p, q}  \|(u_0, u_1)\|_{\mathcal{H}^1\times L^2}. 
 \end{equation} 
 A combination of this inequality with the interior part estimate \eqref{low triangle centre part} finishes the proof. 
\end{proof}

Now we are at the position to prove the main result of this chapter, i.e. Proposition \ref{strichartz estimate}. 

\paragraph{Homogeneous estimates} Let us first consider the homogeneous equation, namely $f\equiv 0$. The Strichartz estimate 
\[
  \|u\|_{L^p L^q (\Rm \times \Rm^d)} \lesssim_{d,p,q} \|(u_0,u_1)\|_{\mathcal{H}^1\times L^2}
\] 
has been verified in Corollary \ref{strichartz low triangle} if $(1/p,1/q)$ satisfies our assumptions and is located on or below the line segment $AD$
 (i.e. in the darker grey region in figure \ref{figure admpair}). If $(1/p,1/q)$ is located above $AD$ instead (i.e. in the lighter grey region), then  we may apply an interpolation among the line segment $AD$ (including $A$ but not $D$) and the points $B$, $C$ to deduce the corresponding Strichartz estimates. Please note that the Strichartz estimates corresponding to the points $B$ and $C$ have been verified in Proposition \ref{special strichartz estimates}. 
 
\paragraph{Inhomogeneous estimates} Now we switch to the contribution of inhomogeneous part. Let $u$ be the solution to 
\[
 \partial_t^2 u - \Delta u + \frac{u}{|x|} = f \in L^1 L^2(\Rm \times \Rm^d)
\]
 with zero initial data. We need to show 
 \[
  \|u\|_{L^p L^q(\Rm \times \Rm^d)} \lesssim_{d,p,q} \|f\|_{L^1 L^2(\Rm \times \Rm^d)}. 
 \]
 In fact we may rewrite ($t>0$)
 \[
  u(t) = \int_0^t \mathbf{S}_{\mathcal{C}}(t-\tau) (0, f(\tau)) {\rm d} \tau = \int_0^\infty \chi_{\tau} (t)  \mathbf{S}_{\mathcal C} (t-\tau) (0,f(\tau)) {\rm d} \tau
 \]
 Here $\mathbf{S}_{\mathcal C}$ is the linear propagation operator; $\chi_\tau (t)$ is the characteristic function of the interval $[\tau,+\infty)$. It immediately follows the Strichartz estimates for the homogeneous linear Coulomb wave equation that 
 \begin{align*}
  \|u\|_{L^p L^q (\Rm^+ \times \Rm^d)} &\leq \int_0^\infty \left\|\chi_{\tau} (t)  \mathbf{S}_{\mathcal C} (t-\tau) (0,f(\tau))\right\|_{L^p L^q (\Rm^+\times \Rm^d)} {\rm d} \tau\\
  & \lesssim_{d,p,q} \int_0^\infty \|f(\tau)\|_{L^2} {\rm d} \tau \\
  & \lesssim_{d,p,q} \|f\|_{L^1 L^2(\Rm^+\times \Rm^d)}. 
 \end{align*}
 This actually finishes our proof because the negative time direction is similar. 
 
 \section{Scaling identity for non-radial estimates}
  Before we conclude this chapter, we give an outline of the proof for the claim given in Remark \ref{nonradial energy level}, i.e. if the Strichartz estimates 
 \[
  \|u\|_{L^p L^q (\Rm \times \Rm^d)} \lesssim \|(u_0,u_1)\|_{\mathcal{H}^1\times L^2}
 \] 
 holds for all free Coulomb waves, possibly with non-radial initial data $(u_0,u_1)$, then we must have
 \[
  \frac{1}{p} + \frac{d}{q} = \frac{d}{2} - 1. 
 \]
 Let $(\varphi, \psi) \in \mathcal{C}_0^\infty(B(0,R_0))$ be initial data and $v$ be the linear free wave $\mathbf{S}_{\mathcal{W}} (\varphi, \psi)$. Given $R_1 > 0$, we consider the solution $u_R$ to 
 \[
  \partial_t^2 u - \Delta u + \frac{u}{|x|} = 0
 \] 
 with initial data $(\varphi(x-R\vec{e}_1), \psi(x-R\vec{e}_1))$. Here $R \gg \max\{R_0,R_1\}$ and $\vec{e}_1=(1,0,\cdots, 0)$. It follows that $w_R (x,t) = u_R (x,t) - v(x-R\vec{e}_1,t)$ solves the linear wave equation 
 \begin{equation} \label{equation w diff}
  \partial_t^2 w_R - \Delta w_R + \frac{w_R}{|x|} = - \frac{v(x-R\vec{e}_1,t)}{|x|} 
 \end{equation} 
with zero initial data. Finite speed of propagation implies that both $w_R$ and $v(x-R\vec{e}_1,t)$ are supported in $B(R \vec{e}_1, |t|+R_0)$. Thus we have for $|t|\leq R_1$ that
\begin{align*}
 \left\|\frac{v(x-R\vec{e}_1,t)}{|x|}\right\|_{L^2} & \lesssim_1 \frac{R_0+R_1}{R}  \left\|\frac{v(x-R\vec{e}_1,t)}{|x-R\vec{e}_1|}\right\|_{L^2} \\
 & \lesssim_d  \frac{R_0+R_1}{R} \|v(t)\|_{\dot{H}^1} \lesssim_d \frac{R_0+R_1}{R} \|(\varphi,\psi)\|_{\dot{H}^1\times L^2}.
\end{align*}
 A similar argument shows that
 \[
  \left\|\frac{w_R (x,t)}{|x|}\right\|_{L^2} \lesssim_d \frac{R_0+R_1}{R} \|w_R (t)\|_{\dot{H}^1}. 
 \]
 We then apply the Strichartz estimates of regular wave equation on \eqref{equation w diff} and obtain 
 \begin{align*}
  \sup_{t\in [-R_1,R_1]} \|w_R (t)\|_{\dot{H}^1} &\leq \int_{-R_1}^{R_1} \left\|\frac{w_R (t)}{|x|} + \frac{v(x-R\vec{e}_1,t)}{|x|}\right\|_{L^2} {\rm d} t \\
  & \lesssim_d \frac{R_0 +R_1}{R} \int_{-R_1}^{R_1} \left(\|(\varphi,\psi)\|_{\dot{H}^1\times L^2} + \|w_R (t)\|_{\dot{H}^1}\right) {\rm d} t \\
  & \lesssim_d \frac{R_1(R_0 + R_1)}{R} \left(\|(\varphi,\psi)\|_{\dot{H}^1\times L^2} + \sup_{t\in [-R_1,R_1]} \|w_R (t)\|_{\dot{H}^1}\right). 
 \end{align*}
 This implies that 
 \[
  \sup_{t\in [-R_1,R_1]} \|w_R (t)\|_{\dot{H}^1}  \rightarrow 0, \quad {\rm as}\; R \rightarrow +\infty. 
 \]
 Thus 
 \[
  \sup_{t\in [-R_1,R_1]} \|u_R (x+R \vec{e}_1,t) - v(x,t)\|_{L^{2^*}(\Rm^d)}  \rightarrow 0, \quad {\rm as}\; R \rightarrow +\infty.
 \]
 By Fatou's lemma we obtain 
 \begin{align*}
  \|v\|_{L^p L^q ([-R_1,R_1]\times \Rm^d)} &\leq \liminf_{R\rightarrow +\infty} \|u_R (x+R \vec{e}_1,t) \|_{L^p L^q ([-R_1,R_1]\times \Rm^d)}\\
  & \leq \liminf_{R\rightarrow +\infty} \|u_R (x,t) \|_{L^p L^q (\Rm \times \Rm^d)}\\
  & \leq C \liminf_{R\rightarrow +\infty} \left\|(\varphi(x-R\vec{e}_1), \psi(x-R\vec{e}_1))\right\|_{\mathcal{H}^1 \times L^2}\\
  & \leq C \left\|(\varphi, \psi)\right\|_{\dot{H}^1 \times L^2}.
 \end{align*}
 Here $C$ is a constant independent of $R_1$ and $(\varphi, \psi)$. We then let $R_1 \rightarrow +\infty$ and obtain 
 \[
  \|v\|_{L^p L^q (\Rm \times \Rm^d)} \leq C\left\|(\varphi, \psi)\right\|_{\dot{H}^1 \times L^2}. 
 \]
 Finally we recall that $C_0^\infty(\Rm^d)$ are dense in the spaces $\dot{H}^1$ and $L^2$, and deduce that the Strichartz estimate
  \[
  \|v\|_{L^p L^q (\Rm \times \Rm^d)} \leq C\left\|(\varphi, \psi)\right\|_{\dot{H}^1 \times L^2}
 \]
 holds for all linear free waves with initial data $(\varphi, \psi) \in \dot{H}^1 \times L^2$. This immediately gives the rescaling identity 
 \[
  \frac{1}{p} + \frac{d}{q} = \frac{d}{2} - 1, 
 \]
 because the linear free wave equation admits a natural rescaling invariance. 
 
\chapter{Scattering of radial solutions}

In this chapter we prove the scattering of radial solutions to defocusing immediate Coulomb wave equation. The proof follows a combination of the radial Strichartz estimates given in the previous chapter and the global decay estimate given by the inward/outward energy theory. We also give a small data scattering theory for possibly focusing nonlinear Coulomb wave equation in the final section. We start by considering the defocusing case. 

\section{Abstract theory}

\begin{lemma} \label{lemma abstract scattering}
 Assume $3\leq d\leq 5$ and that $(p_1,q_1)$ is a pair satisfying the conditions in Proposition \ref{strichartz estimate}. Let $1+\frac{4}{d-1} \leq p \leq 1+\frac{4}{d-2}$, $1\leq p_2, q_2\leq \infty$ and $\eta \in (0,p)$ be constants such that 
\[
 \left(1, \frac{1}{2}\right) = \eta \left(\frac{1}{p_1}, \frac{1}{q_1}\right) + (p-\eta)\left(\frac{1}{p_2}, \frac{1}{q_2}\right).
\]
If $u$ is radial finite-energy solution to defocusing Coulomb wave equation
\[
 \partial_t^2 u - \Delta u + \frac{u}{|x|} = - |u|^{p-1} u
\]
 satisfying
\[
 \lim_{t\rightarrow +\infty}\|u\|_{L^{p_2} L^{q_2}([t,+\infty)\times \Rm^d)} = 0,
\]
then we always have 
 \begin{align*}
  &u \in L^{p_1} L^{q_1} (\Rm^+\times \Rm^d); & &|u|^{p-1} u \in L^1 L^2(\Rm^+ \times \Rm^d).
 \end{align*}
 It immediately follows that $u$ scatters in the positive time direction. 
\end{lemma}
\begin{proof}
 The assumption on the indices implies that the following nonlinear estimate holds 
 \[
  \left\|-|u|^{p-1}u\right\|_{L^1 L^2 (J\times \Rm^d)} \leq \|u\|_{L^{p_1} L^{q_1}(J\times \Rm^d)}^\eta \|u\|_{L^{p_2} L^{q_2}(J\times \Rm^d)}^{p-\eta}.
 \]
 Here $J$ is an arbitrary time interval. We then apply the Strichartz estimate given in Proposition \ref{strichartz estimate} on $u$ and obtain ($0\leq T < T' < +\infty$)
 \begin{align}
  \|u\|_{L^{p_1} L^{q_1} ([T,T']\times \Rm^d)} & \leq C\|(u(T),u_t(T))\|_{\mathcal{H}^1\times L^2} + C \left\|-|u|^{p-1}u\right\|_{L^1 L^2 ([T,T']\times \Rm^d)} \nonumber \\
  & \leq C (2E)^{1/2} + C \|u\|_{L^{p_1} L^{q_1}([T,T']\times \Rm^d)}^\eta \|u\|_{L^{p_2} L^{q_2}([T,T']\times \Rm^d)}^{p-\eta}. \label{recurrence strichartz}
 \end{align}
 According to the local theory, $-|u|^{p-1}u \in L^1 L^2(J \times \Rm^d)$ for all finite time interval $J$. It follows that $u \in L^{p_1} L^{q_1} ([T,T']\times \Rm^d)$ for all $0\leq T < T' < +\infty$.  We may choose a sufficiently small constant $\varepsilon > 0$ so that the inequality 
 \begin{equation} \label{contradiction C}
  2C (2E)^{1/2} > C(2E)^{1/2} + C \varepsilon^{p-\eta} \left[2C (2E)^{1/2}\right]^\eta
 \end{equation}
 holds, and then choose a large time $T$ such that 
 \[
  \|u\|_{L^{p_2} L^{q_2}([T,+\infty)\times \Rm^d)} < \varepsilon. 
 \]
 We claim that 
 \begin{equation} \label{claim C}
  \|u\|_{L^{p_1} L^{q_1}([T,T']\times \Rm^d)} < 2C (2E)^{1/2}, \quad \forall T'>T. 
 \end{equation}
 If this were false, then by continuity there would exist a time $T' > T$, so that 
 \[
  \|u\|_{L^{p_1} L^{q_1}([T,T']\times \Rm^d)} = 2C (2E)^{1/2}.
 \]
 Inserting this into \eqref{recurrence strichartz} gives a contradiction with \eqref{contradiction C}. This verifies \eqref{claim C} and yields
 \[
   \|u\|_{L^{p_1} L^{q_1}([T,+\infty)\times \Rm^d)} \leq 2C (2E)^{1/2}.
 \]
 The estimate of $|u|^{p-1} u$ immediately follows the nonlinear estimate. Finally we let 
 \[
  (\varphi,\psi) = (u_0,u_1) + \int_0^\infty \vec{\mathbf{S}}_{\mathcal{C}}(-\tau) (0, -|u|^{p-1} u(\tau)) {\rm d} \tau \in \mathcal{H}^1 \times L^2. 
 \]
 A direct calculation shows that 
 \begin{align*}
  (u(t),u_t(t)) - \vec{\mathbf{S}}_{\mathcal{C}}(t) (\varphi,\psi) & = \int_t^\infty \vec{\mathbf{S}}_{\mathcal{C}}(t-\tau) (0, |u|^{p-1} u(\tau)) {\rm d} \tau.
 \end{align*}
 It immediately follows that 
 \begin{align*}
  \left\|(u(t),u_t(t)) - \vec{\mathbf{S}}_{\mathcal{C}}(t) (\varphi,\psi)\right\|_{\mathcal{H}^1\times L^2} & \leq \int_t^\infty \left\|\vec{\mathbf{S}}_{\mathcal{C}}(t-\tau) (0, |u|^{p-1} u(\tau))\right\|_{\mathcal{H}^1\times L^2} {\rm d} \tau\\
  & \leq \int_t^\infty \left\||u|^{p-1}u\right\|_{L^2} {\rm d} \tau \rightarrow 0. 
 \end{align*}
 This verifies the scattering of the solution $u$.
\end{proof}

\begin{remark}
 If $p_1<+\infty$, then the continuous dependence of $\|u\|_{L^{p_1} L^{q_1}([T,T']\times \Rm^d)}$ on $T'$ is clear. If $p_1= +\infty$, then we must have $2+\frac{4}{2d-3} \leq q_1 \leq 2 + \frac{4}{d-2}$. A combination of the radial Sobolev embedding
 \begin{equation} \label{radial Sobolev embedding mH1}
  \mathcal{H}^1(\Rm^d) \hookrightarrow L^{q}(\Rm^d), \qquad 2+ \frac{4}{2d-3} \leq q \leq 2 + \frac{4}{d-2} 
 \end{equation}
 given in Remark \ref{radial Sobolev embedding} and the fact $u(t) \in \mathcal{C}(\Rm; \mathcal{H}^1)$ yields
 \[
  u(t) \in \mathcal{C}(\Rm; L^{q_1}(\Rm^d)). 
 \]
 The continuous dependence of $\|u\|_{L^{\infty} L^{q_1}([T,T']\times \Rm^d)}$ on $T'$ then follows. In addition, if $p_2 < +\infty$, the limit 
 \[
  \lim_{t\rightarrow +\infty} \|u\|_{L^{p_2} L^{q_2}([t,+\infty)\times \Rm^d)} = 0
 \]
 holds as long as $\|u\|_{L^{p_2} L^{q_2}([t,+\infty)\times \Rm^d)} < +\infty$ for some time $t\in \Rm$. 
\end{remark}

\section{Scattering in the defocusing case}

\begin{figure}[h]
 \centering
 \includegraphics[scale=1.5]{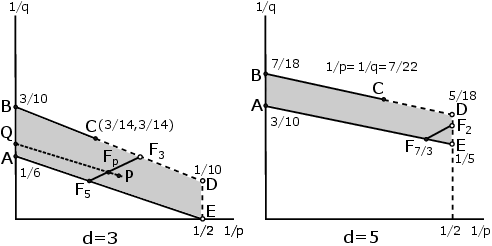}
 \caption{Location of $F_p$ for $d=3,5$} \label{figure admpairp1}
\end{figure}

\begin{figure}[h]
 \centering
 \includegraphics[scale=1.5]{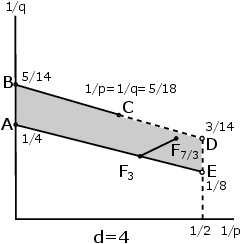}
 \caption{Location of $F_p$ for $d=4$} \label{figure admpairp2}
\end{figure}

Let $u$ be a radial solution to defocusing immediate Coulomb wave equation with initial data $(u_0,u_1)\in \mathcal{H}^1 \times L^2$. By the radial Sobolev embedding \eqref{radial Sobolev embedding mH1}, we have 
\[
 u_0 \in L^{p+1} (\Rm^d), \qquad 1+\frac{4}{d-1} \leq p \leq 1 + \frac{4}{d-2}. 
\]
Thus the solution $u$ comes with a finite energy. The inward/outward energy theory (Proposition \ref{energy distribution summary}) then immediately gives the decay 
\[
 \lim_{t\rightarrow +\infty} \|u(\cdot,t)\|_{L^{p+1}(\Rm^d)} = 0, 
\]
which gives the decay of $L^\infty L^{p+1}$ norm of $u$. 
 In view of this decay estimate and Lemma \ref{lemma abstract scattering}, the remaining work to prove Theorem \ref{main thm scattering} is to show the existence of parameters $p_1,q_1,\eta$ satisfying 
\begin{itemize}
 \item $(p_1,q_1)$ is a radial Coulomb allowed pair as in Proposition \ref{strichartz estimate}; 
 \item The following identity holds
 \[
  \left(1, \frac{1}{2}\right) = \eta \left(\frac{1}{p_1}, \frac{1}{q_1}\right) + (p-\eta)\left(0, \frac{1}{p+1}\right).
 \] 
\end{itemize}
The identity above is actually equivalent to 
\[
 \left(\frac{1}{p}, \frac{1}{2p}\right) = \frac{\eta}{p} \left(\frac{1}{p_1}, \frac{1}{q_1}\right) + \left(1-\frac{\eta}{p}\right)\left(0, \frac{1}{p+1}\right).
\]
From a geometrical perspective, this means that the sufficient and necessary condition for the existence of $\eta\in (0,p)$ satisying the identity above is that the point $F_p = \left(\frac{1}{p},\frac{1}{2p}\right)$ lies on the line segment between $P = \left(\frac{1}{p_1}, \frac{1}{q_1}\right)$ and $Q = \left(0, \frac{1}{p+1}\right)$. For $3\leq d\leq 5$, we consider a few different cases. (Please see figure \ref{figure admpairp1} and \ref{figure admpairp2} for the locations of the point $F_p$.)

\begin{itemize}
 \item[(i)] The energy-critical case $p=1+\frac{4}{d-2}$. In this case the point $F_p$ is located on the line segment $AE$, the bottom part of the boundary of the radial Coulomb allowed region. The point $(0,\frac{1}{p+1}) = (0, \frac{d-2}{2d})$ is exactly the point $A$. Thus we can always choose a point $(\frac{1}{p_1}, \frac{1}{q_1})$ between $F_p$ and $E$, which is still contained in the radial Coulomb allowed region. Clearly the pairs $(\frac{1}{p_1},\frac{1}{q_1})$ and $(0,\frac{1}{p+1})$ satisfy the geometric condition given above, which leads to the scattering of radial finite-energy solutions in this case. 
 \item[(ii)] The immediate case $1+\frac{4}{d-1}<p<1+\frac{4}{d-2}$. In this case the point $F_p$ is contained in the interior of the Coulomb allowed region. Thus we can always pick up a point $P$ on the extension line of the line segment $QF_p$ so that $P$ is still contained in the allowed region. Indeed, $P$ is in the allowed region as long as it is sufficiently close to $F_p$, as shown in the left part of figure \ref{figure admpairp1}. Again this implies the scattering of all radial finite-energy solutions in the immediate case. 
 \item[(iii)] We would also like to consider the conformal case $p = 1 + \frac{4}{d-1}$, although this is not covered by Theorem \ref{main thm scattering}. The situation is slight different in dimensions $d=3,5$ and in dimension $4$. In the case of $d=4$, the point $F_{7/3}$ is still contained in the interior of the Coulomb allowed region. As a result, a similar argument as in (ii) shows that any radial finite-energy solution in the conformal case of $d=4$ also scatters in the energy space. If $d=3$ or $d=5$, however, the point $F_3$ and $F_2$ are on the boundary of the Coulomb allowed region, respectively. 
\end{itemize}

As in the case of the classic nonlinear wave equation $\partial_t^2 u - \Delta u + |u|^{p-1} u$.  We may extend the range of $p$ in the scattering theory by imposing stronger smoothness and decay conditions on the initial data. For example, we have 

\begin{proposition} \label{weighted scattering}
 Let $\kappa > 0$ be a constant and $(u_0,u_1)\in \mathcal{H}^1(\Rm^3) \times L^2(\Rm^3)$ be radial initial data so that 
 \[
  E_\kappa (u_0,u_1) = \int_{\Rm^3} (|x|^\kappa + 1)\left(\frac{1}{2}|\nabla u_0(x)|^2 + \frac{1}{2} |u_1(x)|^2 + \frac{|u_0(x)|^2}{2|x|} + \frac{|u_0(x)|^4}{4}\right) {\rm d} x <+\infty.
 \]
 Then the corresponding solution to the conformal Coulomb wave equation 
 \[
  \partial_t^2 u - \Delta u + \frac{u}{|x|} + |u|^2 u = 0
 \]
 scatters in the energy space in both two time directions. 
\end{proposition}
\begin{proof}
 Again we prove this proposition by applying Lemma \ref{lemma abstract scattering}. We observe that the point $F_3$ is exactly located on the top part of the boundary of the Coulomb allowed region, i.e. the line segment $BD$. Thus if we can find a pair $(1/p_2,1/q_2)$ above the line segment $BD$ such that  
 \[
  \lim_{t\rightarrow \infty} \|u\|_{L^{p_2} L^{q_2}([t,+\infty)\times \Rm^d)} = 0
 \]
 holds for all radial finite-energy solution $u$ satisfying the assumption in Proposition \ref{weighted scattering}, then a similar argument as above immediately gives the scattering of these solutions. Indeed, weighted Morawetz estimate gives the decay (see Corollary \ref{coro weighted Morawetz})
\[
 E_-(t) \in L^{1/\kappa} (\Rm^+),
\]
Here we assume $\kappa>0$ is a small constant, without loss of generality. Combing this with the inequality (see \eqref{GN inequality})
\[
  \int_{\Rm^3} |u(x)|^{10/3} {\rm d} x \lesssim_d \|u\|_{\mathcal{H}^1}^{4/3} \int_{\Rm^3} \frac{|u(x)|^2}{|x|} {\rm d} x,
\]
we obtain 
\[
 \int_{\Rm^3} |u(x)|^{10/3} {\rm d} x \in L^{1/\kappa} (\Rm^+)\quad \Longrightarrow \quad u \in L^{10/3\kappa} L^{10/3} (\Rm^+ \times \Rm^3)
\]
This implies that 
\[
  \lim_{t\rightarrow \infty} \|u\|_{L^{10/3\kappa} L^{10/3} ([t,+\infty) \times \Rm^3)} = 0.
\]
The pair $(3\kappa/10,3/10)$ is located above the line segment $BD$. (see left part of figure \ref{figure admpairp1}) 
\end{proof}

\section{Scattering with small data}

As in the case of classic wave equation, small data scattering theory immediately follows from the corresponding Strichartz estimates. We have

\begin{proposition}\label{small data scattering}
 Let $3\leq d \leq 5$ and $1+\frac{4}{d-1}< p \leq 1+\frac{4}{d-2}$. Assume that the nonlinear function $f(u)$ satisfies 
 \begin{align*}
  &f(0) = 0,& &|f(u)-f(v)|\leq \eta (|u|^{p-1}+|v|^{p-1})|u-v|.
 \end{align*}
 Then there exists a constant $\delta = \delta(d,p,\eta)>0$, such that given any radial initial data $(u_0,u_1)$ satisfying $\|(u_0,u_1)\|_{\mathcal{H}^1\times L^2}<\delta$, the corresponding nonlinear Coulomb wave equation 
 \begin{equation} \label{to solve small data}
  \left\{\begin{array}{l} \partial_t^2 u + \mathbf{H} u = f(u); \\  (u,u_t)|_{t=0} = (u_0,u_1) \end{array}\right.
 \end{equation} 
 admits a unique global solution $u$ with 
 \[
  \|u\|_{L^p L^{2p}(\Rm \times \Rm^d)} \lesssim \|(u_0,u_1)\|_{\mathcal{H}^1\times L^2}. 
 \]
\end{proposition}

\begin{remark}
 Proposition \ref{small data scattering} covers the focusing power-like nonlinearity $f(u)=|u|^{p-1} u$ and the absolute value case $f(u) = |u|^p$. The latter one is frequently used in discussion of minimal existence time of solutions to wave equations. 
\end{remark}

\begin{proof}
 The proof of Proposition \ref{small data scattering} is a classic fixed-point argument. We recall that a radial solution to the linear Coulomb wave equation 
 \[
  \partial_t^2 u + \mathbf{H} u = g
 \]
 satisfies the Strichartz estimate 
 \begin{equation} \label{Strichartz small data}
  \|v\|_{L^p L^{2p}(\Rm \times \Rm^d)} \leq c\left(\|(v(0),v_t(0))\|_{\mathcal{H}^1\times L^2}+\|g\|_{L^1 L^2 (\Rm \times \Rm^d)}\right). 
 \end{equation}
 Here we use the fact that $(p,2p)$ is a Coulomb allowed pair(see figures \ref{figure admpairp1} and \ref{figure admpairp2}) and $c=c(p,d)$ is constant. We consider the complete metric space 
 \[
  X = \{u\in L^p L^{2p}(\Rm \times \Rm^d): u\; \hbox{is radially symmetric},\; \|u\|_{L^p L^{2p}} \leq 2c\|(u_0,u_1)\|_{\mathcal{H}^1\times L^2}\},
 \]
 equipped with the distance
 \[
  d(u,v) = \|u-v\|_{L^p L^{2p}(\Rm \times \Rm^d)};
 \]
 and the map $\mathbf{T}: X \rightarrow L^{p} L^{2p} (\Rm \times \Rm^d)$ defined by
 \[
  \mathbf{T} u = \mathbf{S}_{\mathcal{C}} (t) (u_0,u_1) + \int_0^t \mathbf{S}_{\mathcal{C}}(t-s) (0, f(u(s))) {\rm d} s. 
 \]
 Please note that $\mathbf{T} u$ is actually the solution to the linear Coulomb wave equation 
 \[
  \partial_t^2 v - \mathbf{H} v = f(u)
 \]
 with initial data $(u_0,u_1)$. A combination of the Strichartz estimate \eqref{Strichartz small data} and our assumption on the nonlinear function $f$ implies that 
 \begin{align*}
  \|\mathbf{T} u\|_{L^p L^{2p}} & \leq c\left(\|(u_0,u_1)\|_{\mathcal{H}^1 \times L^2} + \|f(u)\|_{L^1 L^2}\right) \\
  & \leq c\|(u_0,u_1)\|_{\mathcal{H}^1 \times L^2} + c\eta \|u\|_{L^p L^{2p}}^p; \\
  \|\mathbf{T} u - \mathbf{T} v\|_{L^p L^{2p}} & \leq c\|f(u)-f(v)\|_{L^1 L^2} \\
  & \leq c\eta \left(\|u\|_{L^p L^{2p}}^{p-1} + \|v\|_{L^p L^{2p}}^{p-1}\right) \|u-v\|_{L^p L^{2p}}. 
 \end{align*}
 We choose $\delta = \delta(d,p,\eta)$ sufficiently small so that 
 \begin{align*}
  &c\eta (2c\delta)^{p} < c \delta;& &2c\eta (2c\delta)^{p-1} < \frac{1}{2}.
 \end{align*}
 As a result, if $\|(u_0,u_1)\|_{\mathcal{H}^1\times L^2} < \delta$ is sufficiently small, then $\mathbf{T}$ becomes a contraction map from the metric space $X$ to itself. This immediately gives a unique fixed-point, which is exactly the solution to the corresponding nonlinear wave equation \eqref{to solve small data}. The global estimate 
\[
 \|u\|_{L^p L^{2p}} \leq 2c \|(u_0,u_1)\|_{\mathcal{H}^1\times L^2}
\]
then follows the definition of the space $X$.  
\end{proof}

\chapter{Appendix}

In this appendix we give a few results about the solutions to Klein-Gordon equation with good initial data. They are useful in the argument of this work. 

\begin{lemma} \label{fast decay outside parabola} 
 Assume that $u$ is a solution to the one-dimensional Klein-Gordon equation ($m > 0$)
 \[
  u_{tt}-u_{xx} + m^2 u = 0.
 \]
 with initial data $(u_0,u_1)$ in the Schwartz class $\mathcal{S}$. Let $\alpha <1$ and $N$ be two positive constants. Then the following inequality holds for all sufficiently large time $t$:
 \[
  |u(x,t)|+|u_x(x,t)|+|u_t(x,t))|\leq C(m,u_0,u_1,\alpha,N)|t|^{-N}, \qquad  |x|>|t|-|t|^\alpha.
 \]
\end{lemma}
\begin{proof}
 By rescaling it suffices to consider the case $m=1$. We may write the solution explicitly in terms of the Fourier transforms 
 \[
  u(x,t) = c_1 \int_{-\infty}^\infty e^{{\rm i} x \xi} \left(\hat{u}_0(\xi)\cos (t \sqrt{\xi^2 + 1}) + \hat{u}_1(\xi) \frac{\sin (t \sqrt{\xi^2+1})}{\sqrt{\xi^2+1}}\right) {\rm d} \xi
 \]
 Thus it suffices to prove that 
 \[
  v(x,t) = \int_{-\infty}^\infty e^{{\rm i}(x \xi + t \sqrt{\xi^2 + 1})} \varphi(\xi) {\rm d} \xi
 \]
 satisfies the following estimate for any $\varphi \in \mathcal{S}$: 
 \begin{equation} \label{decay to prove} 
  |v(x,t)| \leq C(\varphi,\alpha,N) |t|^{-N}, \qquad |x|>|t|-|t|^\alpha, \; |t|\gg 1. 
 \end{equation} 
 We fix a smooth cut-off function $\phi: \Rm \rightarrow [0,1]$ satisfying 
 \[
  \phi(s) = \left\{\begin{array}{ll} 1, & s\leq 1/2; \\ 0, & s\geq 1; \end{array}\right. 
 \]
 and write 
 \[
  v(x,t) = v_1(x,t) + v_2(x,t);
 \]
 with 
 \begin{align*}
  v_1(x,t) &= \int_{-\infty}^\infty  e^{{\rm i}(x \xi + t \sqrt{\xi^2 + 1})} (1-\phi(|t|^{-\beta}\xi)) \varphi(\xi) {\rm d} \xi; \\
  v_2(x,t) &= \int_{-\infty}^\infty  e^{{\rm i}(x \xi + t \sqrt{\xi^2 + 1})} \phi(|t|^{-\beta}\xi) \varphi(\xi) {\rm d} \xi.
 \end{align*}
 Here $\beta = \beta(\alpha) \in (0,1/4)$ is a small constant satisfying $0<\alpha < 1-2\beta$. The upper bound of $v_1$ can be given by 
 \[
  |v_1(x,t)| \leq \|(1-\phi(|t|^{-\beta}\xi)) \varphi(\xi)\|_{L^1(\Rm)} \leq \|\varphi\|_{L^1(\{\xi: |\xi|>|t|^{\beta}/2\})} \leq C(\varphi, \beta, N) |t|^{-N}. 
 \]
 In order to evaluate the upper bound of $v_2(x,t)$, we need to apply integration by parts. We first utilize the support of $\phi$ and rewrite  
 \[
  v_2(x,t) = \int_{-|t|^\beta}^{|t|^\beta} e^{{\rm i}(x \xi + t \sqrt{\xi^2 + 1})} \phi(|t|^{-\beta}\xi) \varphi(\xi) {\rm d} \xi. 
 \]
 An integration by parts then shows that 
 \begin{align*}
  v_2(x,t) & = \int_{-|t|^\beta}^{|t|^\beta}  \frac{\phi(|t|^{-\beta}\xi) \varphi(\xi)}{{\rm i}\left(x+ t\cdot \frac{\xi}{\sqrt{\xi^2+1}}\right)} {\rm d} e^{{\rm i}(x \xi + t \sqrt{\xi^2 + 1})}\\
   & = \int_{-|t|^\beta}^{|t|^\beta} \left[-\frac{\partial_\xi \left(\phi(|t|^{-\beta}\xi) \varphi(\xi)\right)}{{\rm i}\left(x+ t\cdot \frac{\xi}{\sqrt{\xi^2+1}}\right)} + \frac{\phi(|t|^{-\beta}\xi) \varphi(\xi) \frac{t}{(1+\xi^2)^{3/2}}}{{\rm i}\left(x+ t\cdot \frac{\xi}{\sqrt{\xi^2+1}}\right)^2}\right] e^{{\rm i}(x \xi + t \sqrt{\xi^2 + 1})} {\rm d} \xi.  \end{align*}
We may repeat this process of integration by parts for $M$ times and deduce that 
\begin{equation} \label{expression of v2 integration by parts}
 v_2(x,t) = \int_{-|t|^\beta}^{|t|^\beta} \left(\sum_{k=0}^{M} \frac{A_{M,k,\varphi,\beta} (t, \xi)\cdot t^{k}}{\left(x+ t\cdot \frac{\xi}{\sqrt{\xi^2+1}}\right)^{k+M}}\right) e^{{\rm i}(x \xi + t \sqrt{\xi^2 + 1})} {\rm d}\xi.
\end{equation}
Here $A_{M,k,\varphi,\beta}(t,\xi)$ is a finite linear combination of the following terms 
\begin{align*}
& \partial_\xi^{\ell+1} \left(\phi(|t|^{-\beta}\xi) \varphi(\xi)\right);& &\partial_\xi^{\ell} \left(\phi(|t|^{-\beta}\xi) \varphi(\xi)\right) P(\xi) (1+\xi^2)^{-m/2};
\end{align*}
where $\ell \geq 0$, $m\geq 3$ and $P(\xi)$ is a polynomial of $\xi$ with a degree smaller or equal to $m-3$. Therefore we have 
\begin{equation} \label{L1 AMk}
 \|A_{M,k,\varphi,\beta}(t,\xi)\|_{L_\xi^1(\Rm)} \leq C(M,k,\varphi,\beta)<+\infty.
\end{equation}
Next we observe that if $|t|\gg 1$ and $|\xi|\leq |t|^\beta$, then
\[
 \left| t\cdot \frac{\xi}{\sqrt{1+\xi^2}}\right| \leq \frac{|t|}{\sqrt{1+|\xi|^{-2}}} \leq \frac{|t|}{\sqrt{1+|t|^{-2\beta}}} \leq \frac{|t|}{1+|t|^{-2\beta}/3} \leq |t|(1-|t|^{-2\beta}/4).
\] 
Thus if we also assume $|x| > |t| - |t|^\alpha$, then we always have
\[
 \left|x + t\cdot \frac{\xi}{\sqrt{1+\xi^2}}\right| \geq |x| -\left| t\cdot \frac{\xi}{\sqrt{1+\xi^2}}\right|  \geq |x| - |t| + \frac{|t|^{1-2\beta}}{4} \geq \frac{|t|^{1-2\beta}}{4}  - |t|^\alpha \geq \frac{|t|^{1-2\beta}}{8}.
\]
Here we use the fact $\alpha < 1-2\beta$ and our assumption $|t|\gg 1$. Inserting this lower bound and the upper bound \eqref{L1 AMk} into \eqref{expression of v2 integration by parts}, we obtain 
\begin{align*}
 |v_2(x,t)| \leq C(M,\varphi,\beta) \sum_{k=0}^M \frac{|t|^k}{|t|^{(1-2\beta)(M+k)}} \leq C'(M,\varphi,\beta) |t|^{(4\beta-1)M}, \quad |x|>|t|-|t|^\alpha.
\end{align*}
We then choose a sufficiently large number $M=M(\beta,N)$ such that $(4\beta-1)M\leq -N$ and obtain ($|t|\gg 1$)
\[
 |v_2(x,t)| \leq C(\varphi,\alpha, N) |t|^{-N}, \quad |x|>|t|-|t|^\alpha.
\]
Combining this with the upper bound of $v_1(x,t)$, we finish the proof. 
\end{proof}
\begin{corollary} \label{dense subset K0}
 Let $\mathcal{K}$ and $\mathcal{K}_0$ be spaces defined in Section \ref{sec: construction of operator}. Then $\mathcal{K}_0$ is dense in $\mathcal{K}$. 
\end{corollary}
\begin{proof} 
 Given any $v \in \mathcal{K}$ and $\varepsilon > 0$, since solutions with smooth and compactly supported initial data are dense in the space $\mathcal{K}$, we may first choose such a solution $u$, so that $\|u-v\|_{\mathcal{K}} < \varepsilon /2$. We assume that the initial data are supported in the ball of radius $R$. By finite speed of propagation we have 
 \[
  u(x,t) = 0, \qquad |x|>R+|t|. 
 \]
 By Lemma \ref{fast decay outside parabola}, we may choose a large number $N\gg d$ and obtain for sufficiently large time $t>0$ that
 \[
   |u(x,t)|+|u_x(x,t)|+|u_t(x,t)|\leq C(u_0,u_1,N)|t|^{-N}, \qquad  t-R < |x| < t+R. 
 \]
 This implies that if we fix a smooth cut-off function 
 \[
  \phi(s) = \left\{\begin{array}{ll} 1, & s<-R; \\ 0, & s>0; \end{array}\right. 
 \]
 then the solution $\tilde{u}_{t_0} (t) = \mathbf{S}_\mathcal{K} (t-t_0) (\phi(|x|-t_0) u(x,t_0), \phi(|x|-t_0) u_t(x,t_0)) \in \mathcal{K}_0$  
 satisfies
 \[ 
  \lim_{t_0\rightarrow +\infty} \|\tilde{u}_{t_0} - u\|_{\mathcal{K}} = 0.
 \]
 As a result, we may choose a sufficiently large time $t_0$ so that $\|\tilde{u}_{t_0} - u\|_{\mathcal{K}}< \varepsilon/2$. This implies that $\|\tilde{u}_{t_0} - v\|<\varepsilon$ and finishes the proof. 
\end{proof}

\begin{lemma} \label{general decay of KG} 
 Assume that $u$ is a solution to the one-dimensional Klein-Gordon equation ($c>0$)
 \[
  u_{tt} - u_{xx} + c u = 0. 
 \]
 with initial data $(u_0,u_1)$ in the Schwartz class $\mathcal{S}$. Then the following inequality holds for all sufficiently large time $t$:
 \[
  |u(x,t)|+ |u_x(x,t)| + |u_t(x,t)| \leq C(u_0,u_1,c) |t|^{-1/2}. 
 \]
\end{lemma}
\begin{remark}
 This kind of dispersive estimates is well known. In fact, Lemma \ref{general decay of KG} is a direct consequence of the following dispersive estimate
 \[
  \|u(t)\|_{B_{\infty,2}^\sigma(\Rm)} \leq C |t|^{-1/2} \left(\|u_0\|_{B_{1,2}^{\sigma+3/2}(\Rm)} + \|u_1\|_{B_{1,2}^{\sigma+1/2}(\Rm)}\right), \qquad \sigma\in \Rm
 \]
 and the embedding $B_{\infty,2}^\sigma \hookrightarrow L^\infty$ for $\sigma>0$. For more details on these dispersive estimates, please refer to Brenner \cite{dispersive2} and Ginibre-Velo \cite{dispersive1}. Here we still give a proof of Lemma \ref{general decay of KG} for the reason of completeness. 
\end{remark}

\begin{proof}
By a rescaling, we assume $c=1$ without loss of generality. Again we may write the solution explicitly in terms of the Fourier transforms 
 \[
  u(x,t) = c_1 \int_{-\infty}^\infty e^{{\rm i} x \xi} \left(\hat{u}_0(\xi)\cos (t \sqrt{\xi^2 + 1}) + \hat{u}_1(\xi) \frac{\sin (t \sqrt{\xi^2+1})}{\sqrt{\xi^2+1}}\right) {\rm d} \xi
 \]
 Thus it suffices to prove that 
 \[
  v(x,t) = \int_{-\infty}^\infty e^{{\rm i}(x \xi + t \sqrt{\xi^2 + 1})} \varphi(\xi) {\rm d} \xi
 \]
 satisfies 
 \begin{equation} \label{exp to prove} 
  |v(x,t)| \leq C(\varphi) |t|^{-1/2}, \qquad |t|\gg 1. 
 \end{equation} 
 Here $\varphi \in \mathcal{S}$. Without loss of generality we assume $t > 0$. We split the real line $\Rm$ into two parts:
 \begin{align*}
  &J_1=\left\{\xi\in \Rm: \left|x+ \frac{\xi}{\sqrt{\xi^2+1}} t\right|\geq t^{1/2} \right\}; & &J_2=\left\{\xi\in \Rm: \left|x+ \frac{\xi}{\sqrt{\xi^2+1}} t\right|< t^{1/2} \right\}; &
 \end{align*} 
 and write $v(x,t)= v_1(x,t) + v_2(x,t)$ accordingly. Here $v_1$ and $v_2$ are defined by
 \[
  v_k = \int_{J_k} e^{{\rm i}(x \xi + t \sqrt{\xi^2 + 1})} \varphi(\xi) {\rm d} \xi. 
 \]
 Since the function $x+ \frac{\xi}{\sqrt{\xi^2+1}} t$ is a strictly increasing function of $\xi$, with limits $x-t$, $x+t$ at the negative and positive infinity, respectively, the sets $J_1$ and $J_2$ satisfy either of the following, as long as $t$ is sufficiently large. 
 \begin{itemize}
  \item $J_1 = \Rm$, $J_2 = \varnothing$;
  \item $J_1 = (-\infty, \xi_1]$, $J_2 = (\xi_1, +\infty)$; here $x+ \frac{\xi_1}{\sqrt{\xi_1^2+1}} t = -t^{1/2}$; 
  \item $J_1 = (-\infty, \xi_1]\cup [\xi_2,+\infty)$, $J_2 = (\xi_1, \xi_2)$; here $x+ \frac{\xi_1}{\sqrt{\xi_1^2+1}} t = -t^{1/2}$, $x+ \frac{\xi_2}{\sqrt{\xi_2^2+1}} t = t^{1/2}$;
  \item $J_1 = [\xi_2, +\infty)$, $J_2= (-\infty, \xi_2)$; here $x+ \frac{\xi_2}{\sqrt{\xi_2^2+1}} t = t^{1/2}$. 
 \end{itemize} 
 The integral over $J_1$ can be evaluated via the integration by parts. Let us consider the case $J_1 = [\xi_2, +\infty)$ as an example. In fact, all other cases can be dealt with in the same manner. We have 
 \begin{align*}
  v_1(x,t) & = {-\rm i}\int_{\xi_2}^\infty \left[\frac{\rm d}{{\rm d} \xi} e^{{\rm i}(x \xi + t \sqrt{\xi^2 + 1})}\right] \frac{\varphi(\xi)}{x+ \frac{\xi}{\sqrt{\xi^2+1}}t} {\rm d} \xi \\
  & = {\rm i} e^{{\rm i}(x \xi_2 + t \sqrt{\xi_2^2 + 1})} \frac{\varphi(\xi_2)}{x+ \frac{\xi_2}{\sqrt{\xi_2^2+1}}t} + {\rm i} \int_{\xi_2}^\infty e^{{\rm i}(x \xi + t \sqrt{\xi^2 + 1})} \frac{\varphi_\xi (\xi)}{x+ \frac{\xi}{\sqrt{\xi^2+1}}t} {\rm d} \xi \\
  & \qquad + {\rm i} \int_{\xi_2}^\infty e^{{\rm i}(x \xi + t \sqrt{\xi^2 + 1})} \varphi(\xi) \frac{\rm d}{{\rm d}\xi}\left[\frac{1}{x+ \frac{\xi}{\sqrt{\xi^2+1}}t}\right] {\rm d} \xi.
 \end{align*}
 We recall the definition of $J_1$ and the identity $x+ \frac{\xi_2}{\sqrt{\xi_2^2+1}} t = t^{1/2}$ to obtain 
 \begin{align*}
  |v_1(x,t)| \leq t^{-1/2} \left(\sup |\varphi|\right) + t^{-1/2} \|\varphi_\xi\|_{L^1(\Rm)} + (\sup |\varphi|) \int_{\xi_2}^\infty \left|\frac{\rm d}{{\rm d}\xi}\left[\frac{1}{x+ \frac{\xi}{\sqrt{\xi^2+1}}t}\right] \right| {\rm d} \xi.
 \end{align*}
 A simple calculation shows that 
 \[
  \frac{\rm d}{{\rm d}\xi}\left[\frac{1}{x+ \frac{\xi}{\sqrt{\xi^2+1}}t}\right] = - \left(x+ \frac{\xi}{\sqrt{\xi^2+1}}t\right)^{-2} (1+\xi^2)^{-3/2} t < 0.
 \]
 Thus 
 \[
  \int_{\xi_2}^\infty \left|\frac{\rm d}{{\rm d}\xi}\left[\frac{1}{x+ \frac{\xi}{\sqrt{\xi^2+1}}t}\right] \right| {\rm d} \xi = \frac{1}{x+ \frac{\xi_2}{\sqrt{\xi_2^2+1}}t} - \frac{1}{x+t}\leq t^{-1/2}. 
 \]
 It immediately follows that 
 \[
   v_1(x,t)  \lesssim_1 t^{-1/2} \left(\sup |\varphi| + \|\varphi_\xi\|_{L^1(\Rm)} \right). 
 \]
 Next we consider the integral over $J_2$. If $J_2 = \varnothing$, then $v_2(x,t) = 0$. Otherwise we may write $J_2 = (a,b)$ and 
 \begin{align*}
  |v_2(x,t)| \leq \int_a^b |\varphi(\xi)| {\rm d} \xi &\leq \left(\sup_{\xi\in \Rm} (\xi^2+1)^{3/2} |\varphi(\xi)|\right) t^{-1} \int_a^b t (\xi^2+1)^{-3/2} {\rm d} \xi \\
  & \leq \left(\sup_{\xi\in \Rm} (\xi^2+1)^{3/2} |\varphi(\xi)|\right) t^{-1} \left[\left(x+\frac{b}{\sqrt{b^2+1}} t\right)-\left(x+\frac{a}{\sqrt{a^2+1}} t\right)\right]\\
  & \leq 2\left(\sup_{\xi\in \Rm} (\xi^2+1)^{3/2} |\varphi(\xi)|\right) t^{-1/2}. 
 \end{align*}
 Here we use the fact that the values of $x + \frac{\xi}{\sqrt{\xi^2+1}} t$ at the endpoints $a, b$ of $J_2$ must be contained in the interval $[-t^{1/2}, t^{1/2}]$, which is a direct consequence of the convergence of $x + \frac{\xi}{\sqrt{\xi^2+1}} t$ at the endpoints and the definition of $J_2$. In summary we always have 
 \[ 
  |v(x,t)| \leq C(\varphi) |t|^{-1/2}, \qquad |t|\gg 1.
 \]
 This finishes the proof. 
\end{proof}

\section*{Acknowledgement}
Ruipeng Shen is financially supported by National Natural Science Foundation of China Project 12471230.

\end{document}